\documentclass[12pt,a4paper]{amsart}
\usepackage[a4paper,margin=30mm,headheight=15pt,headsep=7mm,footskip=10mm]{geometry}
\makeatletter
\renewcommand\normalsize{\@setfontsize\normalsize{12}{14.5}%
 \abovedisplayskip 9pt plus 2pt minus 2pt
 \abovedisplayshortskip 4pt plus 2pt
 \belowdisplayshortskip 6pt plus 2pt minus 2pt
 \belowdisplayskip \abovedisplayskip}
\normalsize

\makeatother
\usepackage[T1]{fontenc}

\usepackage[utf8]{inputenc}
\usepackage{lmodern}
\usepackage{mathtools,booktabs,longtable,array,enumitem,microtype,xcolor}
\usepackage{amsmath,amssymb,amsfonts}%
\usepackage{amsthm}%
\usepackage{mathrsfs}%
\usepackage{tikz-cd}
\usetikzlibrary{arrows.meta,calc,positioning,decorations.pathreplacing,decorations.pathmorphing}
\tikzset{%
    symbol/.style={%
        draw=none,
        every to/.append style={%
            edge node={node [sloped, allow upside down, auto=false]{$#1$}}}
    }
}
\usepackage{adjustbox}
\usepackage[cmtip,all]{xy}
\newtheorem{theorem}{Theorem}[section]
\newtheorem{proposition}[theorem]{Proposition}%

\newtheorem{lemma}[theorem]{Lemma}

\newtheorem{example}[theorem]{Example}%
\newtheorem{remark}[theorem]{Remark}%
\newtheorem{corollary}[theorem]{Corollary}

\newcommand{\e}{\operatorname{Ext}}

\newcommand{\CC}{\mathbb{C}}
\renewcommand{\P}{\mathbb{P}}

\newcommand{\I}{\mathcal{F}}

\renewcommand{\O}{\mathcal{O}}
\newcommand{\Z}{\mathbb{Z}}

\newcommand{\Q}{\mathbb{Q}}

\newcommand{\Hom}{\operatorname{Hom}}
\newcommand{\ch}{\mathsf{ch}}
\newcommand{\ext}{\operatorname{Ext}}

\newcommand{\A}{\mathbb{A}}

\newcommand{\cX}{\mathcal{X}}

\usepackage{needspace}
\usepackage{hyperref, cleveref}
\numberwithin{equation}{section}

\renewcommand{\P}{\mathbb{P}}

\renewcommand{\O}{\mathcal{O}}

\newcommand{\R}{\mathrm R}
\newcommand{\LLf}{\mathrm L}
\newcommand{\Perf}{\operatorname{Perf}}
\newcommand{\RPerf}{\mathfrak{Perf}}

\newcommand{\RH}{\operatorname{RHom}}
\newcommand{\cRHom}{\mathcal{R}\!\operatorname{Hom}}
\newcommand{\End}{\operatorname{End}}
\newcommand{\Tor}{\operatorname{Tor}}
\newcommand{\Spec}{\operatorname{Spec}}

\newcommand{\Sym}{\operatorname{Sym}}
\newcommand{\Tot}{\operatorname{Tot}}
\newcommand{\Crit}{\operatorname{Crit}}
\newcommand{\dCrit}{\operatorname{dCrit}}
\newcommand{\MF}{\operatorname{MF}}
\newcommand{\HP}{\operatorname{HP}}

\newcommand{\CyclicChains}{\operatorname{CC}}
\newcommand{\fib}{\operatorname{fib}}

\newcommand{\tr}{\operatorname{tr}}
\newcommand{\id}{\operatorname{id}}
\newcommand{\rk}{\operatorname{rk}}

\newcommand{\At}{\operatorname{At}}

\newcommand{\Diff}{\operatorname{Diff}}

\newcommand{\otL}{\mathbin{\otimes^{\mathbf L}}}
\newcommand{\timesh}{\mathbin{\times^{h}}}

\newcommand{\fM}{\mathfrak M}

\newcommand{\MC}{\operatorname{MC}}

\newcommand{\Gm}{\mathbb G_m}
\newcommand{\red}{\mathrm{red}}
\newcommand{\an}{\mathrm{an}}

\newcommand{\per}{\mathrm{per}}

\newcommand{\RR}{\mathbb R}

\newcommand{\Pic}{\operatorname{Pic}}
\newcommand{\Res}{\operatorname{Res}}

\newcommand{\Bl}{\operatorname{Bl}}
\newcommand{\GL}{\operatorname{GL}}

\newcommand{\sslash}{\mathbin{/\mkern-6mu/}}

\renewcommand{\Bar}{\operatorname{Bar}}

\theoremstyle{definition}
\newtheorem{hypothesis}[theorem]{Hypothesis}

\theoremstyle{remark}
\newtheorem{problem}[theorem]{Problem}
\theoremstyle{plain}
\DeclareFontFamily{U}{rsfs}{\skewchar\font127}
\DeclareFontShape{U}{rsfs}{m}{n}{<-6.5> rsfs5 <6.5-7.5> rsfs7 <7.5-> rsfs10}{}
\setlist[enumerate]{itemsep=3pt,topsep=5pt}
\hypersetup{hidelinks,bookmarksnumbered=true,
 pdftitle={Degeneration of (-1)-shifted symplectic stacks and categorification of Donaldson--Thomas invariants},
 pdfauthor={Artan Sheshmani},
 pdfsubject={Derived Lagrangian intersections, Calabi--Yau duality, and cohomological Donaldson--Thomas specialization}}
\makeatletter
\newcommand{\IHEScopylabel}[2]{%
  \@ifundefined{r@#2}{}{%
    \protected@write\@auxout{}{\string\newlabel{#1}{\csname r@#2\endcsname}}}}

\makeatother

\title[Degeneration of shifted symplectic stacks]{Degeneration of $(-1)$-shifted symplectic stacks and categorification of Donaldson--Thomas invariants}
\author[A. Sheshmani]{Artan Sheshmani}
\address{Artan Sheshmani}
\address{Beijing Institute of Mathematical Sciences and Applications, A6, Room 205
No. 544, Hefangkou Village, Huaibei Town, Huairou District, Beijing 101408,
China}
\address{Massachusetts Institute of Technology (MIT), IAiFi Institute, 182 Memorial Drive, Cambridge, MA 02139, USA}
\address{National Research University, Higher School of Economics, Russian Federation, Laboratory of Mirror Symmetry, NRU HSE, 6 Usacheva Street, Moscow, Russia, 119048}
\email{artan@mit.edu}

\date{October 28, 2026}

\providecommand{\id}{\operatorname{id}}

\begin{document}
\begin{abstract}
We study cohomological Donaldson--Thomas invariants of perfect
complexes in a semistable degeneration of projective Calabi--Yau
threefolds, with smooth components $X_1,X_2$ meeting along a smooth
anticanonical surface $S$.  The generic fiber of degeneration family is a $(-1)$-shifted symplectic stack. Over the special fiber, the derived restriction to $S$ is derived Lagrangian,
and its homotopy fiber product is also $(-1)$-shifted symplectic.  Our goal is to relate the categorical DT invariants of the generic fiber moduli, induced by periodic cyclic homology of matrix factorization category of 0-shifted potential functions over such moduli, to that of the special fiber obtained via derived Lagrangian intersection and geometric quantization. \\ This is achieved via construction of a flat Gauss-Manin connection that relates the generic fiber cohomological DT invariants to that of the special fiber. The
geometric framework, universal trace diagrams, and fourfold
foliation were developed jointly with Jacob Kryczka in \cite{KS}.
Here we construct ordinary critical realizations, oriented
cyclic-to-bar comparisons, and Gauss--Manin operators, and verify
finite specialization in explicit moduli problems. We workout our construction via several explicit examples such as Grassmannians, Relative Hilbert Schemes and moduli space of pure 2 dimensional sheaves with support moving in a fixed linear system of divisors in the ambient threefold. 
\end{abstract}

\keywords{Tyurin Degenerations, Derived Lagrangian Foliations, Categorification of Donaldson--Thomas Invariants}
\maketitle
\tableofcontents

\section{Introduction}\label{sec:intro}
Let $\pi:\cX\to C$ be a projective flat family over a smooth pointed
complex curve $(C,0)$. Assume that $\cX$ is a smooth fourfold,
that the fibers $X_b$ for $b\ne0$ are smooth Calabi--Yau threefolds,
and that the reduced central fiber is a simple normal-crossing union
\begin{equation}\label{eq:intro-family}
 X_0=X_1\cup_S X_2.
\end{equation}
The components $X_i$ and their common divisor $S$ are smooth.
A relative volume form identifies $\omega_{X_i}(S)$ with $\O_{X_i}$
and gives logarithmic forms with opposite residues on $S$.
Section~\ref{sec:geometry} defines this datum, its moduli functors,
and the orientation conventions in detail.

Our objects are perfect complexes. A coherent sheaf on a smooth fiber is an example
concentrated in degree zero; over the special fiber, derived restriction to $S$ may introduce
negative cohomology. We distinguish the category $\Perf(Z)$ from
its derived moduli stack $\RPerf(Z)$ and write
\[
 \mathfrak X_i=\RPerf(X_i),\qquad
 \mathfrak S=\RPerf(S),\qquad
 r_i:\mathfrak X_i\longrightarrow\mathfrak S,
 \quad E_i\longmapsto L\iota_i^*E_i.
\]
\cite{KS}, written jointly with Jacob Kryczka, proves that these maps
are Lagrangian on the specified locally geometric loci, and that
\begin{equation}\label{eq:intro-gluing}
 \mathfrak G=\mathfrak X_1\times^h_{\mathfrak S}\mathfrak X_2
\end{equation}
is $(-1)$-shifted symplectic
\cite[Theorem~2.26 and Corollary~2.27]{KS}.
An object is a triple $(E_1,E_2,\alpha)$, with a chosen equivalence
$\alpha:L\iota_1^*E_1\simeq L\iota_2^*E_2$ and its homotopy data.
The common restriction is a perfect complex, not necessarily a
sheaf. Thus the construction retains the complete boundary
information on the fixed components.

The purpose of this paper is to provide a completely explicit comparison for the cohomological DT
invariants of such moduli problems from $b\ne0$ to $b=0$, and also prove some of the conjectural properties of the Gauss-Manin connection proposed in \cite{KS}.
On a smooth fiber one uses specified ordinary critical charts and
their matrix-factorization categories; on the central fiber one
uses an oriented bar complex for the chosen derived-intersection
locus. The Lagrangian geometry supplies the central space and its
shifted pairing. To compare invariants one must additionally
construct compatible critical structures, orientation isomorphisms,
chain maps, connections, and finite cohomology. We develop these
constructions for perfect-complex moduli problems before verifying
the hypotheses in explicit examples.

\subsection{The general comparison problem}
Let $M_b$ be a finite-type classical moduli locus with a specified
oriented d-critical realization, and let $P_b$ be its perverse
sheaf of vanishing cycles. Our convention is
\[
 \mathrm{DT}_{\mathrm{coh}}^j(M_b)=
       \mathbb H_c^j(M_b^{\an},P_b).
\]
An ordinary critical chart consists of a smooth scheme $U_b$,
a degree-zero function $f_b$, and an identification of the locus
with $\operatorname{Crit}(f_b)$. The orientation is a square root
of the virtual canonical line bundle, with compatible chart
isomorphisms. The full derived
object stack, its classical scalar rigidification, and this
oriented critical realization are distinguished in
Section~\ref{ssec:sc-oriented-moduli}.

Let $t$ be a local coordinate on $C$ at $0$, and let $u$ be the
even formal periodicity parameter, independent of $t$.
On a smooth relative critical chart $U\to C$ with ordinary
potential $F$, put
\[
 K_F=(\Omega^\bullet_{U/C}((u)),D_F),\qquad
 D_F=u\,d_{U/C}-d_{U/C}F\wedge.
\]
The notation $((u))$ denotes Laurent series with a finite lower
bound; the algebraic and completed coefficient conventions are
specified in Section~\ref{sec:v5-formal-transport}.
On a product chart the base operator is
$\nabla_{\partial_t}=\partial_t-u^{-1}\partial_tF$.
Theorem~\ref{thm:v5-Cech-GM} constructs the operator on a smooth
ambient family without a chosen product trivialization and gives
its overlap homotopies. Theorem~\ref{thm:v5-periodic-bar} and
Proposition~\ref{prop:v5-horizontal-bar} construct the corresponding
operator on the periodic bar model.

For differential complexes with connections, a degree-zero chain
comparison $c:(P,d_P)\to(B,d_B)$ is horizontal up to homotopy when
\begin{equation}\label{eq:intro-horizontal}
 d_Bc=cd_P,\qquad
 \nabla_{B,\partial_t}c-c\nabla_{P,\partial_t}
       =d_BH+Hd_P,\qquad |H|=-1.
\end{equation}
It then intertwines the induced connections on cohomology.
The symbol $\mathcal B_{0,\mathrm{DT}}$ denotes the oriented
periodic DT bar complex on the chosen central locus, including
its mixed differential and critical-chart suspension. The
structure-algebra bar computing the derived intersection alone
is a different complex. We write $\mathsf K_b^{\mathrm{DT}}$
for a compatible normalized cyclic realization; on a global
critical chart of dimension $d$ the normalization includes $[d]$.
For nonproper moduli, the support convention is part of the
realization and must agree with the compact-support convention
in the displayed DT groups.

\begin{theorem}[Theorem~\ref{thm:finite-cyclic} and
Theorem~\ref{thm:categorical-degeneration-conditional}]
\label{thm:intro-cyclic}
Let $V\subset C$ be a connected affine neighborhood of $0$, put
$K=\CC((u))$ and $R=\Gamma(V,\O_C)\otimes_\CC K$.
A two-periodic complex of finite projective $R$-modules with a
$K$-connection defined at $0$ and commuting with its differential
has locally free cohomology, a flat induced connection, and
cohomology commuting with base change. Suppose that the
horizontal geometric comparisons of
Hypothesis~\ref{hyp:finite-cyclic} identify its nonzero fibers
with $\mathsf K_b^{\mathrm{DT}}$ and its central fiber with
$\R\Gamma_c(M_0,\mathcal B_{0,\mathrm{DT}})$. These
realizations are then fibers of the same flat bundle.
Over $K[[t]]$, horizontal transport with prescribed initial value
identifies the scalar extension of the central fiber with the
formal generic fiber over $K((t))$.
\end{theorem}

Here a connection defined at $0$ has no pole in the geometric
parameter $t$; powers of $u^{-1}$ belong to the separate periodic
coefficient field. Analytic transport between distinct complex
points requires an analytic realization and a chosen path.
The chain comparisons in the theorem
identify the actual geometric endpoints.

\begin{theorem}[Theorem~\ref{thm:smooth-proper-comparison}]
\label{thm:intro-general-perfect}
Let $p:M\to V$ be a smooth projective relative moduli locus of
simple perfect complexes, of relative dimension $d$. Suppose its
central objects are identified with their effective matching
triples, compatibly with the chosen d-critical realizations.
With the relative orientation $K_{M/V}$ and the associated
zero-critical cyclic and bar models, the fiber DT sheaves are
$\CC_{M_b}[d]$. Their periodic realization has the finite model
\[
 P^{\bar a}=\bigoplus_{j-d\equiv a\ (2)}
       H^j_{\mathrm{dR}}(M/V)\otimes_\CC\CC((u)),
 \qquad d_P=0,
\]
with the relative Gauss--Manin connection. The comparisons with
the oriented bar model are horizontal up to homotopy and commute
with specialization. The objects of $M$ may have arbitrary
finite amplitude; no section-complex presentation is required.
\end{theorem}

The proof constructs the finite model and its horizontal
comparisons from the relative de Rham complex. It separates the
geometric hypotheses on a selected locus from the presentation
of its individual perfect objects. For singular moduli loci,
the general criterion remains Theorem~\ref{thm:intro-cyclic}
with its explicit critical and finite-model hypotheses.
Proposition~\ref{prop:general-cyclic-reduction} constructs a local
reduced formal potential from a specified strictly unital cyclic
model of the full endomorphism algebra, including at rank zero.
Propositions~\ref{prop:general-perfect-fiber} and
\ref{prop:general-equivalence-transport} further show how actual
morphisms of perfect complexes and relative enhanced equivalences
produce compatible matching objects. Stability, properness, and
orientation comparisons are checked for the resulting loci.

\subsection{The comparison conjectures of Paper I}
\cite{KS} asks for a proper completion and a cohomological comparison
in a specified class of divisor-sheaf examples
\cite[Conjectures~1.3 and~1.4]{KS}. We recall their content.
Let $L$ be a line bundle over $\cX$, fix $m>0$, and assume on the fibers
under consideration that
\[
 H^1(X_b,\O_{X_b})=H^2(X_b,\O_{X_b})=0,\qquad
 H^{j}(X_b,L_b^m)=0\quad(j>0),
\]
and that $V_m=\pi_*L^m$ is locally free of rank $r_m>0$ and
commutes with base change, including at $0$.
The smooth fibers have $L_b$ ample and
$\Pic(X_b)=\mathbb Z[L_b]$; the relative volume form is fixed.
Set $M_m=\P_C(V_m)$, using one-dimensional subspaces, and
$d_m=r_m-1$. For an integral twist $k$, consider
\[
 \mathcal E^{m,k}=
 [L^{k-m}\otimes\O_{M_m}(-1)\longrightarrow L^k]
       \quad\text{on }\cX\times_C M_m.
\]

The first conjecture states that this family gives a smooth
proper relative classical locus in the scalar-rigidified moduli
of simple perfect complexes. On a smooth fiber the objects are
$\O_D\otimes L_b^k$, $D\in|L_b^m|$; at $0$ one keeps the
whole perfect complexes and their effective matching triples.
The second conjecture asks for an oriented bar comparison with
$(\Omega^\bullet_{M_m/C}((u)),u\,d)[d_m]$, compatible with
specialization, and a pole-free Gauss--Manin connection whose
horizontal cohomology classes are in terms of the hyperplane class powers.
These are particular moduli problems inside the general framework,
not a description of all perfect objects of the same charge.

\begin{theorem}[Theorem~\ref{thm:all-positive-degrees}]
\label{thm:intro-v6}
Let $h,q_4,q_5$ have degrees $1,4,5$ on $\P^4$, with
$S=(h=q_4=0)$ and $\Gamma=(h=q_4=q_5=0)$ smooth. Use the resolved
family $\cX=\Bl_{(h,t)}\{hq_4=tq_5\}\to C$ over a sufficiently
small Zariski open neighborhood of $0$, and put
$L=a^*\O_{\P^4}(1)$. For every $m>0$ and $k\in\mathbb Z$,
the family $\mathcal E^{m,k}$ verifies
\cite[Conjectures~1.3 and~1.4]{KS}. More explicitly,
\[
 r_m=\binom{m+4}{4}-\binom{m-1}{4},\qquad d_m=r_m-1,
 \qquad \binom{m-1}{4}=0\quad(1\le m\le4).
\]
With orientation $K_{M_m/C}$ the DT sheaf is
$\CC_{M_{m,b}}[d_m]$, and the oriented periodic comparison is
\[
 \mathcal B^{\mathrm{or}}_{M_m,\mathrm{per}}
     \simeq(\Omega^\bullet_{M_m/C}((u)),u\,d)[d_m].
\]
Its cohomology has rank $r_m$ in parity $d_m$ and zero rank in
the other parity. The connection is $\partial_t$ in the
hyperplane-power basis. The normalized ordinary cyclic
realization is $\CyclicChains_{\mathrm{per}}(\MF(M_{m,b},0))[d_m]$.
\end{theorem}

Section~\ref{sec:v6-hyperplanes} supplies the geometric construction,
reconstructs the classical deformation functor, and identifies
the actual central matching locus before calculating its DT
complex. For $m=1,2,3,4,5$ the ranks are $5,15,35,70,125$;
the degree-four example lies in odd parity because $d_4=69$.
The Grassmannian examples of Section~\ref{sec:grassmannian}
and the relative spherical twists of
Section~\ref{sec:complexes-several-sheaves} extend this calculation
to other ranks and to complexes with several cohomology sheaves,
including on smooth fibers.

\subsection{Examples with smooth classical moduli and derived obstruction directions}
\label{ssec:intro-smooth-examples}
The first examples have a smooth \emph{entire classical parameter space}.  This does
not make their reduced derived enhancements smooth.  For a critical
chart $(U,f)$ the tangent complex is
\[
 [T_U|_{\Crit(f)}\xrightarrow{\operatorname{Hess}(f)}
       \Omega_U|_{\Crit(f)}]
\]
in degrees zero and one.  For $f=0$ both terms remain and the
obstruction space is the dual of the tangent space.  The full Artin
stack of objects additionally retains scalar automorphisms and their
Serre-dual degree-three extensions; see Section~\ref{ssec:sc-oriented-moduli}.

\textbf{Complete linear systems and Grassmannians.}
Theorems~\ref{thm:all-positive-degrees}, \ref{thm:grassmannian-moduli}
and~\ref{thm:grassmannian-DT} treat every point of
$\mathbb P_C(V_m)$ and $\operatorname{Gr}_C(q,V_m)$, respectively. Note that in the level of underlying geometry, here, singular, reducible, or nonreduced supporting divisors, common zeros of
several sections, and central complexes with several cohomology sheaves
are all included in our analysis.  They do not however create singularities of these particular
classical object-moduli spaces.  \\The reconstruction of the complete
Artinian deformation functors proves this assertion.  If $d=q(r_m-q)$, the
Grassmannian DT groups on both fibers have a basis indexed by
partitions $\lambda\subset q\times(r_m-q)$, in degrees
$2|\lambda|-d$.  Their periodic rank is $\binom{r_m}{q}$ in parity
$d$.  The projective case has one class in degrees $2a-(r_m-1)$,
$0\leq a<r_m$.  The actual central restrictions, oriented bar maps,
and the horizontal hyperplane or Schubert classes are calculated.
The connection matrix in those classes is zero.  Relative equivalences
transport this calculation to the multiple-cohomology-sheaf examples.
The discussion at the end of Section~\ref{sec:v6-hyperplanes} also
explains why linearly dependent sections do not give a missing singular
stratum of the same simple-object Grassmannian problem.

\textbf{Points and reduced rigid curves.}
On a smooth quintic, the complete one-point and two-point incidence
spaces are projective bundles over $X$ and $X^{[2]}$, including double
points and singular surface supports.  Here by ``\textit{incidence}" we mean a subscheme, say, a point, two points, a point and a line etc, lying scheme theoretically on a divisor moving in a fixed linear system of divisors in the ambient threefold. Their potentials are zero
because their classical moduli spaces are proved smooth.  For $|H|$
the periodic dimensions are $(16,816)$ and $(62172,3672)$; for $|2H|$
they are $(56,2856)$ and $(269412,15912)$.  The nontrivial quintic
Gauss--Manin system occurs in these groups.  Its logarithmic matrices,
not just its curvature, are computed in
Section~\ref{sec:crossing-perfect}.  The central open calculations in
Section~\ref{sec:divisor-sheaves} are explicitly distinguished from a
proper central completion.  The perfect crossing objects and the
nonperfect strata are examined in
Propositions~\ref{prop:crossing-one-perfect},
\ref{prop:crossing-two-perfect}, and~\ref{prop:crossing-perfect-exhaustion}.

For a general quintic and the transverse ``\textit{limiting-curve hypotheses}",
Theorem~\ref{thm:line-full-comparison} gives the complete proper
comparison for the reduced line-ideal component: its $2875$ projective
planes contribute $\mathbb C^{2875}$ in degrees $-2,0,2$ on both
fibers.  Its degree-two version has $2875$ projective spaces of
dimension eleven.  The boundary self-intersection and the K3-pencil
structure-sheaf families are further complete zero-critical
calculations.  In the pencil, singular surface members are retained;
the classical object parameter remains the smooth affine pencil however.

\subsection{Examples with singular classical moduli and nonzero critical structure}
\label{ssec:intro-singular-examples}
Quasi-smoothness is a condition on a derived enhancement, not a claim
that its classical truncation is nonsingular.  The following examples
retain the gradient equations and their higher obstruction terms.
Their nonzero functions represent \textbf{nonzero} d-critical classes, rather
than being quadratic stabilizations of zero functions.

\textbf{Two defects on a rigid line.}
The singular locus of $Q=\operatorname{Quot}_X(I_\Lambda,2)$ is
$\Lambda$, parametrizing $I_\Lambda/\mathfrak m_pI_\Lambda$.
Theorem~\ref{thm:determinant-intrinsic} identifies the intrinsic local
potential with $\det X$ on $\operatorname{Mat}_3\times\mathbb A^s$.
For the surface objects, $s=1$ in the degree-one positivity
neighborhood and $s=10$ in degree two.  The local DT stalks occur in
$-6,-4,-1$ and $-15,-13,-10$, respectively.  Theorem~\ref{thm:determinant-relative}
computes an actual generic-to-central comparison of these germs,
including the central restriction triples and the horizontal bar.

Section~\ref{sec:global-two-defect} goes beyond a stalk calculation.
It constructs a small resolution of the entire proper singular
scheme $Q$, determines the extensions in its perverse DT sheaf, and
computes its global cohomology.  For one rigid line the periodic
dimensions are $(20735,1632)$.
Theorem~\ref{thm:global-quadratic-defect} identifies the oriented DT
sheaf on the full degree-two incidence scheme and gives periodic
dimensions $(16320,207350)$.  This includes all length-two quotients,
not just the smooth or cyclic quotients.  Proposition~\ref{prop:global-defect-GM}
computes the induced flat connection and its logarithmic limit.

\textbf{\textbf{A proper nonreduced component in a non-general quintic family}.}
Section~\ref{sec:cubic-geometry} chooses a coefficient path for which
the quintics are smooth but a selected line-moduli component is a
double point.  Its central component curves meet tangentially in $S$.
Theorems~\ref{thm:cubic-actual-moduli} and~\ref{thm:cubic-comparison}
construct the resulting proper object-moduli component and identify
its actual matching problem with the cubic presentation
$a(t)\lambda^3/3$ over a projective factor, with $a(t)$ a unit.
For $|H|$ the DT groups are $\mathbb C^2$ in degrees $-2,0,2$;
for $|2H|$ they are $\mathbb C^2$ in degrees $-11,-9,\ldots,11$.
The finite connection is
\[
 \partial_t-\frac{a'(t)}{3a(t)}
              \begin{pmatrix}1&0\\0&2\end{pmatrix}
\]
on each hyperplane-power summand.  Here both geometric endpoints
and the central oriented bar are identified; regular specialization
is proved, not inferred from flatness alone.
Proposition~\ref{prop:higher-contact-critical} gives the corresponding
calculation for an identified higher-contact component, with all its
nilpotent directions retained.  Its geometric identification is a
hypothesis of that proposition, not a claim about every special
quintic.

\textbf{\textbf{Boundary strata and the extent of specialization}.}
For a chosen global DT complex, Proposition~\ref{prop:boundary-localization}
retains the boundary groups and their attachment maps, and
Proposition~\ref{prop:specialization-cone-test} retains the cone of
specialization.  These maps are necessary: neither the sum of local
stalk dimensions nor an equality of Euler characteristics computes
the global invariant.  Section~\ref{sec:global-two-defect} illustrates
this concretely by calculating nontrivial extensions.  The known
nonzero monodromy of the full point and defect systems rules out
an unrestricted pole-free comparison.  Their remaining global
central comparison is formulated with nearby cycles, rather than
asserted as a consequence of quasi-smoothness.

\subsection{Organization and conventions}\label{Outlook}
Section~\ref{sec:geometry} develops the geometric setup and the
perfect-complex moduli functors, recalling the foundational
results of \cite{KS} with the sketches needed to specify the maps.
Section~\ref{sec:local} studies general deformation algebras,
normal operations, and ordinary potentials.
Sections~\ref{sec:cyclic} and~\ref{sec:quantization} construct
the Gauss--Manin operators and horizontal bar comparisons.
Section~\ref{sec:categorical-degeneration} proves the general
specialization criteria. The remaining sections give the
resolved-quintic, Grassmannian, divisor-sheaf, and K3-pencil
examples with their explicit calculations. Section~\ref{sec:singular-comparison}
identifies the determinant critical structure. Section~\ref{sec:crossing-perfect}
constructs crossing complexes and the logarithmic matrices. Section~\ref{sec:cubic-geometry}
proves the proper cubic comparison and gives the higher-contact
calculation. Section~\ref{sec:global-two-defect} computes the global
singular two-defect component, including its perverse extensions and
connection. The appendix supplies additional checks of the model
hypotheses.

Complexes are cohomologically graded: $V[1]$ places a vector
space in degree $-1$. We distinguish sheaf-valued $\cRHom$
from global $\RH$, and use $\otimes^{\LLf}$ and $\times^h$
for derived tensor and fiber products. The parameters $t,u,\hbar$
respectively describe geometric deformation, cyclic periodicity,
and quantization. Ordinary matrix-factorization potentials have
degree zero; the fourfold section--cosection potential has degree
$-1$. The main comparison concerns normalized periodic cyclic
and cohomological DT invariants. An equivalence of the
underlying dg categories requires additional categorical data;
where such an equivalence is claimed in an example, the categories
and the inducing isomorphism are specified.

\subsection*{Acknowledgements}
The author thanks Jacob Kryczka for the collaboration underlying
\cite{KS} and the comparison program developed there. Further thanks go to Maxim
Kontsevich, Ludmil Katzarkov, and Vladimir Baranovsky for discussions
during the earlier stages of this project. He is grateful to Oren
Ben-Bassat for correspondence on Lagrangian restriction, and to
Sheldon Katz, Charles Doran, Tony Pantev, Shing-Tung Yau, Martijn
Kool, Marco Robalo, Yukinobu Toda, and Bertrand To\"en for valuable
comments and discussions. The author acknowledges support by grants
from Beijing Institute of Mathematical Sciences and Applications (BIMSA),
the Beijing NSF BJNSF-IS24005, and the China National Science Foundation
(NSFC) NSFC-RFIS program W2432008. He would also like to thank China's
National Program of Overseas High Level Talent for generous support.
Finally, he would like to thank NSF AI Institute for Artificial
Intelligence and Fundamental Interactions at Massachusetts Institute
of Technology (MIT) which is funded by the US NSF grant under
Cooperative Agreement PHY-2019786.

\section{The geometric family and the moduli of perfect complexes}
\label{sec:geometry}
The comparison has three geometric terms: the relative moduli of objects on
$\cX/C$, the moduli of their restrictions on the two central components,
and the moduli of their common restriction on $S$. This section defines
these terms, the maps between them, and the deformation complexes which
will be used throughout the paper. The shifted forms and the fourfold
foliation are results of \cite{KS}. We state them with their hypotheses and
recall the parts of their proofs that fix the maps, shifts, and orientations
needed for the constructions below.

\subsection{An oriented degeneration and its boundary}
\label{ssec:geometric-family}
Work over $\CC$. Let $C$ be a smooth connected curve with a distinguished
point $0$, and write $C^*=C\setminus\{0\}$. In the projective setting our
geometric datum is a projective flat morphism
\[
 \pi:\cX\longrightarrow C
\]
from a smooth fourfold, smooth over $C^*$, with reduced central fiber
$X_0=X_1\cup_S X_2$. The threefolds $X_i$ and the surface
$S=X_1\cap X_2$ are smooth, and the components meet transversely.
Thus, etale locally at a point of $S$, the family has equation
\begin{equation}\label{eq:sc-local-family}
 xy=t,
\end{equation}
where $t$ is a local coordinate on $C$ vanishing at $0$.
An orientation is an isomorphism
$\vartheta:\O_{\cX}\xrightarrow{\sim}\omega_{\cX/C}$.
The dualizing sheaf is an invertible sheaf because this morphism is
Gorenstein. The smooth fibers are projective Calabi--Yau threefolds
with their specified volume forms. This is the datum of
\cite[Definition~2.6]{KS}.

Write $j_i:X_i\hookrightarrow X_0$, $j_S:S\hookrightarrow X_0$,
$\iota_i:S\hookrightarrow X_i$, and $i_b:X_b\hookrightarrow\cX$.
Adjunction identifies
\begin{equation}\label{eq:sc-log-volumes}
 \omega_{X_0}|_{X_i}=\omega_{X_i}(S),\qquad
 \omega_{X_i}(S)\simeq\O_{X_i},\qquad \omega_S\simeq\O_S.
\end{equation}
The two logarithmic volume forms $\eta_i$ induced by $\vartheta$
have opposite residues: $\Res_S\eta_1+\Res_S\eta_2=0$.
Locally these are represented by
$dy/y\wedge dz_1\wedge dz_2$ and
$-dx/x\wedge dz_1\wedge dz_2$ on the two branches of $xy=0$.
Their residues specify the same surface volume after changing the
sign on the second branch. This sign convention will also be used
for the two Lagrangian structures with common target.

The normal line bundles obey
\begin{equation}\label{eq:sc-normal-bundles}
 N_{S/X_1}\otimes N_{S/X_2}\simeq\O_S\otimes_\CC T_{C,0},
 \qquad N_{X_b/\cX}\simeq\O_{X_b}\otimes_\CC T_{C,b}.
\end{equation}
The coordinate $t$ trivializes the first constant factor. The
identities follow by comparing the two normal coordinates in
\eqref{eq:sc-local-family} and by taking the conormal of a fiber.
The constructions do not require the components or the total space
to be Fano. In particular, the anticanonical condition in
\eqref{eq:sc-log-volumes} is a condition on pairs $(X_i,S)$.

A base twist can be retained before choosing $\vartheta$.
Under the connectedness and cohomology-and-base-change assumptions
of \cite[Proposition~2.7]{KS}, the evaluation map gives
\[
 \lambda=\pi_*\omega_{\cX/C},\qquad
 \pi^*\lambda\xrightarrow{\sim}\omega_{\cX/C},\qquad
 \omega_{\cX}\simeq\pi^*(\omega_C\otimes\lambda).
\]
A frame of $\lambda$ gives the relative threefold orientation;
a frame of $\omega_C\otimes\lambda$ gives an absolute fourfold
volume form. These are separate choices. On a coordinate open with
a fixed relative orientation, the latter is the product of that
orientation with $dt$.

\subsection{The relative moduli functor and the selected objects}
\label{ssec:sc-perfect-moduli}
For a scheme or derived scheme $Z$, $\Perf(Z)$ is the stable
category of perfect complexes: its objects are locally quasi-isomorphic
to bounded complexes of finite-rank locally free sheaves. The derived
stack $\mathfrak{Perf}(Z)$ is its moduli stack. More precisely, a
$T$-point is a perfect complex on $Z\times^{\mathrm h}T$, with
isomorphisms interpreted as equivalences and with their higher
homotopies retained. For the family define
\begin{equation}\label{eq:sc-relative-functor}
 \mathfrak M(T)=\mathfrak{Perf}(\cX/C)(T)
 =\{(\beta,E_T):\beta:T\to C,\ E_T\in\Perf(\cX\times_C^{\mathrm h}T)\}^{\simeq}.
\end{equation}
The superscript denotes the maximal infinity-groupoid. A map
$T'\to T$ acts by derived pullback. The fiber at $b$ is
$\mathfrak{Perf}(X_b)$, including $b=0$; perfectness on the
singular fiber is part of this definition.

We use compatible locally geometric, locally finitely presented
derived Artin substacks of these functors, with perfect cotangent
complexes on the loci under discussion. Their existence is the
mapping-stack representability input of
\cite{ToenVaquie,PTVV13,Lieblich2006}; a finite-type moduli problem
also requires a bounded selection of objects. For example, one may
fix a numerical class, a bound on Tor amplitude, and an open
stability condition wherever such a condition has been specified.
A fixed charge alone is not a boundedness or properness assertion.
The numerical class over $C^*$ is understood in a chosen locally
constant system of classes, with the required monodromy invariance.

There is no fixed presentation imposed on the perfect complexes in
\eqref{eq:sc-relative-functor}. A flat coherent sheaf on a smooth
fiber is a perfect object with cohomology in degree zero. The
complexes in the examples may instead have several nonzero
cohomology sheaves. Given a coherent sheaf $F$ on a smooth
component and an effective Cartier inclusion $i:D\hookrightarrow Z$,
\begin{equation}\label{eq:sc-derived-cartier}
 i_*Li^*F\simeq[F(-D)\longrightarrow F],\qquad
 \mathcal H^{-1}(Li^*F)=\Tor_1^{\O_Z}(F,\O_D),\quad
 \mathcal H^0(Li^*F)=F\otimes\O_D.
\end{equation}
The first expression is a complex of sheaves on $Z$ and the last
two sheaves are on $D$. This is \cite[Lemma~2.2]{KS}: tensor the
Cartier resolution of $i_*\O_D$ with $F$. For a general perfect
complex the derived tensor product is taken in all degrees and its
cohomology need not be confined to these two degrees.

Two elementary local cases clarify these conditions. If $D=(x=0)$
in a smooth affine component and $F=\O_D$, then the Cartier
multiplication map is zero, so $Li^*F=\O_D[1]\oplus\O_D$.
This is an admissible perfect boundary object; its negative
cohomology must be included in the matching equivalence.
In contrast, near the central crossing put
$A=\CC[x,y,z,w]$ and $R=A/(xy)$. The sheaf $A/(x)$ has a
length-one free resolution over $A$, but as an $R$-module it has
an infinite resolution with alternating multiplication maps
$x,y$. Tensoring that resolution with the residue field at the
crossing leaves a nonzero term in every homological degree.
Thus this sheaf is perfect on the smooth total space but is not
perfect on the singular fiber. These are the distinctions in
\cite[Examples~2.3 and~2.19]{KS}. Perfectness on $X_0$ in
\eqref{eq:sc-relative-functor} cannot be replaced by perfectness
of its pushforward to $\cX$.

On a simple universally gluable locus we require
\[
 \operatorname{Ext}^{a}(E,E)=0\quad(a<0),\qquad
 \End(E)=\CC\,\id_E
\]
at geometric points. The vanishing of negative self-extensions is a
condition on the object, not on the number of its cohomology sheaves.
The classical simple stack has scalar automorphisms; its scalar
rigidification, when used below, is a classical algebraic space on
the representable loci of \cite{Lieblich2006}. We distinguish that
operation from any reduction of a derived obstruction theory.
The full derived stacks remain the spaces on which the restriction
and shifted-symplectic constructions are first made.

\subsection{Universal complexes, trace, and the relative shifted form}
\label{ssec:sc-relative-trace}
On $\cX\times_C^{\mathrm h}\mathfrak M$ let $\mathcal E$ be the
universal perfect complex and let $q$ be projection to $\mathfrak M$.
A universal object exists on the universal stack; on a coarse
algebraic space it may require a gerbe or compatible local universal
families. The endomorphism complexes and maps below are formed on
that stack or on those specified local families.

For a proper perfect morphism $p:Z\to T$, an orientation of relative
dimension $d$ supplies the trace
\[
 Rp_*\O_Z\longrightarrow Rp_*\omega_{Z/T}
       \longrightarrow\O_T[-d].
\]
Perfect pushforward and derived base change are the finiteness
properties needed to use this trace in a family. For perfect
$E,F$, relative duality identifies the adjoint of the composition--trace
pairing with
\[
 Rp_*\cRHom(E,F)\simeq
       (Rp_*\cRHom(F,E))^\vee[-d].
\]
These are the orientation and duality of
\cite[Definition~2.14 and Proposition~2.15]{KS}.

\begin{theorem}[{\cite[Theorem~3.7]{KS}}]\label{thm:sc-relative-symplectic}
Under the oriented family hypotheses above, on the specified
locally geometric loci the relative stack $\mathfrak M\to C$
has a relative $(-1)$-shifted symplectic form. Put
$\mathcal A=Rq_*\cRHom(\mathcal E,\mathcal E)$. Then
\begin{equation}\label{eq:sc-relative-tangent}
 T_{\mathfrak M/C}\simeq\mathcal A[1],\qquad
 L_{\mathfrak M/C}\simeq\mathcal A^\vee[-1],\qquad
 T_{\mathfrak M/C}\xrightarrow{\sim}L_{\mathfrak M/C}[-1].
\end{equation}
All maps commute with derived base change. Before trivializing
$\lambda$, the last equivalence is
$T_{\mathfrak M/C}\simeq L_{\mathfrak M/C}[-1]\otimes p^*\lambda^\vee$.
\end{theorem}
\begin{proof}(Sketch of proof of \cite[Theorem~3.7]{KS}).
The universal Atiyah class identifies the tangent complex with
$\mathcal A[1]$. Evaluation of the universal $2$-shifted symplectic
form on $\Perf$ followed by integration against the relative
threefold trace produces the closed form, by \cite{PTVV13}.
Its underlying pairing is
\[
 \mathcal A[1]\otimes\mathcal A[1]
 \xrightarrow{\mathrm{composition},\,\tr}Rq_*\O[2]
 \longrightarrow\O_{\mathfrak M}[-1].
\]
Relative Serre duality gives $\mathcal A^\vee\simeq\mathcal A[3]$,
so its adjoint is precisely the equivalence in
\eqref{eq:sc-relative-tangent}. Integration supplies closedness,
not merely a pointwise perfect pairing. Without a frame of
$\lambda$, the same duality is
$\mathcal A^\vee\simeq\mathcal A[3]\otimes p^*\lambda$.
Each arrow is evaluation, composition, trace, or proper perfect
pushforward, and hence is compatible with the stated derived base
changes. This also treats the Gorenstein central fiber.
\end{proof}

The absolute tangent of $\mathfrak M$ fits into
$T_{\mathfrak M/C}\to T_{\mathfrak M}\to p^*T_C$.
Adding this base direction does not change the degree of the
relative form to $-2$. The fourfold moduli problem considered
below is a different functor, with endomorphisms computed on
$\cX$ rather than on $X_b$.

\subsection{Lagrangian restriction to the common divisor}
\label{ssec:sc-lagrangian}
Write
$\mathfrak X_i=\mathfrak{Perf}(X_i)$ and
$\mathfrak S=\mathfrak{Perf}(S)$. The two derived restriction maps
are defined on every test scheme by
\begin{equation}\label{eq:sc-restriction-maps}
 r_i:\mathfrak X_i\longrightarrow\mathfrak S,\qquad
 E_T\longmapsto L(\iota_i\times1_T)^*E_T.
\end{equation}
Derived pullback preserves perfectness, whether or not a cohomology
sheaf meets $S$ transversely. The orientation of $S$ gives
$\mathfrak S$ its $0$-shifted symplectic form.
An isotropic structure on a morphism $r:\mathfrak U\to\mathfrak S$
is a path from $0$ to $r^*\omega_{\mathfrak S}$ in the space of
closed two-forms. It is Lagrangian when the induced map
$T_r\to L_{\mathfrak U}[-1]$ is an equivalence. Thus the assertion
includes a closed-form homotopy as well as a relative duality.

\begin{theorem}[{\cite[Theorem~2.26 and Corollary~2.27]{KS}}]\label{thm:sc-lagrangian-restriction}
For each smooth projective anticanonical pair $(X_i,S)$ with the
chosen logarithmic volume form, \eqref{eq:sc-restriction-maps} is
Lagrangian. At $E_i$ its relative tangent and its Lagrangian map are
\begin{equation}\label{eq:sc-restriction-tangent}
 T_{r_i,E_i}=\RH_{X_i}(E_i,E_i(-S))[1]
       \xrightarrow{\sim}L_{\mathfrak X_i,E_i}[-1].
\end{equation}
The result restricts to any specified open derived substack of the
source. The target remains the full stack of restricted perfect
complexes. With the common-target residue convention, the derived
intersection
\begin{equation}\label{eq:sc-matching-stack}
 \mathfrak G=\mathfrak X_1\times_{\mathfrak S}^{\mathrm h}\mathfrak X_2
\end{equation}
is $(-1)$-shifted symplectic.
\end{theorem}
\begin{proof}[Sketch of proof of \cite{KS}, Theorem~2.26 and Corollary~2.27]
Let $\mathcal E_i$ be universal on $X_i\times\mathfrak X_i$.
Write $p_i$ for projection, $\mathcal E_{S,i}$ for its restriction,
and $p_{S,i}$ for the boundary projection. Set
\[
 \begin{aligned}
 \mathcal A_i&=Rp_{i*}\cRHom(\mathcal E_i,\mathcal E_i),\\
 \mathcal A_i^-&=Rp_{i*}\cRHom(\mathcal E_i,\mathcal E_i(-S)),\\
 \mathcal D_i&=Rp_{S,i*}\cRHom(\mathcal E_{S,i},\mathcal E_{S,i}).
 \end{aligned}
\]
Tensoring the divisor sequence with the universal endomorphism
complex and applying pushforward gives
\begin{equation}\label{eq:sc-universal-restriction}
 \mathcal A_i^-\longrightarrow\mathcal A_i\longrightarrow
 \mathcal D_i\xrightarrow{+1}.
\end{equation}
The logarithmic residue sequence identifies the trace on $S$ with
the composite of its boundary map and the trace on $X_i$.
The consecutive maps in the divisor triangle have a specified
nullhomotopy. Integration sends it to the isotropic path of closed
forms, as in \cite{CalaqueBoundary}. Its induced relative pairing is
\[
 \mathcal A_i^-[1]\otimes\mathcal A_i[1]
    \longrightarrow Rp_{i*}\omega_{X_i}[2]
    \longrightarrow\O_{\mathfrak X_i}[-1].
\]
Because $\O_{X_i}(-S)=\omega_{X_i}$, Serre duality identifies
$\mathcal A_i^-[1]$ with $\mathcal A_i^\vee[-2]$.
The compatibility of the residue boundary with the trace identifies
this equivalence with the actual map of the isotropic path.
The universal diagrams proving that identification are in
\cite[Section~2.5.2--2.5.4]{KS}.

Pulling the two isotropic paths to \eqref{eq:sc-matching-stack}
and composing them using its universal matching equivalence gives
a loop of closed degree-zero forms, hence a closed degree-$-1$
form. Differentiating the fiber product and using
\eqref{eq:sc-universal-restriction} proves nondegeneracy.
This is the derived-intersection construction of \cite{PTVV13}.
The residue sign change is made in the second isotropic path, not
in the restriction functor.
\end{proof}

\subsection{Matching triples, mapping complexes, and effectivity}
\label{ssec:sc-effectivity}
An object over $T$ of \eqref{eq:sc-matching-stack} is
\[
 (E_{1,T},E_{2,T},\alpha_T),\qquad
 \alpha_T:L(\iota_1\times1)^*E_{1,T}
                  \xrightarrow{\sim}L(\iota_2\times1)^*E_{2,T}.
\]
A morphism consists of component morphisms together with a homotopy
between their restrictions; higher homotopies are part of the datum.
The commutative diagram of stacks is
\[
 \begin{tikzcd}[column sep=large,row sep=large]
 \mathfrak G\ar[r,"p_2"]\ar[d,"p_1"']&\mathfrak X_2\ar[d,"r_2"]\\
 \mathfrak X_1\ar[r,"r_1"']&\mathfrak S .
 \end{tikzcd}
\]
With $E_S$ denoting the identified restriction, its tangent is
\begin{equation}\label{eq:recalled-gluing-tangent}
 T_{\mathfrak G,(E_1,E_2,\alpha)}
 \simeq\fib\!\left(
 \bigoplus_{i=1}^2\RH_{X_i}(E_i,E_i)
       \xrightarrow{r_1-r_2}\RH_S(E_S,E_S)\right)[1].
\end{equation}
Conjugation by $\alpha$ is understood in the second restriction.
The kernel and boundary terms of this triangle record both
component deformations and deformations of the matching equivalence.

There is a canonical restriction morphism
\[
 \rho:\mathfrak{Perf}(X_0)\longrightarrow\mathfrak G.
\]
A triple is \emph{effective} when it is induced by a perfect complex
on $X_0$, in families as well as at a point.
For perfect $E,F$ on $X_0$, the normalization sequence gives
\begin{equation}\label{eq:sc-mapping-triangle}
 \RH_{X_0}(E,F)\longrightarrow
 \bigoplus_i\RH_{X_i}(Lj_i^*E,Lj_i^*F)
 \longrightarrow\RH_S(Lj_S^*E,Lj_S^*F)\xrightarrow{+1}.
\end{equation}
This is \cite[Proposition~3.3]{KS}. To recover the proof, tensor
$\O_{X_0}\to j_{1*}\O_{X_1}\oplus j_{2*}\O_{X_2}
\to j_{S*}\O_S$ with $E^\vee\otimes^{\mathbf L}F$.
Projection formula and global sections give the triangle.
The homotopy fiber of its last map is exactly the mapping complex
of triples. Thus restriction is fully faithful, and its map on
endomorphism dg algebras gives an equivalence of the corresponding
formal deformation problems at an effective object. This argument
commutes with derived base change.

We use the following sufficient effectivity test from
\cite[Proposition~3.4]{KS}. If locally on $X_0$ the component
objects have bounded vector-bundle representatives $P_i^\bullet$
and $\alpha$ is represented by a degreewise chain isomorphism over
$S$, then the triple descends to a perfect complex.
At a boundary point put $R=A_1\times_D A_2$, with $A_i\twoheadrightarrow D$.
After local changes of bases, the matching matrices become identities;
the equalizers
\[
 P^a=\{(v_1,v_2):\bar v_1=\bar v_2\}\simeq R^{r_a}
\]
are free and the component differential matrices define a differential
on $P^\bullet$. Full faithfulness in \eqref{eq:sc-mapping-triangle}
lifts overlap equivalences and their homotopies, so descent glues
these local complexes. This is a construction of the descended
object, not a properness theorem for a selected moduli locus.

In a DT comparison one also fixes compatibility of the closed forms
and their orientations along $\rho$. At the level of pairings it is
expressed by
\begin{equation}\label{eq:sc-symplectic-comparison}
 \begin{tikzcd}[column sep=large,row sep=large]
 T_{\mathfrak{Perf}(X_0)}\ar[r,"d\rho","\sim"']\ar[d,"\omega_0^\flat"']&
 \rho^*T_{\mathfrak G}\ar[d,"\rho^*\omega_{\mathfrak G}^\flat"]\\
 L_{\mathfrak{Perf}(X_0)}[-1]\ar[r,"\sim"']&\rho^*L_{\mathfrak G}[-1].
 \end{tikzcd}
\end{equation}
The closed-form homotopy and the isomorphism of orientation line
bundles are retained in addition to this diagram; see
\cite[Remark~3.6]{KS}. Equality of tangent dimensions would not
supply either of these comparisons.

\subsection{From the derived stack to an oriented DT realization}
\label{ssec:sc-oriented-moduli}

A derived scheme is quasi-smooth when its cotangent complex is perfect
of Tor amplitude $[-1,0]$.  In a critical presentation this is the
complex $[T_U|_R\xrightarrow{\operatorname{Hess}f}\Omega_U|_R]$
in cotangent degrees $-1,0$; its dual is the tangent complex in
degrees $0,1$.  Its classical truncation $R$ may be singular or
nonreduced.  A perfect obstruction theory on a classical space,
by itself, is not a specified shifted-symplectic derived enhancement
or a choice of orientation.

The unreduced stack $\RPerf(X_b)$ is an Artin stack, not such a
two-term derived scheme: at a simple object its tangent cohomology
contains $\operatorname{Ext}^0(E,E)$ in degree $-1$ and
$\operatorname{Ext}^3(E,E)$ in degree $2$, in addition to
$\operatorname{Ext}^{1,2}$.  Classical scalar rigidification removes
classical scalar automorphisms, not their dual higher derived
term automatically.  The reduced cyclic critical realization used
here is the separately specified one in Proposition~\ref{prop:general-cyclic-reduction}.
No smoothness assumption on its classical critical scheme is made
in the local connection and bar theorems.

For a chosen finite-type classical moduli space $M$, a critical
chart is a closed embedding $R\hookrightarrow U$ into a smooth
scheme, together with a regular function $f:U\to\A^1$, such that
$R=\Crit(f)$ and $f|_{R^{\red}}=0$. Compatible classes of the
functions modulo the squares of the chart ideals give a d-critical
structure $s$. The Darboux and d-critical constructions
\cite{BravBussiJoyceDarboux,BenBassatBravBussiJoyce2015ShiftedDarboux,Joyce}
provide such structures from the appropriate $(-1)$-shifted
symplectic derived enhancements. For a derived Artin stack the
construction is stack-theoretic. When passing to a classical
algebraic space of simple objects, the compatible reduced critical
realization must be specified as well.

The virtual canonical line bundle $K_{M,s}$ is locally identified
with $K_U^{\otimes2}|_{R^{\red}}$. An orientation is a line bundle
$\mathcal L^{\mathrm{or}}$ with an isomorphism
$(\mathcal L^{\mathrm{or}})^{\otimes2}\simeq K_{M,s}$, satisfying
the overlap compatibility of the critical charts. The oriented
vanishing-cycle construction \cite{BBDJS} then defines
$\mathcal P_{M,s}$, locally the perverse vanishing cycles of
$\CC_U[\dim U]$ with the orientation comparison. We set
\begin{equation}\label{eq:sc-coh-DT}
 \mathrm{DT}_{\mathrm{coh}}^a(M,s)
       =\mathbb H_c^a(M^{\an},\mathcal P_{M,s}).
\end{equation}
For proper $M$ compact and ordinary hypercohomology agree.
For a smooth $d$-dimensional $M$ with zero critical structure,
$K_{M,0}=K_M^{\otimes2}$; choosing $K_M$ as the orientation gives
$\mathcal P_{M,0}=\CC_M[d]$. No square root of the ordinary
canonical line bundle is needed in this example.

There are two reductions to keep separate at rank zero.
For a perfect object $E$, the composite of scalar unit and trace
\[
 \O\longrightarrow\cRHom(E,E)\xrightarrow{\tr}\O
\]
is multiplication by the virtual rank. For a fixed line bundle
$\Lambda$, the determinant moduli problem consists of pairs
$(E,\alpha:\det E\xrightarrow{\sim}\Lambda)$; its tangent is
\[
 \fib\bigl(\RH(E,E)[1]\xrightarrow{\tr[1]}R\Gamma(\O)[1]\bigr).
\]
If $\rk E=0$, a scalar automorphism acts trivially on $\det E$
and preserves $\alpha$. Hence this condition does not remove the
scalar automorphism. These are \cite[Proposition~2.20 and Section~2.4.3]{KS}.
On a Calabi--Yau threefold scalar $\operatorname{Ext}^0$ is paired with $\operatorname{Ext}^3$;
removing only one of them is not a symmetric derived reduction.
Our finite-type DT statements therefore specify the classical
moduli problem, its critical realization, and its orientation.
They apply to perfect objects of arbitrary amplitude satisfying
these conditions, not just to the examples presented by two line
bundles.

\subsection{The two-sided bar complexes and their differentials}
\label{ssec:sc-two-bars}
Take compatible affine derived charts $\Spec A_i\to\Spec A$ of
\eqref{eq:sc-matching-stack}. Its structure algebra is
$A_1\otimes_A^{\mathbf L}A_2$. After flat or cofibrant replacements,
the unnormalized two-sided bar has simplicial degree $r$ term
\begin{equation}\label{eq:sc-structure-bar}
 \Bar_r(A_1,A,A_2)=A_1\otimes_\CC A^{\otimes_\CC r}\otimes_\CC A_2.
\end{equation}
The outside face maps are the module actions, the inside faces
multiply adjacent factors, and the degeneracies insert the unit.
Normalization and the internal dg differentials give a total
complex; simplicial degree $r$ has cohomological degree $-r$.
This recalls \cite[Section~2.1.3]{KS} and fixes the structure-algebra
model used below.

A critical chart has a different, DT realization of that geometry.
On smooth coordinates $x_1,\ldots,x_d$ its Koszul algebra is
\[
 K_f=\O_U[\xi_1,\ldots,\xi_d],\quad |\xi_i|=-1,\qquad
 \delta_f=\sum_i (\partial_i f)\partial_{\xi_i}.
\]
The oriented periodic model has, in these coordinates, mixed
operator $-\delta_f+u\Delta$ with
$\Delta=\sum_i\partial_{x_i}\partial_{\xi_i}$ and the critical-chart
suspension. Tensoring with the orientation coefficient and
contracting polyvectors with top forms compares it with
\begin{equation}\label{eq:sc-twisted-de-rham}
 (\Omega_U^\bullet((u)),u\,d-df\wedge)[d].
\end{equation}
The signs and coordinate changes are fixed by the explicit
construction in Section~\ref{sec:v5-bar-bridge}. In particular,
$\mathcal B_{0,\mathrm{DT}}$ always denotes the oriented bar/BV
realization with this extra differential, not
\eqref{eq:sc-structure-bar} with only its structure differential.

The categorical input on a chart is $\MF(U,f)$ for an ordinary
degree-zero potential. Its periodic cyclic realization is compared
with twisted de Rham cohomology by \cite{Efimov2018CyclicMF}.
A collection of charts is used globally only with its specified
category comparisons and coherent descent. The perverse-sheaf
construction alone is not a proof that an arbitrary collection of
matrix-factorization categories has that descent. Likewise a
quasi-isomorphism of periodic complexes is a statement about a
realization, not an equivalence of the underlying dg categories.

\subsection{Fourfold duality and the normal deformation triangle}
\label{sec:CY}
For a finite subset $Z\subset C$ let $\Perf_{\pi^{-1}Z}(\cX)$
consist of perfect complexes whose cohomology sheaves have
set-theoretic support in $\pi^{-1}Z$. These objects may live on
infinitesimal thickenings of the fibers. Their endomorphisms are
computed on the smooth fourfold, not by
\eqref{eq:sc-relative-functor}.

The following is \cite[Proposition~3.9 and Corollary~3.11]{KS}.
If $\omega_{\cX}\simeq\pi^*B$ and a frame of $B$ is fixed on
an open neighborhood of $Z$, the category of these objects has
functorial perfect pairings
\begin{equation}\label{eq:sc-fourfold-pairing}
 \RH_{\cX}(E,F)\otimes\RH_{\cX}(F,E)\longrightarrow\CC[-4].
\end{equation}
The corresponding moduli of perfect complexes with proper support
has a closed $(-2)$-shifted symplectic form on its locally
geometric finite-presentation loci. Indeed proper-support Serre
duality first gives
$\RH(E,F)^\vee\simeq\RH(F,E\otimes\omega_{\cX})[4]$.
The chosen frame acts on a whole neighborhood of the support and
therefore on every infinitesimal thickening. It trivializes the
canonical twist in this formula. The fourfold volume form has a
closed negative-cyclic representative; its moduli transgression
supplies the closed form, as in \cite{BD2}. Shifting the two
endomorphism factors by $[1]$ changes the target of the pairing
to degree $-2$.

For the comparison with threefold deformations, let
$i:X\hookrightarrow Y$ be an effective Cartier divisor with normal
line bundle $N$, and let $E\in\Perf(X)$. The natural exact triangle
is
\begin{equation}\label{eq:normal-deformation}
 \begin{split}
 \RH_X(E,E)&\longrightarrow\RH_Y(i_*E,i_*E)\\
 &\longrightarrow\RH_X(E,E\otimes N)[-1]
          \longrightarrow\RH_X(E,E)[1].
 \end{split}
\end{equation}
We recall the map, not just the three terms, from
\cite[Proposition~4.19]{KS}. The self-intersection kernel has a
filtration giving
\[
 E\otimes N^\vee[1]\longrightarrow Li^*i_*E\longrightarrow E
       \xrightarrow{\alpha_E}E\otimes N^\vee[2].
\]
Apply $\RH_X(-,E)$ and derived adjunction. The extension class is
\begin{equation}\label{eq:sc-atiyah-conormal}
 \alpha_E:
 E\xrightarrow{\At(E)}E\otimes L_X[1]
   \xrightarrow{1\otimes\kappa_i[1]}E\otimes N^\vee[2],
 \qquad \kappa_i:L_X\to N^\vee[1].
\end{equation}
For a sheaf, the corresponding part of the long exact sequence is
\[
 0\longrightarrow\operatorname{Ext}_X^1(E,E)\longrightarrow
 \operatorname{Ext}_Y^1(i_*E,i_*E)\longrightarrow\Hom_X(E,E\otimes N)
       \xrightarrow{\partial_E}\operatorname{Ext}_X^2(E,E).
\]
Here $\partial_E(a)=(a\otimes1_{N^\vee})[2]\circ\alpha_E$,
with the connecting-arrow convention above. In a fiber embedding,
$a=\id_E\otimes v$ gives $\partial_E(a)=\At(E)\cdot\kappa(v)$
for $v\in T_{C,b}$. This is the Atiyah--Kodaira--Spencer obstruction
\cite{HT}. The full triangle, including negative self-extensions
when present, is the formula for a general perfect object.

In families, apply the same construction to the universal perfect
complex. The relative Atiyah class in the moduli direction is
\[
 \mathcal E\longrightarrow\mathcal E\otimes q^*L_{\mathfrak M/C}[1].
\]
Adjunction identifies its induced morphism
$(Rq_*\cRHom(\mathcal E,\mathcal E))^\vee[-1]
\to L_{\mathfrak M/C}$ with the cotangent identification of
\eqref{eq:sc-relative-tangent}. Thus the deformation triangle,
trace pairing, and moduli differential refer to the same universal
object; the full diagrams are \cite[Appendix~A]{KS}.

\subsection{The fourfold foliation used in the potential construction}
\label{ssec:sc-foliation}
The fourfold foliation has a specific presentation hypothesis.
Let $Z$ be a smooth Calabi--Yau fourfold with chosen volume form,
contained in a smooth projective compactification $\overline Z$.
Fix a polarization, a Hilbert polynomial, and a bounded stable
sheaf locus on $\overline Z$ whose supports avoid
$\overline Z\setminus Z$. Choose a uniform regularity integer
$r$, a framing $\CC^N\otimes\O(-r)\twoheadrightarrow F$, and
a finite-type derived Quot presentation with the invariant product
charts and minimal Artin reductions of
\cite[Hypothesis~4.7]{KS}. Its quotient by $\GL_N$ presents that
sheaf locus; it is the quotient, rather than the framed atlas,
which carries the fourfold shifted form.

For these presentation data,
\cite[Theorem~4.11 and Corollary~4.12]{KS} give an invariant purely
derived real foliation, Lagrangian for the imaginary part of the
form and negative definite for its real part in degree one.
The argument first identifies the Quot quotient pairing with
\eqref{eq:sc-fourfold-pairing}, then passes through the invariant
Artin charts and applies the strictification and softness results
of \cite{BKSY2024,BKSY2025}. Its overlap maps compare representatives
of the same distribution and transport their isotropic homotopies;
integrability is not obtained by averaging brackets.
These hypotheses concern this sheaf presentation. The general
perfect-complex comparison below uses critical realizations and
does not assume that every perfect complex is a point of a
coherent-sheaf Quot scheme.

The local output needed in the present paper is a vector bundle
$V$ on a smooth scheme or real smooth chart $T$, a section $s$
and a cosection $\tau$ satisfying $\tau(s)=0$.
The foliation quotient has Koszul algebra
\begin{equation}\label{eq:recalled-shifted-potential}
 B=(R[\theta_1,\ldots,\theta_r],\delta),\quad
 |\theta_a|=-1,\quad\delta\theta_a=s_a,\quad
 \nu=\sum_a\tau_a\theta_a,\quad\delta\nu=\sum_a\tau_as_a=0.
\end{equation}
The shifted critical reconstruction from
\cite[Proposition~4.17]{KS} adjoins $\eta_a$ of degree $-1$ and
$p_i$ of degree $-2$, with
\begin{equation}\label{eq:sc-fourfold-critical}
 \delta\eta_a=\tau_a,\qquad
 \delta p_i=-\sum_a\bigl((\partial_i s_a)\eta_a+
                                  (\partial_i\tau_a)\theta_a\bigr).
\end{equation}
The equation $\delta^2p_i=-\partial_i(\tau(s))=0$ explains the
closed-cosection condition. Cotangent coordinate changes glue
these formulas and the canonical degree-$-2$ form.
The function $\nu$ has degree $-1$ in a dg algebra. The ordinary
potential $f$ used in matrix factorizations has degree zero on a
smooth ambient scheme. Section~\ref{sec:local} computes their
possible comparison from the actual normal deformation operations.

\subsection{Proper supports and an additional anticanonical boundary}
\label{ssec:sc-open-boundary}
The second geometric example will use a nonprojective Calabi--Yau
open set inside a projective pair. We record a form of the restriction
construction which specifies exactly the allowed objects.

\begin{proposition}\label{prop:sc-two-boundaries}
Let $V$ be a smooth projective threefold, with disjoint smooth
(possibly disconnected) divisors $D$ and $B$ and a chosen
isomorphism $\omega_V(D+B)\simeq\O_V$.
The restriction of $\mathfrak{Perf}(V)$ to the open substack of
objects whose cohomology supports avoid $B$ has a Lagrangian map
to $\mathfrak{Perf}(D)$. The relative tangent at such an object is
$\RH_V(E,E(-D))[1]\simeq L_{\mathfrak M,E}[-1]$.
The same conclusion holds on finite-type loci of perfect objects
with proper support in $V\setminus B$, by extension by zero.
\end{proposition}
\begin{proof}
The union $D\amalg B$ is a smooth anticanonical divisor in $V$.
Theorem~\ref{thm:sc-lagrangian-restriction} gives a Lagrangian map
\[
 \mathfrak{Perf}(V)\longrightarrow
          \mathfrak{Perf}(D)\times\mathfrak{Perf}(B).
\]
The zero object is an open substack of $\mathfrak{Perf}(B)$:
for a perfect family over a test scheme, the locus on which all
cohomology vanishes on the proper $B$ is open. Its tangent complex
is zero. Restricting the target to
$\mathfrak{Perf}(D)\times\{0\}$ and pulling back the source therefore
preserves the Lagrangian structure. The source condition is
$L i_B^*E=0$. For a perfect complex this is equivalent to its
cohomology supports avoiding $B$, by derived Nakayama at every
point of $B$.

On a neighborhood of such a support the canonical section of
$\O_V(B)$ is invertible. Hence
$\RH_V(E,E(-D))=\RH_V(E,E(-D-B))$, and the trace proof above gives
the asserted tangent equivalence. A perfect family with support
closed and proper over its parameter in $V\setminus B$ extends
by zero to a perfect family on $V$: it is unchanged near its
support and zero near $B$. This identifies the indicated moduli
problems and their maps.
\end{proof}

This proposition will allow the ruled component in
Section~\ref{sec:second-degeneration} to contribute with its correct
extra boundary. It also illustrates why specifying support
conditions on perfect complexes is part of the geometric
construction, rather than a consequence of a numerical charge.

\section{Normal operations and ordinary potentials}\label{sec:local}
The geometric fourfold foliation and its shifted critical algebra are
inputs from \cite{KS}. Our task here is different: compute the normal
operations of a given object, identify the coefficients that an
ordinary potential would have to use, and test a proposed product
factorization, as explicitly as possible. These calculations specify when the shifted presentation
can supply the ordinary function used by the cyclic connection.

\subsection{The formal deformation algebra of a perfect object}
\label{ssec:sc-formal-moduli}
Let $E\in\Perf(X_b)$ be any of the objects under consideration.
Choose an enhanced endomorphism dg algebra representing
$\RH_{X_b}(E,E)$ and give it the commutator dg Lie bracket.
Write this Lie algebra as $\mathfrak g_E$.
For a local Artinian connective dg algebra $A$ with residue
field $\CC$ and maximal ideal $\mathfrak m_A$, the formal
deformations with an identification of the residue object are
modeled by
\[
 \operatorname{Def}_E(A)=
       \MC_\bullet(\mathfrak g_E\otimes\mathfrak m_A).
\]
This is the characteristic-zero formal-moduli description used
in \cite[Section~4.5.1]{KS}; cyclic homotopy transfer is treated
in \cite{KontsevichSoibelman2006}. A local change of differential
$d_E+a$ satisfies $(d_E+a)^2=0$ exactly when
$da+\tfrac12[a,a]=0$. Gauge transformations and their higher
homotopies account for changes of representative and for the
isomorphisms of deformations. Using an enhanced Cech resolution
gives the global endomorphism algebra, rather than only the
endomorphisms of one local vector-bundle complex.

For a minimal $L_\infty$ model $L_E$ its complete algebra of formal
functions is
\begin{equation}\label{eq:sc-full-CE}
 \widehat{\mathrm{CE}}(L_E)=
       (\widehat\Sym(L_E^\vee[-1]),d_{\mathrm{CE}}).
\end{equation}
The coefficients of the differential are all its transferred
brackets $\ell_r$. A class in $L_E^i=\operatorname{Ext}^i(E,E)$ gives a dual
function generator of degree $1-i$. In particular a surviving
$\operatorname{Ext}^3$ gives degree-$-2$ functions, while a scalar $\operatorname{Ext}^0$
records an automorphism direction. This explains why the full
Artin deformation problem and a two-term reduced critical
presentation are different objects even for a coherent sheaf.
For a simple universally gluable perfect object on a
Calabi--Yau threefold, duality makes $\operatorname{Ext}^{>3}=0$ and pairs
$\operatorname{Ext}^1$ with $\operatorname{Ext}^2$ and $\operatorname{Ext}^0$ with $\operatorname{Ext}^3$.
These conclusions depend on self-extensions, not on the
amplitude of the cohomology sheaves of $E$.

When a reduced cyclic deformation problem is used below, its
reduction and its comparison with the selected geometric
critical charts are specified. On the smooth loci of our
examples the critical structure is zero, so this comparison
can be checked directly. For singular loci the local ordinary
potential retains the full transferred operations on the
reduced model. Completion in \eqref{eq:sc-full-CE} is essential:
a formal potential is not assumed to be an algebraic or
convergent function without an additional realization.

\subsection{A cyclic reduction for simple perfect objects}
\label{ssec:sc-cyclic-reduction}
The following local calculation explains when the two-term
\emph{deformation algebra} can be obtained from the full object
algebra. The cyclic model is specified in the statement; its comparison
with algebraic critical charts is a further geometric step.

\begin{proposition}\label{prop:general-cyclic-reduction}
Let $E$ be a simple universally gluable perfect object on an
oriented proper Calabi--Yau threefold. Suppose its enhanced
endomorphism algebra has a finite-dimensional, strictly unital,
minimal uncurved cyclic $A_\infty$ model $A$ on
$A^i=\operatorname{Ext}^i(E,E)$, whose cyclic pairing is the
threefold trace pairing. Let $L$ be its associated minimal
commutator $L_\infty$ algebra. Then the only nonzero higher
brackets of $L$ are
\[
 \ell_r:\operatorname{Sym}^r A^1\longrightarrow A^2,
 \qquad r\ge2,
\]
with the symmetric convention on degree-one inputs.
The summands $A^0=\CC\,1$ and $A^3=\CC\,\epsilon$, where
$\langle1,\epsilon\rangle=1$, are central $L_\infty$ summands.
The subspace $L_{\mathrm{red}}=A^1\oplus A^2$ is a cyclic
$L_\infty$ algebra, with $A^2\simeq(A^1)^\vee$.
On ordinary local Artinian bases, the scalar-rigidified formal
deformation functor of $E$ is the critical functor of the
formal series
\begin{equation}\label{eq:sc-general-potential}
 f_E(x)=\sum_{r\ge2}\frac1{(r+1)!}
              \langle\ell_r(x,\ldots,x),x\rangle.
\end{equation}
This statement applies equally to rank-zero objects.
\end{proposition}
\begin{proof}
Universal gluability and threefold duality give $A^i=0$ for
$i<0$ and $i>3$, and simplicity identifies $A^0$ with the
strict unit. Duality then identifies $A^3$ with its dual and
pairs $A^1$ perfectly with $A^2$. Minimality makes $\ell_1=0$.
Strict unitality implies that every higher commutator bracket
with an input $1$ vanishes: for $r=2$ this is the commutator
with the identity, and for $r>2$ every contributing $A_\infty$
operation with a unit vanishes.

For positive-degree inputs, $\ell_r$ has degree $2-r$.
If all inputs have degree one, its output has degree two.
If exactly one input has degree two, the output has degree
three; cyclicity pairs that output with $1$ by a bracket
having $1$ as an input. This is zero by the preceding paragraph.
The pairing between $A^3$ and $A^0$ is nondegenerate, so that
output is zero. Two degree-two inputs, or any degree-three
input with the remaining inputs positive, would give degree at
least four, and hence also vanish. Together with the unit
calculation these exhaust all inputs. They prove the asserted
bracket description and the centrality of $A^0$ and $A^3$.
In particular $L_{\mathrm{red}}$ inherits the cyclic structure.

For an ordinary Artinian base with maximal ideal $\mathfrak m$,
a Maurer--Cartan element lies in $A^1\otimes\mathfrak m$ and
satisfies
\[
 \sum_{r\ge2}\frac1{r!}\ell_r(x,\ldots,x)=0
                     \quad\text{in }A^2\otimes\mathfrak m.
\]
These sums are finite on that base. Cyclicity makes the
$(r+1)$ differentiated terms of the $r$th summand of
\eqref{eq:sc-general-potential} equal, and therefore
\[
 (df_E)_x(v)=\sum_{r\ge2}\frac1{r!}
                   \langle\ell_r(x,\ldots,x),v\rangle.
\]
The perfect pairing identifies this critical equation with the
Maurer--Cartan equation. The degree-zero gauge algebra is the
central scalar algebra and acts trivially; scalar rigidification
removes its automorphisms. This proves the claim about the
classical formal functor. The reduced complete
Chevalley--Eilenberg algebra is the corresponding Koszul algebra
of $df_E$, as in \cite[Proposition~4.22]{KS}.

For a rank-zero object the degree-zero trace of $1$ may vanish.
The reduction above uses instead the nondegenerate pairing
$A^0\otimes A^3\to\CC$ and the strict unit. It does not split
the degree-zero scalar--trace composite by division by the rank.
Hence the same argument applies at rank zero.
\end{proof}

The proposition constructs a local reduced formal model from a
specified full cyclic model. For dg Artinian bases, the central
degree-zero and degree-three summands still contribute to the
full derived deformation functor. They are not discarded in
$\mathfrak{Perf}(X_b)$ or in its derived restriction maps.
An algebraic realization of \eqref{eq:sc-general-potential},
its orientation, and coherent comparisons on overlaps are the
critical-chart data used in the following sections. In particular,
the proposition proves neither convergence of the formal series
nor automatic global descent of a matrix-factorization category.

\subsection{The reduced cyclic input}
\label{ssec: Local models and shifted symplectic structures}
We now specify the reduction which permits an ordinary potential.
The preceding complete deformation algebra retains all tangent degrees.
A degree-zero critical presentation uses a reduced cyclic model;
its identification with the chosen classical moduli problem is part
of the input.
We use the following input from \cite[Hypothesis~4.21]{KS}.
\begin{hypothesis}\label{hyp:cyclic-model}
The chosen reduced formal deformation problem has a minimal cyclic
$L_\infty$ model with $L^1=U$, $L^2=U^\vee$, and no other terms.
In the suspended convention its cyclic tensors
$\langle\ell_r(x_1,\ldots,x_r),x_0\rangle$ are symmetric.
The reduction and its pairing are specified as part of the data.
\end{hypothesis}
Under this hypothesis, \cite[Proposition~4.22]{KS} identifies the
ordinary formal potential and the critical algebra as
\begin{equation}\label{eq:formal-potential}
 f(x)=\sum_{r\ge2}\frac{\langle\ell_r(x,\ldots,x),x\rangle}{(r+1)!},
 \qquad
 \widehat{\mathrm{CE}}(L)
  =\bigl(\widehat\Sym(U^\vee)\otimes\Lambda^\bullet U,\iota_{df}\bigr).
\end{equation}
The coefficient $1/(r+1)!$ fixes the convention used in the transfer
calculations below. Algebraic or analytic matrix-factorization charts
require the corresponding geometric critical presentation, in
addition to this formal series. The smooth divisor examples in
Section~\ref{sec:v6-hyperplanes} have zero reduced potential; the
verification of their model and orientation is in
Proposition~\ref{prop:v15-three-hypotheses}.

\subsection{The normal complex of the object and the local coordinates}
\label{sssec:normal-coordinates}
Fix a smooth fiber $X=X_b$, a simple perfect object $E$ on $X$ with
$\ext^{<0}_X(E,E)=0$, and a Cartier embedding $i:X\hookrightarrow Z$ in a smooth fourfold.
The two applications of this formula are $Z=\cX$ and
$Z=Y=\P^4$; the corresponding normal line bundles will be
distinguished explicitly. Put
\[
 \begin{aligned}
 A_E&=\RH_X(E,E),& M_E&=\RH_X(E,E\otimes N),\\
 N&=N_{X/Z},& U&=\ext^1_X(E,E).
 \end{aligned}
\]
The chosen threefold volume form identifies $\ext^2_X(E,E)$ with
$U^\vee$. The symbols used for the normal terms are
\begin{equation}\label{eq:r13-normal-spaces}
 W_1=H^0(M_E)=\Hom_X(E,E\otimes N),\qquad
 W_2=H^1(M_E)=\ext^1_X(E,E\otimes N).
\end{equation}
They occur in consecutive degrees of $M_E$, and in degrees one and
two of the shifted complex $M_E[-1]$ in
\eqref{eq:normal-deformation}.
For the fiber embedding, write $\mathfrak M_{\cX}$ for the
moduli stack of perfect complexes with proper support on $\cX$.
The differentiated pushforward gives the diagram of exact triangles
\begin{equation}\label{eq:r13-moduli-tangent-diagram}
\begin{tikzcd}[column sep=small,row sep=large]
 \mathbb T_{\mathfrak X_b,E}\arrow[r,"{d(i_{b*})}"]\arrow[d,"{\sim}"']
 &\mathbb T_{\mathfrak M_{\cX},i_{b*}E}\arrow[r]\arrow[d,"{\sim}"']
 &M_E\arrow[d,equal]\\
 A_E[1]\arrow[r]
 &\RH_{\cX}(i_{b*}E,i_{b*}E)[1]\arrow[r]&M_E .
\end{tikzcd}
\end{equation}
The next term in both rows is $A_E[2]$. Thus the occurrences of
$W_1$ and $W_2$ below are the degree-zero and degree-one cohomology
of the normal term in the actual tangent triangle of the moduli stacks.
 For a sheaf $E$ one has $H^{-1}(M_E)=0$. The same vanishing holds
for a fiber embedding under the assumed negative self-Ext vanishing,
since $N=\O_X\otimes T_{C,b}$. Under this vanishing the portion of
the long exact sequence needed here is
\begin{equation}\label{eq:r13-normal-long-exact}
\begin{aligned}
 0\longrightarrow U&\longrightarrow\ext^1_Z(i_*E,i_*E)
    \longrightarrow W_1\xrightarrow{\partial_{E,0}}U^\vee,\\
 U^\vee&\longrightarrow\ext^2_Z(i_*E,i_*E)
    \longrightarrow W_2\xrightarrow{\partial_{E,1}}\ext^3_X(E,E).
\end{aligned}
\end{equation}
The second row continues the first row. The maps $\partial_{E,j}$
are induced by the same extension class $\alpha_E$ in
\cite[equation~(4.26)]{KS}. Thus the normal quotient of the first
order deformation space is $\ker\partial_{E,0}$, not an arbitrarily
chosen copy of $W_1$.

If $\alpha_E=0$, choose splittings of these sequences. They give
\begin{equation}\label{eq:r13-cohom-splittings}
 \ext^1_Z(i_*E,i_*E)\simeq U\oplus W_1,\qquad
 \ext^2_Z(i_*E,i_*E)\simeq U^\vee\oplus W_2.
\end{equation}
A basis $u_i$ of $U$ and a basis $w_\alpha$ of $W_1$ then determine
first order parameters $x=\sum_i x_i u_i$ and
$a=\sum_\alpha a_\alpha w_\alpha$. They also determine degree-zero
function generators $x_i,a_\alpha$ in a formal model. The generators
dual to $U^\vee\oplus W_2$ have degree $-1$. Generators dual to
$\ext^3_Z(i_*E,i_*E)$ have degree $-2$ and remain in the fourfold
complex. More precisely, the degrees of the dual function spaces are
\[
\begin{array}{c|c}
 \text{cohomology of the deformation algebra}&
                 \text{degree of its dual functions}\\ \hline
 U\oplus W_1&0\\
 U^\vee\oplus W_2&-1\\
 \ext^3_Z(i_*E,i_*E)&-2.
\end{array}
\]
The splittings in \eqref{eq:r13-cohom-splittings} specify the linear
coordinates. The higher terms of the differential are the transferred
$L_\infty$ operations; these must be computed with the same choices.
If $\alpha_E\ne0$, \eqref{eq:r13-normal-long-exact} is used before
selecting the minimal coordinates on the ambient moduli problem.

For the embedding of a fiber, the normal line bundle is
$\O_X\otimes T_{C,b}$. Consequently
\begin{equation}\label{eq:r13-fiber-normal-spaces}
 W_1=\CC\,\id_E\otimes T_{C,b},\qquad
 W_2=U\otimes T_{C,b}.
\end{equation}
Choose a coordinate $t$ on $C$ near $b$, and set
$v=\partial/\partial t|_b$. Then $a\,\id_E\otimes v$ is the
candidate first order displacement of the support in the base
direction. The condition that it lift to a deformation of $i_*E$ is
\begin{equation}\label{eq:r13-base-obstruction}
 \partial_{E,0}(a\,\id_E\otimes v)
       =a\,\At(E)\cdot\kappa(v)=0.
\end{equation}
This is the actual map controlling transverse motion in the family,
as proved in \cite[Proposition~4.19]{KS}.
The identifications \eqref{eq:r13-fiber-normal-spaces} do not identify
$W_1$ with $W_2$: the first has dimension one, while the second has
dimension $\dim U$.

These spaces can be computed in the hyperplane-section application.
Let $E=\O_D$ for $D=X\cap H$ on a smooth quintic. For
$X\hookrightarrow\cX$, Proposition~\ref{prop:Ext-example} gives
\begin{equation}\label{eq:r13-hyperplane-dimensions}
 \dim U=4,\qquad \dim W_1=1,\qquad\dim W_2=4.
\end{equation}
For $X\hookrightarrow\P^4$, instead $N=\O_X(5)$ and
\[
 W_1=H^0(D,\O_D(5)),\qquad W_2=H^0(D,\O_D(6)).
\]
Indeed, applying sheaf Hom to
$[\O_X(-1)\to\O_X]$ with target $\O_D(5)$ gives the complex
$[\O_D(5)\xrightarrow{0}\O_D(6)]$. The intermediate cohomology
of these line bundles on the quintic surface vanishes, and
\[
 \dim W_1=\binom{8}{3}-1=55,\qquad
 \dim W_2=\binom{9}{3}-4=80.
\]
Thus the fixed projective fourfold and the total space of the
fibration give different normal deformation problems for the same
$E$. In both embeddings $E$ extends as the two-term complex
$[L^{-1}\xrightarrow{\ell}\O]$ to a neighborhood of $X$.
Its restriction to the first infinitesimal neighborhood supplies a
lift, so $\alpha_E=0$. The dimensions just computed concern the
splittings of \eqref{eq:r13-normal-long-exact} for this object;
they do not identify the formal neighborhood of every quintic with
its normal line bundle.

\subsection{The mixed operations on the normal complex}
\label{sssec:normal-operations}
To calculate the normal equations, choose a dg model of $A_E$ and
its dg module $M_E$. On resolutions of $E$, the action is induced by
composition of morphisms:
\begin{equation}\label{eq:r13-module-action}
 [z,m]_N=(z\otimes\id_N)\circ m-(-1)^{|z||m|}m\circ z,
 \qquad z\in A_E,\ m\in M_E.
\end{equation}
The same formula is used in a totalized resolution computing global
derived morphisms. If $\xi$ is a Maurer--Cartan element describing
a deformation of $E$ on the fixed $X$, the differential on its
normal morphism complex is
\begin{equation}\label{eq:r13-twisted-module}
 d_{M,\xi}m=d_Mm+[\xi,m]_N,\qquad
 d_{M,\xi}^{\,2}m=
 [d_A\xi+\tfrac12[\xi,\xi],m]_N=0.
\end{equation}
This calculation describes deformations on the fixed variety $X$.
The displacement of $X$ itself is governed by
\eqref{eq:r13-base-obstruction} and the corresponding higher
operations in the ambient deformation algebra.

Here is a chain calculation of the normal Taylor coefficients.
Write $\mathcal W^\bullet=H^\bullet(M_E)$, including all its
degrees, and choose a contraction of complexes of $\CC$-vector spaces
\[
 (\mathcal W^\bullet,0)
 \underset{p_M}{\stackrel{i_M}{\rightleftarrows}}(M_E,d_M),
 \qquad d_Mh_M+h_Md_M=1-i_Mp_M,
\]
with $p_Mi_M=1$ and $p_Mh_M=h_Mi_M=h_M^2=0$.
Over a local Artinian coefficient algebra, put
$T_\xi=[\xi,-]_N$. Its coefficients lie in the maximal ideal.
The following series are therefore finite there, and converge in
the maximal-ideal completion:
\begin{equation}\label{eq:r13-transfer-formulas}
 i_\xi=(1+h_MT_\xi)^{-1}i_M,
 \qquad d_{\mathcal W,\xi}=p_MT_\xi i_\xi.
\end{equation}
They satisfy
\begin{equation}\label{eq:r13-transfer-identities}
 (d_M+T_\xi)i_\xi=i_\xi d_{\mathcal W,\xi},
 \qquad p_Mi_\xi=1,
 \qquad d_{\mathcal W,\xi}^{\,2}=0.
\end{equation}
For completeness, the first identity follows by substituting
$i_\xi=i_M-h_MT_\xi i_\xi$:
\begin{align*}
 (d_M+T_\xi)i_\xi
 &=i_Mp_MT_\xi i_\xi+h_Md_MT_\xi i_\xi\\
 &=i_Md_{\mathcal W,\xi}
       -h_MT_\xi(d_M+T_\xi)i_\xi.
\end{align*}
The second equality uses $(d_M+T_\xi)^2=0$.
Multiplication by $(1+h_MT_\xi)^{-1}$ proves the first identity
in \eqref{eq:r13-transfer-identities}. The second follows from
$p_Mh_M=0$. Applying $d_M+T_\xi$ once more, and then $p_M$,
proves the third. The map $i_\xi$ is a quasi-isomorphism: filter
both complexes by powers of the maximal ideal of the Artinian
coefficient algebra. On each associated graded quotient it is
$i_M$ tensored with that quotient, hence a quasi-isomorphism.
The finite filtration proves the assertion for $i_\xi$.
Thus the part from $W_1=\mathcal W^0$ to
$W_2=\mathcal W^1$ is calculated by
\begin{equation}\label{eq:r13-Phi-transfer}
 \Phi_\xi=
 \left.p_MT_\xi i_M-p_MT_\xi h_MT_\xi i_M
       +p_MT_\xi h_MT_\xi h_MT_\xi i_M-\cdots\right|_{W_1}.
\end{equation}
This expresses the coefficients in terms of the original
composition maps and the chosen contraction.

Equivalently, homotopy transfer gives mixed operations
\[
 \mu_n:\Sym^nU\otimes W_1\longrightarrow W_2,
 \qquad n\geq1,
\]
with the suspended convention of
\cite{KontsevichSoibelman2006}. Their series, on the chosen formal
model, is
\begin{equation}\label{eq:r13-Phi-mixed}
 \Phi(x)a=\sum_{n\geq1}\frac1{n!}\mu_n(x,\ldots,x;a).
\end{equation}
On an Artinian deformation satisfying the Maurer--Cartan equation,
this is the $W_1\to W_2$ component of
\eqref{eq:r13-transfer-formulas}. Values of the series away from
the Maurer--Cartan locus belong to the specified transferred model.
In bases $u_i,w_\alpha,v_j$ of $U,W_1,W_2$, respectively,
\begin{equation}\label{eq:r13-Phi-coefficients}
 \Phi_{j\alpha}(x)
 =\sum_{n\geq1}\frac1{n!}
   \sum_{i_1,\ldots,i_n}
   \bigl(\mu_n(u_{i_1},\ldots,u_{i_n};w_\alpha)\bigr)_j
                      x_{i_1}\cdots x_{i_n}.
\end{equation}
In particular $\Phi(0)=0$ in the minimal normal model. The equality
$\Phi(x)a=0$ means that the class $a$ lifts to a degree-zero cocycle
of the deformed normal morphism complex. Its obstruction lies in
$W_2$. The other degrees of $\mathcal W^\bullet$ and their mixed
operations remain part of the full deformation complex.

There is a useful additional check for a simple object in a fiber.
After the choice of $v$ in \eqref{eq:r13-fiber-normal-spaces},
$M_E=A_E\otimes\CC v$. Choose the contraction with
$i_M(\id_E\otimes v)=\id_E\otimes v$. Then
\[
 T_\xi(\id_E\otimes v)=0,
 \qquad i_\xi(\id_E\otimes v)=\id_E\otimes v,
 \qquad\Phi_\xi(\id_E\otimes v)=0.
\]
The first equality follows directly from
\eqref{eq:r13-module-action}; the others follow from the series.
Thus this scalar normal morphism has no obstruction while $E$
deforms on a fixed $X$. The nonconstant variation of the varieties
$X_b$ enters through the Kodaira--Spencer class and the
$t$-dependence of the cyclic operations, not through an invertible
map $W_1\to W_2$.

\subsection{The cyclic extension of the specified normal equations}
\label{sssec:r13-cyclic-normal}
Choose the normal model and mixed operations just described. The
part of the ambient obstruction equation that is linear in the
normal parameter has the form $\Phi(x)a\in W_2$ when the chosen
comparison with $M_E[-1]$ has that presentation. More generally the
normal component may contain terms of degree at least two in $a$,
and the tangential component contains the connecting map
$\partial_{E,0}$ and its higher analogues. The next proposition
constructs the cyclic extension for the stated equation
$\Phi(x)a=0$; its input is this presentation, rather than the vector
spaces $W_1,W_2$ alone.

For the construction, $\eta\in W_2^\vee$ is an additional variable
of degree zero. Write $a=\sum_\alpha a_\alpha w_\alpha$,
$\eta=\sum_j\eta_jv^j$ in the bases of
\eqref{eq:r13-Phi-coefficients}. Define
\begin{equation}\label{eq:r13-h-normal}
 h_N(x,a,\eta)=\langle\eta,\Phi(x)a\rangle
  =\sum_{n\geq1}\frac1{n!}
                 \langle\eta,\mu_n(x,\ldots,x;a)\rangle.
\end{equation}
Thus $h_N$ is the scalar contraction of the normal mixed operations,
while $f$ is the cyclic primitive of the fiber obstruction map in
\eqref{eq:formal-potential}. Both have degree zero. The coordinates
$x_i,a_\alpha,\eta_j$ are all degree zero. Degree-$-1$ generators
are introduced only when forming their derived critical algebra.
We use $g$ exclusively for the degree-$-1$ factor of a shifted
potential in Section~\ref{sec:graded-factorization}.

\begin{proposition}\label{prop:additive-potential}
Let $U,W_1,W_2$ be finite-dimensional vector spaces, let $f$ be a degree-zero formal function on $U$, and let
\[
 \Phi\in\Hom(W_1,W_2)\otimes\widehat\Sym(U^\vee).
\]
Put $V=U\oplus W_1\oplus W_2^\vee$, with coordinates $(x,a,\eta)$, and define
\begin{equation}\label{eq:additive-W}
 F(x,a,\eta)=f(x)+\langle\eta,\Phi(x)a\rangle.
\end{equation}
The differential of $F$ is
\[
 \begin{aligned}
 d_UF(v)&=df(x)(v)+\langle\eta,(d\Phi)_x(v)a\rangle,\\
 d_{W_1}F&=\Phi(x)^\vee\eta,\\
 d_{W_2^\vee}F&=\Phi(x)a.
 \end{aligned}
\]
Consequently the completed Koszul algebra
\[
 \bigl(\widehat\Sym(V^\vee)\otimes\Lambda^\bullet V,
                       \iota_{dF}\bigr),
 \qquad\deg(\Lambda^1 V)=-1,
\]
represents the derived zero locus of these three groups of equations. It is the formal derived critical locus of the ordinary degree-zero function $F$.
\end{proposition}
\begin{proof}
Choose bases of $W_1,W_2$ and write $\Phi(x)=(\Phi_{ji}(x))$. Then
\[
 F=f+\sum_{i,j}\eta_j\Phi_{ji}(x)a_i.
\]
Differentiation in $\eta_j$ gives $\sum_i\Phi_{ji}(x)a_i$; differentiation in $a_i$ gives $\sum_j\eta_j\Phi_{ji}(x)$. Differentiation in $x$ gives the first displayed component. The Koszul differential sends each degree-$-1$ generator to the corresponding partial derivative and is zero on coordinate functions. Its square is zero, and it is precisely the defining differential of the derived critical algebra. The Hessian symmetry supplies the cyclic pairing of this two-term model.
\end{proof}

The ordinary potential $F=f+h_N$ records both the normal equation
and its dual. The next corollary controls a different operation:
the addition or removal of an acyclic two-term summand in a chosen
resolution. Denote the terms of that summand by $P^0,P^1$, so they
are distinguished from the cohomology groups $W_1,W_2$ of $M_E$.
This operation accounts for auxiliary variables in equivalent
critical presentations.

\begin{corollary}\label{cor:quadratic}
Let $P^0,P^1$ be finite-dimensional vector spaces and let
$d_P(x):P^0\to P^1$ be a formal family of linear maps with
$d_P(0)$ an isomorphism. For $z\in P^0$ and $\lambda\in(P^1)^\vee$,
put
\[
 F_P(x,z,\lambda)=f(x)+\langle\lambda,d_P(x)z\rangle.
\]
The formal change $z'=d_P(x)z$ identifies $F_P$ with
$f(x)+\langle\lambda,z'\rangle$. The derived critical locus of
$F_P$ is equivalent to that of $f$.
\end{corollary}
\begin{proof}
The constant term of $\det d_P(x)$ is nonzero. It is a unit in the
complete local ring, so the indicated coordinate change is
invertible. Write $r=\dim P^0=\dim P^1$. In the new coordinates the
critical algebra is the completed tensor product of the critical
algebra of $f$ with
\[
 Q=\CC[[z'_1,\ldots,z'_r,\lambda_1,\ldots,\lambda_r]]
      [\alpha_1,\ldots,\alpha_r,\beta_1,\ldots,\beta_r],
 \quad |\alpha_i|=|\beta_i|=-1,
\]
where
\[
 d_Q\alpha_i=\lambda_i,\qquad d_Q\beta_i=z'_i,
 \qquad d_Qz'_i=d_Q\lambda_i=0.
\]
Its augmentation sends all generators to zero. Define the
odd operator
\[
 H_0=\sum_i\left(\beta_i\partial_{z'_i}
                        +\alpha_i\partial_{\lambda_i}\right).
\]
The graded commutator satisfies
\[
 d_QH_0+H_0d_Q
 =\sum_i\left(z'_i\partial_{z'_i}+\lambda_i\partial_{\lambda_i}
                 +\alpha_i\partial_{\alpha_i}+\beta_i\partial_{\beta_i}\right).
\]
On each monomial of positive total degree $m$, set $H=m^{-1}H_0$,
and set $H(1)=0$. This is continuous on the completed algebra and
gives $d_QH+Hd_Q=1-\iota\epsilon$, where $\epsilon:Q\to\CC$ and
$\iota:\CC\to Q$ are augmentation and inclusion. Tensoring this
contraction with the critical algebra of $f$ proves the asserted
quasi-isomorphism.
\end{proof}

For the minimal transfer onto $H^\bullet(M_E)$ used above,
\eqref{eq:r13-Phi-mixed} gives $\Phi(0)=0$ on the cohomology
groups $W_j$ of \eqref{eq:r13-normal-spaces}. These groups are therefore not the acyclic
summand of the corollary. In the fiberwise hyperplane example,
$\dim W_1=1$ and $\dim W_2=4$ by
\eqref{eq:r13-hyperplane-dimensions}; no isomorphism between them
is imposed anywhere in that application.

\begin{remark}\label{rem:foliation}The real foliation is constructed from the $(-2)$-shifted symplectic form of the specified fourfold moduli presentation. The local potentials of the cyclic two-term model give $(-1)$-shifted derived critical loci. The parameter calculation below concerns these ordinary local potentials.\end{remark}

\begin{remark}\label{rem:odd-potential}
The local functions of Proposition~\ref{prop:additive-potential} have degree zero. Their critical equations and Koszul differentials are the input for the ordinary matrix-factorization calculations in the next section.
\end{remark}

\subsection{The coefficients of the shifted potential at the moduli point}
\label{sssec:r13-odd-normal}
We next identify the degree-$-1$ generators in a fourfold chart.
Return to $i_b:X_b\hookrightarrow\cX$ and fix the volume form
from the chosen base coordinate and relative volume. Suppose the
normal triangle for $E$ splits, and choose compatible linear
identifications in the reduced fourfold chart. Then
\[
 \ext^1_{\cX}(i_{b*}E,i_{b*}E)=U\oplus W_1,
 \qquad
 \ext^2_{\cX}(i_{b*}E,i_{b*}E)=U^\vee\oplus W_2,
 \qquad W_2\simeq U.
\]
The last isomorphism uses $v=\partial_t|_b$ in
\eqref{eq:r13-fiber-normal-spaces}. A splitting compatible with
the fourfold pairing identifies its middle quadratic space with
\begin{equation}\label{eq:r13-hyperbolic-pairing}
 Q((\alpha,w),(\alpha',w'))=\alpha(w')+\alpha'(w),
 \qquad(\alpha,w),(\alpha',w')\in U^\vee\oplus W_2.
\end{equation}
We specify this compatibility when using the following coordinate
presentation. A splitting of the vector spaces alone specifies
neither the higher operations nor the foliation quotient.

Let $z=(x,a)$ be the degree-zero coordinates and let
$\ell_n^{\cX}$ denote the transferred operations of the ambient
moduli problem in the chosen chart. Its obstruction section is
\begin{equation}\label{eq:r13-ambient-obstruction}
 \sigma(z)=\sum_{n\geq1}\frac1{n!}
                   \ell_n^{\cX}(\widetilde z,\ldots,\widetilde z)
       =(K(z),s(z))\in (U^\vee\oplus W_2)\widehat\otimes R.
\end{equation}
Here $R$ is the complete degree-zero coordinate ring and
$\widetilde z$ is the element represented by $z$ in the chosen
deformation model. In a minimal model $\ell_1^{\cX}=0$.
The maps $K$ and $s$ are obtained by projecting this series onto
the two indicated obstruction spaces. A comparison respecting the
fiber deformation model is required to identify its restriction as
\begin{equation}\label{eq:r13-fiber-compatibility}
 K(x,0)=df(x),\qquad s(x,0)=0.
\end{equation}
More generally the comparison includes the change of obstruction
frame and the higher components of the morphism of deformation
algebras induced by $i_{b*}$. Equation
\eqref{eq:r13-fiber-compatibility} is the normalization used when
that comparison has been constructed. It is included among the
hypotheses on the chart in the next proposition.

The first derivatives in $a$ are the mixed deformation operations.
For example, when the model identifies their normal component with
\eqref{eq:r13-Phi-mixed},
\begin{equation}\label{eq:r13-normal-Taylor}
 s(x,a)=\Phi(x)a+\sum_{m\geq2}\frac1{m!}s_m(x)(a,\ldots,a).
\end{equation}
Thus $h_N$ of \eqref{eq:r13-h-normal} records the linear-in-$a$
part of $s$ after pairing with an additional even variable.
The odd potential instead uses the complementary obstruction
coefficients $K$ and odd generators, as follows.

\begin{proposition}\label{prop:r13-complex-coefficients}
Assume a reduced complex fourfold Darboux presentation at $i_{b*}E$
has the identifications \eqref{eq:r13-hyperbolic-pairing}--
\eqref{eq:r13-fiber-compatibility}, with obstruction section
$\sigma=(K,s)$ and $Q(\sigma,\sigma)=0$. In dual bases write
$K=(K_1,\ldots,K_n)$ and $s=(s_1,\ldots,s_n)$.
For the complex polarization retaining the $W_2$ obstruction terms,
the quotient and its shifted potential are
\begin{equation}\label{eq:r13-complex-quotient}
 B=\bigl(R[\theta_1,\ldots,\theta_n],\delta\bigr),
 \quad |\theta_j|=-1,\quad\delta\theta_j=s_j,
 \qquad\nu=\sum_jK_j\theta_j.
\end{equation}
The generators $\theta_j$ belong to $W_2^\vee[1]$ and are dual to
normal obstruction directions. They are not the degree-zero
coordinates of a variable $\eta\in W_2^\vee$ in
\eqref{eq:r13-h-normal}.

Let $f\in R$ be a specified extension of the ordinary fiber
potential. If multiplication by $f$ is injective, a regular factor
$g$ in this presentation exists precisely when $K_j=f\beta_j$ in
$R$ for every $j$. In that case
\begin{equation}\label{eq:r13-g-normal}
 g=\sum_j\beta_j(x,a)\theta_j,
 \qquad\sum_j\beta_js_j=0,
 \qquad \nu=fg.
\end{equation}
The coefficient column of $g$ is therefore the complementary
projection of \eqref{eq:r13-ambient-obstruction}, divided by $f$
when this division is regular.
\end{proposition}
\begin{proof}
The quadratic identity is
$\sum_jK_js_j=0$, so \eqref{eq:r13-complex-quotient} gives
$\delta\nu=0$. \cite[Proposition~4.17]{KS}
reconstructs the full reduced critical algebra by adjoining
$\zeta_j$ of degree $-1$ and $p_\mu$ of degree $-2$, with
\begin{equation}\label{eq:r13-fourfold-differential}
 \begin{aligned}
  \delta\theta_j&=s_j,&\delta\zeta_j&=K_j,\\
  \delta p_\mu&=-\sum_j
       \left((\partial_{z_\mu}s_j)\zeta_j
                     +(\partial_{z_\mu}K_j)\theta_j\right).
 \end{aligned}
\end{equation}
These are respectively the normal and tangential obstruction
coefficients and their derivatives. The square on $p_\mu$ is
$-\partial_{z_\mu}\sum_jK_js_j=0$; the other squares vanish
immediately. Under the specified Darboux identification these are
the chart's differentials.

The degree-$-1$ part of $B$ is free over $R$ on the $\theta_j$.
Comparing coefficients of $\nu=fg$ gives precisely
$K_j=f\beta_j$. Injectivity of $f$ gives uniqueness and
\[
 0=\delta\nu=f\delta g=f\sum_j\beta_js_j,
\]
which proves closedness. This also proves the coefficient
interpretation in \eqref{eq:r13-g-normal}.

For an invertible change of frame $\theta'=A\theta$, the new
columns are $s'=As$, $K'=A^{-T}K$ and
$\beta'=A^{-T}\beta$. Hence
\[
 K'^T\theta'=K^T\theta,\quad
 \beta'^T\theta'=\beta^T\theta,\quad
 K'^Ts'=K^Ts.
\]
The quotient, potential, closedness relation and factorization
are invariant under these frame changes. The remaining momenta
transform by the cotangent changes of coordinates of
\cite[Proposition~4.17]{KS}.
\end{proof}

\begin{proposition}\label{prop:v15-fiber-product-obstruction}
Let $R=\CC[[x_1,\ldots,x_n,a_1,\ldots,a_r]]$ and let
$\overline R=R/(a_1,\ldots,a_r)$. Suppose the coefficients of
\eqref{eq:r13-complex-quotient} satisfy
\[
 K_j\bmod(a_1,\ldots,a_r)=\partial_{x_j}f_0,
 \qquad f_0\in(x_1,\ldots,x_n)^2\subset\overline R.
\]
Let $f\in R$ restrict to $f_0$. If $f_0\ne0$, then there is no
regular degree-$-1$ element $g=\sum_j\beta_j\theta_j$ with
$\nu=fg$ in this presentation. If also $s_j\bmod(a_1,\ldots,a_r)=0$
for every $j$, no identity $\nu=fg+\delta h$ with $h\in B^{-2}$
exists in this Koszul algebra either.
\end{proposition}
\begin{proof}
Freeness of the degree-$-1$ module implies
$K_j=f\beta_j$ for every $j$. Reducing modulo the base-direction
ideal gives
\begin{equation}\label{eq:v15-gradient-divisibility}
 \partial_{x_j}f_0=f_0\overline\beta_j
          \quad\text{in }\CC[[x_1,\ldots,x_n]].
\end{equation}
Write $m=\operatorname{ord}(f_0)\ge2$ and let $f_{0,m}$ be its
nonzero homogeneous term of degree $m$. Euler's identity
\[
 \sum_j x_j\partial_{x_j}f_{0,m}=m f_{0,m}
\]
shows that $\partial_{x_j}f_{0,m}\ne0$ for some $j$. For this index,
$\operatorname{ord}(\partial_{x_j}f_0)=m-1$. If
$\overline\beta_j=0$, the right side of
\eqref{eq:v15-gradient-divisibility} is zero. Otherwise its order is
$m+\operatorname{ord}(\overline\beta_j)\ge m$. Both alternatives
contradict the left side. Thus the assumed coefficient factorization
cannot hold. The argument uses only restriction to the fiber and
regularity of the coefficients; it does not depend on a choice of
higher terms in $f$ away from that fiber.

For the last assertion, $\delta$ is zero on $R$ and sends
$\theta_j$ to $s_j$. If all $s_j$ restrict to zero, the induced
differential on $B\otimes_R\overline R$ is zero. Thus
$(\delta h)\bmod(a_1,\ldots,a_r)=0$, and restriction of the
putative homotopy factorization gives exactly
\eqref{eq:v15-gradient-divisibility}. The preceding order argument
applies without change.
\end{proof}

For the zero-potential divisor loci the premise $f_0\ne0$ does not
hold. For a nonzero fiber potential, the proposition applies to the
specified complex polarization satisfying
\eqref{eq:r13-fiber-compatibility}. It does not identify a real BKSY
quotient with this complex polarization. An application using a
different quotient must specify its function algebra and its
comparison with the fiber deformation algebra; the coefficient
factorization cannot be inferred from the vector spaces $U,W_1,W_2$.

This proposition applies to the stated complex polarization.
The real negative-definite foliation used in the BKSY construction
has a different half of the real obstruction space. Its coefficients
can also be related to $K$ and $s$, rather than denoted by unrelated
functions. The next calculation supplies that relation in a standard
split presentation. It is used after the chart and foliation
representative have been chosen as in
\cite[Hypothesis~4.7]{KS}.

\begin{proposition}\label{prop:r13-real-coefficients}
In the quadratic space of \eqref{eq:r13-hyperbolic-pairing}, choose
a positive Hermitian form $h$ on $W_2\simeq U$, linear in its first
argument, and the antilinear isomorphism
\[
 J:W_2\longrightarrow U^\vee,\qquad (Jw)(z)=h(z,w).
\]
Set
\[
 E^-_{\RR}=\{(-Jw,w):w\in W_2\},\qquad
 E^+_{\RR}=\{(Jw,w):w\in W_2\}.
\]
Both are real subspaces. The first is negative definite for
$\operatorname{Re}Q$ and isotropic for $\operatorname{Im}Q$.
For a local representative with degree-one foliation subbundle
$E^-_{\RR}$, use $E^+_{\RR}$ to represent the quotient. Write
\begin{equation}\label{eq:r13-real-projections}
 v_+=\tfrac12(s+J^{-1}K),\quad v_-=\tfrac12(s-J^{-1}K),
 \qquad\sigma_+=(Jv_+,v_+),\quad\sigma_-=(-Jv_-,v_-).
\end{equation}
The section and cosection of the corresponding real
bundle-and-cosection presentation are
\begin{equation}\label{eq:r13-real-section-cosection}
 s_+=\sigma_+,\qquad
 \tau_+=\operatorname{Im}Q(\sigma_-,-)|_{E^+_{\RR}},
 \qquad\tau_+(s_+)=0.
\end{equation}
In a unitary frame, with real basis $(Je_j,e_j),(J(ie_j),ie_j)$
of $E^+_{\RR}$ and its dual odd generators
$\theta_j^{\mathrm R},\theta_j^{\mathrm I}$, these formulas are
\begin{equation}\label{eq:r13-real-odd-potential}
 \begin{aligned}
 \delta\theta_j^{\mathrm R}&=\operatorname{Re}v_{+,j},&
 \delta\theta_j^{\mathrm I}&=\operatorname{Im}v_{+,j},\\
 \nu_{\RR}&=\sum_j\left(
       2\operatorname{Im}v_{-,j}\,\theta_j^{\mathrm R}
       -2\operatorname{Re}v_{-,j}\,\theta_j^{\mathrm I}\right).
 \end{aligned}
\end{equation}
Thus any proposed factor $\nu_{\RR}=fg$, after the specified
complexification when necessary, has coefficient equations
\[
 f\beta_j^{\mathrm R}=2\operatorname{Im}v_{-,j},\qquad
 f\beta_j^{\mathrm I}=-2\operatorname{Re}v_{-,j},\qquad
 g=\sum_j(\beta_j^{\mathrm R}\theta_j^{\mathrm R}
                  +\beta_j^{\mathrm I}\theta_j^{\mathrm I}).
\]
These are equations in the function algebra of the chosen real
derived chart; regularity and holomorphicity are stated in that
algebra and its specified holomorphic comparisons.
\end{proposition}
\begin{proof}
The Hermitian symmetry gives
\[
 Q((-Jw,w),(-Jz,z))=-h(z,w)-h(w,z)
                              =-2\operatorname{Re}h(w,z).
\]
This is real, and is negative on nonzero diagonal vectors.
Changing both signs to plus gives the positive real form on
$E^+_{\RR}$. The imaginary pairing between the two summands is
nondegenerate: if $w\ne0$, evaluating it on the vector corresponding
to $iw$ gives a nonzero multiple of $h(w,w)$.
Consequently it identifies $E^-_{\RR}$ with $(E^+_{\RR})^\vee$.

Equations \eqref{eq:r13-real-projections} give
$\sigma_-+\sigma_+=\sigma$. The two diagonal values of $Q$ on
these summands are real. Since $Q(\sigma,\sigma)=0$,
\[
 0=\operatorname{Im}Q(\sigma,\sigma)
   =2\operatorname{Im}Q(\sigma_-,\sigma_+)
   =2\tau_+(s_+).
\]
This proves closedness of the cosection. In a unitary frame,
$Jw=\overline w$ and
\[
 Q((-\overline v_-,v_-),(\overline w,w))
          =-\overline v_-\cdot w+\overline w\cdot v_-.
\]
Substituting $w=e_j$ and $w=ie_j$ gives, respectively,
$2\operatorname{Im}v_{-,j}$ and $-2\operatorname{Re}v_{-,j}$.
These are exactly the coefficients in
\eqref{eq:r13-real-odd-potential}. Its internal differential is
$\tau_+(s_+)=0$. Finally, freeness of the degree-$-1$ module makes
factorization equivalent to the displayed coefficient equations.
\end{proof}

In particular, the real quotient generators are dual to a real
quotient of $U^\vee\oplus W_2$. They are related to the original
obstruction complex by \eqref{eq:r13-real-projections}; they are
not copies of the normal parameter space $W_1$. The negative
subspace and its global integrable representative are furnished
by the BKSY construction \cite{BKSY2025,BKSY2024}; the preceding
proposition is the coefficient calculation on a representative
with the stated splitting. If that representative is described
in another frame, the section and cosection transform by dual
matrices as in Proposition~\ref{prop:r13-complex-coefficients}.

\subsection{Variation along the base and the hyperplane-section example}
\label{sssec:r13-base-variation}
Suppose a compatible relative critical presentation has been
constructed near $E$ with degree-zero variables $(x,t)$ and
ordinary potential $f(x,t)$. Let $a=t-b$. In a frame in which
the cyclic pairing is constant, the fiberwise formula is
\begin{equation}\label{eq:r13-cyclic-base-variation}
 \begin{aligned}
 f(x,t)&=\sum_{n\geq2}\frac1{(n+1)!}
          \langle\ell_{n,t}(x,\ldots,x),x\rangle,\\
 \partial_tf(x,t)&=\sum_{n\geq2}\frac1{(n+1)!}
          \langle(\partial_t\ell_{n,t})(x,\ldots,x),x\rangle.
 \end{aligned}
\end{equation}
The first equation assumes the centered minimal fiberwise models
with these operations. In a general nonminimal presentation one
uses its given $f(x,t)$, with its linear and quadratic terms as
well. In either presentation,
\[
 d_xf(x,b+a)=d_xf(x,b)+a\,d_x\partial_tf(x,b)+O(a^2).
\]
The coefficient of $a$ is therefore a tangential obstruction term,
obtained from the variation of the original deformation operations.
For the differentiated fourfold model of
Proposition~\ref{prop:v5-cotangent-completion}, the corresponding
coefficient occurs in the actual differential
\begin{equation}\label{eq:r13-mixed-transverse-differential}
 \nu_f=\sum_j\partial_{x_j}f(x,t)\theta_j,
 \qquad
 \delta p_a=-\sum_j\partial_t\partial_{x_j}f(x,t)\theta_j.
\end{equation}
Here $a$ is dual to $\id_E\otimes v\in W_1$, whereas the
$\theta_j$ are the odd functions dual to the $W_2$ obstruction
directions in the specified complex polarization. This is the
coordinate relation between normal motion and the fourfold
differential. If this model is identified with a geometric chart,
that identification must intertwine all the displayed differentials.

When a regular factorization has instead been verified in that
chart, \eqref{eq:r13-g-normal} gives $g$ from its actual
coefficients. Choosing $\kappa(g)=1$ then gives
\[
 \partial_tf=\kappa(\partial_t\nu)-f\kappa(\partial_tg),
 \qquad\nabla_{\partial_t}=\partial_t-u^{-1}\partial_tf.
\]
Thus the factor, when available, determines the coefficient of the
ordinary Gauss--Manin operator. A normal parameter in $W_1$, the
even function $h_N$, and the odd function $g$ enter three distinct
displayed formulas; none is identified with another by its name.

For the actual hyperplane-section locus of
Section~\ref{sec:v6-hyperplanes}, the universal two-term complex
extends in both the hyperplane and base directions. Equations
\eqref{eq:r13-hyperplane-dimensions} and
\eqref{eq:r13-base-obstruction} give the normal spaces and the
vanishing first order obstruction. The reconstruction of Artinian
deformations in Proposition~\ref{prop:v6-projective-locus} proves
that the classical moduli locus is $\P^4\times C$, and
Theorem~\ref{thm:v6-gluing-locus} supplies its zero critical
presentation. Hence $f(x,t)=0$, the reduced obstruction section
in that presentation is zero, and both the complex and real
coefficient formulas above give $\nu=0$. A product $0=0\cdot g$
does not determine $g$. The regular comparison in
Theorem~\ref{thm:v6-comparison} is induced by the constructed
zero-potential bar maps and $\nabla_{\partial_t}=\partial_t$.

\subsection{Graded factorization of the shifted potential}\label{sec:graded-factorization}
The potential of the fourfold foliation and the potential of a
threefold critical chart have different degrees. A proposed relation
between them must therefore be an identity in a specified differential
graded algebra. In this subsection we formulate the identity, compute its
obstruction, and examine the shifted critical algebra when it holds.
We also construct a cotangent-completion model in which the relation to
the ordinary potential is differentiation rather than multiplication.
All statements concern the displayed algebraic or formal charts.

Fix a characteristic-zero field $k$ and a regular local $k$-algebra $R$,
or a power-series algebra over $k$. A parameter $t$ on the geometric base
curve can be included among the degree-zero variables. Let
\begin{equation}\label{eq:v5-Koszul-base}
 B=(R[\theta_1,\ldots,\theta_r],\delta),\qquad
 |\theta_a|=-1,\qquad \delta\theta_a=s_a\in R,\quad \delta R=0.
\end{equation}
Brackets in this expression mean the free graded-commutative algebra;
in particular, the $\theta_a$ anticommute. The shifted potential is a
closed degree-$-1$ element
\begin{equation}\label{eq:v5-nu-coefficients}
 \nu=\sum_{a=1}^r\tau_a\theta_a,\qquad
 \delta\nu=\sum_a\tau_as_a=0.
\end{equation}
In the complex presentation of
Proposition~\ref{prop:r13-complex-coefficients}, $r=\dim W_2$,
the generators span $W_2^\vee[1]$, and $\tau_a=K_a$ is the
complementary component of the transferred fourfold obstruction map.
For the real BKSY representative in
Proposition~\ref{prop:r13-real-coefficients}, the generators are dual
to $E^+_{\RR}$ and the $\tau_a$ are the explicit real projections
in \eqref{eq:r13-real-odd-potential}. In both cases the quotient
has the displayed form. The internal differential is $\delta$;
$d_{\mathrm{dR}}$ differentiates its coefficient functions.

Let $f\in R$ be the chosen ordinary potential. In an application, a
comparison of charts must first place the fiberwise function $f_t$ in
this $R$; its existence is not implied merely by the degrees of $f$ and
$\nu$. We consider a degree-$-1$ factor
\[
 g=\sum_ab_a\theta_a\in B^{-1}.
\]
The equality $\nu=fg$ has the required degree. The following test says
exactly when it is an equality of regular functions on this derived
chart, as opposed to a meromorphic expression.

\begin{proposition}\label{prop:v5-factor-test}
Assume that multiplication by $f$ is injective on $R$.
\begin{enumerate}
\item A strict factorization $\nu=fg$ exists if and only if
$\tau_a\in fR$ for every $a$. The coefficients $b_a$ are then unique,
and $\delta g=0$.
\item The obstruction to strict factorization is the image
\[
 \overline\tau=(\tau_1\bmod f,\ldots,\tau_r\bmod f)
                  \in(R/fR)^r.
\]
Its vanishing is independent of an invertible change of frame of the
bundle of degree-$-1$ generators.
\item A factorization up to internal homotopy, in the precise sense
\begin{equation}\label{eq:v5-homotopy-product}
 \nu=fg+\delta h,\qquad
 g\in B^{-1},\ \delta g=0,\ h\in B^{-2},
\end{equation}
exists if and only if the image of $[\nu]$ vanishes in
$H^{-1}(B/fB)$.
\end{enumerate}
\end{proposition}
\begin{proof}
The degree-$-1$ part of $B$ is the free $R$-module with basis
$\theta_1,\ldots,\theta_r$. Comparing coefficients in
$\sum\tau_a\theta_a=f\sum b_a\theta_a$ gives
$\tau_a=fb_a$ for every $a$. These equations are soluble precisely under
the divisibility condition. Injectivity of $f$ gives uniqueness.
Moreover,
\[
 0=\delta\nu=\delta(fg)=f\delta g,
\]
and $B^0=R$ has no $f$-torsion. Therefore $\delta g=0$.
If a frame change replaces the coefficient column $\tau$ by
$M\tau$, with $M\in\mathrm{GL}_r(R)$, its reduction modulo $f$ is
invertible. Thus $M\overline\tau=0$ exactly when
$\overline\tau=0$.

Every graded piece of $B$ is free over $R$. Consequently
\begin{equation}\label{eq:v5-product-exact-sequence}
 0\longrightarrow B\xrightarrow{\ f\ }B
       \longrightarrow B/fB\longrightarrow0
\end{equation}
is a short exact sequence of complexes. Its long exact sequence includes
\[
 H^{-1}(B)\xrightarrow{\ f\ }H^{-1}(B)
                \longrightarrow H^{-1}(B/fB).
\]
This already proves the last equivalence. More explicitly, if the image
of $\nu$ is a boundary, choose $\overline h\in(B/fB)^{-2}$ with
$\overline\nu=\delta\overline h$ and lift it to $h\in B^{-2}$.
Then $\nu-\delta h=fg$ for some $g\in B^{-1}$.
Applying $\delta$ gives $f\delta g=0$; injectivity again proves that
$g$ is closed. The converse follows by reducing
\eqref{eq:v5-homotopy-product} modulo $f$.
\end{proof}

The homotopy test is weaker than the coefficient test: one may change a
representative of $\nu$ by a boundary. It is also the appropriate test
under quasi-isomorphisms preserving the specified ordinary function and
shifted-potential class. If $f$ is not injective on the underlying
complex, the quotient in the homotopy test must be replaced by the
cofiber of multiplication by $f$, or equivalently the derived zero
locus of $f$. We will use the free, $f$-torsion-free presentation above.
A change to a different foliation quotient is additional geometric data;
it is not an invertible change of frame in this proposition.

\begin{corollary}\label{cor:v5-order-obstruction}
Suppose $R=k[[x_1,\ldots,x_n]]$, $\mathfrak m=(x_1,\ldots,x_n)$,
and $0\ne f\in\mathfrak m$. If one coefficient $\tau_a$ has
$\mathfrak m$-adic order less than the order of $f$, then
$\nu$ has no strict factorization $fg$ in $B$. In particular, a
quadratic leading term in a coefficient of $\nu$ cannot be obtained
by multiplying a cubic-leading ordinary potential by a regular
homogeneous element of degree $-1$.
\end{corollary}
\begin{proof}
If $\tau_a=fb_a$ with $b_a\ne0$, multiplication of leading nonzero
homogeneous terms gives
$\operatorname{ord}_{\mathfrak m}(\tau_a)
=\operatorname{ord}_{\mathfrak m}(f)+
\operatorname{ord}_{\mathfrak m}(b_a)
\geq\operatorname{ord}_{\mathfrak m}(f)$.
A zero coefficient does not affect the test. This proves both claims.
\end{proof}

\subsubsection{The critical equations for a product of graded functions}
Suppose that the strict factorization has been established, so
$\tau_a=fb_a$ and $\sum b_as_a=0$. Substitution in
\cite[equation~(4.24)]{KS} gives the full differential
\begin{equation}\label{eq:v5-product-critical}
\begin{aligned}
 \delta x_i&=0,&\delta\theta_a&=s_a,&
 \delta\eta_a&=fb_a,\\
 \delta p_i&=-\sum_a(\partial_i s_a)\eta_a
             -(\partial_i f)\sum_ab_a\theta_a
             -f\sum_a(\partial_i b_a)\theta_a.
\end{aligned}
\end{equation}
In particular, its degree-zero cohomology is
\begin{equation}\label{eq:v5-product-truncation}
 H^0(\mathcal A)=R/(s_1,\ldots,s_r,fb_1,\ldots,fb_r).
\end{equation}
Indeed, the boundaries in degree zero are precisely the ideal generated
by the differentials of the degree-$-1$ generators; products involving
$p_i$ have lower degree and contribute no new degree-zero boundaries.
This calculation must be compared with the actual truncation of the
ambient moduli chart.

Call $g$ \emph{primitive} at a point if the vector $(b_1,\ldots,b_r)$
has a unit entry in the local ring. Equivalently its coefficients
locally generate the unit ideal. Replacing the generators by an
invertible $R$-linear change of frame makes $\theta'_1=g$. Since
$\delta g=0$, this generator has zero internal differential.
In these coordinates the coefficients of $\nu$ are $(f,0,\ldots,0)$,
and \eqref{eq:v5-product-truncation} becomes
$R/(s'_2,\ldots,s'_r,f)$.
Thus a primitive odd factor imposes the equation $f=0$ among the
classical equations. It does not by itself impose only the gradient
equations of the fiberwise critical locus.

\subsubsection{A differentiated potential and its cotangent completion}
There is a useful explicit model that retains all degree-$-2$
generators. Let
\[
 R=k[t,x_1,\ldots,x_n]\quad\text{or its local/formal version},
 \qquad f=f(x,t),\qquad f_i=\partial f/\partial x_i.
\]
The relative critical algebra is
\begin{equation}\label{eq:v5-relative-critical-algebra}
 K_f=(R[\eta_1,\ldots,\eta_n],\delta),\qquad
 |\eta_i|=-1,\quad \delta\eta_i=f_i.
\end{equation}
Here differentiation is in the $x$-variables, not in $t$.
Thus its ordinary equations are the equations for critical points on
the fibers of the family.

\begin{proposition}\label{prop:v5-cotangent-completion}
On $B_0=R[\theta_1,\ldots,\theta_n]$, with $|\theta_i|=-1$ and
zero internal differential, put
\begin{equation}\label{eq:v5-differentiated-potential}
 \nu_f=\sum_{i=1}^n f_i(x,t)\theta_i.
\end{equation}
Its shifted critical algebra is
\[
 \mathcal A_f=R[\theta_i,\eta_i,p_1,\ldots,p_n,p_t],
 \qquad |p_i|=|p_t|=-2,
\]
with
\begin{equation}\label{eq:v5-cotangent-differential}
\begin{aligned}
 \delta\theta_i&=0,&\delta\eta_i&=f_i,\\
 \delta p_j&=-\sum_i f_{ij}\theta_i,&
 \delta p_t&=-\sum_i f_{it}\theta_i.
\end{aligned}
\end{equation}
This is also the algebra of the absolute shifted cotangent completion
$T^*[-2]\operatorname{Spec}K_f$. In particular
$H^0(\mathcal A_f)=R/(f_1,\ldots,f_n)$, and the projection to the
relative critical locus is induced by $K_f\hookrightarrow\mathcal A_f$.
\end{proposition}
\begin{proof}
In \cite[Proposition~4.17]{KS} set $s_i=0$ and
$\tau_i=f_i$, and include $t$ among the degree-zero coordinates.
That gives all the differentials in
\eqref{eq:v5-cotangent-differential}; their squares vanish since
$\delta f_i=\delta\theta_i=0$. The Hamiltonian is
$\sum_i f_i\theta_i$ for the bracket pairing $x_j$ with $p_j$,
$t$ with $p_t$, and $\eta_i$ with $\theta_i$.
Its self-bracket is zero, since it involves no $p$ or $\eta$
variables. The nondegenerate pairing on the free coordinate tangent
module is the one in
\cite[Proposition~4.17]{KS}.

For the second description start with $K_f$. The cotangent lift of
its differential $\delta\eta_i=f_i$ has Hamiltonian
$\sum_i f_i\theta_i$, where $\theta_i$ is cotangent-dual to
$\eta_i$. Differentiation in $x_j$ and $t$ gives the two displayed
negative derivatives on $p_j,p_t$. This proves the identification
by generators, differential, and symplectic pairing, rather than
only by equality of tangent dimensions.
The inclusion $K_f\hookrightarrow\mathcal A_f$ is a chain map.
Sending $\theta_i,p_j,p_t$ to zero defines its retraction as an
algebra map, and the zero-degree boundaries are generated by the $f_i$.
\end{proof}

\begin{proposition}\label{prop:v5-gradient-not-divisible}
Let $0\ne f\in\mathfrak m^2\subset k[[x_1,\ldots,x_n]]$.
For the differentiated model \eqref{eq:v5-differentiated-potential},
$\nu_f$ is not of the form $fg$ with a regular homogeneous
$g$ of degree $-1$.
\end{proposition}
\begin{proof}
Let $m=\operatorname{ord}_{\mathfrak m}(f)$ and write its first
nonzero homogeneous term as $f_m$. In characteristic zero the
nonconstant polynomial $f_m$ has a nonzero partial derivative:
otherwise Euler's identity $\sum_i x_i\partial_i f_m=mf_m$
would give $f_m=0$. For such an index $i$,
$\operatorname{ord}_{\mathfrak m}(f_i)=m-1$.
Corollary~\ref{cor:v5-order-obstruction} now applies.
\end{proof}

The differentiated construction has a geometric origin in the split
normal-direction model. If a minimal cyclic threefold deformation
algebra has degree-one part $U$, degree-two part $U^\vee$, and
potential $f$, its cotangent completion adds the dual obstruction
generators. The differential on those generators is the derivative
of the obstruction map, namely the Hessian of $f$.
This is exactly \eqref{eq:v5-cotangent-differential} when $f$ is
independent of $t$. In a family, the additional term
$-\sum_i f_{it}\theta_i$ records its variation in the base direction.
The proposition constructs this completion; identifying it with a
particular vertical sheaf-moduli chart requires a comparison of its
deformation differential with that of the chart.

One can also inspect the tangent differential directly. At a classical
point with $f_i=0$ and all negative-degree generators set to zero,
a tangent vector $(v_1,\ldots,v_n,v_t)$ maps to the linearization
\[
 \bigl(\sum_jf_{ij}v_j+f_{it}v_t\bigr)_{i=1}^n
\]
in the obstruction block dual to the $\eta_i$. The following map is
its signed transpose in the symplectic pairing, by the displayed
formulas for $\delta p_j$ and $\delta p_t$. Thus the mixed partial
derivatives are actual entries of the fourfold deformation complex;
they are not additional degree-zero factors inserted into $\nu_f$.

\subsubsection{Recovering an ordinary potential from a primitive odd factor}
When a product does exist, it gives a precise local extraction formula.
The extraction and variation formulas are as follows.

\begin{proposition}\label{prop:v5-extraction}
Suppose $\nu=fg$ in \eqref{eq:v5-Koszul-base} and that $g$ is
primitive. Choose functions $c_a\in R$ with $\sum_ac_ab_a=1$, and
let
\[
 \kappa=\sum_ac_a\partial_{\theta_a}
\]
be the $R$-linear graded derivation of degree $+1$, where odd derivatives
are left derivatives. Then
\begin{equation}\label{eq:v5-extract-f}
 \kappa(g)=1,\qquad f=\kappa(\nu).
\end{equation}
For a base derivation $\partial_t$ acting on the coefficients and fixing
the chosen $\theta_a$,
\begin{equation}\label{eq:v5-extract-variation}
 \partial_t f=\kappa(\partial_t\nu)
                     -f\,\kappa(\partial_t g).
\end{equation}
The right side is independent of the choice of $\kappa$ satisfying
$\kappa(g)=1$.
\end{proposition}
\begin{proof}
Because $f$ has degree zero and belongs to $R$, $\kappa(f)=0$.
The graded Leibniz rule therefore gives
$\kappa(fg)=f\kappa(g)=f$. Differentiating the product first gives
\[
 \partial_t\nu=(\partial_tf)g+f\partial_tg.
\]
Applying $\kappa$ proves \eqref{eq:v5-extract-variation}.
If $\kappa'$ is another choice, both right sides equal the same
$\partial_tf$, which proves independence. On the open set where
$b_j$ is a unit, one may take $\kappa=b_j^{-1}\partial_{\theta_j}$;
then the formula reads
\[
 \partial_tf=b_j^{-1}\partial_t\tau_j
                         -f b_j^{-1}\partial_tb_j.
\]
No division by an odd element occurs: $g^2=0$, so $g$ is not invertible.
\end{proof}

If $\nu$ and the primitive $g$ are global sections on a fixed smooth
family of such charts, local extracted functions agree on overlaps.
Indeed $(f_\alpha-f_\beta)g=0$, and a unit coefficient of $g$ implies
$f_\alpha=f_\beta$. If instead one changes the chosen factor to
$g_\beta=a_{\alpha\beta}g_\alpha$, the extracted function changes to
$f_\beta=a_{\alpha\beta}^{-1}f_\alpha$.
Thus a global shifted potential by itself does not fix the transition
rule of the ordinary potentials; the factor and its transition data
are part of the extraction construction. Section~\ref{sec:v5-global-GM}
constructs the connection once an ordinary function on a smooth family
has been obtained in this way or otherwise.

The coefficient formulas \eqref{eq:r13-g-normal} and
\eqref{eq:r13-real-odd-potential} specify the possible odd factor
in terms of the original fourfold obstruction section. The mixed
base derivative is calculated in
\eqref{eq:r13-cyclic-base-variation} and
\eqref{eq:r13-mixed-transverse-differential}. Once the primitive
factorization and the comparison with ordinary critical presentations
are given, the extraction formula determines the coefficient of
$\partial_t-u^{-1}\partial_tf$. We now prove that this operator
glues on the stated ordinary-potential family, and establish the
finiteness and specialization conditions for its cohomology.

\section{Gauss--Manin connection on fiberwise periodic cyclic homology}
\label{sec:cyclic}
We now construct the connection on a specified ordinary critical
realization and prove the finite-complex specialization criterion.
The coordinate $t$ is on the degeneration curve; the even variable
$u$ records periodicity. The point of the construction is the
compatibility of the connection with the differential and the overlap
homotopies, before cohomology is taken.

\subsection{Orientations and the chosen cyclic realization}
\label{ssec:orientation-input}
Let $M$ be a finite-type classical moduli space of simple perfect
objects after scalar rigidification, with compatible critical charts
$(R_\alpha,U_\alpha,f_\alpha)$ and their d-critical structure $s$.
These data are the geometric input recalled in
\cite[Section~5.1]{KS}; in particular rigidification and a critical
presentation are separate operations.
\begin{hypothesis}\label{hyp:DT}
Choose a line bundle $\mathcal O_M^{\mathrm{or}}$ on $M^{\red}$
and an isomorphism
$(\mathcal O_M^{\mathrm{or}})^{\otimes2}\simeq K_{M,s}$,
compatible with the critical-chart identifications
$K_{M,s}|_{R_\alpha^{\red}}\simeq K_{U_\alpha}^{\otimes2}|_{R_\alpha^{\red}}$.
In a family these identifications include their overlap compatibility.
\end{hypothesis}
The oriented vanishing-cycle sheaf $\mathcal P_{M,s}$ is the one
constructed in \cite{BBDJS,Joyce}. We use
\begin{equation}\label{eq:DT-convention}
 \mathrm{DT}^j_{\mathrm{coh}}(M,s)
   =\mathbb H_c^j(M^{\an},\mathcal P_{M,s}).
\end{equation}
On a smooth zero critical locus of dimension $d$, the orientation
$K_M$ gives $\mathcal P_{M,0}=\CC_M[d]$; other orientation choices
can introduce a local system.

For a smooth critical chart $(U,f)$, $\MF(U,f)$ denotes the coherent
matrix-factorization category with odd differential squaring to
$f$. We use its localization and periodic cyclic complex in the
sense of \cite{Efimov2018CyclicMF}. The comparison in that reference,
for $U$ smooth separated of finite type and
$\Crit(f)\subset f^{-1}(0)$, is
\begin{equation}\label{eq:cyclic-DR-input}
 \HP_\bullet(\MF(U,f))
  \simeq\mathbb H^\bullet
   \bigl(U,(\Omega_U^\bullet((u)),u\,d-df\wedge)\bigr).
\end{equation}
The normalization used for the DT bar includes the shift $[d]$.
Where several critical charts are used, their category comparisons,
orientations, and coherent descent are part of the chosen realization.
The global cyclic comparison below is asserted only for such specified
data. The explicit projective and Grassmannian examples use the
single global zero critical presentation.

For clarity, coherent descent here includes quasi-equivalences between
the chosen local categorical presentations on double overlaps,
specified homotopies between their composites on triple overlaps,
and the higher compatibility homotopies on higher intersections.
On their normalized periodic complexes one also specifies the
comparison with the oriented critical-chart models. The resulting
totalization and its support convention define
$\mathsf K^{\mathrm{DT}}$. The identification of that totalization
with the periodic cyclic complex of a global category, when used,
is an additional comparison in this datum; we do not commute
periodic cyclic homology with an arbitrary limit of categories.
The bar maps constructed below are required to respect these same
overlap identifications. This formulation makes explicit the
geometric input needed in addition to the local operator identity.

\subsection{Local operators and the geometric parameter}
\label{ssec: Periodic cyclic homology of curved DG-categories}
For a curved category the periodic differential has the mixed form
$b+uB$ on Laurent-completed Hochschild chains. On the ordinary
function chart it is represented by $D_f=u\,d-df\wedge$ in
\eqref{eq:cyclic-DR-input}. The connection in $u$ and the connection
in $t$ have different meanings. With form-degree operator
$\Gamma|_{\Omega^p}=-p/2$, the former is represented by
\begin{equation}\label{eq:u-connection}
 \nabla_{\partial_u}=\partial_u+\frac{f}{u^2}+\frac{\Gamma}{u}.
\end{equation}
It satisfies $[\nabla_{\partial_u},D_f]=D_f/(2u)$ in this
normalization \cite{Efimov2018CyclicMF}. The base operator instead
commutes with $D_f$ itself. The one-parameter product formula is
recorded in \cite[Section~5.2]{KS}; the following version also records
its mixed-parameter curvature, which is useful before descent.
\begin{proposition}\label{prop:parameter-GM}
For a smooth $U$, a smooth base $T$ with local coordinates $t_i$,
and $F\in\Gamma(U\times T,\O)$, set
\[
 D=u\,d_U-d_UF\wedge,\qquad
 \nabla_i=\partial_{t_i}-u^{-1}F_{t_i}.
\]
Then $[\nabla_i,D]=0$ and $[\nabla_i,\nabla_j]=0$.
\end{proposition}
\begin{proof}
Writing $M_g$ for multiplication by $g$, one has
$[\partial_{t_i},D]=-d_UF_{t_i}\wedge$ and
$[M_{-F_{t_i}/u},u\,d_U]=d_UF_{t_i}\wedge$.
The remaining multiplication operators commute. The curvature is
$-u^{-1}M_{F_{t_it_j}}+u^{-1}M_{F_{t_jt_i}}=0$.
\end{proof}
These identities do not imply finiteness or base change of cohomology.
For a nonproduct ambient family the lifts of $\partial_t$ differ on
overlaps. The next calculation constructs the correcting homotopies,
in agreement with the cyclic Gauss--Manin formalism of \cite{Getzler1993}.
\subsection{A global chain connection for a specified family of ordinary potentials}\label{sec:v5-global-GM}
\subsubsection{The connection for the fiberwise moduli problem}
Fix the geometric family $\pi:\cX\to C$ and an affine coordinate
open $V\subset C$. At an object $E\in\Perf(X_b)$, the complexes
being compared are
\[
 A_E=\RH_{X_b}(E,E),\qquad
 B_E=\RH_{\cX}(i_{b*}E,i_{b*}E),
\]
\[
 A_E\longrightarrow B_E\longrightarrow
 A_E\otimes T_{C,b}[-1]\xrightarrow{+1}.
\]
The connecting map is the Atiyah--Kodaira--Spencer map of
\cite[Proposition~4.19]{KS}. In particular, for
$v=\partial_t|_b$ the candidate displacement of the support is
$\mathrm{id}_E\otimes v\in W_1$, and its obstruction is
$\operatorname{At}(E)\cdot\kappa(v)\in\ext^2_{X_b}(E,E)$.
The connection below acts on the cohomological invariant of this
\emph{fiberwise} moduli problem as $b$ varies. Its base variable
$t$ is the coordinate of $C$, not a coordinate on the fourfold
moduli stack.

Choose a relative critical presentation for a locus of the moduli
problem over $V$. This means a smooth ambient family
$q:\mathscr U\to V$, a regular degree-zero function $F$ on
$\mathscr U$, and the specified comparison of critical charts
\begin{equation}\label{eq:r14-relative-moduli-chart}
\begin{tikzcd}[column sep=large,row sep=large]
 \mathfrak M^{\mathrm{crit}}\arrow[r,"{\simeq}"]\arrow[dr]
   &\dCrit_{\mathscr U/V}(F)\arrow[r]\arrow[d]
   &\mathscr U\arrow[dl,"q"]\\
   &V.&
\end{tikzcd}
\end{equation}
Here $t_0\mathfrak M^{\mathrm{crit}}$ is the chosen classical
moduli locus, with the critical structure and orientation specified
in Hypothesis~\ref{hyp:DT}. A presentation of this kind can also
be made on an open locus. Its comparison is part of the data when
applying the next theorem to a moduli problem.

The relation to the deformation algebra is the following. At
$E\in M_b$ choose the reduced cyclic model
$L_{E,b}^{\mathrm{red}}$ of Hypothesis~\ref{hyp:cyclic-model}.
In minimal formal coordinates $x_1,\ldots,x_n$, the fiber of
\eqref{eq:r14-relative-moduli-chart} identifies
\begin{equation}\label{eq:r14-local-CE-identification}
 \widehat{\operatorname{CE}}(L_{E,b}^{\mathrm{red}})
 \simeq
 \left(\CC[[x_1,\ldots,x_n]][\xi_1,\ldots,\xi_n],
       \delta\xi_i=\partial_{x_i}f_b\right),
 \quad f_b=F|_{\mathscr U_b},\quad |\xi_i|=-1.
\end{equation}
The $x_i$ are dual to $U=\ext^1_{X_b}(E,E)$ in the chosen
reduction; the $\xi_i$ are dual to its obstruction space
$U^\vee$. In a cyclic trivialization their potential is
\[
 f_b(x)=\sum_{r\ge2}\frac1{(r+1)!}
        \langle\ell_{r,b}(x,\ldots,x),x\rangle.
\]
On a family of cyclic trivializations with fixed pairing,
$\partial_tF$ differentiates the higher-operation coefficients, as in
\eqref{eq:r13-cyclic-base-variation}. Thus the operator
$\partial_t-u^{-1}\partial_tF$ compares the complexes of the
chosen fiberwise deformation problems, rather than cyclic
homology of $\Perf(\cX)$ itself.

The fourfold quotient supplies another possible construction of
this same $F$. On a specified comparison with its complex
presentation, write
\begin{equation}\label{eq:r14-fourfold-to-fiber-function}
 B_\alpha=(\mathcal O_{\mathscr U_\alpha}[\theta_a],\delta),
 \quad\delta\theta_a=s_a,\quad
 \nu_\alpha=\sum_a\tau_a\theta_a,
 \quad g_\alpha=\sum_a\beta_a\theta_a,
 \quad\tau_a=F_\alpha\beta_a.
\end{equation}
In Proposition~\ref{prop:r13-complex-coefficients}, the $\tau_a$
are the complementary coefficients $K_a$ of the ambient
obstruction section and the $\theta_a$ belong to $W_2^\vee[1]$.
For the real BKSY representative they are the real coefficient
functions of Proposition~\ref{prop:r13-real-coefficients}.
Use a holomorphic comparison before treating $F_\alpha$ as an
ordinary holomorphic potential. In addition to the coefficient
equations in \eqref{eq:r14-fourfold-to-fiber-function}, the
extracted $F_\alpha$ must satisfy the moduli-chart identification
\eqref{eq:r14-local-CE-identification}. These are the two geometric
comparisons used when deriving $F$ from the fourfold foliation.

If $g_\alpha$ is primitive, choose
$\kappa_\alpha=\sum_ac_a\partial_{\theta_a}$ with
$\sum_ac_a\beta_a=1$. Then
\[
 F_\alpha=\kappa_\alpha(\nu_\alpha),\qquad
 \partial_tF_\alpha
   =\kappa_\alpha(\partial_t\nu_\alpha)
       -F_\alpha\kappa_\alpha(\partial_tg_\alpha).
\]
Consequently the operator of the next theorem is exactly
\[
 \nabla_{\partial_t}
  =\partial_t-\frac1u\left(
      \kappa_\alpha(\partial_t\nu_\alpha)
       -F_\alpha\kappa_\alpha(\partial_tg_\alpha)\right)
\]
in this presentation. Corollary~\ref{cor:v5-product-GM} proves
independence of $\kappa_\alpha$ and gluing when the product data
are compatible. This coefficient calculation gives the precise
role of the degree-$-1$ function in the connection.

For general comparisons of local periodic complexes
$\mathcal K_\alpha$, the required overlap maps obey
\[
 c_{\alpha\beta}:\mathcal K_\beta\xrightarrow{\sim}\mathcal K_\alpha,
 \qquad D_\alpha c_{\alpha\beta}=c_{\alpha\beta}D_\beta,
\]
\begin{equation}\label{eq:r14-connection-descent-data}
 \nabla_\alpha c_{\alpha\beta}
       -c_{\alpha\beta}\nabla_\beta
   =D_\alpha H_{\alpha\beta}+H_{\alpha\beta}D_\beta,
 \qquad |H_{\alpha\beta}|=-1.
\end{equation}
The maps and their homotopies are required to satisfy the higher
overlap compatibilities. The theorem below constructs them for
one smooth ambient family with a global $F$, using restriction
of relative forms and differences of horizontal lifts. It does
not identify arbitrary local critical charts by restriction of
functions. In the divisor-section loci of
Theorem~\ref{thm:r14-divisor-examples}, the ambient family is
$M_m\times V$, $F=0$, and the coordinate and bar comparison maps
are independent of $t$. Those loci therefore satisfy this
part of the connection construction globally. A nonzero product
on a prescribed fourfold quotient is subject to the separate
coefficient test of Proposition~\ref{prop:v5-factor-test}.

We now fix $q:\mathscr U\to V$ and $F$ as in
\eqref{eq:r14-relative-moduli-chart}. The following calculation
also applies to a specified smooth family with an ordinary
potential independently of a moduli interpretation.

The relative twisted de Rham differential is
\begin{equation}\label{eq:v5-relative-DR-differential}
 D=u\,d_{\mathscr U/V}-d_{\mathscr U/V}F\wedge
        \quad\hbox{on}\quad\Omega^\bullet_{\mathscr U/V}((u)).
\end{equation}
The variable $u$ is even after periodization. A Laurent series in
$u$ has a lower bound on its exponents; this bound is required for
all formulas below. The operators we use lower the $u$-order by at
most a fixed finite amount. The ordinary-function comparison of
\cite{Efimov2018CyclicMF} applies fiberwise when its hypotheses,
including the condition on the critical values, are satisfied.
The Gauss--Manin calculations along the base curve are made directly on
\eqref{eq:v5-relative-DR-differential}.

\subsubsection{Horizontal lifts and Cartan homotopies}
A lift of $\partial_t$ need not exist globally on the smooth family.
We compare two local lifts first: the contraction by their difference
will provide the correction on a Cech overlap.

A lift $v$ of $\partial_t$ is a derivation of $\O_{\mathscr U}$
with $v(q^*a)=q^*(\partial_ta)$ for base functions $a$.
Its Lie derivative on relative forms is characterized by
\[
 \mathcal L_v(a\,d b_1\wedge\cdots\wedge d b_p)
 =v(a)\,d b_1\wedge\cdots\wedge d b_p
  +\sum_j a\,d b_1\wedge\cdots\wedge d(vb_j)
                      \wedge\cdots\wedge d b_p.
\]
Every $d$ in this formula is relative to $V$. Define
\begin{equation}\label{eq:v5-lift-connection}
 \nabla_v=\mathcal L_v-u^{-1}v(F).
\end{equation}
The final term is multiplication, not contraction by a vector field.
The Leibniz rule is
$\nabla_v(a\alpha)=(\partial_ta)\alpha+a\nabla_v\alpha$
for a function $a$ on $V$.

\begin{lemma}\label{lem:v5-Cartan}
For a lift $v$ and a vertical vector field $w$, the following identities
hold as graded operators on relative forms:
\begin{align}
 D^2&=0,\label{eq:v5-D-square}\\
 [\nabla_v,D]&=0,\label{eq:v5-GM-commutator}\\
 [D,u^{-1}\iota_w]&=\mathcal L_w-u^{-1}w(F),
                      \label{eq:v5-Cartan-homotopy}\\
 \nabla_{v'}-\nabla_v&=[D,u^{-1}\iota_{v'-v}]
           \quad\text{when }v',v\text{ lift the same base derivation}.
                      \label{eq:v5-change-lift}
\end{align}
If $v_a,v_b$ lift commuting base derivations in a higher-dimensional
base, their commutator on twisted forms is
\[
 [\nabla_{v_a},\nabla_{v_b}]
       =[D,u^{-1}\iota_{[v_a,v_b]}].
\]
\end{lemma}
\begin{proof}
Put $\beta=d_{\mathscr U/V}F$. Since $d\beta=0$ and
$\beta\wedge\beta=0$, expansion gives
\[
 (u d-\beta\wedge)^2
  =u^2d^2-u(d\beta)\wedge+\beta\wedge\beta\wedge=0.
\]
The Lie derivative commutes with $d$ and sends $dF$ to $d(vF)$.
Consequently
\[
 [\mathcal L_v,D]=-d(vF)\wedge.
\]
For a function $a$, the rule $d(a\alpha)=da\wedge\alpha+a d\alpha$
gives $[m_a,u d]=-u\,da\wedge$, where $m_a$ is multiplication by
$a$. Taking $a=-u^{-1}vF$ gives
$[m_a,D]=d(vF)\wedge$. This cancels the preceding term and proves
\eqref{eq:v5-GM-commutator}.

Both $D$ and $\iota_w$ are odd. On a form $\alpha$,
\[
 \iota_w(\beta\wedge\alpha)
       =w(F)\alpha-\beta\wedge\iota_w\alpha.
\]
Together with $d\iota_w+\iota_wd=\mathcal L_w$, this yields
\[
 D\iota_w+\iota_wD=u\mathcal L_w-w(F),
\]
proving \eqref{eq:v5-Cartan-homotopy}. Subtracting the two
expressions \eqref{eq:v5-lift-connection} proves
\eqref{eq:v5-change-lift}.

Finally, multiplication operators commute, so
\[
 [\mathcal L_{v_a}-u^{-1}v_aF,
   \mathcal L_{v_b}-u^{-1}v_bF]
   =\mathcal L_{[v_a,v_b]}-u^{-1}[v_a,v_b]F.
\]
When the base derivations commute the vector field $[v_a,v_b]$ is
vertical. Applying \eqref{eq:v5-Cartan-homotopy} proves the curvature
identity.
\end{proof}

The homotopies \eqref{eq:v5-change-lift} define the overlap
correction in the following Cech total complex.

\subsubsection{The connection on a finite Cech complex}
The Cartan identity gives a homotopy on each double overlap. We now
insert those homotopies into the total complex so that the connection
commutes with both the internal and Cech differentials.

Choose a finite ordered affine cover $\{U_\alpha\}$ of $\mathscr U$
and a lift $v_\alpha$ on each member. Smoothness permits these local
lifts. Refine the cover if necessary to have affine multiple
intersections. The Cech model used here is obtained at finite orders
in $u$ from the relative differential forms, followed by the derived
$u$-adic inverse limit and inversion of $u$. On affine intersections
the coefficient sheaves at every finite order are acyclic and the
transition maps on sections are surjective; hence this model computes
the corresponding completed derived global sections.

Write $c^{p,q}$ for a Cech cochain of Cech degree $p$ and relative
form degree $q$, with Laurent coefficients after that completion.
Its total parity is $p+q$. Let $\check\delta$ be the usual alternating
Cech differential, and put
\begin{equation}\label{eq:v5-Cech-D}
 \mathbb D c^{p,q}=\check\delta c^{p,q}+(-1)^pDc^{p,q}.
\end{equation}
The two terms square to zero and their cross terms cancel, because
restriction of forms commutes with $D$.
On overlaps put
\[
 h_{\alpha\beta}=u^{-1}\iota_{v_\beta-v_\alpha}.
\]
Then
\begin{equation}\label{eq:v5-Cech-h-identities}
 [D,h_{\alpha\beta}]=\nabla_{v_\beta}-\nabla_{v_\alpha},\qquad
 h_{\beta\gamma}-h_{\alpha\gamma}+h_{\alpha\beta}=0.
\end{equation}
The second equation follows from addition of vector fields and
linearity of contraction.

Define operators $N_0$ and $N_1$ by
\begin{align}
 (N_0c)_{\alpha_0\ldots\alpha_p}
     &=\nabla_{v_{\alpha_0}}c_{\alpha_0\ldots\alpha_p},
                             \label{eq:v5-Cech-N0}\\
 (N_1c)_{\alpha_0\ldots\alpha_{p+1}}
     &=(-1)^p h_{\alpha_0\alpha_1}
                         c_{\alpha_1\ldots\alpha_{p+1}}.
                             \label{eq:v5-Cech-N1}
\end{align}
Here $N_1$ raises Cech degree by one and lowers form degree by one,
so it has total degree zero.

\begin{theorem}\label{thm:v5-Cech-GM}
Let $q:\mathscr U\to V\subset C$ be smooth and separated of finite
type, and let $F\in\Gamma(\mathscr U,\O_{\mathscr U})$.
Choose the affine cover and local lifts of $\partial_t$ above,
and define $N_0,N_1$ by \eqref{eq:v5-Cech-N0} and
\eqref{eq:v5-Cech-N1}. Then the operator
\begin{equation}\label{eq:v5-Cech-connection}
 \mathbb\nabla_{\partial_t}=N_0+N_1
\end{equation}
is a connection on the Cech total complex along $V$ and satisfies
\[
 [\mathbb\nabla_{\partial_t},\mathbb D]=0.
\]
Its action on cohomology is independent of the chosen horizontal lifts
and is compatible with restriction in the base. These actions define
a connection along the smooth curve $C$ on the corresponding
cohomology sheaves. Its curvature is zero. The assertion does not
include finite generation of those sheaves.
\end{theorem}
\begin{proof}
The operator $N_0$ satisfies the connection Leibniz rule, and $N_1$
is linear over base functions. Their sum is therefore a connection.
We check the chain equation on a homogeneous Cech cochain.
The internal commutator with $N_0$ is zero by
\eqref{eq:v5-GM-commutator}. In
$\check\delta N_0-N_0\check\delta$, all terms except the term
omitting the first index cancel. The remaining term is
\[
 (\nabla_{v_{\alpha_1}}-\nabla_{v_{\alpha_0}})
                    c_{\alpha_1\ldots\alpha_{p+1}}
 =[D,h_{\alpha_0\alpha_1}]
                    c_{\alpha_1\ldots\alpha_{p+1}}.
\]
The commutator of the internal part of $\mathbb D$ with $N_1$ is
\[
 (-1)^{p+1}D((-1)^phc)
       -(-1)^ph((-1)^pDc)=-(Dh+hD)c.
\]
It cancels the preceding expression. Lastly,
\[
 (\check\delta N_1-N_1\check\delta)c
      =(-1)^p(\check\delta h)\smile c=0
\]
by the second identity in \eqref{eq:v5-Cech-h-identities}.
This proves the chain equation with all totalization signs included.

To compare lifts, write $v'_\alpha-v_\alpha=w_\alpha$ and put
\[
 (Jc)_{\alpha_0\ldots\alpha_p}
      =(-1)^p u^{-1}\iota_{w_{\alpha_0}}
                              c_{\alpha_0\ldots\alpha_p}.
\]
The operator $J$ has total degree $-1$. Its internal commutator
with $\mathbb D$ gives
$\nabla_{v'_{\alpha_0}}-\nabla_{v_{\alpha_0}}$ by
\eqref{eq:v5-Cartan-homotopy}. Its Cech commutator gives
\[
 (-1)^pu^{-1}\iota_{w_{\alpha_1}-w_{\alpha_0}}
                         c_{\alpha_1\ldots\alpha_{p+1}},
\]
which is the change in $N_1$. Thus
\begin{equation}\label{eq:v5-Cech-independence}
 \mathbb\nabla'-\mathbb\nabla=\mathbb D J+J\mathbb D.
\end{equation}
This identity proves equality on cohomology. The formulas commute
with restriction, and a common refinement compares two covers through
the same identities. Multiplying a base derivation by a base function
multiplies the induced covariant derivative by that function.
Consequently the local operators transform as a connection under
changes of the coordinate on the base. Since a smooth curve has
$\Omega_C^2=0$, the curvature of this connection is zero. This proves
flatness on cohomology, without assuming that the cohomology is a
finite-rank vector bundle.
\end{proof}

\begin{corollary}\label{cor:v5-product-GM}
Let $q:\mathfrak U\to V\subset C$ be a smooth separated family
of finite type. Suppose it carries a sheaf of dg algebras locally
of the form
$B=\O_{\mathfrak U}[\theta_1,\ldots,\theta_r]$,
$|\theta_a|=-1$, with internal differential zero on
$\O_{\mathfrak U}$. Let $\nu,g\in B^{-1}$ be global closed
sections, and suppose $g$ is primitive. Assume on a cover that
$\nu=f_\alpha g$ for ordinary functions $f_\alpha$.
Then the $f_\alpha$ glue to a global function $F$ on
$\mathfrak U$. For $\kappa(g)=1$ as in
Proposition~\ref{prop:v5-extraction}, the operator on a smooth
ambient presentation $U\times V\to V$
\begin{equation}\label{eq:v5-product-connection}
 \nabla_{\partial_t}
   =\partial_t-\frac1u
       \bigl(\kappa(\partial_t\nu)-F\kappa(\partial_tg)\bigr)
\end{equation}
is independent of $\kappa$, commutes with
$D_F=u\,d_x-d_xF\wedge$, and has the global Cech extension
\eqref{eq:v5-Cech-connection}. Its action on cohomology is a
flat Gauss--Manin connection along $V$.
\end{corollary}
\begin{proof}
On an overlap,
$(f_\alpha-f_\beta)g=0$.
Write $g=\sum_ab_a\theta_a$ in a local frame and choose $c_a$
with $\sum_ac_ab_a=1$. Comparing coefficients and multiplying
by $c_a$ gives
\[
 f_\alpha-f_\beta
   =\sum_ac_a(f_\alpha-f_\beta)b_a=0.
\]
Thus the $f_\alpha$ define $F$.
The local operator $\kappa=\sum_ac_a\partial_{\theta_a}$ satisfies
$\kappa F=0$ and $\kappa g=1$. Differentiating $\nu=Fg$ gives
\[
 \partial_t\nu=(\partial_tF)g+F\partial_tg,
 \qquad
 \kappa(\partial_t\nu)-F\kappa(\partial_tg)=\partial_tF.
\]
If $\kappa'$ is another contraction with $\kappa'g=1$,
subtracting the two formulas gives zero. The expression is
therefore the same function $\partial_tF$ in every local frame.
A change of frame depending on $t$ changes both differentiated
odd expressions, but their combination is still $\partial_tF$
by this identity.

Consequently \eqref{eq:v5-product-connection} is
$\partial_t-u^{-1}\partial_tF$.
On relative forms its commutator is
\[
 \begin{aligned}
 [\partial_t,D_F]&=-d_x(\partial_tF)\wedge,\\
 [-u^{-1}\partial_tF,D_F]&=d_x(\partial_tF)\wedge,
 \end{aligned}
\]
so their sum vanishes. For two commuting base coordinates the
curvature vanishes by Proposition~\ref{prop:parameter-GM}.
On a general smooth ambient family $q:\mathfrak U\to V$, choose
local lifts $v_\alpha$ of $\partial_t$ and put
$h_{\alpha\beta}=u^{-1}\iota_{v_\beta-v_\alpha}$.
The local differences are
\[
 \nabla_\beta-\nabla_\alpha=[D_F,h_{\alpha\beta}],
 \qquad
 h_{\beta\gamma}-h_{\alpha\gamma}+h_{\alpha\beta}=0.
\]
Theorem~\ref{thm:v5-Cech-GM}, with these specified homotopies,
constructs $\mathbb\nabla=N_0+N_1$ and proves
$[\mathbb\nabla,\mathbb D]=0$.
It also proves independence of the lifts by the explicit
homotopy $\mathbb\nabla'-\mathbb\nabla=\mathbb DJ+J\mathbb D$.
The resulting sheaf connection on the curve has curvature in
$\Omega^2_V=0$. This proves its flatness globally.
\end{proof}

This corollary constructs the flat connection when the stated local product data are given. Proposition~\ref{prop:v5-factor-test} is the independent test
for existence of those data on a fourfold quotient. Finite rank and
the identification of the geometric special fiber are not assumptions
of the chain calculation and do not follow from its flatness alone.
Theorem~\ref{thm:finite-cyclic} deals with finite models, and
Section~\ref{sec:v5-bar-bridge} constructs the central bar comparison
in an explicit relative critical-chart family.

\subsection{Local freeness of the cohomology of a finite periodic complex}
The global construction gives a connection on cohomology. To use it
at the distinguished fiber we require a finite model with a connection
regular at $0$. The next lemma explains the role of regularity; the
following theorem proves the base-change assertion used in
Theorem~\ref{thm:intro-cyclic}.

\begin{lemma}\label{lem:connection-free}
Let $R$ be the coordinate ring of a smooth curve $D=\Spec R$ over a characteristic-zero field $K$. In the geometric application $D=V_K$, for $V\subset C$ and $K=\CC((u))$ as in Theorem~\ref{thm:intro-cyclic}. A finite $R$-module with a $K$-connection is locally free.
\end{lemma}
\begin{proof}
It suffices to work at a closed point, where the local ring is a discrete valuation ring. Let $M_{\mathrm{tor}}$ be the torsion submodule. If $a m=0$ for a nonzero $a$, the Leibniz identity gives $a\nabla m=-da\otimes m$, and hence $a^2\nabla m=0$. Thus torsion is preserved by the connection.

Choose an etale local coordinate $z$ at the point and use the vector field $\partial_z$, which sends the uniformizer $z$ to $1$ after an etale localization if necessary. The annihilator ideal $I$ of the finite torsion module is preserved by $\partial_z$: for $a\in I$ and $m\in M_{\mathrm{tor}}$,
\[
 0=\nabla_{\partial_z}(am)
 =\partial_z(a)m+a\nabla_{\partial_z}m,
\]
and the last term is zero. If the torsion module were nonzero, its annihilator would be $(z^N)$ for some $N>0$. But $\partial_z(z^N)=Nz^{N-1}\notin(z^N)$ in characteristic zero, a contradiction. Thus $M$ is torsion-free and hence free over the discrete valuation ring. Etale descent gives local freeness on $C$.
\end{proof}

\begin{theorem}\label{thm:finite-cyclic}
Let $R$ be as in Lemma~\ref{lem:connection-free}, and let
\[
 P^{\overline0}\xrightarrow{d_0}P^{\overline1}
 \xrightarrow{d_1}P^{\overline0},\qquad d_1d_0=d_0d_1=0,
\]
be a $2$-periodic complex of finite projective $R$-modules. Suppose it carries a degree-preserving regular $K$-connection $\nabla$ satisfying $[\nabla,d]=0$. Then:

\noindent (a) $H^{\overline0}(P)$ and $H^{\overline1}(P)$ are finite locally free.\par
\noindent (b) For every point $q\in D=\Spec R$ with residue field $\kappa(q)$, the canonical maps
\[
 H^{\bar i}(P)\otimes_R \kappa(q)\xrightarrow{\sim}
 H^{\bar i}(P\otimes_R \kappa(q))
\]
are isomorphisms.\par
\noindent (c) The induced connections are flat and their ranks are constant on each connected component.\par
\noindent (d) If $P\otimes_R \kappa(q)$ is identified with $\CyclicChains_{\per}(\mathscr C_q)$ for specified fiber categories, these conclusions apply to $\HP(\mathscr C_q)$.
\end{theorem}
\begin{proof}
The two cohomology modules are finite because $R$ is noetherian and the terms are finite. The chain connection descends to them, so Lemma~\ref{lem:connection-free} proves (a). Work locally on the curve. Images and kernels of the differentials are finite torsion-free submodules of free modules, so are free. The sequences
\[
 0\to\ker d_i\to P^{\bar i}\to\operatorname{im}d_i\to0,
 \qquad
 0\to\operatorname{im}d_{i-1}\to\ker d_i\to H^{\bar i}(P)\to0
\]
split locally. Tensoring by any residue field preserves them, proving (b). Curvature takes values in $\Omega^2_{R/K}=0$, so the connections on a smooth curve are flat. A locally free sheaf has constant rank on connected components, giving (c). Part (d) is the given fiberwise chain comparison combined with (b).
\end{proof}

\begin{hypothesis}\label{hyp:finite-cyclic}
Let $\pi:\cX\to C$ be the geometric degeneration, let $V\subset C$ be a connected affine neighborhood of $0$, and put
\[
 K=\CC((u)),\qquad R=\Gamma(V,\O_C)\otimes_{\CC}K.
\]
For each closed point $b\in V\subset C$, let $K_b$ denote evaluation at $b$ as an $R$-module. Fix a $2$-periodic complex $P$ of finite projective $R$-modules with a regular chain connection along $V_K$.

For $b\ne0$, let $\mathscr C_b$ be the specified category for the
smooth fiber $X_b$, with its local matrix-factorization presentations
and compatible descent data. Let $\mathsf K_b^{\mathrm{DT}}$
be its chosen normalized cyclic realization, including the
critical-chart suspensions and orientation comparisons. On a
global critical chart of dimension $d$ it is
$\CyclicChains_{\per}(\mathscr C_b)[d]$. For nonproper loci a
compact-support realization and its comparison are required when
compactly supported DT is used. We assume quasi-isomorphisms
\[
 \iota_b:P_b=P\otimes_R K_b
       \xrightarrow{\sim}\mathsf K_b^{\mathrm{DT}}.
\]
Fix the oriented gluing locus $M_0$ on the special fiber and its periodic DT complex $\mathcal B_{0,\mathrm{DT}}$. In the smooth quantization setting this is the periodic bar complex defined in \cite[equation~(6.3)]{KS}. Put
\[
 \mathsf H_0^{\mathrm{DT}}
      =\R\Gamma_c(M_0,\mathcal B_{0,\mathrm{DT}}).
\]
For the chosen locus in $\fM_1\timesh_{\RPerf(S)}\fM_2$, where $\fM_i$ is a moduli stack on $X_i$, we assume a central-fiber quasi-isomorphism
\begin{equation}\label{eq:central-bar-identification}
 \iota_0:P\otimes_R K_0
       \xrightarrow{\sim}\mathsf H_0^{\mathrm{DT}}.
\end{equation}
In the relative critical-chart model, Theorem~\ref{thm:v5-periodic-bar} constructs its local counterpart explicitly as $\mathcal I\Pi_0$. A global identification in the present hypothesis also includes the specified presentation and descent maps for the moduli stacks on $X_1$ and $X_2$. The smooth-fiber comparisons are compatible with the Gauss--Manin connection on the family of periodic cyclic complexes. Thus $P$ is one complex over the base with the indicated descriptions of its fibers, including \eqref{eq:central-bar-identification}. The regularity assumption means that its connection is defined at $0$. In the formal version replace $R$ by $K[[t]]$, where $t$ is a local parameter at $0\in C$, and use a continuous connection.
\end{hypothesis}
Theorem~\ref{thm:r14-divisor-examples} and
Proposition~\ref{prop:v15-three-hypotheses} verify these data for
$E^{m,k}_{s,b}$, $m=1,2,3$: the finite complex is concentrated in even
parity, with ranks $5,15,35$, zero differential, and connection
$\partial_t$. The central map is the explicit bar lift of the scaled
hyperplane-power cocycles. The parity convention in this hypothesis
must be respected when using another critical chart: on an ambient
chart of dimension $n$, the normalized DT bar complex corresponds to
$\CyclicChains_{\mathrm{per}}(\MF(U,f))[n]$, as in
\cite[Proposition~6.6]{KS}. In the three displayed divisor examples
$n=d_m$ is even. An odd-dimensional presentation requires the
corresponding suspension of the periodic complex in a normalized
comparison; it cannot be identified with the unshifted cyclic complex
by an even chain map solely by choosing an orientation.

\subsubsection{Verification for sheaves on divisors in a fixed linear system}
\label{sssec:r14-hypothesis-examples}
The central comparison \eqref{eq:central-bar-identification} is a
statement about the chosen moduli locus $M_0$ and its DT complex. We give
explicit families of the two-dimensional sheaves considered in
this paper for which all its terms and maps can be constructed.

Use the resolved quintic family and the line bundle
$L=a^*\O_{\P^4}(1)$ of \cite[Proposition~2.10]{KS}.
For $m=1,2,3$ and an integer $k$, set
\[
 V_m=H^0(\P^4,\O_{\P^4}(m)),\quad M_m=\P(V_m),
 \quad r_m=\dim V_m=\binom{m+4}{4},\quad d_m=r_m-1.
\]
At a nonzero section $s\in V_m$ define the perfect complex
\begin{equation}\label{eq:r14-divisor-family}
 E^{m,k}_{s,b}=
       [L_b^{k-m}\xrightarrow{s}L_b^k],\qquad |L_b^{k-m}|=-1.
\end{equation}
For $b\ne0$, the threefold $X_b$ is smooth and integral, so a
nonzero $s$ is a non-zero-divisor. Writing $D_{s,b}=(s=0)$ gives
\[
 E^{m,k}_{s,b}\simeq j_{D*}(L_b^k|_{D_{s,b}}),
 \qquad D_{s,b}\in|L_b^m|.
\]
This is a pure two-dimensional sheaf. On the open locus of
integral divisors it is Gieseker stable: a nonzero subsheaf has
generic rank one on its integral support, and a proper quotient
has smaller dimension and positive Hilbert polynomial. The
full projective family \eqref{eq:r14-divisor-family} is taken
in the moduli of simple perfect complexes. 
Theorem~\ref{thm:r14-divisor-examples} proves that the classical
scalar-rigidified locus is $M_m\times C$. Its critical structure
is zero, with orientation $K_{M_m}$. The examples and the finite
ranks obtained there are
\[
\begin{array}{c|c|c|c}
 m& D_{s,b}\text{ in }X_b& M_m &
          \operatorname{rank}_{\CC((u))}\HP(\MF(M_m,0))\\ \hline
 1&\text{hyperplane section}&\P^4&5\\
 2&\text{quadric section}&\P^{14}&15\\
 3&\text{cubic section}&\P^{34}&35.
\end{array}
\]
The last row concerns degree-three divisors of the quintic; it
is not an auxiliary polynomial singularity example. Tensoring
by $L_b^k$ gives the indicated line-bundle sheaves on the same
divisors. Tensoring by the extending line bundle $L^k$ is an
exact equivalence on the ambient and fiberwise perfect-complex
categories. On endomorphism complexes it sends
$a\mapsto a\otimes\mathrm{id}_{L^k}$ and preserves composition,
trace, and derived restriction. Hence it also identifies the
cyclic deformation data and the chosen orientations.

Here is the central map in these examples. Let
$\mathcal B_{m,\mathrm{per}}^{\mathrm{or}}$ be the periodic
oriented bar complex of the zero critical presentation on $M_m$.
The maps in
Proposition~\ref{prop:v6-global-bar} give
\begin{equation}\label{eq:r14-actual-central-diagram}
\begin{tikzcd}[column sep=small,row sep=large]
 R\Gamma(M_m,\mathcal B_{m,\mathrm{per}}^{\mathrm{or}})
     \arrow[r,"\Pi"]&
 R\Gamma(M_m,\mathcal K_m^{\mathrm{or}}((u)))
     \arrow[r,"\mathcal I"]&
 R\Gamma(M_m,\Omega_{M_m}^{\bullet}((u))[d_m]).
\end{tikzcd}
\end{equation}
The first arrow is a quasi-isomorphism and has the chain section
$I$; the second is contraction with top forms, with the signs
of \eqref{eq:v6-contraction} in dimension $d_m$. The inverse
arrow in the derived category is represented by
$I\mathcal I^{-1}$.

Fix an affine open $V\subset C$ and retain the ring $R$ from
Hypothesis~\ref{hyp:finite-cyclic}. All $d_m$ in the table are
even. Take the finite periodic complex
\begin{equation}\label{eq:r14-explicit-finite-P}
 P^{\overline0}=\bigoplus_{j=0}^{d_m}R e_j,
 \qquad P^{\overline1}=0,\qquad d_P=0,\qquad
 \nabla e_j=0.
\end{equation}
Let $z_j$ be the Cech--de Rham cocycle for $H^j$, where
$H=c_1(\O_{M_m}(1))$, and put $Z_j=u^jz_j$ for the differential
$u\,d$. Equation \eqref{eq:central-bar-identification} is represented
in these models by
\begin{equation}\label{eq:r14-iota0-basis}
 \iota_0(e_j)=(I\mathcal I^{-1})(Z_j)|_{b=0}.
\end{equation}
These are cocycles, their classes are a basis, and the defining
maps are independent of $t$. The finite-contraction and
base-change proof is given in
Theorem~\ref{thm:r14-divisor-examples}. In a nonzero fiber,
Efimov's comparison for the chosen presentation $(M_m,0)$
identifies the same complex with
$\CyclicChains_{\mathrm{per}}(\MF(M_m,0))$. A morphism such as
$\iota_b$ may be represented by a zigzag of chain
quasi-isomorphisms, or by a single map after replacing the target
by that quasi-isomorphic model. All references to these categories
mean the categories of these specified global critical
presentations.

The coefficient ring $\Gamma(V,\O_C)[[u]][u^{-1}]$ used in
completed de Rham forms is not identified with $R$.
The finite complex \eqref{eq:r14-explicit-finite-P} is already
defined over $R$; extension of its scalars gives the completed
finite representative. Both constructions specialize to
$K_b=\CC((u))$ with the same basis and maps. Properness of $M_m$
identifies ordinary and compactly supported hypercohomology.

\subsubsection{The scope for other sheaves supported on divisors}
For a general pure two-dimensional sheaf $G$ with support
$D\in|L_b^m|$, the two-term resolution in
\eqref{eq:r14-divisor-family} is not part of the definition of
$G$. Thus the projective-space reconstruction is a verified
case of Hypothesis~\ref{hyp:finite-cyclic}, not a consequence of
purity or of the support class alone. To apply that hypothesis
to another locus, one specifies its critical family, proves
finiteness and base change for its periodic complex, and
constructs \eqref{eq:central-bar-identification} for its own
matching triples on $X_1,X_2$.

For example, if $D$ is smooth and $G=j_{D*}A$ for a line bundle
$A$ on $D$, deforming the pair $(D,A)$ in a direction $v$ has
line-bundle obstruction
\[
 \operatorname{At}(A)\cdot\kappa_D(v)
                      \in H^2(D,\O_D).
\]
Its image under pushforward contributes to the obstruction
$\operatorname{At}(G)\cdot\kappa_{X_b}(v)$ in the normal
triangle \cite{HT}. An arbitrary $A$ need not be a restriction
of a line bundle on the total family. For
$A=L_b^k|_D$, the explicit family
\eqref{eq:r14-divisor-family} supplies a lift, and its
Atiyah--Kodaira--Spencer obstruction vanishes. This is a geometric
reason why the examples above satisfy the lifting requirement;
no such vanishing is deduced for every sheaf on $D$.

\begin{corollary}\label{cor:categorical-specialization}
Under Hypothesis~\ref{hyp:finite-cyclic}, the cohomology of $P$ is a finite-rank bundle with regular flat connection on $V_K$. Its fibers are
\[
 H^{\bar i}(P)\otimes_R K_b\simeq
 \begin{cases}
 H^{\bar i}(\mathsf K_b^{\mathrm{DT}}),&b\ne0,\\
 \mathbb H_c^{\bar i}(M_0,\mathcal B_{0,\mathrm{DT}}),&b=0.
 \end{cases}
\]
For a proper $M_0$, ordinary hypercohomology gives the same central-fiber group. In the formal version, horizontal trivialization with prescribed initial value gives
\[
 H^{\bar i}(P)\simeq
 \mathbb H_c^{\bar i}(M_0,\mathcal B_{0,\mathrm{DT}})
       \otimes_K K[[t]].
\]
\end{corollary}

\begin{proof}
Use Theorem~\ref{thm:finite-cyclic} for the locally free cohomology and base-change maps. Composing them with $\iota_b$ gives the displayed smooth-fiber groups; composing at $0$ with \eqref{eq:central-bar-identification} gives the central-fiber group. The same torsion argument applies over $K[[t]]$. Choose a basis of the resulting free module, so the connection is $\partial_t+A(t)$. Solve $G'=-A(t)G$, $G(0)=I$, coefficient by coefficient. If $A=\sum A_jt^j$ and $G=\sum G_nt^n$, the recursion
\[
 (n+1)G_{n+1}=-\sum_{j=0}^{n}A_jG_{n-j}
\]
has a unique solution in characteristic zero. The constant term is invertible, so $G$ is invertible and gives the horizontal frame.
\end{proof}

The consequence of construction of regular connection is a derived-geometric realization of
deformation invariance \emph{for the specified categorified DT
problem}. Over the formal disk, every horizontal cohomology class is
uniquely determined by its special-fiber value. In particular,
transport gives an equivalence of the periodic cohomological vector
spaces provided that the quasi-isomorphisms $\iota_b$ and $\iota_0$
identify the two geometric endpoints. Flatness of the chain
connection does not itself manufacture the special-fiber DT
complex or establish the descent maps. This distinction is necessary
when moving beyond the global zero-potential hyperplane-section
example to arbitrary stable coherent sheaves. It is the input for
the categorical degeneration identity in
Section~\ref{sec:categorical-degeneration}.

\begin{remark}\label{rem:completion-order}
The horizontal solution above is taken in $K[[t]]=\CC((u))[[t]]$.
The order of the two completions and the exponential transport formula
are specified in Section~\ref{sec:v5-formal-transport}. For an analytic
finite-rank local system, parallel transport instead uses a chosen
path; the formal initial-value construction fixes the comparison
near $0$.
\end{remark}

For the DT application the theorem is used with the chain maps
$\iota_b$ and $\iota_0$ of Hypothesis~\ref{hyp:finite-cyclic}. Their
construction is therefore the remaining geometric input:
Section~\ref{sec:v5-bar-bridge} supplies the local horizontal maps, and
Section~\ref{sec:v6-hyperplanes} supplies their global realization for
the hyperplane-section locus.

\subsubsection{Compatibility with arbitrary base change}
The local splitting argument in Theorem~\ref{thm:finite-cyclic} actually shows that, locally on the curve, the finite periodic complex is the direct sum of its cohomology with contractible two-term summands. Indeed split the image and kernel sequences in the proof; the differential is then an isomorphism from a chosen complement of its kernel onto its image. It follows that its cohomology commutes with \emph{any} scalar extension on that local open, not only residue-field extensions. The splitting need not preserve the connection, but this does not affect the derived base-change conclusion.

This observation can be used when comparing two finite cyclic models. A chain quasi-isomorphism between them induces an isomorphism of the locally free cohomology bundles. If it intertwines the connections up to a chain homotopy, the induced isomorphism on cohomology intertwines the connections exactly. Thus the finite-model theorem is invariant under suitable homotopical changes of model; it does not depend on a particular free resolution.

\subsection{Formal transport and the order of completion}\label{sec:v5-formal-transport}
Corollary~\ref{cor:categorical-specialization} gives horizontal
specialization for the finite model. At the level of a potential one
can also write a formal transport operator explicitly. We prove that
formula with its coefficient ring. The geometric applications use the
finite comparison maps of Theorem~\ref{thm:r14-divisor-examples}.

There is an explicit horizontal transport operator, but its coefficient
ring must be stated. Put $K=k((u))$ and let $A$ be the coordinate
algebra of a fixed smooth affine chart over $k$. Write
$\Omega_A^\bullet((u))$ for the ordinary Laurent-completed complex.
For a series $F(t)\in A[[t]]$ put $F_0=F(0)$ and consider
\begin{equation}\label{eq:v5-t-completed-complex}
 \widehat K_t=
       \bigl(\Omega_A^\bullet((u))[[t]],
                    u d-dF(t)\wedge\bigr).
\end{equation}
Thus the $t$-adic completion is taken \emph{after} forming Laurent
series in $u$. The lower bound on $u$-exponents may depend on the
coefficient of $t^j$.

\begin{theorem}\label{thm:v5-formal-transport}
Let $a(t)=F(t)-F_0\in tA[[t]]$. Multiplication by
\begin{equation}\label{eq:v5-exponential-transport}
 T(t)=\exp(a(t)/u)=\sum_{m\geq0}\frac{a(t)^m}{m!u^m}
\end{equation}
is an invertible operator on the completed coefficient module in
\eqref{eq:v5-t-completed-complex}. It satisfies
\begin{equation}\label{eq:v5-transport-identities}
 D_tT=TD_0,\qquad
 (\partial_t-u^{-1}\partial_tF)T=T\partial_t,\qquad T(0)=1.
\end{equation}
In particular it induces a horizontal isomorphism
\[
 H^\bullet(\widehat K_t)
      \simeq H^\bullet(\Omega_A^\bullet((u)),D_0)[[t]].
\]
If the latter $K$-vector space is finite dimensional, this is a finite
free $K[[t]]$-module. The assertion concerns the indicated completed
family, not an interchange of the $u$- and $t$-adic completions.
\end{theorem}
\begin{proof}
Since $a(t)^m\in t^mA[[t]]$, for each power $t^j$ only the summands
$m\leq j$ occur. Each coefficient of $t^j$ is therefore a finite
sum of Laurent monomials in $u$ with coefficients in $A$.
Multiplication by this series is well defined coefficientwise in
$t$, even on an input whose coefficients already belong to
$\Omega_A^\bullet((u))$. Replacing $a$ by $-a$ gives its inverse.
Termwise differentiation is valid in this completion and gives
\[
 dT=u^{-1}(dF-dF_0)T,\qquad
 \partial_tT=u^{-1}(\partial_tF)T.
\]
For a relative form $\alpha$ we compute
\begin{align*}
 D_t(T\alpha)
 &=u(dT)\wedge\alpha+uT d\alpha-dF\wedge T\alpha\\
 &=T\bigl(u d\alpha-dF_0\wedge\alpha\bigr)=TD_0\alpha.
\end{align*}
The second differentiation identity gives the horizontal equation
in \eqref{eq:v5-transport-identities}. The constant term is one.
Thus $T$ identifies $\widehat K_t$ with the coefficientwise constant
complex $K_0[[t]]$. Kernels in that complex are coefficientwise
kernels. Boundaries are also coefficientwise boundaries: for a
series of boundaries choose a preimage of each coefficient in the
$K$-vector-space complex $K_0$. This proves
$H(K_0[[t]])=H(K_0)[[t]]$. A basis of a finite-dimensional $H(K_0)$
then gives the asserted finite free module and horizontal frame.
\end{proof}

The horizontal equation also has a useful recursive form. If
$F(t)=F_0+\sum_{m\geq1}t^mF_m$ and
$T(t)=\sum_{j\geq0}T_jt^j$, then
\begin{equation}\label{eq:v5-transport-recursion}
 T_0=1,\qquad
 (j+1)T_{j+1}=u^{-1}\sum_{m=1}^{j+1}mF_mT_{j+1-m}.
\end{equation}
Each $T_j$ is computed by finitely many additions and multiplications.
The recursion explains both the existence of formal transport and
its possible unbounded poles in $u$ as $j$ grows.

\subsubsection{The two completions and specialization}
The two rings
\[
 A[[t]]((u))\quad\text{and}\quad A((u))[[t]]
\]
have different elements. In the first, a series has a common lower
bound for the powers of $u$, independent of its $t$-degree. In the
second this bound may tend to $-\infty$ with the $t$-degree.
There is a natural coefficientwise inclusion from the first ring to
the second, but it is not surjective. This distinction is material
for the proposed comparison with a geometric generic fiber.

The connection and the finite-specialization theorem now give a precise
target for the special-fiber calculation. It remains to represent
$R\Gamma_c(M_0,\mathcal B_{0,\mathrm{DT}})$ by chain maps compatible with
the ordinary-potential complexes. Section~\ref{sec:quantization}
constructs these maps on the stated Lagrangian and relative critical
presentations. The actual global identification is carried out for the
quintic hyperplane-section locus in Section~\ref{sec:v6-hyperplanes}.

\section{Horizontal bar comparisons and finite models}\label{sec:quantization}
The geometric bar records matching along the divisor. To compute a
DT realization it must carry the orientation and the appropriate
periodic or quantized differential. \cite{KS} constructs the complete
quantization lattices and their finite-order spectral sequences.
Here we use those results as inputs, compute the additional operator
filtration, and construct explicit relative cyclic-to-bar chain maps.

\subsection{Quantization input and its two filtrations}
\label{ssec: Deformation quantization approach}
On a smooth holomorphic symplectic chart with embedded oriented
Lagrangians, the quantization comparison of Gunningham--Safronov
is recalled in \cite[Theorem~6.1]{KS}. Its right and left
quantization modules, denoted $\mathcal N_\hbar$ and
$\mathcal M_\hbar$, compute the DT complex in the specified
Riemann--Hilbert realization after inverting $\hbar$.
For integral calculations we use exactly the lattice hypotheses of
\cite[Hypothesis~6.2]{KS}: the balanced completed tensor product has
a bounded complete $\hbar$-torsion-free representative $C_\hbar$,
and its derived reduction is identified with the classical tensor
product of the two reduced modules. Thus
\begin{equation}\label{eq:DQ-basechange}
 C_\hbar\otL_{\CC[[\hbar]]}\CC
 \simeq\mathcal N_0\otL_{\O_Z}\mathcal M_0.
\end{equation}
For a fixed $N$, \cite[Theorem~6.4]{KS} gives the finite spectral
sequence
\[
 E_1^{p,q}=\mathcal H^{p+q}(C_\hbar/\hbar C_\hbar)\,\hbar^p,
 \quad0\le p<N,
 \qquad E_1\Longrightarrow\mathcal H^{p+q}(C_\hbar/\hbar^N).
\]
Completion is taken by the derived inverse limit. The local Spencer
resolution, including the side of each module action and the signs,
is \cite[Proposition~6.3]{KS}. These results fix the input needed
below; they do not identify a general derived-Lagrangian morphism
with a smooth Lagrangian embedding.
\subsubsection{Conormal symbols and the order filtration}
The $\hbar$-filtration determines successive approximations to the
complex. To write the operators in those approximations, we also need
the action of the ideal of a Lagrangian. The following calculation
identifies that action with differential operators under a local PBW
hypothesis; this is a different filtration from the powers of $\hbar$.

For a line bundle $P$ on a smooth variety $L$, differential operators have the symbol sequences
\begin{equation}\label{eq:Diff-symbol}
 0\longrightarrow\Diff^{\leq k-1}(P)\longrightarrow\Diff^{\leq k}(P)
 \longrightarrow\Sym^kT_L\longrightarrow0 \quad(k\geq1).
\end{equation}
For $k=1$ this is the Atiyah-algebra extension
\begin{equation}\label{eq:Atiyah-line}
 0\longrightarrow\O_L\longrightarrow\Diff^{\leq1}(P)
 \longrightarrow T_L\longrightarrow0.
\end{equation}
Trivializing $P$ in coordinates writes an operator as a finite sum $\sum a_\alpha\partial^\alpha$. Its top coefficients give its principal symbol; those of lower order form the kernel. A change of trivialization does not change the principal symbol, so the sequences glue. They need not split globally.

For the smooth closed embedding $j:L\hookrightarrow Z$ used here, with classical ideal sheaf $I$, the conormal sheaf $I/I^2$ is locally free. Indeed, smoothness of $L$ and $Z$ makes $j$ a regular immersion, so a local regular sequence cutting out $L$ gives a basis of $I/I^2$; see \cite[Tag 0638]{Stacks}. Equivalently, in holomorphic coordinates in which $L$ is defined by $p_1=\cdots=p_n=0$, the classes of the $p_i$ are this basis. The conormal exact sequence is
\[
 0\longrightarrow I/I^2\longrightarrow j^*\Omega_Z^1
         \longrightarrow\Omega_L^1\longrightarrow0.
\]
The symplectic form identifies the normal bundle $N_{L/Z}$ with $T_L^*$ because $L$ is Lagrangian. Dualizing gives $I/I^2=N_{L/Z}^{\vee}\simeq T_L$. These are vector-bundle statements for the smooth closed embedding in this section; a general derived Lagrangian morphism is described by its relative cotangent complex. Let $r:\mathcal A_\hbar\to\O_Z$ be reduction modulo $\hbar$ and set
\[
 J=r^{-1}(I),\qquad \hbar\mathcal A_\hbar\subset J.
\]
The quantization action defines the divided action of $J/J^2$ on $P$.
The class of $\hbar$ supplies the scalar operator in
\eqref{eq:Atiyah-line}; the remaining symbols lift tangent vector fields.
Proposition~\ref{prop:PBW-operators} calculates this action.

\begin{hypothesis}\label{hyp:PBW}
Let $\mathcal A_\hbar$ be a quantization near $L$ and let $\mathcal M_\hbar$ be a flat lattice reducing to a line bundle $P$ on $L$. Let $J$ be the preimage of $I$ under reduction modulo $\hbar$. Assume locally along $L$ that, after completing in $J$, the algebra and its lattice have coordinates and a trivialization with ordered monomials in $p_1,\ldots,p_n,\hbar$, coefficients in $\O_L$, and relations
\[
 [p_i,a]=\hbar\partial_{x_i}a,\qquad
 [p_i,p_j]=0,
 \qquad p_i\cdot m=\hbar\partial_{x_i}m,
\]
with $\hbar$ central and $J=(p_1,\ldots,p_n,\hbar)$. The ordered expansion is unique. Coordinate and line-bundle changes are required to identify the given lattice action. This is a specified local normal-form and lattice hypothesis, not an assertion for every quantization period.
\end{hypothesis}
For $L=M_m$ and $P=\O_{M_m}$, the completed Rees algebra and
its standard left lattice satisfy these formulas intrinsically:
$p_v=\hbar v$ and $p_v\cdot a=\hbar v(a)$. The ordered expansion,
convergence in $J$, and compatibility under coordinate changes are
proved in Proposition~\ref{prop:v15-global-Rees}.

\begin{proposition}\label{prop:PBW-operators}
Under Hypothesis~\ref{hyp:PBW}, division of the action by $\hbar^k$ followed by reduction modulo $\hbar$ induces an isomorphism of left $\O_L$-modules
\begin{equation}\label{eq:PBW-operators}
 J^k/J^{k+1}\xrightarrow{\sim}\Diff^{\leq k}(P)
 \quad(k\geq0).
\end{equation}
These isomorphisms intertwine multiplication in the associated $J$-graded algebra with composition of filtered differential operators. For degree one, the bracket $[\ ,\ ]/\hbar$ corresponds to the commutator of first-order differential operators. Multiplication by the class of $\hbar$ corresponds to the inclusion of operators of order at most $k$ into those of order at most $k+1$.
\end{proposition}
\begin{proof}
The local model gives $J^k\mathcal M_\hbar\subset\hbar^k\mathcal M_\hbar$. Flatness makes multiplication by $\hbar^k$ injective, so division on this image is well defined. An element of $J^{k+1}$ acts with image in $\hbar^{k+1}\mathcal M_\hbar$, hence disappears after division and reduction. Changing a lift of a section of $P$ by a multiple of $\hbar$ also changes the result by a multiple of $\hbar$. This defines the map intrinsically.

By the ordered normal form, the degree-$k$ quotient has a basis of the monomials
\[
 \hbar^{k-|\alpha|}p_1^{\alpha_1}\cdots p_n^{\alpha_n},
 \qquad |\alpha|\leq k.
\]
Their divided action is exactly $\partial^\alpha$. These form a local basis of $\Diff^{\leq k}(P)$, proving bijectivity. There are $\binom{n+k}{k}$ such monomials, which agrees with the rank of the filtered operator bundle, rather than just the rank of its top-symbol quotient.

If $a\in J^k$ and $b\in J^\ell$, their divided actions compose as the divided action of $ab$ because $\hbar$ is central. Reducing gives multiplicativity. For $a,b\in J$, the local relations give $[J,J]\subset\hbar J$. Changing $a$ or $b$ by $J^2$ changes $[a,b]/\hbar$ by $J^2$, so the bracket is defined on $J/J^2$. Its divided action is the commutator of the two first-order operators, as claimed. Finally the extra factor $\hbar$ cancels the extra division, giving the order-filtration inclusion.

The construction used only the intrinsic action on the given lattice. Therefore on overlaps it changes by conjugation under the transition functions of $P$, precisely as differential operators on a line bundle do. The local isomorphisms glue whenever the lattice data in Hypothesis~\ref{hyp:PBW} glue.
\end{proof}

The associated graded here need not be commutative over its degree-zero coefficients: the relation $[p_i,a]=\hbar\partial_i a$ remains homogeneous for the $J$-degree. It is the Rees-type algebra of the \emph{filtered} differential operators, with $\hbar$ recording the order-filtration inclusion. Confusing this with $\Sym(T_L)$ discards the scalar direction and the anchor action.

\subsubsection{Two filtrations, two different spectral sequences}
The $\hbar$-adic filtration on the tensor-product complex in \cite[Theorem~6.4]{KS} has an especially transparent associated graded: every graded piece is a copy of the classical derived intersection, shifted by a power of $\hbar$. The ideal-adic filtration by the quantized ideals $J_i$ is different. In that filtration one assigns degrees to the preimages $J_i$ of the two Lagrangian ideals in the quantization algebra and sums these degrees across the bar factors.

The Poincaré–Birkhoff–Witt (PBW) identification in Proposition~\ref{prop:PBW-operators} gives a Rees algebra of differential operators, including their action on functions. The bar differential uses this multiplication and action. For the spectral-sequence calculation we use the two-sided bar resolution and its $\hbar$-adic filtration from \cite[Theorem~6.4]{KS}; the operator calculation describes the distinct filtration by powers of the quantized ideal.

A fully specified ideal-filtered construction can still be studied. At finite filtration order, form the quotient bar complex, compute its associated graded using the actual multiplication and lattice actions, and take cohomology to obtain its first page. An identification with a Chevalley--Eilenberg or Koszul complex is a separate local quasi-isomorphism to be checked. Globalization must preserve that filtration. The resulting filtration and quasi-isomorphisms are required to be compatible with the transition maps of the quantization modules.

The filtration of \cite[Theorem~6.4]{KS} is
$F^pC=\hbar^pC$. Its first page is determined by
\eqref{eq:DQ-basechange}, and the higher differentials are calculated
by the successive lifting equations below.

\subsubsection{Computing the successive differentials}
Suppose a locally split finite-order model is written
\[
 d=d_0+\hbar d_1+\hbar^2d_2+\cdots
 \pmod{\hbar^N}.
\]
The equation $d^2=0$ gives
\[
 d_0^2=0,\quad d_0d_1+d_1d_0=0,\quad
 d_0d_2+d_2d_0+d_1^2=0,
\]
and the analogous relations in every order. A class $[z_0]\in H(d_0)$ survives the first differential when there is $z_1$ satisfying $d_1z_0+d_0z_1=0$. Its next obstruction is represented by
\begin{equation}\label{eq:second-differential}
 d_2z_0+d_1z_1.
\end{equation}
This is a $d_0$-cocycle because of the displayed equations; changing the lift changes it by the equivalence relation defining the second page. Thus the second differential contains both the second-order deformation and the correction needed to lift the first-order cocycle.

The same recursion continues: to lift $z_0$ through order $r$, solve for $z_1,\ldots,z_r$ successively. The obstruction at each stage is the next spectral-sequence differential. Merely proving that $d_1^2$ is nullhomotopic does not make all those later obstructions vanish. The successive lifting equations specify the higher obstructions.

\begin{proposition}\label{prop:first-order-class}
For a flat completed complex $C$ as in \cite[Hypothesis~6.2]{KS}, the connecting morphism of
\[
 0\to C/\hbar C\xrightarrow{\hbar}C/\hbar^2C
 \to C/\hbar C\to0
\]
is the first Bockstein differential on $H(C/\hbar C)$. Under a graded splitting $d=d_0+\hbar d_1$ modulo $\hbar^2$, it is induced by $d_1$.
\end{proposition}
\begin{proof}
Lift a $d_0$-cocycle $z$ to the middle complex. Its differential lies in the image of the left term and is $\hbar d_1z$ modulo $\hbar^2$. Dividing by $\hbar$ gives the connecting class. Replacing the lift changes this by a $d_0$-boundary. This is both the standard connecting morphism and the first filtered differential.
\end{proof}

Here the first-order deformation is that of the \emph{entire} balanced tensor-product complex. It can involve the algebra multiplication, the right-module action, and the left-module action. An extension class for only one lattice gives the whole Bockstein only when the other data are fixed compatibly and the resulting identification is proved.

The completed and hypercohomological versions use the derived inverse
limit and the support hypotheses of \cite[Section~6.1]{KS}. The explicit
comparisons below are constructed on the completed complexes themselves;
no degeneration of an infinite spectral sequence is needed for those
chain identities.

\subsection{The mixed Koszul model}\label{ssec: SpecialCase}
For a relative function $f$ in coordinates $x_1,\ldots,x_n$, put
\[
 K_f=R[\eta_1,\ldots,\eta_n],\quad |\eta_i|=-1,\qquad
 \delta_f=\sum_i f_i\partial_{\eta_i},\qquad
 \Delta=\sum_i\partial_{x_i}\partial_{\eta_i}.
\]
The identities $\delta_f^2=\Delta^2=\delta_f\Delta+\Delta\delta_f=0$
and the contraction-with-volume comparison are proved in
\cite[Proposition~6.6]{KS}. We fix the resulting periodic differential
$-\delta_f+u\Delta$. The new point in the following construction is
to realize this differential on an actual two-sided bar resolution,
with explicit projection, inclusion, homotopy, and base connection.
The normalization by $[n]$ and its signs are kept in those maps.
\subsection{An explicit cyclic-to-bar comparison on a relative critical chart}\label{sec:v5-bar-bridge}
We now construct the local map required at the central fiber, including
its differential and connection. The proof follows the order
\[
 \begin{aligned}
 \text{bar resolution}&\longleftrightarrow\text{Koszul resolution},\\
 \text{periodic BV complex}&\longleftrightarrow\text{twisted de Rham complex}.
 \end{aligned}
\]
We first build the retraction, then transfer the periodic differential,
and finally correct the bar connection so that both comparison maps
are horizontal. All maps are defined over the base and therefore
specialize together.

Let $R_0$ be the coordinate ring of a smooth affine open $V\subset C$
with a derivation $\partial_t$. Put
\[
 R=R_0[x_1,\ldots,x_n],\qquad A=R[p_1,\ldots,p_n],
          \qquad f\in R,\quad f_i=\partial_{x_i}f.
\]
The scheme $\operatorname{Spec}A$ is the relative cotangent bundle of
$\operatorname{Spec}R$ over $V$. Its zero section and the graph of the
relative differential of $f$ give the two maps
\begin{equation}\label{eq:v5-graph-maps}
 \epsilon_0:A\to R,\quad p_i\mapsto0,
 \qquad
 \epsilon_f:A\to R,\quad p_i\mapsto f_i.
\end{equation}
Use $R_f$ as a right $A$-module via $\epsilon_f$, and $R_0^{\rm sec}$
as a left $A$-module via $\epsilon_0$; the superscript distinguishes
the zero-section module from the base ring $R_0$.
Their derived tensor product represents the relative critical locus:
\begin{equation}\label{eq:v5-graph-intersection}
 R_f\otimes_A^{\mathrm L}R_0^{\rm sec}
       \simeq(R[\eta_1,\ldots,\eta_n],\delta_f),\qquad
 \delta_f=\sum_if_i\partial_{\eta_i}.
\end{equation}
We prove this equivalence together with the chain maps needed after
adding the BV differential.

\subsubsection{The two resolutions and an explicit projection}
Let $\overline A=(p_1,\ldots,p_n)\subset A$ be the augmentation ideal.
The normalized one-sided bar resolution of $R_0^{\rm sec}$ has terms
\[
 C_j=A\otimes_R\overline A^{\otimes_R j},\qquad j\geq0.
\]
Write an element as $a_0[a_1|\cdots|a_j]$. Its differential is
\begin{align}\label{eq:v5-bar-resolution-diff}
 b(a_0[a_1|\cdots|a_j])
 &=a_0a_1[a_2|\cdots|a_j]\notag\\
 &\quad+\sum_{l=1}^{j-1}(-1)^l
       a_0[a_1|\cdots|a_la_{l+1}|\cdots|a_j].
\end{align}
The usual last term is zero since $\epsilon_0(a_j)=0$.
Its augmentation in degree zero is $\epsilon_0:A\to R$.
The Koszul resolution is
\[
 L_j=A\otimes_R\Lambda^j(R^n),\qquad
 d_L=\sum_ip_i\partial_{\eta_i},\qquad |\eta_i|_{\mathrm{hom}}=1.
\]
In this paragraph the indices are homological. Later $C_j,L_j$
are placed in cohomological degree $-j$.

Define the $A$-linear inclusion
\begin{equation}\label{eq:v5-bar-inclusion}
 I(\eta_{i_1}\cdots\eta_{i_j})
  =\sum_{\sigma\in S_j}\operatorname{sgn}(\sigma)
                  [p_{i_{\sigma(1)}}|\cdots|p_{i_{\sigma(j)}}].
\end{equation}
No division by $j!$ occurs in this formula.
For a polynomial $a\in A$ and an auxiliary scalar $s$, write
\[
 da(sp)=\sum_i (\partial_{p_i}a)(sp)\eta_i.
\]
The notation means evaluation of the partial derivative at $sp$;
there is no additional factor $s$. It differentiates in the fiber
variables $p_i$, not in $x_i$ or $t$.
Let $\Sigma_j=\{1\geq s_1\geq\cdots\geq s_j\geq0\}$.
Define
\begin{equation}\label{eq:v5-simplex-projection}
 \Pi(a_0[a_1|\cdots|a_j])
   =a_0\int_{\Sigma_j}
         da_1(s_1p)\wedge\cdots\wedge da_j(s_jp)
                                    \,ds_1\cdots ds_j.
\end{equation}
The integrals mean integration of polynomials in the auxiliary
variables. Their coefficients are rational numbers, so the expression
is defined over every characteristic-zero field. In degree zero,
$\Pi:A\to A$ is the identity. For $j>n$ the exterior product is zero.

\begin{proposition}\label{prop:v5-bar-Koszul}
The maps just defined satisfy
\begin{equation}\label{eq:v5-bar-chain-maps}
 bI=Id_L,\qquad \Pi b=d_L\Pi,\qquad \Pi I=1.
\end{equation}
Both resolutions are free $A$-module resolutions of $R_0^{\rm sec}$.
There is an explicitly recursive $A$-linear homotopy $h:C_j\to C_{j+1}$
which may be normalized to satisfy
\begin{equation}\label{eq:v5-SDR}
 bh+hb=1-I\Pi,\qquad hI=0,\quad\Pi h=0,\quad h^2=0.
\end{equation}
After tensoring over $A$ with $R_f$, the same maps and identities
hold for the two-sided bar complex and the Koszul complex in
\eqref{eq:v5-graph-intersection}.
\end{proposition}
\begin{proof}
In $bI$, the interior products in the antisymmetrized sum cancel in
pairs, because $p_ip_j=p_jp_i$. The first terms are
\[
 \sum_{l=1}^j(-1)^{l-1}p_{i_l}
                   I(\eta_{i_1}\cdots\widehat\eta_{i_l}\cdots\eta_{i_j}),
\]
which is $Id_L$. This proves the first identity.

For the projection, contraction of its integrand by the radial
vector $(p_1,\ldots,p_n)$ gives a sum whose $l$th term is
\[
 (-1)^{l-1}\partial_{s_l}a_l(s_lp)\,
       da_1(s_1p)\wedge\cdots\widehat{da_l(s_lp)}
                                  \cdots\wedge da_j(s_jp).
\]
Integration in $s_l$ gives its upper and lower faces. The face
$s_1=1$ gives $a_1\Pi[a_2|\cdots|a_j]$. The face $s_j=0$ gives
$(-1)^j\epsilon_0(a_j)\Pi[a_1|\cdots|a_{j-1}]=0$.
At the common face $s_l=s_{l+1}$, the two adjacent contributions
combine, with sign $(-1)^l$, as
\[
 a_l\,da_{l+1}+da_l\,a_{l+1}=d(a_la_{l+1}).
\]
They give precisely the $l$th interior term in
\eqref{eq:v5-bar-resolution-diff}. These are all faces of the ordered
simplex. Hence $d_L\Pi=\Pi b$.
For an ordered product of $j$ generators, each summand of $I$ is
sent by $\Pi$ to $1/j!$ times the same ordered exterior product:
the volume of $\Sigma_j$ is $1/j!$, and the two permutation signs
cancel. Summing its $j!$ terms proves $\Pi I=1$.

For exactness of $C$, the augmented complex has the following
$R$-linear extra degeneracy:
\[
 s(a_0[a_1|\cdots|a_j])
       =1[a_0-\epsilon_0(a_0)|a_1|\cdots|a_j].
\]
A direct substitution in \eqref{eq:v5-bar-resolution-diff} gives
$bs+sb=1$ in positive homological degrees and
$bs=1-\iota_R\epsilon_0$ in degree zero. The map is not $A$-linear,
but proves exactness. All $C_j$ are free over $A$, since monomials
$p^\mu$ with $|\mu|>0$ form an $R$-basis of $\overline A$.
The elements $p_1,\ldots,p_n$ form a regular sequence in the
polynomial algebra over $R$. Successively taking the two-term
resolutions of these variables proves exactness of $L$.

We now give an $A$-linear comparison homotopy rather than using
the non-$A$-linear map $s$ after tensoring. Put $h_{-1}=h_0=0$.
For $j\geq1$, choose an $A$-basis element
$e=[p^{\mu_1}|\cdots|p^{\mu_j}]$ of $C_j$. Having constructed
$h$ in degrees below $j$, set
\begin{equation}\label{eq:v5-homotopy-recursion}
 r_e=e-I\Pi(e)-h_{j-1}b(e),\qquad h_j(e)=s(r_e),
\end{equation}
and extend this assignment $A$-linearly. The chain-map identities
and induction give
\[
 b r_e=b e-I\Pi b e-(1-I\Pi-hb)b e=0.
\]
Since $r_e$ is a positive-degree cycle, $bs(r_e)=r_e$.
Thus $bh_j(e)+h_{j-1}b(e)=e-I\Pi(e)$. Every polynomial and every
bar word appearing at each step is finite, so this is a well-defined
recursive construction of the homotopy on every basis element.

To impose the side conditions, put $Q=1-I\Pi$ and
$k=QhQ$. Then $bk+kb=Q$, $kI=0$, and $\Pi k=0$.
Replace $k$ by $h'=kbk$. We have
\[
 bh'+h'b=bk+kb=Q.
\]
The identity $bk^2=k^2b$, obtained by multiplying $bk+kb=Q$
on the two sides by $k$, gives
$(h')^2=kbk^2bk=k^3b^2k=0$.
The other side conditions remain true. Rename $h'$ as $h$.
This proves \eqref{eq:v5-SDR}.
All identities are now $A$-linear. Tensoring with $R_f$ preserves
them without any flatness assumption on $R_f$ over $A$.
Since we tensor resolutions, the resulting complex computes the
derived tensor product. This proves the last statement and
\eqref{eq:v5-graph-intersection}.
\end{proof}

For example, with one variable $p$, the projection is
\[
 \Pi[p^m]=p^{m-1}\eta\quad(m\geq1),\qquad
 I(\eta)=[p].
\]
The initial homotopy gives $h[p^2]=-[p|p]$, and indeed
$b(-[p|p])=[p^2]-p[p]$. On the graph one replaces the coefficient
$p$ by $f_x$, obtaining $[p^2]-f_x[p]$ as required. This verifies
the role of the graph specialization in the comparison; the bar
slots themselves are not replaced by scalars.

\subsubsection{Putting the periodic differential on the bar complex}
The retraction computes the underlying derived intersection. DT
cohomology requires the additional BV operator. We transfer it by the
constructed maps $I,\Pi$ and verify that the resulting periodic bar
differential squares to zero.

From now on write $B_f$ for the two-sided bar complex after the graph
specialization and $K_f$ for its Koszul retract. Use cohomological
degree $-j$ for a bar word of length $j$ or an exterior monomial
of degree $j$. Keep the same letters $I,\Pi,h$ for the specialized
maps. Fix the relative volume $dx_1\wedge\cdots\wedge dx_n$ and put
\[
 \Delta=\sum_i\partial_{x_i}\partial_{\eta_i},\qquad
 d_f=-\delta_f+u\Delta\quad\hbox{on }K_f((u)).
\]
\cite[Proposition~6.6]{KS} proves $d_f^2=0$ and identifies this
complex, with its parity shift, with the twisted de Rham complex.
Define on the bar complex
\begin{equation}\label{eq:v5-transferred-BV}
 \Delta_B=I\Delta\Pi,\qquad
 d_B=-b+u\Delta_B,\qquad H=-h.
\end{equation}
The operator $\Delta_B$ is the square-zero operator obtained by
transferring $\Delta$ along the chosen contraction. It is not being identified with the Connes
operator of an unspecified category.

\begin{theorem}\label{thm:v5-periodic-bar}
The operators and comparison maps satisfy
\begin{align}
 d_B^2&=0,&d_B I&=I d_f,&\Pi d_B&=d_f\Pi,
                          \label{eq:v5-periodic-bar-maps}\\
 \Pi I&=1,&d_BH+Hd_B&=1-I\Pi,&H I&=\Pi H=H^2=0.
                          \label{eq:v5-periodic-bar-homotopy}
\end{align}
They give the chain-homotopy equivalence
\begin{equation}\label{eq:v5-bar-DR-diagram}
\begin{tikzcd}[column sep=large,row sep=large]
 (B_f((u)),d_B)\arrow[r,shift left=0.65ex,"\Pi"]
 &(K_f((u)),d_f)\arrow[l,shift left=0.65ex,"I"]
                  \arrow[r,"{\mathcal I}","\sim"']
 &(\Omega^\bullet_{R/R_0}((u)),u d-df\wedge)[n],
\end{tikzcd}
\end{equation}
where $\mathcal I$ is the signed contraction to the shifted complex
of \cite[Proposition~6.6]{KS}. The differential of the rightmost
object is $(-1)^n(u\,d-df\wedge)$.
After specialization at any point $b\in V$, the diagram is the
same constructed diagram for $f_b$. If the critical locus of $f_b$
is contained in its zero fiber, its rightmost cohomology in parity $\bar i$ is
$\HP_{\overline{i+n}}(\MF(\operatorname{Spec}k(b)[x],f_b))$ by the
ordinary-function comparison theorem. This is the periodic shift
required by the normalized DT complex.
\end{theorem}
\begin{proof}
Using $bI=I\delta_f$, $\Pi b=\delta_f\Pi$, and $\Pi I=1$,
we obtain
\begin{align*}
 d_B^2
 &=-uI(\delta_f\Delta+\Delta\delta_f)\Pi
                         +u^2I\Delta(\Pi I)\Delta\Pi\\
 &=0.
\end{align*}
The same identities give $d_BI=I(-\delta_f+u\Delta)$ and
$\Pi d_B=(-\delta_f+u\Delta)\Pi$.
Since $hI=\Pi h=0$, the perturbation terms vanish in the homotopy:
\[
 d_B(-h)+(-h)d_B=bh+hb=1-I\Pi.
\]
Thus the stated homotopy identities are exact, with no infinite
homological-perturbation series. The map $\mathcal I$ and the
identity $\mathcal I d_f=(-1)^n(u d-df\wedge)\mathcal I$ are those of
\cite[Proposition~6.6]{KS}.

The inclusion, projection, and recursive homotopy are polynomial
in the coefficients $f_i$ before Laurent completion. Their formulas
therefore commute with base evaluation. Their $u$-linear extensions
and the operator $u\Delta_B$ are continuous on Laurent series with
bounded-below exponents. They commute with evaluation coefficientwise,
giving the specialized diagrams. The last assertion uses
the comparison theorem of \cite{Efimov2018CyclicMF}; all matrix factorizations in it have the
ordinary degree-zero potential $f_b$.
\end{proof}

Different normalized comparison contractions give homotopy-equivalent
periodic bar models: between two such models use $I'\Pi$ and $I\Pi'$.
Their two composites are $I'\Pi'$ and $I\Pi$, homotopic to the
respective identities by \eqref{eq:v5-periodic-bar-homotopy}.
Consequently the periodic model just constructed does not depend,
up to this explicit homotopy equivalence, on the recursion chosen
to normalize the comparison homotopy. The volume change is handled
by \cite[equation~(6.14)]{KS}.

\subsubsection{A connection on the bar complex itself}
A quasi-isomorphism of complexes is not yet the horizontal comparison
needed in Theorem~\ref{thm:intro-cyclic}. We now construct the bar
connection with correction term $HE$ and verify that $I$ and $\Pi$
intertwine it with the ordinary-potential Gauss--Manin operator.

The next calculation makes the arrows in
\eqref{eq:v5-bar-DR-diagram} horizontal. It is not sufficient merely
to transfer a connection on cohomology.
On $K_f((u))$ put
\[
 \nabla_K=\partial_t-u^{-1}\partial_t f.
\]
As in Lemma~\ref{lem:v5-Cartan},
$[\nabla_K,d_f]=0$. The underlying bar module has a fixed basis of
words in the monomials $p^\mu$, so it has the coefficientwise
operator $\partial_t$. The maps $\Pi,h$ may depend on $t$ through
$f_i$. Write $Q=1-I\Pi$ and define
\begin{equation}\label{eq:v5-bar-connection-construction}
 \nabla^0_B=I\nabla_K\Pi+Q\partial_tQ,\qquad
 E=[\nabla^0_B,d_B],\qquad
 \nabla_B=\nabla^0_B+H E.
\end{equation}
Composition with a derivation means composition as an operator;
thus the first two terms include derivatives of $\Pi$ and $Q$.

\begin{proposition}\label{prop:v5-horizontal-bar}
The operator $\nabla_B$ is a connection and satisfies
\begin{equation}\label{eq:v5-bar-horizontal}
 [\nabla_B,d_B]=0,\qquad
 \nabla_BI=I\nabla_K,\qquad
 \Pi\nabla_B=\nabla_K\Pi.
\end{equation}
It is regular in $t$ whenever the displayed data are regular in $t$.
Its cohomology connection is flat along the curve.
\end{proposition}
\begin{proof}
The Leibniz rule for $\nabla_B^0$ follows from $I\Pi+Q=1$ and
$Q^2=Q$. The commutator $E$ is linear over base functions; the
Leibniz terms cancel because $d_B$ is base-linear. Thus $HE$ is
base-linear and $\nabla_B$ remains a connection.
The identities $QI=0$ and $\Pi Q=0$ give
\[
 \nabla_B^0I=I\nabla_K,\qquad \Pi\nabla_B^0=\nabla_K\Pi.
\]
For example $Q\partial_tQI=0$ as an operator because $QI=0$
is an identity for every $t$; the derivative acts on the entire
zero operator. Since $I,\Pi$ are chain maps and $\nabla_K$
commutes with $d_f$, we obtain $EI=0$, $\Pi E=0$, and hence
$QE=EQ=E$.

The relation $d_B^2=0$ gives $d_BE+Ed_B=0$. Using
$d_BH+Hd_B=Q$,
\begin{align*}
 [HE,d_B]
 &=HEd_B-d_BHE\\
 &=HEd_B-(Q-Hd_B)E\\
 &=-E+H(Ed_B+d_BE)=-E.
\end{align*}
Adding this identity to the definition of $E$ proves the first
formula in \eqref{eq:v5-bar-horizontal}. The remaining two follow
from $EI=0$ and $\Pi H=0$. All operators in
\eqref{eq:v5-bar-connection-construction} involve only polynomial
coefficient operations, their derivatives, and finitely many negative
powers of $u$ on each operator; they do not require division by $t$.
This proves regularity. Flatness on cohomology follows because the
base has dimension one.
\end{proof}

We have therefore constructed the commutative diagram of complexes
with base connections
\begin{equation}\label{eq:v5-horizontal-comparison-diagram}
\begin{tikzcd}[column sep=large,row sep=large]
 (B_f((u)),d_B)\arrow[r,"\Pi"]\arrow[d,"\nabla_B"']
 &(K_f((u)),d_f)\arrow[r,"\mathcal I"]\arrow[d,"\nabla_K"]
 &(\Omega^\bullet((u)),D_f)[n]\arrow[d,"{\partial_t-u^{-1}\partial_t f}"]\\
 (B_f((u)),d_B)\arrow[r,"\Pi"']
 &(K_f((u)),d_f)\arrow[r,"\mathcal I"']
 &(\Omega^\bullet((u)),D_f)[n].
\end{tikzcd}
\end{equation}
All vertical maps are covariant derivatives along the same base
vector field $\partial_t$. In the chosen coordinate volume,
$\mathcal I$ is independent of $t$ and commutes with multiplication
by $\partial_t f$, giving the last square.

\subsubsection{The central comparison and formal horizontal transport}
At $t=0$, the left complex in
\eqref{eq:v5-horizontal-comparison-diagram} is the explicitly
constructed periodic bar complex of the zero section and
$\operatorname{graph}(df_0)$. Thus the central chain comparison is
$\mathcal I\Pi_0$, not an additional assumed quasi-isomorphism in
this chart. The application to the moduli stacks on $X_1$ and $X_2$ requires
identifying the chosen critical chart of their derived intersection, with its orientation and
transition maps, with this presentation.

There is also a completely explicit formal transport on bar complexes.
Use exactly the $t$-after-$u$ completion of
Section~\ref{sec:v5-formal-transport}. Denote the comparisons at
$t$ by $I_t,\Pi_t,H_t$ and at zero by $I_0,\Pi_0,H_0$.
Multiplication by $T(t)=\exp((f(t)-f_0)/u)$ acts on the Koszul
complex as well as on differential forms. Define
\begin{equation}\label{eq:v5-bar-formal-transport}
 \mathcal T_t=I_tT(t)\Pi_0,
 \qquad \mathcal S_t=I_0T(t)^{-1}\Pi_t.
\end{equation}

\begin{theorem}\label{thm:v5-local-bridge}
In the stated completion, $\mathcal T_t$ and $\mathcal S_t$ are
chain-homotopy inverses between the constant central periodic bar
complex and the family periodic bar complex. They satisfy
\begin{align}
 d_{B,t}\mathcal T_t&=\mathcal T_t d_{B,0},&
 \nabla_{B,t}\mathcal T_t&=\mathcal T_t\partial_t,
                           \label{eq:v5-bar-transport-chain}\\
 \mathcal S_t\mathcal T_t&=I_0\Pi_0,&
 \mathcal T_t\mathcal S_t&=I_t\Pi_t.
                           \label{eq:v5-bar-transport-inverse}
\end{align}
The two products are homotopic to the respective identities by
$H_0$ and $H_t$. Thus central and family cohomology are explicitly
identified after $t$-adic completion in this local model. On ordinary
geometric fibers satisfying the ordinary-function comparison theorem,
the same bar model computes their periodic cyclic homology shifted by
$[n]$; formal
transport to those fibers additionally requires compatibility with
the order of completion specified in Section~\ref{sec:v5-formal-transport}.
\end{theorem}
\begin{proof}
Theorem~\ref{thm:v5-formal-transport}, transported through
$\mathcal I$, gives $d_{f,t}T=Td_{f,0}$ and
$\nabla_{K,t}T=T\partial_t$. Compose these identities with the
horizontal chain maps $I_t$ and $\Pi_0$ of
Proposition~\ref{prop:v5-horizontal-bar}. The zero-fiber complex
in this calculation is extended constantly in $t$, so its maps
$I_0,\Pi_0$ commute with $\partial_t$. This proves
\eqref{eq:v5-bar-transport-chain}.
Using $\Pi_t I_t=1$ and $T^{-1}T=1$ gives
\[
 \mathcal S_t\mathcal T_t
  =I_0T^{-1}(\Pi_t I_t)T\Pi_0=I_0\Pi_0.
\]
The other composition is the same calculation with $0,t$
interchanged. The homotopy identities are exactly
\eqref{eq:v5-periodic-bar-homotopy}, extended coefficientwise in
$t$. They prove the cohomological comparison. The statements about
ordinary fibers and completion follow from
Theorem~\ref{thm:v5-periodic-bar} and
Theorem~\ref{thm:v5-formal-transport}, with their specified coefficient
rings.
\end{proof}

For the degeneration moduli in the paper, the construction to be
compared with this local result is
\[
 \mathfrak M_1\times^h_{\mathfrak{Perf}(S)}\mathfrak M_2,
\]
using the full derived restrictions, including their negative
cohomology. A compatible relative critical presentation gives the
local bar comparison above. A global application requires the
presentation maps and their homotopies on overlaps, together with
finiteness and the chosen support convention. In a single global
ordinary-potential presentation, Theorem~\ref{thm:v5-Cech-GM}
constructs its Cech connection. Hypothesis~\ref{hyp:finite-cyclic}
specifies the global comparisons for the moduli stacks on $X_1$ and $X_2$.
The application below satisfies those comparisons on the stated
hyperplane-section locus.

\subsection{Finite cohomology for families with finite Jacobian algebras}\label{sec:v5-finite-critical}
The preceding maps identify the fiber complexes. To obtain finite
cohomology and regular specialization, we impose a condition that can
be checked on the gradient Koszul complex itself. The contraction
constructed below produces the finite representative and its
connection matrix, and proves Theorem~\ref{thm:v5-finite-Jacobian}.

Let $R_0=\Gamma(V,\O_C)$ as above,
$R=R_0[x_1,\ldots,x_n]$, and $f\in R$. Suppose that the Koszul
complex of $f_1,\ldots,f_n$ is exact except in its last form degree
and that
\begin{equation}\label{eq:v5-finite-Jacobian}
 J=R/(f_1,\ldots,f_n)
\end{equation}
is a finite projective $R_0$-module. A regular sequence of partial
derivatives with finite projective quotient is a sufficient condition.
After shrinking $V$ to trivialize $J$, choose a basis $e_1,\ldots,e_\mu$
and polynomial lifts $\sigma_j\in R$. Put $\omega=dx_1\wedge\cdots\wedge dx_n$.

The coefficient ring of the family in this statement is
\[
 \mathcal R=R_0[[u]][u^{-1}].
\]
It permits a power series in $u$ whose individual coefficients are
regular functions on $V$, with a common lower bound on $u$-powers.
It is not silently identified with $R_0\otimes_k k((u))$.
Evaluation at $b\in V$ gives $\mathcal R\to k(b)((u))$.

\begin{theorem}\label{thm:v5-finite-Jacobian}
Under these hypotheses, the complex
\[
 (\Omega^\bullet_{R/R_0}[[u]],D_f),\qquad
 D_f=-df\wedge+u\,d_{R/R_0},
\]
has cohomology only in form degree $n$. The splittings and maps constructed in the proof identify it
with $J[[u]]$.
The classes $[\sigma_j\omega]$ form a basis wherever $J$ is free.
After inverting $u$ its cohomology is finite projective over
$\mathcal R$, in parity $n$, and its formation commutes with evaluation
at every $b\in V$, including $0$.

The operator $\partial_t-(\partial_tf)/u$ induces a connection on
this module. In the displayed basis its matrix has entries in
$u^{-1}R_0[[u]]$, so it is regular at every point of $V$ in the
geometric variable. It is flat along the curve. If the critical locus
of each $f_b$ is contained in $f_b^{-1}(0)$, its fibers are the ordinary
periodic cyclic homology groups of $\MF(\operatorname{Spec}k(b)[x],f_b)$.
The constructed periodic bar complex of
Theorem~\ref{thm:v5-periodic-bar} has the same finite rank and
connection after the shift $[n]$: its nonzero cohomology has even
parity.
\end{theorem}
\begin{proof}
Write $d_0=-df\wedge$ and consider first the bounded complex of
$R_0$-modules $(\Omega^\bullet_{R/R_0},d_0)$. Its terms are free over
$R_0$, possibly of infinite rank. Its last cohomology is $J\omega$,
which is projective by hypothesis. The quotient map
$\Omega^n\to J\omega$ therefore has an $R_0$-linear section $i$,
chosen to send $e_j$ to $\sigma_j\omega$ on a free open.
The image $B^n=\operatorname{im}(d_0:\Omega^{n-1}\to\Omega^n)$
is a projective direct summand of $\Omega^n$. The surjection
$\Omega^{n-1}\to B^n$ splits. Its kernel is projective and, by
exactness below degree $n$, equals $B^{n-1}$. Repeating downward
constructs decompositions
\[
 \Omega^q=B^q\oplus L^q\quad(q<n),\qquad
 \Omega^n=B^n\oplus i(J\omega),
 \quad d_0:L^q\xrightarrow{\sim}B^{q+1}.
\]
Let $p$ be the projection to $J\omega$ in degree $n$ and zero in
other degrees. Define $h$ on $B^{q+1}$ to be the inverse of the last
displayed isomorphism and set $h=0$ on $L^{q+1}$ and on $i(J\omega)$.
This specifies every component of $h$. We obtain
\begin{equation}\label{eq:v5-Jacobian-contraction}
 d_0h+hd_0=1-ip,\quad pi=1,\quad hi=ph=h^2=0,
 \quad d_0i=pd_0=0.
\end{equation}
These identities follow on $B^q$, $L^q$, and $i(J\omega)$ separately.

Set $Q=u\,d_{R/R_0}$ and $M=1+Qh$. The operator $Qh$ raises
$u$-adic order by one and preserves form degree, so
\[
 M^{-1}=\sum_{a\geq0}(-Qh)^a
\]
is a well-defined $R_0[[u]]$-linear operator on formal series.
No bound on polynomial degree is needed: each coefficient is obtained
by finitely many polynomial operations and applications of $h$.
We have $Q^2=0$, $d_0Q+Qd_0=0$, and $Qi=0$, the last identity
because $i$ takes values in top-degree forms. Consequently
\begin{align*}
 D_fM
 &=(d_0+Q)(1+Qh)
   =d_0+Q-Qd_0h\\
 &=d_0+Q-Q(1-ip-hd_0)
   =(1+Qh)d_0=M d_0.
\end{align*}
Thus $D_f=M d_0M^{-1}$. In particular the required deformation
retract is given explicitly by
\begin{equation}\label{eq:v5-Jacobian-perturbed-contraction}
 I_u=i,\qquad P_u=pM^{-1},\qquad H_u=hM^{-1}.
\end{equation}
Indeed $Mi=i$ and $Mh=h$, by $hi=h^2=0$. Conjugating
\eqref{eq:v5-Jacobian-contraction} by $M$ gives
\[
 D_fH_u+H_uD_f=1-I_uP_u,\qquad P_uI_u=1,
 \qquad D_fI_u=P_uD_f=0.
\]
This proves that the cohomology is $J[[u]]$ in degree $n$ and that
its basis is represented by $\sigma_j\omega$. Localization at $u$
is exact, giving $J[[u]][u^{-1}]$.

All the unperturbed splittings are $R_0$-linear. Their identities
therefore remain true after any base evaluation, even though the
individual terms have infinite rank. The series defining $M^{-1}$
is evaluated coefficientwise, and at every fixed $u$-order only
finitely many terms occur. Hence the same contraction gives the
cohomology of the evaluated complex. Moreover evaluation of series
satisfies
\[
 R_0[[u]][u^{-1}]/(t-b)\simeq k(b)((u))
\]
in a local coordinate neighborhood of a closed point: coefficientwise
division by $t-b$ does not change the lower bound in $u$.
The corresponding assertion for a point ideal follows by choosing
its local parameter. Since $J$ is finite projective, this proves
\[
 H^n(D_f)\otimes_{\mathcal R}k(b)((u))
    \xrightarrow{\sim} H^n(D_{f_b}).
\]

It remains to compute the connection, not merely the dimensions.
Let $\nabla=\partial_t-u^{-1}\partial_t f$. Lemma~\ref{lem:v5-Cartan} proves
$[\nabla,D_f]=0$. The connection matrix is the following explicit
formal series of operators applied to the basis lifts:
\begin{equation}\label{eq:v5-Jacobian-connection-matrix}
 \sum_i A_{ij}(t,u)e_i\omega
 =p\sum_{a\geq0}(-u\,d_{R/R_0}h)^a
       \left((\partial_t\sigma_j)\omega
                   -u^{-1}(\partial_tf)\sigma_j\omega\right).
\end{equation}
Only the second term in parentheses has a negative $u$-power,
namely $u^{-1}$; all powers supplied by the sum are nonnegative.
Every coefficient belongs to $R_0$, since the chosen splittings
are $R_0$-linear and the input forms have coefficients in $R$.
Thus $A_{ij}\in u^{-1}R_0[[u]]$. This is regular in the base
coordinate throughout the chosen free open. A change of basis of
$J$ has its usual connection transformation law, so these formulas
define a connection on the finite projective cohomology module.
Its curvature along the curve is zero. The last two assertions
follow respectively from the comparison theorem of \cite{Efimov2018CyclicMF} and the horizontal
homotopy equivalence in Proposition~\ref{prop:v5-horizontal-bar}.
\end{proof}

Proposition~\ref{prop:v15-rigid-stabilization} realizes the hypotheses
of this theorem on the hyperbolic stabilization $f=zw$ of the reduced
moduli point $\O_{X_b}$ in the same quintic family. For the
positive-dimensional divisor loci with $f=0$, finiteness instead
comes from proper global cohomology, as in
Theorem~\ref{thm:r14-divisor-examples}.

The local comparison maps and finite-cohomology criteria are now
available. We next apply the construction to a family of actual
objects on the resolved quintic. Its critical locus is projective
and positive-dimensional, so its finiteness is proved by the global
bar complex and projective-space cohomology rather than by the
finite-Jacobian hypothesis.

\section{Cohomological degeneration and its logarithmic extension}
\label{sec:categorical-degeneration}
The preceding constructions compare normalized periodic realizations
of the chosen categorical DT theory. The result is an isomorphism of
cohomology with connection, subject to the geometric comparison
hypotheses. An equivalence between the underlying dg categories is
a stronger assertion and is not inferred from periodic cyclic homology.
We first give a geometric criterion that does not restrict the
amplitude or presentation of the perfect objects.

\subsection{A sufficient criterion for arbitrary perfect-complex loci}
Let $V\subset C$ be a smooth affine coordinate neighborhood, put
$R_0=\Gamma(V,\O_C)$ and $K=\CC((u))$. In this subsection the
algebraic coefficient ring is $R_K=R_0\otimes_\CC K$.
Periodization is applied to the proper de Rham pushforward and then
extended to $K$. The $u$-adically completed coefficient convention
instead gives $R_0[[u]][u^{-1}]$; the comparison with that convention
is obtained by extending coefficients, not by identifying these two
rings. This distinction is the one used in
Section~\ref{sec:v5-formal-transport}.

\begin{theorem}\label{thm:smooth-proper-comparison}
Suppose $p:M\to V$ is smooth and projective of relative dimension
$d$ and is identified with a classical locus of simple perfect
complexes on the fibers of $\pi$. Assume that its central objects
are identified with their effective derived matching triples and
that the chosen d-critical realizations agree under these
identifications. Equip their zero d-critical structures with the
relative orientation $K_{M/V}$. Use the corresponding zero critical
matrix-factorization and oriented bar realizations.
Then their cohomological DT sheaves on the fibers are
$\CC_{M_b}[d]$, and their normalized periodic cohomology has finite
locally free algebraic model
\begin{equation}\label{eq:general-smooth-model}
 P^{\bar a}=\bigoplus_{j-d\equiv a\ (2)}
       H^j_{\mathrm{dR}}(M/V)\otimes_\CC K,
       \qquad d_P=0,
\end{equation}
where $H^j_{\mathrm{dR}}(M/V)=\mathcal H^j(Rp_*\Omega^\bullet_{M/V})$.
This model has the ordinary relative Gauss--Manin connection,
defined without poles on $V$. It has chain comparisons to the
chosen oriented bar realization, horizontal up to homotopy,
and commutes with every fiber specialization. Consequently it
satisfies the finite regular comparison conditions for these
realizations. The perfect objects parametrized by $M$ are not
required to be sheaves or two-term section complexes.
\end{theorem}
\begin{proof}
On each smooth $M_b$ the d-critical structure is zero and its
virtual canonical bundle is $K_{M_b}^{\otimes2}$. The specified
orientation therefore gives the constant sheaf shifted by $d$.
On a smooth affine chart of $M/V$, take the two zero sections of
its relative cotangent bundle. The relative Koszul-to-bar maps
of Section~\ref{sec:v5-bar-bridge}, tensored with $K_{M/V}$ and
contracted with top forms, give
\begin{equation}\label{eq:general-smooth-bar}
 \mathcal B^{\mathrm{or}}_{M,\mathrm{per}}
    \longrightarrow(\Omega^\bullet_{M/V}\otimes_\CC K,u\,d)[d].
\end{equation}
The formulas are tensorial under coordinate changes, and the
orientation Jacobians are the transition functions of $K_{M/V}$.
Thus the maps glue. Their local homotopy-equivalence property
proves that the global map is a quasi-isomorphism. The base
connection is the one constructed by horizontal lifts in
Theorem~\ref{thm:v5-Cech-GM}, with $F=0$; the horizontal bar
identities of Proposition~\ref{prop:v5-horizontal-bar} identify
its induced connection with that on relative de Rham cohomology.

For a smooth proper morphism in characteristic zero, relative
de Rham cohomology is finite locally free, commutes with base
change, and carries the Gauss--Manin connection constructed from
the absolute de Rham filtration \cite{KO}. The connection is
defined on $V$, including $0$, because $p$ is smooth there.
The same operator is obtained from the local horizontal-lift
construction: both lift a base derivation to absolute forms and
use Cartan homotopy on intersections. Since $V$ is a curve this
connection is flat; its usual integrability also follows directly
from the filtered de Rham construction.

Choose a bounded finite-projective representative $C^\bullet$
of $Rp_*\Omega^\bullet_{M/V}$ over the affine ring $R_0$.
Its cohomology modules are projective. Splitting successively
from the top degree shows that boundaries and cycles are also
projective and gives a contraction onto
$\bigoplus_jH^j(C)[-j]$. To see this explicitly, at the top degree
split $C^N\twoheadrightarrow H^N(C)$; its kernel is the top
boundary module. Split $C^{N-1}\twoheadrightarrow B^N$, then
$Z^{N-1}\twoheadrightarrow H^{N-1}(C)$, and continue downward.
These splittings produce maps $i,p,h$ with
$pi=1$ and $1-ip=dh+hd$. Periodization and extension to $R_K$
give \eqref{eq:general-smooth-model}.

Multiplication by $u^r$ on relative $r$-forms intertwines $d$
with $u\,d$. Only $0\le r\le d$ occur, so this changes Laurent
bounds by a uniformly bounded amount. The contraction therefore
also gives the comparison in the completed coefficient convention.
The splittings express $C^\bullet$ as its cohomology together with
two-term contractible complexes $B^j\xrightarrow{1}B^j$.\\
Choose algebraic connections on these projective boundary modules;
such connections exist over the smooth affine base, since the
Atiyah extension of a projective module splits. Use the same
connection on both terms of each contractible pair, and use the
Gauss--Manin connection on the cohomology summands. This gives
a chain connection on $C^\bullet$ commuting with its differential.
The projection and inclusion of the cohomology summands are
horizontal for this connection. \\
Comparison with the geometric
connection induces the same map on de Rham cohomology, so the
difference between the two routes of a comparison is an
$R_K$-linear chain map inducing zero on cohomology. Using the displayed
contraction, such a map is nullhomotopic: its component on
cohomology is zero and all remaining components factor through
contractible summands. This proves horizontality up to homotopy.
Every comparison is preserved after derived evaluation at $b$,
by perfectness and proper de Rham base change. Finally the zero
critical instance of \eqref{eq:cyclic-DR-input} gives the specified
normalized cyclic realization on each fiber. These are all the
assertions of the theorem.
\end{proof}

The hypotheses of this criterion are geometric: smoothness,
properness, the actual central matching identification, and the
compatible orientation must be established for the selected moduli
locus. The criterion removes the special form of the objects from
the comparison once those hypotheses are known. The projective
and Grassmannian loci constructed in Section~\ref{sec:v6-hyperplanes} satisfy them. The
one- and two-point incidence families have nontrivial monodromy
on the punctured base, and therefore cannot admit a smooth proper
extension of this kind with the same cohomological variation.

\subsection{Charges and the central matching complex}
Let $E_0$ be a perfect complex on $X_0$, with derived restrictions
$E_i$ and $E_S$. The normalization triangle used in
\cite[Proposition~3.3]{KS}, now tensored with $E_0$, gives
\begin{equation}\label{eq:Kgluing-MV}
 [E_0]=[j_{1*}E_1]+[j_{2*}E_2]-[j_{S*}E_S]
             \quad\text{in }G_0(X_0).
\end{equation}
Here $G_0$ is the Grothendieck group of coherent sheaves. Pushforward
from a component need not be perfect on the singular union, so
this relation is stated in $G_0$, not asserted in the perfect
$K$-group. It includes the boundary contribution, which must be
retained when indexing a degeneration by charges.

For a chosen charge $\gamma$, write
\begin{equation}\label{eq:def-periodic-invariants}
 \mathsf D_{b,\gamma}=H^{\mathbb Z/2}(\mathsf K_{b,\gamma}^{\mathrm{DT}}),\qquad
 \mathsf D_{0,\gamma}=\mathbb H_c^{\mathbb Z/2}
       (M_{0,\gamma},\mathcal B_{0,\gamma,\mathrm{DT}}).
\end{equation}
The cyclic realization and bar complex are those specified in
Hypothesis~\ref{hyp:finite-cyclic}. The critical-chart suspension
and support convention are included; when the suspension is odd,
these parity groups differ from unshifted ordinary periodic cyclic
homology by that parity shift. They are not determined by $\gamma$ alone.
A decomposition by component charges also requires a defined
specialization of those charges and a finite or suitably locally
finite collection of triples satisfying \eqref{eq:Kgluing-MV}.

\begin{theorem}\label{thm:categorical-degeneration-conditional}
Suppose Hypothesis~\ref{hyp:finite-cyclic} holds for the chosen
charge $\gamma$. Then the two parity ranks of its finite regular
model are constant and its fiber cohomology commutes with base
change. In particular, for every closed $b\in V$, the finite
$K$-supervector spaces in \eqref{eq:def-periodic-invariants} have
the same even and odd dimensions. Their Grothendieck classes agree.
On the completed local base $K[[t]]$ there is a horizontal
identification with the constant central fiber; after inverting $t$
this gives
\begin{equation}\label{eq:coefficientwise-cat-degeneration}
 H(P_\gamma\otimes_{R_K}K((t)))
   \simeq \mathsf D_{0,\gamma}\otimes_KK((t)).
\end{equation}
Analytic transport between distinct closed fibers additionally
uses an analytic realization and a chosen path. If the central
complex has a specified charge-compatible finite direct-sum
decomposition, its cohomology and Grothendieck class decompose
accordingly. For a monoid of charges and a completion in which
these sums are defined, this yields the coefficientwise identity
\begin{equation}\label{eq:categorical-partition}
 \sum_\gamma q^\gamma[\mathsf D_{b,\gamma}]
       =\sum_\gamma q^\gamma[\mathsf D_{0,\gamma}]
       \quad\text{in }K_0(\mathrm{sVect}_K)[[q^\Gamma]].
\end{equation}
The group $K_0(\mathrm{sVect}_K)$ here retains both parity ranks;
it is the Grothendieck group of the abelian category of
finite-dimensional supervector spaces.
\end{theorem}
\begin{proof}
Theorem~\ref{thm:finite-cyclic} gives local freeness and base
change of $H(P_\gamma)$. On a connected base its parity ranks
are constant, proving the assertion for closed fibers and their
Grothendieck classes without choosing a transport path. The
geometric quasi-isomorphisms in Hypothesis~\ref{hyp:finite-cyclic}
identify those fibers with \eqref{eq:def-periodic-invariants}.
Over $K[[t]]$, write the connection matrix as $A(t)$. The recursion
\[
 G_0=1,\qquad (r+1)G_{r+1}
       =-\sum_{i=0}^r A_iG_{r-i}
\]
constructs the unique invertible solution to $G'=-AG$ with
$G(0)=1$. It gives the formal horizontal identification and
\eqref{eq:coefficientwise-cat-degeneration}. This is an identity
over the formal Laurent field, not an evaluation of $G(t)$ at
a nonzero complex number. An assumed direct-sum comparison
passes to cohomology, and taking Grothendieck classes and summing
in the stated completion proves \eqref{eq:categorical-partition}.
\end{proof}

\begin{corollary}\label{cor:all-degree-specialization}
For the resolved quintic and the divisor-section loci of
Theorem~\ref{thm:all-positive-degrees}, the comparisons are constant
in the hyperplane-power basis for every $m>0$ and $k\in\mathbb Z$.
Their common coefficient is $[K^{r_m}]$ in parity $d_m=r_m-1$.
Thus the geometric completion and the cohomological comparison
asked for in \cite[Conjectures~1.3 and~1.4]{KS} hold on these loci.
The Grassmannian and transported perfect-complex loci have the
corresponding statement with their horizontal Schubert bases.
\end{corollary}
\begin{proof}
The chain comparisons and constant finite models are constructed
in Theorems~\ref{thm:all-positive-degrees},
\ref{thm:grassmannian-DT}, and Proposition~\ref{prop:transport-equivalences}.
Their central fiber identifications are part of those constructions.
Apply Theorem~\ref{thm:categorical-degeneration-conditional}.
\end{proof}

\subsection{Logarithmic comparison for general divisor charges}
Theorem~\ref{thm:incidence-monodromy} and
Theorem~\ref{thm:two-point-incidence} exhibit charge-fixed
surface-sheaf loci with nonzero unipotent monodromy. A connection
extending without a pole over a disk has a single-valued local
horizontal frame, and hence trivial monodromy on the punctured
disk. These examples therefore exclude an unrestricted pole-free
version of \eqref{eq:coefficientwise-cat-degeneration} for all
charges $(0,H,\beta,n)$.

The remaining generalization is the logarithmic problem formulated
in \cite[Problem~5.2]{KS}. Its categorical input and its limiting
cohomological assertion can be stated separately as follows.
\begin{problem}\label{prob:logarithmic-categorical-comparison}
For a specified bounded family of oriented moduli of perfect
complexes on $C^*=C\setminus\{0\}$, construct compatible cyclic
realizations with finite cohomology and a logarithmic extension
\[
 \nabla:\mathcal H\longrightarrow
       \mathcal H\otimes\Omega_C^1(\log\{0\}).
\]
Identify its residue and limiting fiber with the monodromy and
nearby-cycle realization of the actual relative DT complex.
Construct, for a defined oriented central matching problem, a
specialization morphism from its DT bar complex to that
nearby-cycle realization, and determine its cone. Require the
maps to preserve the charge relation \eqref{eq:Kgluing-MV} and
the orientation comparisons.
\end{problem}

Nearby cycles require a constructible relative DT complex on the
punctured family, not just a collection of vector spaces. Even
when a logarithmic extension of the cohomology local system is
available, its limiting fiber is not automatically the DT
cohomology of an ordinary central moduli space. The specialization
morphism in the problem is the additional comparison that would
justify such an identification. For the one-point incidence
family in Proposition~\ref{prop:incidence-log-extension} we describe
its logarithmic extension; the two-point calculation determines
its new monodromy contribution explicitly. The boundary self-intersection and the two distinct degeneration
families will be constructed directly in Sections~\ref{sec:divisor-sheaves}
and~\ref{sec:second-degeneration}. Their geometric realizations,
rather than numerical gluing identities, determine the comparisons.

\subsection{Exact constructions and change of perfect-complex moduli problem}
\label{ssec:general-perfect-functors}
The geometry and the comparison criteria above refer to moduli
of perfect objects, not to a fixed number of terms in a chosen
resolution. Two functorial constructions make this generality
explicit. In both statements the orientation and critical
realization remain part of the DT datum.

\begin{proposition}\label{prop:general-perfect-fiber}
Let $T\to C$ and let $a_T:A_T\to Q_T$ be a morphism of perfect
complexes on $\cX_T=\cX\times_C^{\mathrm h}T$. Its homotopy fiber
$E_T=\fib(a_T)$ is perfect. This construction is functorial in the
infinity-category of arrows and commutes with arbitrary derived
base changes. On the central fiber its two component restrictions
and their common boundary restriction form the homotopy-commutative
diagram of exact triangles
\begin{equation}\label{eq:general-fiber-restriction}
 \begin{tikzcd}[column sep=large,row sep=large]
 Lj_i^*E_0\ar[r]\ar[d]&Lj_i^*A_0\ar[r]\ar[d]&Lj_i^*Q_0\ar[d]\\
 Lj_S^*E_0\ar[r]&Lj_S^*A_0\ar[r]&Lj_S^*Q_0 .
 \end{tikzcd}
\end{equation}
The two lower rows are identified by associativity of derived
pullback. Thus the resulting central triple is effective.
\end{proposition}
\begin{proof}
The category of perfect complexes is a stable full subcategory:
locally the cone of a map between finite locally free resolutions
is again a bounded finite locally free complex. A derived pullback
is an exact functor of stable categories. It therefore sends the
fiber sequence defining $E_T$ to the fiber sequence of the
pulled-back arrow, with its natural equivalence. The associativity
constraints of pullback give the compatibility at $S$ and their
higher homotopies. This proves the assertion on families, rather
than only the equality of their classes in a Grothendieck group.
\end{proof}

This proposition defines a morphism from the derived moduli of
arrows of perfect complexes to the moduli of their fibers. To use
an image as a DT moduli locus one still proves simplicity,
boundedness, and the required representability or properness for
that image. None follows just from closure under cones.
Proposition~\ref{prop:derived-point-modification} verifies the
perfectness input for finite point modifications in a geometric
family. The spherical-twist examples verify a stronger condition:
the construction is an equivalence of entire relative moduli
problems.

\begin{proposition}\label{prop:general-equivalence-transport}
Let $\cX/C$ and $\cX'/C$ be families of the kind used above.
Suppose enhanced equivalences
$\Phi_T:\Perf(\cX_T)\xrightarrow{\sim}\Perf(\cX'_T)$
and their inverses commute coherently with derived base change.
Let $M\to C$ be a specified moduli locus and $M'\to C$ its image.
Transport the critical realization, orientation, and cyclic/bar
comparisons along the induced isomorphism of moduli functors.
Then the cohomological DT complexes and their finite Gauss--Manin
comparisons are identified under $M\simeq M'$. If the target
critical realization is constructed independently, the same
conclusion holds when a compatible isomorphism with the transported
oriented realization is given.
\end{proposition}
\begin{proof}
On every test scheme the equivalence and its inverse identify the
object groupoids, isomorphisms, and higher homotopies. They
preserve the endomorphism dg algebras, and hence the simplicity
condition and full deformation problems. This gives the indicated
isomorphism of the selected relative moduli functors.
Base change carries an actual object on $X_0$ to one on $X'_0$;
restriction of the latter gives its effective central triple.
No equivalence between individual component categories is needed
for this assertion.

An isomorphism of oriented critical realizations pulls back the
vanishing-cycle sheaf and its overlap maps, so also its compactly
supported hypercohomology. Pulling back the specified cyclic and
bar chain comparisons gives the target comparisons. If
$\nabla_Bc-c\nabla_P=d_BH+Hd_P$ is the original horizontal
identity, pullback preserves that identity, its differential,
and its homotopy. Thus the finite model and the induced
Gauss--Manin connection are identified. This proves the claim
without assuming that an independently selected orientation or
stability condition is preserved automatically.
\end{proof}

\subsection{The additional specialization map and a test for its being an isomorphism}
\label{ssec:actual-specialization-test}
A connection differentiates sections of a cohomology bundle; it is
not itself a map from one fiber to another. Such a map is obtained
by horizontal transport along a path, or by the formal initial-value
construction of Theorem~\ref{thm:finite-cyclic}.
On a punctured disk, a flat connection can have nontrivial monodromy.
A logarithmic extension retains this monodromy in its residue. Its
fiber at the origin is a limiting realization, which must still be
compared with the DT complex on the chosen central moduli space.
The following criterion isolates this latter comparison.

\begin{proposition}\label{prop:specialization-cone-test}
Let $p:\mathcal M\to\Delta$ be proper over a complex analytic disk,
let $i:\mathcal M_0\hookrightarrow\mathcal M$ and
$j:\mathcal M_{\Delta^*}\hookrightarrow\mathcal M$, and let
$\mathcal P$ be a bounded constructible complex on $\mathcal M$.
Suppose that its restrictions to nonzero fibers are the specified
oriented DT complexes, with a fixed common grading convention, and
that an identification
\[
 i^*\mathcal P\simeq\mathcal P^{\rm bar}_0
\]
with the oriented central DT realization has been constructed.
There is an exact triangle
\begin{equation}\label{eq:actual-specialization-triangle}
 R\Gamma(\mathcal M_0,\mathcal P^{\rm bar}_0)
 \xrightarrow{\operatorname{sp}}
 R\Gamma(\mathcal M_0,R\psi_p\mathcal P)
 \longrightarrow R\Gamma(\mathcal M_0,R\phi_p\mathcal P)
 \xrightarrow{+1}.
\end{equation}
Here $R\phi_p$ denotes the cone of the specialization morphism
without a perverse shift. The morphism $\operatorname{sp}$ is a
quasi-isomorphism if and only if its third term is acyclic. The
same assertion holds after a specified exact periodization.
In particular, local acyclicity of $\mathcal P$ with respect to $p$
is sufficient. If the smooth-fiber realization is identified with a
finite regular cyclic model as in Hypothesis~\ref{hyp:finite-cyclic},
compatibly with this specialization, the horizontal comparison of
Theorem~\ref{thm:categorical-degeneration-conditional} identifies both endpoints.
\end{proposition}
\begin{proof}
The defining nearby-cycle triangle is
$i^*\mathcal P\to R\psi_p\mathcal P\to R\phi_p\mathcal P\to i^*\mathcal P[1]$.
Apply derived global sections and the stipulated central
identification. An exact triangle has its first arrow a
quasi-isomorphism precisely when its cone is acyclic. This proves
the criterion. Proper base change identifies the middle term with
the limiting cohomology of the nearby fibers, including its
monodromy action. Local acyclicity makes the vanishing-cycle
complex zero. Exact periodization preserves the triangle and gives
the corresponding criterion for the periodic realization. Finally
the finite regular model supplies horizontal transport and base
change by Theorem~\ref{thm:finite-cyclic}; the
assumed compatibility identifies that transport with the indicated
geometric endpoints.
\end{proof}

The existence of $\mathcal P$ and the identification of $i^*\mathcal P$
are genuine hypotheses. A symmetric obstruction theory or an equality
of Euler characteristics does not construct either one. For a
nonproper moduli problem, the analogous compact-support statement
requires a chosen compactification and a compatible proper-support
base-change assertion; it is not obtained by interchanging $R p_!$
and nearby cycles without that verification. In the explicit
nonproper central loci below the comparison is constructed inside
a displayed proper parameter space, so this issue is specified.

For the complete section-complex loci, the moduli family is a
projective bundle and the relative oriented DT complex is the
shifted constant complex. The central and nearby identifications
are constructed directly in Section~\ref{sec:v6-hyperplanes}.
For the point and curve ideals of Section~\ref{sec:divisor-sheaves},
we separately determine the generic invariant, the central locus
that has actually been constructed, and the specialization map.
This order prevents flatness on $\Delta^*$ from being substituted
for the missing central comparison.

\subsection{Retaining the singular boundary and its attachment maps}
\begin{proposition}\label{prop:boundary-localization}
Let $M$ be a finite-type oriented d-critical space with DT complex
$\mathcal P$, let $j:U\hookrightarrow M$ be an open stratum and
$i:Z\hookrightarrow M$ its closed complement.  The full DT
cohomology is computed by the triangle
\[
 R\Gamma_c(U,j^*\mathcal P)\longrightarrow R\Gamma_c(M,\mathcal P)
 \longrightarrow R\Gamma_c(Z,i^*\mathcal P)\xrightarrow{+1}.
\]
Thus the groups of the two strata determine the full groups only
after the connecting maps have been determined.  For a finite
closed filtration the repeated triangles give a finite spectral
sequence with the compactly supported cohomology of the restricted
DT complexes on the strata as its first-page terms.
\end{proposition}
\begin{proof}
Apply $R\Gamma_c$ to
$j_!j^*\mathcal P\to\mathcal P\to i_*i^*\mathcal P\to$.
Iteration for the closed filtration gives the asserted finite
filtered complex and its spectral sequence.  Its filtration is
finite, so convergence involves no inverse-limit assertion.
The attachment differentials are those of these triangles.
\end{proof}

In particular, the DT complex on a boundary is $i^*\mathcal P$.
It is not automatically the DT complex of a separately chosen
critical structure on that boundary.  Similarly, a sum of local
Milnor numbers is not a computation of $R\Gamma_c(M,\mathcal P)$.
Section~\ref{sec:global-two-defect} determines the attachment
extensions and computes a global singular example explicitly.
For a degeneration, Proposition~\ref{prop:specialization-cone-test}
applies in addition to these spatial attachment maps: it tests the
variation in the geometric parameter.

\section{Cohomological comparison for hyperplane-section sheaves and their limits as perfect complexes}\label{sec:v6-hyperplanes}

We begin the proof of Theorem~\ref{thm:intro-v6} with the hyperplane
case, constructing every part of its comparison. The extension to all
positive degrees is proved in Section~\ref{ssec:all-positive-degrees}. The proof first identifies
the scalar-rigidified moduli locus as $\P^4\times C$ and reconstructs
its Artinian deformations. It then identifies the central derived
restrictions and their orientation, constructs the global bar/BV
comparison, and calculates its horizontal cohomology. The normalized
ordinary potential is zero on this locus, so the connection is constant;
the main work is to prove that the same complex represents both
geometric endpoints, including the limiting perfect complex at $0$.

Fix the resolved quintic family of \cite[Proposition~2.10]{KS}, over a
connected affine open $C\subset\A^1$ containing $0$. We shrink $C$ so that all
fibers other than $X_0$ are smooth. Let
\[
 a:\cX\longrightarrow\P^4,\qquad
 L=a^*\O_{\P^4}(1),\qquad V=H^0(\P^4,\O_{\P^4}(1)).
\]
Thus $\dim V=5$. The letter $h\in V$ denotes the equation of the central hyperplane, whereas $\ell\in V\setminus\{0\}$ denotes a variable
hyperplane equation. Write $M=\P(V)$ for the projective space of
one-dimensional subspaces of $V$, and $\mathscr M=M\times C$.
The symbol $M$ in this section denotes this particular moduli locus, not the
moduli space of all objects of a prescribed Chern character.

The objects that complete the family at $0$ are
\begin{equation}\label{eq:v6-cones}
 E_{\ell,b}=\bigl[L_b^{-1}\xrightarrow{\ell}\O_{X_b}\bigr],
 \qquad L_b=L|_{X_b},
\end{equation}
in cohomological degrees $-1,0$. A change of nonzero scalar in $\ell$ gives
an isomorphic complex. On a smooth fiber this is the structure sheaf of its
hyperplane divisor. On $X_0$, the complex is retained even if multiplication
by $\ell$ is not injective. This convention will be essential at $\ell=h$.
All restrictions to $X_1,X_2,S$ are derived restrictions of this displayed
bounded vector-bundle complex.

\textbf{Smooth classical parameters and non-smooth derived moduli.}
A statement that $M_b=\mathbb P^4$ is smooth refers to its classical
scheme after scalar rigidification. It does not say that its reduced
derived DT enhancement is smooth. In minimal coordinates
$x_1,\ldots,x_d$ on a reduced critical chart that enhancement has
\[
 K_f=(\mathbb C[[x_1,\ldots,x_d]][\xi_1,\ldots,\xi_d],
             \delta\xi_i=\partial_{x_i}f),\qquad |\xi_i|=-1.
\]
When $f=0$, its tangent complex at a point is
$[\mathbb C^d\xrightarrow{0}(\mathbb C^d)^\vee]$ in degrees $0,1$.
In particular the obstruction space is not removed. The full Artin
object stack also retains scalar automorphisms and their dual degree-three
self-extensions, as explained in Section~\ref{ssec:sc-oriented-moduli}.
Quasi-smoothness of the reduced enhancement is compatible with either a
smooth or a singular classical truncation.

Here the vanishing of the potential is a conclusion, not a smoothness
hypothesis imposed on all perfect-complex moduli. The reconstruction below
proves that the classical completed local ring is regular of dimension
$d=\dim\operatorname{Ext}^1(E,E)$. In a minimal critical chart it is
$\mathbb C[[x_1,\ldots,x_d]]/(\partial f)$. A nonzero ideal in this
regular local domain lowers the dimension, whereas the quotient here
still has dimension $d$. Thus $(\partial f)=0$. In characteristic zero
$f$ is constant, and normalization on the critical locus gives $f=0$.
Exactly this argument applies to the Grassmannian loci after the explicit
matrix calculation in Theorem~\ref{thm:grassmannian-moduli}. It does not
apply to the singular examples in Sections~\ref{sec:singular-comparison}
and~\ref{sec:cubic-geometry}, where nonzero critical potentials are
computed instead. Adding a nondegenerate quadratic stabilization to $0$
would not create a nonzero d-critical structure.

\textbf{The five steps of this example.}
The classical moduli on a smooth fiber is $M_b=\P^4$ and its reduced
potential is $f_b=0$, proved by reconstructing its infinitesimal
deformations below. With orientation $K_{\P^4}$ its DT groups have
one generator in degrees $-4,-2,0,2,4$.
At $0$ we restrict the actual universal cone to $X_1,X_2,S$;
its matching locus is again $\P^4$, including the exceptional
non-sheaf restriction at $[h]$. The oriented two-sided bar of
that matching locus is compared with
$(\Omega_{\P^4}^{\bullet}((u)),u d)[4]$ by constructed chain maps.
Its cohomology therefore has the same five generators. Finally
all maps are defined on $\P^4\times C$, and the connection is
$\nabla_{\partial_t}=\partial_t$ in the basis of hyperplane powers.
The proof of these assertions is Theorem~\ref{thm:v6-comparison};
it verifies the central identification as well as flatness.
The later point examples are different moduli problems and do
not inherit this constant connection merely because their
local potential is also zero.

\subsection{Construction of the resolved family and its polarization}
\label{ssec:quintic-construction}
Choose general forms $h,q_4,q_5$ of degrees $1,4,5$ on $\P^4$.
Assume that the quartic threefold $(q_4=0)$, the quartic surface
$S=(h=q_4=0)$, and the complete-intersection curve
$\Gamma=(h=q_4=q_5=0)$ are smooth. The family before resolution is
\begin{equation}\label{eq:sc-quintic-pencil}
 \cX'=\{hq_4=tq_5\}\subset\P^4\times C,
\end{equation}
where $C\subset\A^1_t$ is a Zariski open neighborhood of $0$ on
which the other fibers are smooth. Blowing up the ideal $(h,t)$
and taking the strict transform gives
\[
 \pi:\cX\longrightarrow C,\qquad
 X_1=(q_4=0),\qquad X_2=\Bl_\Gamma\P^3,
 \quad \P^3=(h=0).
\]
This is the construction of \cite[Proposition~2.10]{KS}.
We recall its local verification because the central restrictions
of the universal perfect complexes use these actual components.
At $\Gamma$, coordinates $a=h$, $b=q_4$, $c=q_5$ and a coordinate
along the curve give $ab=tc$. The two charts of the small resolution
are
\begin{equation}\label{eq:sc-quintic-charts}
 t=a\tau,\ b=\tau c,
 \qquad\hbox{and}\qquad
 a=t\alpha,\ c=\alpha b.
\end{equation}
Both charts are smooth. The first has base map $t=a\tau$; its
central branches identify with the quartic and the blowup of
$\P^3$. Away from $\Gamma$ the modification is an isomorphism.
The exceptional locus has codimension two in the fourfold, so
the relative volume form on the anticanonical hypersurface
\eqref{eq:sc-quintic-pencil} extends over the smooth resolution.
Consequently
\begin{equation}\label{eq:sc-quintic-volume}
 \omega_{\cX/C}\simeq\O_{\cX},\qquad
 \omega_{\cX}\simeq\pi^*\omega_C.
\end{equation}

Write $a:\cX\to\P^4$ and $L=a^*\O_{\P^4}(1)$.
On $X_2$ use $H_2$ for the hyperplane pullback and $E_\Gamma$
for the exceptional divisor. The strict transform of $S$ has
class $4H_2-E_\Gamma$ and
\[
 K_{X_2}=-4H_2+E_\Gamma=-S,\qquad K_{X_1}=-H_1=-S.
\]
Because $\Gamma$ is a Cartier divisor of class $5h_S$ on $S$,
the surface itself is unchanged by this blowup. Its normal line
bundles are $\O_S(h_S)$ and $\O_S(-h_S)$.
Here $h_S=H_1|_S=H_2|_S$, with $h_S^2=4$.

An extending polarization useful for the surface-sheaf problem is
\begin{equation}\label{eq:sc-extending-polarization}
 \begin{gathered}
 \mathcal A=L^{\otimes2}\otimes\O_{\cX}(X_1),\qquad
 \mathcal A|_{X_b}=\O_{X_b}(2),\\
 \mathcal A|_{X_1}=\O_{X_1}(H_1),\qquad
 \mathcal A|_{X_2}=\O_{X_2}(6H_2-E_\Gamma).
 \end{gathered}
\end{equation}
Indeed $\O_{\cX}(X_1)|_{X_1}=\O_{X_1}(-S)$ and its
restriction to $X_2$ is $\O_{X_2}(S)$, which gives the formulas.
On the common intersection surface $S=X_1\cap X_2$, both restrictions
are $\O_S(h_S)$. The surface $S$ is a smooth quartic K3 surface, of
complex dimension two; it is a reduced divisor in each threefold $X_i$.
In particular, $S$ is not a nonreduced surface of multiplicity two.
It is the common boundary on which the two component objects are matched.
The bundle $\O(6H_2-E_\Gamma)$ is ample: the ideal of the
$(4,5)$ curve twisted by $5$ is generated by $q_5$ and the
linear multiples of $q_4$, giving a projective embedding of the
blowup over $\P^3$; the tensor product of its tautological bundle
with the pullback of $\O_{\P^3}(1)$ is the indicated bundle.
The restriction on $X_1$ is ample as well. 

\subsection{Numerical classes of component objects and the boundary}
\label{ssec:quintic-component-polynomials}
The sheaf problem starts with a pure two-dimensional sheaf $F_b$ on
a smooth fiber $X_b$, scheme-theoretically supported on a divisor
$D_b\in|H|$. A central limit in the chosen perfect-complex moduli
problem is denoted $E_0\in\Perf(X_0)$. Its component and common-boundary
objects are
\[
 E_i=Lj_i^*E_0\quad\hbox{on }X_i,\qquad
 E_S=L\iota_i^*E_i=Lj_S^*E_0\quad\hbox{on }S.
\]
On the locus where $E_i$ is concentrated in degree zero, write $E_i=F_i$.
Thus this locus consists of component sheaves $F_1,F_2$ with an
equivalence of their \emph{derived} restrictions to $S$. The boundary
complex may have both $\operatorname{Tor}_1^{\O_{X_i}}(F_i,\O_S)$ in
degree $-1$ and $F_i\otimes\O_S$ in degree zero. Neither term is discarded.
The proper section-complex completion also has points at which an
$E_i$ is not a sheaf; Proposition~\ref{prop:v6-exceptional-cone} describes such
a point explicitly. A statement about the component-sheaf locus and
a statement about this larger completion therefore have different scopes.

For the charge calculation assume
\[
 \ch(E_i)=(0,H_i,\beta_i,n_i[\mathrm{pt}]),\qquad
 \ch_0(E_S)=0,\quad \ch_1(E_S)=h_S,\quad
 c_S:=\int_S\ch_2(E_S).
\]
The symbol $[\mathrm{pt}]$ denotes the cohomology class of a point,
normalized to have integral one on the variety in that expression.
It lies in $H^6(X_i,\mathbb Q)$ for a threefold and in
$H^4(S,\mathbb Q)$ for the surface. Thus one could write
$\ch(E_S)=(0,h_S,c_S[\mathrm{pt}])$; the definition by the integral
specifies the degree without an extra subscript. On a smooth fiber
$\ch(F_b)=(0,H,\beta,n[\mathrm{pt}])$. For a complex the Hilbert
polynomial is the alternating Euler characteristic. The formulas
below apply to both the component-sheaf locus and its indicated
perfect-complex completion. Restriction of the Chern character gives
\begin{equation}\label{eq:sc-boundary-charge}
 \beta_1H_1=\beta_2(4H_2-E_\Gamma)=c_S.
\end{equation}
Riemann--Roch gives
\begin{align}
 P_{2H}(E_b;m)&=10m^2+2(H\beta)m+n+\tfrac{25}{6},\label{eq:sc-component-polynomials}\\
 P_{H_1}(E_1;m)&=2m^2+(\beta_1H_1+2)m
                  +n_1+\tfrac12\beta_1H_1+\tfrac73,\notag\\
 P_{6H_2-E_\Gamma}(E_2;m)&=8m^2+
       (\beta_2(6H_2-E_\Gamma)+2)m\notag\\
 &\hspace{9mm}+n_2+\beta_2(2H_2-\tfrac12E_\Gamma)+\tfrac{11}{6},\notag\\
 P_{h_S}(E_S;m)&=4m+c_S.\notag
\end{align}
Here is a direct verification of the constants and leading terms.
For the quintic, $H^3=5$ and $c_2=10H^2$. For the quartic
threefold, $H_1^3=4$, $c_1=H_1$ and $c_2=6H_1^2$.
On $X_2$ one has
$H_2^3=1$, $H_2^2E_\Gamma=0$, and
$H_2E_\Gamma^2=-20$. A generic member
$B\in|H_2|$ is $\Bl_{20}\P^2$, so
$c_1(X_2)|_B=4h-\sum e_a$ and $c_2(X_2)|_B=26[\mathrm{pt}]_B$.
It follows that
$\int_B\operatorname{td}_2(X_2)=((-4)+26)/12=11/6$.
Also $H_2(6H_2-E_\Gamma)^2/2=8$ and
$H_2(6H_2-E_\Gamma)c_1(X_2)/2=2$.
These identities prove every term of
\eqref{eq:sc-component-polynomials}.

The normalization triangle implies
\[
 P_{\mathcal A}(E_0)=P_{H_1}(E_1)+
          P_{6H_2-E_\Gamma}(E_2)-P_{h_S}(E_S).
\]
For an extending family this equals $P_{2H}(E_b)$.
Comparing coefficients and using \eqref{eq:sc-boundary-charge}
gives
\begin{equation}\label{eq:sc-charge-conservation}
 H\beta=\beta_1H_1+\beta_2H_2,\qquad n=n_1+n_2.
\end{equation}
The boundary correction is present in the polynomial identity
even though it cancels the displayed degree-six corrections in
this particular numerical expression. Neither this identity nor
\eqref{eq:sc-charge-conservation} identifies the DT complexes;
that comparison is constructed on the moduli loci below.

\subsection{Cohomology on all fibers and the universal family}
Fix a point $b\in C$, which may equal $0$. An infinitesimal base
change at $b$ means a local Artinian $\CC$-algebra $A$ with residue
field $\CC$, together with a morphism $\operatorname{Spec}A\to C$
whose closed point maps to $b$. We write
\[
 X_A=\cX\times_C\operatorname{Spec}A,\qquad
 L_A=L|_{X_A},\qquad \O_{X_A}=\text{the structure sheaf of }X_A.
\]
The ordinary fiber product agrees with the derived one here because
$\pi$ is flat. A perfect deformation $E_A$ of $E_{\ell,b}$ is a perfect
complex on $X_A$, with a specified equivalence
$E_A\otimes_A^L\CC\simeq E_{\ell,b}$. If $A=\CC[\epsilon]/(\epsilon^2)$,
the base map may send $t$ to $b+a\epsilon$. It can therefore vary both
the object and its threefold fiber. At $b=0$ its residue object lies
on the singular union $X_0$; at $b\ne0$ it lies on the smooth quintic
$X_b$. The letter $A$ here is not a component of $X_0$.

The reconstruction will lift the distinguished section to a map
$\O_{X_A}\to E_A$ and identify its homotopy fiber with $L_A^{-1}$.
We first compute the cohomology and self-extensions on the
\emph{whole fiber} $X_b$. In particular, when $b=0$ every entry of
\eqref{eq:v6-cohom-table} is cohomology on $X_1\cup_S X_2$,
not cohomology on a single component.

\begin{lemma}\label{lem:v6-line-cohom}
With the chosen relative volume form, the nonzero cohomology groups on every
fiber, including $X_0$, are
\begin{equation}\label{eq:v6-cohom-table}
 \begin{array}{c|cccc}
 &H^0(X_b,-)&H^1(X_b,-)&H^2(X_b,-)&H^3(X_b,-)\\ \hline
 \O_{X_b}&\CC&0&0&\CC\\
 L_b&V&0&0&0\\
 L_b^{-1}&0&0&0&V^\vee .
 \end{array}
\end{equation}
The restriction map $V\to H^0(X_b,L_b)$ is the indicated isomorphism. The
corresponding relative cohomology sheaves are locally free, commute with base
change, and have these constant descriptions over $C$.
\end{lemma}
\begin{proof}
For $b\ne0$, tensor the hypersurface sequence in $\P^4$ by
$\O_{\P^4}(m)$, for $m=0,1,-1$:
\[
 0\longrightarrow\O_{\P^4}(m-5)
 \longrightarrow\O_{\P^4}(m)
 \longrightarrow\O_{X_b}(m)\longrightarrow0.
\]
The intermediate cohomology of line bundles on $\P^4$ vanishes.
For $m=0$, the remaining top cohomology is
$H^4(\P^4,\O(-5))=\CC$; for $m=1$ there is no top cohomology;
for $m=-1$ it is $H^4(\P^4,\O(-6))=V^\vee$. This proves the table
on the smooth fibers, including the asserted restriction of sections.

For the central fiber use the normalization sequence, tensoring by
$L_0^m$:
\begin{equation}\label{eq:v6-normalized-L}
 0\longrightarrow L_0^m\longrightarrow
 j_{1*}L_1^m\oplus j_{2*}L_2^m
 \xrightarrow{\mathrm{res}_1-\mathrm{res}_2}
 j_{S*}L_S^m\longrightarrow0.
\end{equation}
Here $L_1=\O_{X_1}(1)$, $L_2=\beta^*\O_{\P^3}(1)$ for
$\beta:X_2=\Bl_\Gamma\P^3\to\P^3$, and $L_S=\O_S(1)$.
For comparison, the component and boundary groups used in this
sequence are the following; a tuple lists $H^0,H^1,H^2,H^3$ on
the variety in the first column:
\[
\begin{array}{c|ccc}
 &\O&L&L^{-1}\\\hline
 X_1 &(\CC,0,0,0)&(\CC^5,0,0,0)&(0,0,0,\CC)\\
 X_2 &(\CC,0,0,0)&(\CC^4,0,0,0)&(0,0,0,0)\\
 S &(\CC,0,\CC,0)&(\CC^4,0,0,0)&(0,0,\CC^4,0).
\end{array}
\]
The whole-fiber table follows by taking the long exact cohomology
sequence, including its connecting maps. It is not the direct sum
of the two component tables.

Here $\R\beta_*\O_{X_2}=\O_{\P^3}$ can be checked on the blowup
charts. On an affine neighborhood of the smooth codimension-two center
write its ideal as $(a,b)$ in a regular ring $A$. The two charts and
their intersection have coordinate rings
\[
 A[v]/(av-b),\qquad A[w]/(bw-a),\qquad
 A[v,v^{-1}]/(av-b),\quad w=v^{-1}.
\]
Every Laurent polynomial on the intersection is the sum of its
nonnegative powers of $v$, coming from the first chart, and its negative
powers, coming from the second. Thus the degree-one Cech cokernel is
zero. There is no higher Cech term, since these are two affine charts
with affine intersection. Their common regular functions are $A$:
regarded in the fraction field they are regular outside the
codimension-two center, and normality of $A$ extends them across that
center. This proves the claimed derived pushforward on each affine
neighborhood. Projection formula now computes the cohomology of
$L_2^m$ on $\P^3$.

The quartic hypersurface sequence gives
\[
 H^{>0}(X_1,\O_{X_1})=0,\qquad H^{>0}(X_1,L_1)=0.
\]
The same vanishings hold on $X_2$. On the quartic surface $S$, its hypersurface
sequence in $\P^3$ gives
\[
 H^0(\O_S)=\CC,\quad H^1(\O_S)=0,\quad H^2(\O_S)=\CC,
 \qquad H^0(L_S)=\CC^4,\quad H^{>0}(L_S)=0.
\]
For $m=0$, the map on degree-zero cohomology in
\eqref{eq:v6-normalized-L} is $(c_1,c_2)\mapsto c_1-c_2$.
It is surjective and its kernel is $\CC$. The only remaining connecting
isomorphism is
\[
 H^2(S,\O_S)\xrightarrow{\sim}H^3(X_0,\O_{X_0}).
\]
For $m=1$, the map on global sections is
\[
 V\oplus V/\CC h\longrightarrow V/\CC h,
 \qquad(v,\bar w)\longmapsto\bar v-\bar w.
\]
Its kernel is $V$ via $v\mapsto(v,\bar v)$ and its cokernel is zero.
This gives the middle row of \eqref{eq:v6-cohom-table}. Finally,
$\omega_{X_0}\simeq\O_{X_0}$ and Gorenstein Serre duality applied to the
line bundle $L_0$ gives the last row. This uses perfectness of $L_0$,
not smoothness of $X_0$.

All dimensions in the table are constant. Cohomology and base change for the
proper flat family and its invertible sheaves therefore identifies their
relative cohomology sheaves with vector bundles. The map
$V\otimes\O_C\to\pi_*L$ is an isomorphism on every fiber and hence an
isomorphism. The relative trace trivializes $R^3\pi_*\O_{\cX}$, and its
pairing with $\pi_*L$ identifies $R^3\pi_*L^{-1}$ with
$V^\vee\otimes\O_C$. These maps give the asserted descriptions and commute
with base change.
\end{proof}

On $\cX\times M$, the tautological inclusion
$\O_M(-1)\hookrightarrow V\otimes\O_M$ defines the universal perfect complex
\begin{equation}\label{eq:v6-universal}
 \mathcal E=
 \bigl[L^{-1}\boxtimes\O_M(-1)\longrightarrow\O_{\cX\times M}\bigr].
\end{equation}
The structure morphism of this product to $\mathscr M=M\times C$ is
$\pi\times\id_M$, with the factors ordered in the evident way.
Each term of \eqref{eq:v6-universal} is locally free on $\cX\times M$ and
flat over $\mathscr M$. Consequently it is a relative perfect family, even
at points where its differential has a kernel after restriction to a fiber.
The support of each fiber complex is proper since $X_b$ is projective.

The determinant of this universal family also records an important
moduli distinction. Applying the alternating determinant formula in \cite[Section~2.4.3]{KS} to its two
terms gives
\begin{equation}\label{eq:v7-universal-determinant}
 \det(\mathcal E)
    =\det(\O_{\cX\times M})\otimes
       \det\bigl(L^{-1}\boxtimes\O_M(-1)\bigr)^{-1}
    =L\boxtimes\O_M(1).
\end{equation}
At a fixed $[\ell]$, its restriction is isomorphic to $L_b$ after
choosing a basis of the one-dimensional parameter factor. Over all of
$M$, the factor $\O_M(1)$ is retained. Thus the construction of
$\mathscr M$ uses scalar-rigidified simple objects, not a globally
chosen isomorphism $\det(\mathcal E)\simeq p_{\cX}^*L$.
Indeed, restricting such an isomorphism to $\{x\}\times M$ would
trivialize $\O_M(1)$, whose restriction to a projective line has degree
one. This is a direct check that the determinant constraint and the
scalar rigidification used below are not interchangeable.

\begin{lemma}\label{lem:v6-Ext}
For every nonzero $\ell$ and every $b$, one has
\begin{equation}\label{eq:v6-Ext-table}
 \begin{aligned}
 \ext^0(E_{\ell,b},E_{\ell,b})&=\CC, &
 \ext^1(E_{\ell,b},E_{\ell,b})&=V/\CC\ell,\\
 \ext^2(E_{\ell,b},E_{\ell,b})&=\{\lambda\in V^\vee:\lambda(\ell)=0\}, &
 \ext^3(E_{\ell,b},E_{\ell,b})&=\CC.
 \end{aligned}
\end{equation}
The middle identifications use a choice of generator of $\CC\ell$.
Intrinsically the middle groups are the tangent and cotangent spaces of
$M$ at $[\ell]$. All other self-Ext groups vanish. In addition,
\begin{equation}\label{eq:v6-Hom-from-O}
 H^0(\R\Gamma E_{\ell,b})=\CC,
 \qquad H^1(\R\Gamma E_{\ell,b})=0,
 \qquad H^2(\R\Gamma E_{\ell,b})\simeq\CC^4,
\end{equation}
and no other hypercohomology group is nonzero.
\end{lemma}
\begin{proof}
Write $P^{-1}=L_b^{-1}$, $P^0=\O_{X_b}$ and $d_P=\ell$.
The sheaf endomorphism complex of $P$ is
\begin{equation}\label{eq:v6-sheaf-End}
 L_b^{-1}\xrightarrow{w\mapsto(w\ell,\ell w)}
 \O_{X_b}\oplus\O_{X_b}
 \xrightarrow{(a,c)\mapsto(a-c)\ell}L_b,
\end{equation}
in degrees $-1,0,1$. These are exactly the signs in
$d_{\End}(f)=d_Pf-(-1)^{|f|}fd_P$.
Use the spectral sequence with
$E_1^{p,q}=H^q(X_b,\mathcal Hom^p(P,P))$.
By \eqref{eq:v6-cohom-table} its only rows are
\begin{equation}\label{eq:v6-two-rows}
 \begin{array}{c|ccc}
 q=3&V^\vee&\xrightarrow{\lambda\mapsto(\lambda(\ell),\lambda(\ell))}
                       \CC^2&0\\
 q=0&0&\CC^2\xrightarrow{(a,c)\mapsto(a-c)\ell}&V.
 \end{array}
\end{equation}
The first map has rank one since $\ell\ne0$, as does the second map.
After their cohomology is taken, no higher differential can connect the
two rows: such a differential would have to move from $q=3$ to $q=0$ and
increase $p$ by four, whereas the columns range only from $-1$ to $1$.
Thus \eqref{eq:v6-Ext-table} follows. The diagonal kernel in degree zero
is the identity endomorphism; the degree-three quotient is paired with it
by supertrace $(a,c)\mapsto c-a$.

In families, the lower row of \eqref{eq:v6-two-rows} is the Euler sequence
\[
 0\longrightarrow\O_M\longrightarrow V\otimes\O_M(1)
 \longrightarrow T_M\longrightarrow0,
\]
after separating the diagonal scalar summand. The upper row gives its dual
\[
 0\longrightarrow\Omega_M^1\longrightarrow V^\vee\otimes\O_M(-1)
 \longrightarrow\O_M\longrightarrow0.
\]
This proves the intrinsic identifications, including the perfect tangent--
obstruction pairing. They describe the cohomology sheaves; a global splitting
of the complete endomorphism complex is not used.

Apply $\R\Gamma$ instead to the triangle
$L_b^{-1}\xrightarrow{\ell}\O_{X_b}\to E_{\ell,b}$.
The map on third cohomology is the surjection
$V^\vee\to\CC$, $\lambda\mapsto\lambda(\ell)$, by Serre duality.
The degree-zero map starts with $0\to\CC$. The long exact sequence
therefore gives exactly \eqref{eq:v6-Hom-from-O}.
\end{proof}

\textbf{Singular supports are included in the cohomology calculation.}
The vanishing table above is on $X_b$, not on its supporting divisor.
Neither it nor the self-endomorphism spectral sequence assumes that
the divisor $(\ell=0)$ is nonsingular.  The argument below therefore
includes singular hyperplanes and, in higher degree, reducible and
nonreduced divisors.  At $b=0$ it also includes sections whose
restriction to a component is zero.  The distinction is between
singularity of the support and singularity of the classical moduli
of its perfect complex.

\subsection{Every classical deformation remains a cone}
The uniform Ext dimensions do not yet identify a moduli space. The
next proof recovers the defining section from an arbitrary Artinian
deformation and thereby proves the open immersion of
$\P(V)\times C$ into the scalar-rigidified moduli space.

We use the algebraic space obtained by
rigidifying scalar automorphisms of simple universally gluable perfect
complexes. Its existence follows from Lieblich's representability and
rigidification results \cite[Theorem 4.2.1, Corollary 4.3.3]{Lieblich2006}.
Here universal gluability means vanishing of negative global self-Ext; it does
not mean that the cohomology sheaves of the complex vanish in negative degrees.
Lemma~\ref{lem:v6-Ext} verifies the former condition at every point of our
family.

\begin{proposition}\label{prop:v6-projective-locus}
The family \eqref{eq:v6-universal} induces an open immersion
\[
 j:\mathscr M=M\times C\longrightarrow
 \mathfrak{Perf}(\cX/C)^{\mathrm{simp,cl}}\sslash\Gm
\]
into the scalar-rigidified classical simple-perfect-complex moduli space.
In particular its image is a smooth proper moduli locus over $C$, with
fiber $\P^4$. Every infinitesimal deformation of one of its objects,
including deformations of the point of $C$, is given by deforming its
nonzero section $\ell$, modulo multiplication by a unit.
\end{proposition}
\begin{proof}
We give the deformation argument over local Artinian rings in detail.
Let $A$ be such a ring with residue field $\CC$, let
$\Spec A\to C$ be a base deformation of $b$, and let $E_A$ be a perfect
deformation of $E_{\ell,b}$ on $X_A=\cX\times_C\Spec A$.
Proper perfect pushforward makes $\R\Gamma(X_A,E_A)$ a perfect
$A$-complex. By \eqref{eq:v6-Hom-from-O}, its derived residue fiber has
cohomology only in degrees $0,2$, of dimensions $1,4$.
A minimal finite free representative over the local ring $A$ consequently
has one generator in degree zero, four in degree two, and no other
terms. Indeed, in a bounded finite free representative cancel all
invertible entries of the differential; its remaining reduction has zero
differential and its ranks are the residue cohomology dimensions.
There are now no adjacent nonzero degrees, so the remaining differential
itself is zero. In particular
\begin{equation}\label{eq:v6-H0-base}
 H^0(X_A,E_A)\simeq A,
 \qquad H^1(X_A,E_A)=0.
\end{equation}
Choose a generator lifting the canonical map on the residue fiber. It gives
$q_A:\O_{X_A}\to E_A$. Define its homotopy fiber $K_A$ by the triangle
\begin{equation}\label{eq:v6-lift-triangle}
 K_A\xrightarrow{s_A}\O_{X_A}\xrightarrow{q_A}E_A
 \longrightarrow K_A[1].
\end{equation}
The derived residue fiber of $K_A$ is $L_b^{-1}$.
A perfect complex whose derived fiber at every residue point is a vector
space of dimension one in degree zero is locally free of rank one in degree
zero: at a local ring cancel the invertible differential entries in a
finite free representative, and reduce to its residue field. Only one
rank-one term remains. Applied to $K_A$, this shows that $K_A$ is a line
bundle on $X_A$.

It remains to identify that line bundle. For a small extension
$A'\twoheadrightarrow A$ with kernel $J$ killed by the maximal ideal,
the multiplicative sequence is
\[
 1\longrightarrow 1+J\O_{X_b}\longrightarrow\O_{X_{A'}}^*
 \longrightarrow\O_{X_A}^*\longrightarrow1.
\]
The first sheaf is additively $J\otimes_\CC\O_{X_b}$.
Since $H^1(X_b,\O_{X_b})=0$, the induced map
$\Pic(X_{A'})\to\Pic(X_A)$ is injective. The known bundle
$L_{A'}^{-1}$ supplies a lift. Induction on the length of $A$ therefore
identifies $K_A\simeq L_A^{-1}$. With this choice,
\eqref{eq:v6-lift-triangle} writes $E_A$ as
\[
 [L_A^{-1}\xrightarrow{\ell_A}\O_{X_A}],
 \qquad\ell_A\in V\otimes_\CC A,
 \qquad\ell_A\bmod\mathfrak m_A=\ell.
\]
The vector $\ell_A$ generates a direct summand of rank one since one
of its coordinates is a unit.

The choices are exactly projective rescalings. A different generator in
\eqref{eq:v6-H0-base} multiplies $q_A$ by a unit, and a different
identification $K_A\simeq L_A^{-1}$ also multiplies by a unit because
$H^0(X_A,\O_{X_A})=A$. Conversely, an isomorphism between two such
complexes sends their generators $\O_{X_A}\to E_A$ to scalar multiples.
Taking homotopy fibers gives a diagram of triangles
\begin{equation}\label{eq:v6-cone-isom-diagram}
\begin{tikzcd}[column sep=small]
 L_A^{-1}\arrow[r,"\ell_A"]\arrow[d,"u"']&
 \O_{X_A}\arrow[r]\arrow[d,"v"]&E_{\ell_A,A}\arrow[r]\arrow[d,"\simeq"]&
 L_A^{-1}[1]\arrow[d,"{u[1]}"]\\
 L_A^{-1}\arrow[r,"\ell'_A"']&\O_{X_A}\arrow[r]&E_{\ell'_A,A}\arrow[r]&L_A^{-1}[1]
\end{tikzcd}
\end{equation}
with $u,v\in A^*$ and $v\ell_A=u\ell'_A$.
Thus the two sections give the same point of $\P(V)(A)$.
This proves the equivalence of the infinitesimal deformation functors
after scalar rigidification, not just the equality of tangent dimensions.

The family map is locally of finite presentation into the cited algebraic
space. The infinitesimal calculation makes it formally etale at every
point and therefore etale. The same triangle argument over any
algebraically closed extension of $\CC$ shows that distinct projective
sections give distinct objects. Hence this etale morphism is universally
injective and is an open immersion. Its source is $\P(V)\times C$, which
is smooth and proper over $C$. This proves all assertions.
\end{proof}

For $b\ne0$, every hyperplane divisor of the smooth quintic is integral.
Indeed the Grothendieck--Lefschetz theorem \cite{SGA2} gives
$\Pic(X_b)=\Z[\O_{X_b}(1)]$. Each nonzero effective divisor has positive
class in this group, so a divisor of class $\O_{X_b}(1)$ can neither split
nor have multiplicity greater than one. Its Cartier structure is
Cohen--Macaulay, hence has no embedded components. A nonzero proper subsheaf
of $\O_D$ has generic rank one on this integral divisor; its quotient has
smaller dimension and positive Hilbert polynomial. Its reduced Hilbert
polynomial is therefore strictly smaller. Thus these smooth-fiber objects
are Gieseker-stable sheaves. The completion specified in
\eqref{eq:v6-cones} is a completion by simple perfect complexes; no
particular polynomial stability condition on its central objects is imposed.

\subsection{Derived restrictions of the hyperplane-section complexes}
Having identified the family of moduli spaces, we determine its point
in the derived intersection on $X_0$. The complex with $\ell=h$ is
particularly important: it supplies a limit in the proper family,
but its restriction to $X_2$ has two nonzero terms. The following
calculation includes that point and the matching equivalence on $S$.

\begin{proposition}\label{prop:v6-exceptional-cone}
At $b=0$, the section $h$ gives the perfect complex
$E_{h,0}=[L_0^{-1}\xrightarrow{h}\O_{X_0}]$. Its derived restrictions are
\begin{equation}\label{eq:v6-exceptional-restrictions}
 \begin{aligned}
 E_{h,0}|_{X_1}^{\mathbf L}&\simeq\O_S,\\
 E_{h,0}|_{X_2}^{\mathbf L}&\simeq\O_{X_2}\oplus L_2^{-1}[1],\\
 E_{h,0}|_S^{\mathbf L}&\simeq\O_S\oplus L_S^{-1}[1].
 \end{aligned}
\end{equation}
The two restrictions to $S$ are identified by the identity on their
bounded two-term representatives. On $X_0$ one has
\begin{equation}\label{eq:v6-negative-central}
 \mathcal H^{-1}(E_{h,0})
     \simeq j_{2*}\bigl(L_2^{-1}(-S)\bigr)\ne0.
\end{equation}
Nevertheless $E_{h,0}$ is simple and has no negative global self-Ext.
\end{proposition}
\begin{proof}
On $X_1$, $h$ is the defining section of $S$ and $L_1=\O_{X_1}(S)$.
Its two-term complex is the usual resolution of $\O_S$.
On $X_2$, which maps into the hyperplane $\{h=0\}$, the
section $h$ is identically zero. Thus the restricted complex has zero
differential and splits as asserted. The same holds on $S$.
Because the displayed complex consists of line bundles, termwise
restriction computes each derived restriction. In particular the
nonzero degree $-1$ summand on $S$ is retained.

To compute the kernel of multiplication by $h$ on $X_0$, use
\eqref{eq:v6-normalized-L} for $m=-1$. A local section is a pair
$(u_1,u_2)$ with equal restrictions to $S$. If $h(u_1,u_2)=0$,
then $u_1=0$, since $h$ is a non-zero-divisor on the smooth integral
$X_1$. The matching condition becomes $u_2|_S=0$, namely
$u_2\in L_2^{-1}(-S)$. This proves
\eqref{eq:v6-negative-central}. Simplicity and the vanishing of
negative global self-Ext were computed in Lemma~\ref{lem:v6-Ext} without
any assumption that $\ell$ is a non-zero-divisor.
\end{proof}

The second object in \eqref{eq:v6-exceptional-restrictions} has two
nonzero direct summands. It is not a stable sheaf or a simple complex
on $X_2$. It is, however, a valid point of $\mathfrak{Perf}(X_2)$.
Thus the projective completion just constructed uses the full component
stacks in \eqref{eq:intro-gluing}. Restricting both component stacks
to stable sheaves would remove this explicitly identified point and would
be a different moduli problem.

\begin{theorem}\label{thm:v6-gluing-locus}
Derived restriction of \eqref{eq:v6-universal} at $0$ gives a family of
matching triples in
\[
 \mathfrak G=\mathfrak{Perf}(X_1)
       \times^h_{\mathfrak{Perf}(S)}\mathfrak{Perf}(X_2).
\]
The corresponding scalar-rigidified classical locus is $M=\P(V)$.
Its formal deformation problem at each point is the actual deformation
problem of $E_{\ell,0}$; in particular the locus is smooth of dimension
four. The oriented d-critical structure on this locus is $(M,0)$, with
orientation $K_M$. It agrees with the oriented smooth-fiber d-critical
locus through the same zero critical charts on $M\times C$.
\end{theorem}
\begin{proof}
There is a strict commutative diagram of complexes
\begin{equation}\label{eq:v6-restriction-square}
\begin{tikzcd}[column sep=small]
 \bigl[L_1^{-1}\boxtimes\O_M(-1)\to\O_{X_1\times M}\bigr]
  \arrow[r,"\mathbf L\iota_1^*"]&
 \bigl[L_S^{-1}\boxtimes\O_M(-1)\to\O_{S\times M}\bigr]\\
 \bigl[L_2^{-1}\boxtimes\O_M(-1)\to\O_{X_2\times M}\bigr]
  \arrow[r,"\mathbf L\iota_2^*"']&
 \bigl[L_S^{-1}\boxtimes\O_M(-1)\to\O_{S\times M}\bigr]
  \arrow[u,equal]
\end{tikzcd}
\end{equation}
The arrows are termwise restriction and the right vertical map is the
identity. The displayed complexes are defined on the respective spaces;
the diagram indicates their restriction functors. They give the
matching family and all its parameter-base changes.

For the deformation comparison let $P=E_{\ell,0}$ be its two-term
representative. Tensoring the normalization sequence by the complex
$\mathcal Hom^\bullet(P,P)$, whose terms are vector bundles, gives the
termwise exact sequence of complexes
\begin{equation}\label{eq:v6-End-normalization}
 \begin{aligned}
 0\longrightarrow\mathcal Hom^\bullet(P,P)
 &\longrightarrow j_{1*}\mathcal Hom^\bullet(P_1,P_1)
       \oplus j_{2*}\mathcal Hom^\bullet(P_2,P_2)\\
 &\longrightarrow j_{S*}\mathcal Hom^\bullet(P_S,P_S)
 \longrightarrow0.
 \end{aligned}
\end{equation}
The last map is the difference of restrictions. Its differential is
$[d_P,-]$, so the sequence is compatible with the deformation differentials.
After derived global sections, the left complex is quasi-isomorphic
to the homotopy fiber of the right-hand restriction map. For explicit
cone conventions, if $r:C_1\oplus C_2\to C_S$ is that map, the homotopy
fiber has terms
\[
 F^n=C_1^n\oplus C_2^n\oplus C_S^{n-1},
 \quad d(x_1,x_2,z)=(dx_1,dx_2,r(x_1,x_2)-dz).
\]
The natural map from the global endomorphism complex sends $x$ to
$(x|_{X_1},x|_{X_2},0)$ in a compatible resolution. Sequence
\eqref{eq:v6-End-normalization} proves that it is a quasi-isomorphism.
The homotopy-fiber deformation Lie algebra is taken with its induced
homotopy Lie operations, as in \cite[Proposition~3.3]{KS};
the cone display here describes its underlying complex. Restriction is
a morphism of the deformation problems and induces this quasi-isomorphism,
so their characteristic-zero formal deformation problems agree.

Proposition~\ref{prop:v6-projective-locus} now identifies their classical
formal neighborhoods with those of $\P(V)$. Restriction is fully faithful,
and the map is locally of finite presentation on the loci under
consideration; equivalently one may apply the infinitesimal criterion
point by point on a finite-presentation atlas. It gives an etale
monomorphism on the simple scalar-rigidified loci. Thus the indicated
classical locus in the gluing stack is $M$ with its actual smooth
scheme structure, not merely a family mapping into a potentially thicker
moduli space.

The universal complex also fixes the treatment of scalars. On this
classical simple locus the scalar gerbe is neutral: its objects over a
test scheme are the pullbacks of $\mathcal E$ tensored with a line bundle
from the test scheme. Thus the gerbe is $M\times B\Gm$, and
$M\to M\times B\Gm$ supplied by $\mathcal E$ is a smooth atlas of
relative dimension one. We use the d-critical structure on this atlas
and its scalar-rigidified DT normalization. In the stack formulation,
pullback of the stack perverse sheaf is shifted by one for this smooth
atlas, as in \cite[Theorem 4.8]{BenBassatBravBussiJoyce2015ShiftedDarboux}.
The cohomology groups below are those on $M$, not the cohomology of
$M\times B\Gm$.

The boundary forms and the derived Lagrangian intersection give its
d-critical structure by \cite{PTVV13,BenBassatBravBussiJoyce2015ShiftedDarboux}.
A smooth scheme has $\mathcal S_M^0=0$: using its identity embedding,
the defining sequence identifies $\mathcal S_M$ with the kernel of
$d:\O_M\to\Omega_M^1$, and the summand with zero restriction to
$M^{\mathrm{red}}=M$ is zero. Hence this d-critical structure is
necessarily zero \cite[Example 2.15]{Joyce}.
The identity critical chart $(U,U,0,\id)$ is therefore available on
every smooth coordinate open $U\subset M$.
Its virtual canonical bundle is $K_U^{\otimes2}$, and these
identifications glue to $K_{M,0}=K_M^{\otimes2}$.
Choose the square root $K_M$, or relatively $K_{\mathscr M/C}$.
Exactly the same reasoning applies on each smooth fiber of the moduli
family. Thus the chosen orientation and the zero critical charts agree
throughout $M\times C$. This concludes the comparison at the level of
oriented d-critical loci.
\end{proof}

\subsection{Explicit critical charts, coordinate changes, and fourfold tangents}
The central matching triples are now identified with the same smooth
classical locus as the generic-fiber objects. We next give its critical
charts and orientation transitions. Their explicit compatibility is
what allows the bar/BV maps to be defined globally, rather than chosen
independently on each fiber.

For completeness, the zero potential also follows from the actual
minimal critical equations on each fixed fiber. Use the minimal Artin
Darboux atlas of \cite[Theorem 2.10, Corollary 2.11]{BenBassatBravBussiJoyce2015ShiftedDarboux}.
The simple stack has a one-dimensional scalar stabilizer. Its minimal
smooth atlas has four classical tangent coordinates; the additional
negative generator paired with the scalar is removed by the Darboux
reduction, not by asserting a trace splitting for rank-zero objects.
The resulting reduced formal critical chart at $[\ell]$ is
\[
 \widehat A_b=\CC[[x_1,\ldots,x_4]][\xi_1,\ldots,\xi_4],
 \qquad\delta\xi_i=\partial_{x_i}f_b.
\]
Its classical ring is smooth of dimension four by the deformation proof.
The surjection $\CC[[x_1,\ldots,x_4]]\to H^0(\widehat A_b)$ induces an
isomorphism on tangent spaces and has regular target of the same dimension.
It is an isomorphism: after choosing formal coordinates on the target,
its composite with that coordinate isomorphism has invertible linear part,
and formal inversion gives zero kernel. Thus every $\partial_{x_i}f_b$
is zero. In characteristic zero $f_b$ is constant; the critical-value
normalization makes it zero. This identifies the normalized fiberwise
minimal potential from the actual deformation problem.

For the family of DT charts we use the smooth relative d-critical
presentation already established, namely $(M\times C,0)$ relative to $C$.
Its explicit critical algebra on a coordinate open is
\[
 \O_C[x_1,\ldots,x_4][\xi_1,\ldots,\xi_4],\qquad\delta=0.
\]
These charts present the oriented d-critical locus and its DT complex.
The following transition maps define the oriented DT-chart family.
The scalar automorphisms belong to the full derived moduli stack,
whereas these charts are on its scalar-rigidified classical locus.

Let $s_0,\ldots,s_4$ be the coordinates on $V$ in a fixed basis.
On $U_i=\{s_i\ne0\}$, put $x_j=s_j/s_i$ for $j\ne i$.
On $U_i\cap U_k$, the other coordinates are
\begin{equation}\label{eq:v6-coordinate-change}
 y_i=x_k^{-1},\qquad
 y_j=x_j/x_k\quad(j\ne i,k).
\end{equation}
Let $J_{ki}=(\partial y_a/\partial x_b)$. The degree-$-1$
generators transform by
\begin{equation}\label{eq:v6-odd-change}
 \xi_a^{(k)}=\sum_b(J_{ki}^{-1})_{ba}\xi_b^{(i)}.
\end{equation}
These are the transition maps of $T_M[1]$. They preserve the zero
differential, and the chain rule gives
$J_{li}=J_{lk}J_{ki}$ on triple overlaps. Thus the maps in
\eqref{eq:v6-odd-change} satisfy the exact cocycle identity.
For the signed coordinate volume forms
\[
 \omega_i=(-1)^i\bigwedge_{j\ne i} dx_j,
\]
with the factors in increasing index order, the determinant calculation is
\begin{equation}\label{eq:v6-volume-change}
 \omega_k=x_k^{-5}\omega_i.
\end{equation}
Indeed the unsigned Jacobian is
$\det J_{ki}=(-1)^{i+k}x_k^{-5}$.
Squaring \eqref{eq:v6-volume-change} gives precisely the transition
function of $K_{M,0}$. Using \eqref{eq:v6-volume-change} itself as its
square root therefore fixes the orientation transition, with no
unspecified sign torsor. All these formulas are independent of $t$.

\begin{proposition}\label{prop:v6-ambient}
For $i_b:X_b\hookrightarrow\cX$, the normal-direction obstruction of
$E_{\ell,b}$ vanishes. The dimensions of
$\ext^k_{\cX}(i_{b*}E_{\ell,b},i_{b*}E_{\ell,b})$ are
\begin{equation}\label{eq:v6-ambient-dim}
 (1,5,8,5,1),\qquad k=0,1,2,3,4,
\end{equation}
and all other groups vanish. The classical ambient simple-complex
moduli problem is formally smooth of dimension five along this family.
On a minimal scalar-reduced $(-2)$-shifted Darboux chart its differential
is zero. A local shifted potential can be chosen to be zero.
\end{proposition}
\begin{proof}
The complex $[L^{-1}\xrightarrow{\ell}\O_{\cX}]$ extends
$E_{\ell,b}$ to the whole total space. More explicitly, on the open with
parameter $t$ there is a morphism of endomorphism complexes arising from
the globally defined parameter Koszul resolution
\[
 [\O_{\cX}\xrightarrow{t-b}\O_{\cX}].
\]
Tensor this with the extending two-term complex. When derived Hom into
$i_{b*}E_{\ell,b}$ is taken, multiplication by $t-b$ is zero. Adjunction
therefore gives, at this object and with this chosen extension,
\begin{equation}\label{eq:v6-ambient-split}
 \RH_{\cX}(i_{b*}E_{\ell,b},i_{b*}E_{\ell,b})
 \simeq\RH_{X_b}(E_{\ell,b},E_{\ell,b})
       \oplus\RH_{X_b}(E_{\ell,b},E_{\ell,b})[-1].
\end{equation}
This is a splitting for the specified extending object. It is not a
splitting assertion for every complex on a fiber. It also verifies
directly that its Atiyah--Kodaira--Spencer obstruction is zero.
Equation \eqref{eq:v6-ambient-dim} follows from
\eqref{eq:v6-Ext-table}.

The tangent sequence in degree one is
\[
 0\longrightarrow\ext^1_{X_b}(E,E)
 \longrightarrow\ext^1_{\cX}(i_{b*}E,i_{b*}E)
 \longrightarrow\End_{X_b}(E)\longrightarrow0.
\]
The derivative of the $M$-direction is the first summand. The derivative
of the base displacement $b\mapsto b+\epsilon$ maps to $\id_E$ in
the last summand: differentiate the Koszul differential $t-b$ to see
this edge map. Therefore the family $M\times C$ has an isomorphism on
five-dimensional classical tangent spaces to the ambient simple-complex
moduli space.

Here is the formal smoothness argument, including its use of the actual
family. Let $A$ be a complete local ring of the scalar-rigidified ambient
space at this object. It is noetherian and has embedding dimension five,
so there is a surjection $\CC[[z_1,\ldots,z_5]]\twoheadrightarrow A$.
The family gives a map $A\to\CC[[w_1,\ldots,w_5]]$ inducing an isomorphism
on tangent spaces. The composite power-series map has invertible linear
part and hence is an isomorphism, by recursively solving its homogeneous
terms. Its kernel is zero; the kernel of the first surjection is therefore
zero as well. Thus $A$ is the smooth power-series ring and the family
identifies the formal neighborhoods.

The total space has a chosen fourfold volume form, so its proper-support
perfect-object moduli stack is $(-2)$-shifted symplectic
\cite[Corollary 6.2]{BD2}. Apply the Artin Darboux reduction, retaining
and then removing the scalar degree and its paired dual as in
Section~\ref{sec:local}. The reduced minimal chart has five degree-zero
coordinates $z_j$, eight degree-$-1$ coordinates $y_a$, and five
degree-$-2$ coordinates $p_j$. In strong Darboux form its equations have
\[
 \delta y_a=s_a(z),\qquad
 \delta p_j=\sum_a c_a y_a\partial_{z_j}s_a(z),
\]
with fixed nonzero normalization constants $c_a$ after the quadratic
form has been diagonalized. The master equation supplies the usual
quadratic relation among the $s_a$. Since its classical local ring is
already smooth of dimension five in minimal tangent coordinates, the
ideal $(s_1,\ldots,s_8)$ is zero. Hence every $s_a$ is zero, and both
sets of differentials vanish. Over $\CC$, split the rank-eight
nondegenerate symmetric obstruction space into two dual rank-four
subspaces. The resulting algebra is
\[
 \CC[[z_1,\ldots,z_5]][\theta_1,\ldots,\theta_4,
 \eta_1,\ldots,\eta_4,p_1,\ldots,p_5],\qquad\delta=0,
\]
with $|\theta|=|\eta|=-1$, $|p|=-2$. It is the shifted critical algebra
of the zero function of degree $-1$ on
$\CC[[z_1,\ldots,z_5]][\theta_1,\ldots,\theta_4]$.
This proves the stated local assertion.
\end{proof}

For this locus, the normal spaces have dimensions one and four,
as calculated in \eqref{eq:r13-hyperplane-dimensions}. The
coefficient interpretation of the real quotient is given in
Proposition~\ref{prop:r13-real-coefficients}; its section and
cosection vanish in the present zero-equation chart.
The real negative-definite subspace in this chart can also
be written explicitly. Diagonalize the complex obstruction form as
$q(y)=\sum_{a=1}^8 y_a^2$ and write $y=u+iv$. Then
\[
 \Re q=\sum_a(u_a^2-v_a^2),\qquad
 \Im q=2\sum_a u_av_a.
\]
The real subspace $u=0$ is negative definite for $\Re q$ and is
Lagrangian for the bilinear form underlying $\Im q$. Its coordinate
generators commute and the internal differential is zero. This is
the local negative-definite choice used in the BKSY construction. The
holomorphic polarization in the preceding proof and this real negative
subspace are different choices, and the equations in both local constructions vanish
in this chart. The comparison below uses the explicit ordinary zero
critical charts; it does not require identifying a prescribed global
BKSY quotient or extending the coherent-sheaf Quot theorem to the
non-sheaf object \eqref{eq:v6-negative-central}.

\subsection{An oriented global two-sided bar complex}
The zero critical potential and the orientation $K_M$ reduce the local
models to the zero section of $T^*M$. We now globalize their two-sided
bar resolution. The projection to symmetric degree one and
antisymmetrization define $\Pi$ and $I$ without a coordinate choice;
contraction with $K_M$ supplies the differential. These maps give the
central identification used in the comparison theorem.

We construct the DT bar complex of the smooth d-critical locus just
identified. Set $T=T_M$, $K_M=\bigwedge^4\Omega_M^1$, and
\[
 \mathcal A=\Sym_{\O_M}T,
 \qquad\epsilon:\mathcal A\longrightarrow\O_M,
 \qquad\overline{\mathcal A}=\ker\epsilon.
\]
The relative spectrum of $\mathcal A$ is $T^*M$ and $\epsilon$ is its
zero section. The two equal augmentations give the derived
self-intersection, which is the zero critical presentation of $(M,0)$.
Its normalized two-sided bar complex over $\O_M$ has
\begin{equation}\label{eq:v6-bar}
 \mathcal B^{-q}=\overline{\mathcal A}^{\otimes_{\O_M}q},
 \qquad
 b[a_1|\cdots|a_q]=
 \sum_{j=1}^{q-1}(-1)^j
 [a_1|\cdots|a_ja_{j+1}|\cdots|a_q].
\end{equation}
The terms at the two ends vanish because every $a_j$ is in the
augmentation ideal. Tensor coefficients with $K_M$ and write
$\mathcal B^{\mathrm{or}}=\mathcal B\otimes_{\O_M}K_M$.
Its finite Koszul comparison is
\[
 \mathcal K^{\mathrm{or}}=\bigoplus_{q=0}^4
        \bigwedge^q T\otimes K_M,
\]
in cohomological degrees $-q$, initially with zero differential.
The tensor with $K_M$ implements the chosen orientation of the critical
locus; it does not require a square root of $K_M$ itself.

\begin{proposition}\label{prop:v6-global-bar}
There are globally defined comparison maps
\[
 \mathcal K^{\mathrm{or}}\xrightarrow{I}\mathcal B^{\mathrm{or}}
          \xrightarrow{\Pi}\mathcal K^{\mathrm{or}},
 \qquad\Pi I=1,
\]
which are quasi-isomorphisms for the differentials $0,b$.
Let $\Delta$ be the operator on $\mathcal K^{\mathrm{or}}$
defined by signed contraction with top forms and the exterior derivative.
Then
\begin{equation}\label{eq:v6-periodic-bar}
 \mathcal B^{\mathrm{or}}_{\mathrm{per}}
   =\bigl(\mathcal B^{\mathrm{or}}((u)),
                  D_B=-b+uI\Delta\Pi\bigr)
\end{equation}
is a well-defined periodic complex, and there is a global quasi-isomorphism
\begin{equation}\label{eq:v6-global-I-Pi}
 \mathcal I\Pi:\mathcal B^{\mathrm{or}}_{\mathrm{per}}
 \xrightarrow{\sim}
 \bigl(\Omega_M^\bullet((u))[4],u\,d\bigr).
\end{equation}
The maps and differentials are compatible with the projective-chart
coordinate changes and pull back to $M\times C$ independently of $t$.
\end{proposition}
\begin{proof}
For local vector fields $v_1,\ldots,v_q$, define
\begin{equation}\label{eq:v6-I-Pi-formulas}
 \begin{aligned}
 I(v_1\wedge\cdots\wedge v_q)
     &=\sum_{\sigma\in\mathfrak S_q}\mathrm{sgn}(\sigma)
                 [v_{\sigma(1)}|\cdots|v_{\sigma(q)}],\\
 \Pi[a_1|\cdots|a_q]
     &=\frac1{q!}\pi_1(a_1)\wedge\cdots\wedge\pi_1(a_q),
 \end{aligned}
\end{equation}
where $\pi_1:\overline{\mathcal A}\to T$ is projection onto
symmetric degree one. Both maps tensor with the identity on $K_M$.
They are intrinsic to the symmetric algebra of the vector bundle $T$.
The commutativity of its multiplication cancels the terms of $bI$ in
pairs. In each summand of $b$, a product of two positive symmetric degrees
has degree at least two; hence $\Pi b=0$. Antisymmetrizing and then
wedging multiplies by $q!$, proving $\Pi I=1$.

To verify that these are quasi-isomorphisms, locally write
$\mathcal A=\O_M[p_1,\ldots,p_4]$.
Its zero-section module has the free Koszul resolution
\[
 \bigl(\mathcal A[\eta_1,\ldots,\eta_4],
                d_K=\sum_i p_i\partial_{\eta_i}\bigr).
\]
Exactness follows directly from the Euler contraction: on a monomial
of polynomial degree $r$ and exterior degree $q$, with $r+q>0$, put
\[
 H_K=\frac1{r+q}\sum_i\eta_i\partial_{p_i}.
\]
The graded anticommutator is $d_KH_K+H_Kd_K=1$ in those positive
weights. In weight zero the augmentation remains. The one-sided
bar resolution has the usual augmentation contraction and is likewise
a flat resolution. The antisymmetrization map from the Koszul resolution
to the one-sided bar resolution lifts the identity of the augmentation;
its first bar face is exactly $d_K$, and the internal faces cancel in
pairs. Thus it is a comparison of flat resolutions. Derived tensoring
with the second zero-section module gives the map $I$ in
\eqref{eq:v6-I-Pi-formulas}, proving that it is a quasi-isomorphism.
Since $\Pi I=1$, $\Pi$ is a quasi-isomorphism as well.

For a local coordinate volume form $\omega$ define
\begin{equation}\label{eq:v6-contraction}
 \mathcal I_q(v_1\wedge\cdots\wedge v_q\otimes\omega)
 =(-1)^{q(q-1)/2}\iota_{v_q}\cdots\iota_{v_1}\omega.
\end{equation}
It maps degree $-q$ to degree $(4-q)-4$ in $\Omega_M^\bullet[4]$.
Contraction of tensor fields is natural, so these maps glue to a
bundle isomorphism $\mathcal I$ globally. Define
$\Delta=\mathcal I^{-1}d\mathcal I$. In coordinates, for
$P(x,\xi)\otimes\omega$, it is the explicit operator
\begin{equation}\label{eq:v6-Delta}
 \Delta(P\otimes\omega)
    =\left(\sum_i\partial_{x_i}\partial_{\xi_i}P\right)\otimes\omega,
\end{equation}
where left odd derivatives are used and the coefficient of $\omega$
is included in $P$. For verification on a monomial of exterior degree
$q$, the usual unsigned contraction obeys
$d\mathcal I_q^{\mathrm{uns}}=(-1)^{q-1}
\mathcal I_{q-1}^{\mathrm{uns}}\sum_i\partial_{x_i}\partial_{\xi_i}$.
The factor in \eqref{eq:v6-contraction} removes exactly this sign.
This proves \eqref{eq:v6-Delta} and $\Delta^2=0$.

Put $B_{\mathrm{BV}}=I\Delta\Pi$. The preceding identities give
\[
 b^2=0,\quad bB_{\mathrm{BV}}=B_{\mathrm{BV}}b=0,
 \quad B_{\mathrm{BV}}^2=I\Delta(\Pi I)\Delta\Pi=0.
\]
They prove $D_B^2=0$, $D_BI=Iu\Delta$, and
$\Pi D_B=u\Delta\Pi$. In particular \eqref{eq:v6-global-I-Pi} is a
chain map. To check it is a quasi-isomorphism after Laurent-series
completion, work on a coordinate open and use the bar--Koszul homotopy
of Proposition~\ref{prop:v5-bar-Koszul} with both augmentations zero.
There is a normalized homotopy $h$ satisfying
\[
 bh+hb=1-I\Pi,\qquad hI=\Pi h=h^2=0.
\]
Its recursive construction on the free bar basis is given there before
base change. In the present case the projection in that construction
is precisely \eqref{eq:v6-I-Pi-formulas}: each polynomial factor of
degree greater than one has zero first derivative at the zero section,
and the ordered-simplex integral of $1$ is $1/q!$.
The same equations imply
\[
 D_B(-h)+(-h)D_B=1-I\Pi
\]
because the two perturbation terms vanish. Extend this homotopy
coefficientwise in $u$. It neither introduces new negative powers of
$u$ nor changes their lower bound. Thus it proves the local completed
quasi-isomorphism, and hence the global sheaf quasi-isomorphism.

Finally \eqref{eq:v6-odd-change} is the transformation of $T$, and
\eqref{eq:v6-volume-change} is the transformation of $K_M$.
Naturality of contraction gives
$\mathcal I\Theta_{ki}=\varphi_{ki}^*\mathcal I$ on overlaps;
naturality of symmetric degree and antisymmetrization gives the same
identity for $I,\Pi$. Thus $\Delta$ and $D_B$ are intertwined.
The chain rule and the volume determinant identity give the exact
triple-overlap condition. Since none of these maps contains $t$,
their pullbacks to the family commute with $\partial_t$.
\end{proof}

In particular, the critical-chart choices, their orientation, and the
periodic bar comparison are all global in this example. The complex
\eqref{eq:v6-periodic-bar} is an oriented bar/BV presentation of the DT
complex of the locus in Theorem~\ref{thm:v6-gluing-locus}.
It is not being identified with the unmodified structure-sheaf Tor
complex: the operator $uI\Delta\Pi$ is included in its differential.

\subsection{Finiteness, horizontal specialization, and the computed invariant}
The bar complex has now been compared with de Rham forms by global
maps independent of $b$. It remains to prove finiteness and calculate
specialization. The Euler sequences on $M=\P^4$ give five cohomology
classes, and the same representatives on $M\times C$ are horizontal.
This completes the comparison by an explicit basis, including at $0$.

The perverse-sheaf normalization on a four-dimensional smooth critical
locus is $\CC_M[4]$. Indeed for the zero function the nearby-cycle
complex is zero (there are no nonzero fibers), and the defining
vanishing-cycle triangle gives the perverse-normalized vanishing
cycles as $\CC_M[4]$. The orientation $K_M$
just constructed makes its transition maps the ordinary transition maps
of $\CC_M[4]$, by the oriented vanishing-cycle theorem
\cite[Theorem 6.9]{BBDJS}. Write this DT sheaf as $\mathcal P_b$ on
the copy of $M$ labeling the fiber over $b$.

In the categorical statement below we use the category of the specified
global critical presentation,
\[
 \mathscr C_b^{\mathrm{crit}}=\MF(M,0).
\]
The superscript records this choice. The cohomological DT complex is
fixed by the oriented d-critical structure $(M,0,K_M)$, and the
comparison computes its periodic realization. Every use of the category in
the following theorem refers to this specified global critical presentation;
the geometric comparison of the actual moduli locus is the oriented
d-critical comparison of Theorem~\ref{thm:v6-gluing-locus}.

For the periodic de Rham convention, multiplication by $u^p$ on
$\Omega_M^p$ is the chain isomorphism from the differential $d$ to
$u d$. It commutes with the Cech differential and with $\partial_t$.
The analytic de Rham comparison therefore identifies
\eqref{eq:v6-global-I-Pi} with the periodic realization of
$\CC_M[4]$. We use the full algebraic derived global sections on the
proper $M$ for the calculation below; no nonproper critical locus or
formal exchange of two unrelated completions is involved.

\begin{theorem}\label{thm:v6-comparison}
For the moduli locus $\mathscr M=\P(V)\times C$ of
Proposition~\ref{prop:v6-projective-locus}, with orientation
$K_{\mathscr M/C}$, the following assertions hold.
\begin{enumerate}[label=(\roman*)]
\item The special fiber is the actual locus of derived matching triples
of \eqref{eq:v6-restriction-square}, including the object at $[h]$.
The cohomological DT groups are
\begin{equation}\label{eq:v6-DT-groups}
 \mathbb H^k(M,\mathcal P_b)=H^{k+4}(\P^4,\CC)
   =\begin{cases}
        \CC,&k=-4,-2,0,2,4,\\
        0,&\text{otherwise}.
     \end{cases}
\end{equation}
They are finite for every $b\in C$.
\item For the global ordinary critical presentation $(M,0)$, define
$\mathscr C_b^{\mathrm{crit}}:=\MF(M,0)$. This is the category of the
chosen critical presentation. Put $K=\CC((u))$. There is a specified
horizontal comparison
\begin{equation}\label{eq:v6-main-HP}
 \begin{aligned}
 \HP_{\bar i}(\MF(M,0))
 &\xrightarrow{\sim}
 \mathbb H^{\bar i}(M,\mathcal B^{\mathrm{or}}_{\mathrm{per}})\\
 &\simeq\bigoplus_{k\equiv\bar i\ (2)}
                 \mathbb H^k(M,\mathcal P_0)\otimes_\CC K.
 \end{aligned}
\end{equation}
In particular its even rank is five and its odd rank is zero.
\item Over every affine open $W\subset C$, these groups form the
free module of rank five over
$\Gamma(W,\O_C)[[u]][u^{-1}]$ in even parity. It has the regular
flat connection $\nabla_{\partial_t}=\partial_t$ in the basis of
hyperplane powers. The central-fiber comparison is obtained by
specializing the global chain maps in
\eqref{eq:v6-global-I-Pi}; it is the specialization of those same chain maps.
\end{enumerate}
\end{theorem}
\begin{proof}
The identification of the parameter spaces, their d-critical structures,
and the central matching objects was proved in
Theorem~\ref{thm:v6-gluing-locus}. The preceding vanishing-cycle
calculation gives $\mathcal P_b=\CC_M[4]$. We now compute all its
cohomology, and the corresponding periodic cohomology, explicitly.

The Euler sequence and its exterior powers on $\P(V)$ give
\begin{equation}\label{eq:v6-Euler-p}
 0\longrightarrow\Omega_M^p\longrightarrow
 \bigwedge^p V^\vee\otimes\O_M(-p)
 \longrightarrow\Omega_M^{p-1}\longrightarrow0,
 \qquad 1\le p\le4.
\end{equation}
For $1\le p\le4$ the line bundle $\O_M(-p)$ has no cohomology.
Starting with $H^0(M,\O_M)=\CC$ and $H^{>0}(M,\O_M)=0$,
the long exact sequences in \eqref{eq:v6-Euler-p} give inductively
\begin{equation}\label{eq:v6-Hodge}
 H^q(M,\Omega_M^p)=
   \begin{cases}\CC,&q=p,\\0,&q\ne p.\end{cases}
\end{equation}
The hypercohomology spectral sequence of the de Rham complex has
$E_1^{p,q}=H^q(M,\Omega_M^p)$. Every possibly nonzero term lies on
the diagonal $p=q$. A differential of bidegree $(r,1-r)$ could
connect two such positions only if $r=1-r$, which is impossible for
an integer $r\ge1$. The sequence therefore degenerates and gives
one cohomology class in each degree $0,2,4,6,8$. This proves
\eqref{eq:v6-DT-groups} with the shift by four.

The classes can be chosen independently of the base. On the standard
affine cover of $M$, the Cech one-cocycle with values in one-forms
\[
 \alpha_{ij}=d\log(s_j/s_i)
\]
represents the de Rham first Chern class $H$ of $\O_M(1)$, with a
fixed sign convention for its transition functions.
The identity $\alpha_{ij}+\alpha_{jk}=\alpha_{ik}$ follows from
multiplication of transition functions, and $d\alpha_{ij}=0$.
Its cup powers represent $1,H,H^2,H^3,H^4$. To check nonvanishing
explicitly, restrict to the coordinate $\P^p\subset\P^4$. On its
intersection $U_0\cap\cdots\cap U_p$, put $x_j=s_j/s_0$.
The ordered cup representative is
\[
 \alpha_{01}\wedge\alpha_{12}\wedge\cdots\wedge\alpha_{p-1,p}
 =\frac{dx_1}{x_1}\wedge\cdots\wedge\frac{dx_p}{x_p}.
\]
The equality follows by expanding $d\log(x_j/x_{j-1})$ successively;
all terms with a repeated factor vanish. The Cech trace on
$H^p(\P^p,\Omega_{\P^p}^p)$ takes the coefficient of
$x_1^{-1}\cdots x_p^{-1}dx_1\wedge\cdots\wedge dx_p$.
A representative extending across any omitted coordinate chart has zero
such coefficient, so this functional annihilates the Cech coboundaries.
Its value on the displayed cocycle is one. Thus every indicated cup
power is nonzero and spans the corresponding one-dimensional summand.
For the differential $u d$ multiply a representative of form degree
$p$ by $u^p$. This produces explicit cocycles, still independent of
$t$, for the five periodic classes. Pulling them back through
$\mathcal I^{-1}$ and $I$ gives their bar/BV representatives by
\eqref{eq:v6-I-Pi-formulas} and \eqref{eq:v6-contraction}.

We justify the coefficient completion as well as the finite ranks.
Let $C^\bullet$ be the total Cech--de Rham complex on the finite standard
affine cover of $M$, initially with differential $d_{\mathrm{tot}}$
using $d$ rather than $u d$. It is a bounded complex of $\CC$-vector
spaces. Its terms may be infinite dimensional; its cohomology has the
five bases just computed. Let $B^n=\operatorname{im}d_{\mathrm{tot}}^{n-1}$,
choose the indicated cocycles as a complement $H^n$ to $B^n$ in the
cycles, and choose a complement $W^n$ to the cycles in $C^n$. Then
\[
 C^n=B^n\oplus H^n\oplus W^n,\qquad
 d_{\mathrm{tot}}:W^n\xrightarrow{\sim}B^{n+1}.
\]
Define $i$ by the chosen cocycles, $p$ by projection onto $H^n$, and
$h$ on $B^n$ by the inverse of the displayed differential, with $h=0$
on $H^n\oplus W^n$. Directly on the three summands these maps satisfy
\[
 pi=1,\qquad d_{\mathrm{tot}}h+hd_{\mathrm{tot}}=1-ip.
\]
They are fixed $\CC$-linear maps, independent of the parameter.
All sheaves and maps on $M\times W$ are obtained by extending their
coefficients to $R_W=\Gamma(W,\O_C)$; the same identities hold there.
The maps can be applied coefficientwise to power series in $u$.
They preserve the lower Laurent bound. Conjugating by the scaling
$u^p$ on form degree $p$ gives the contraction for $u d$.
Since $0\le p\le4$, conjugation changes Laurent exponents by a bounded
amount. It therefore defines maps on the specified Laurent-series
completion, not only on a larger iterated completion.

The result is a finite representative with five free even generators,
zero odd term, and zero differential over $R_W[[u]][u^{-1}]$.
Its inclusion is given by the scaled hyperplane-power cocycles above.
All contraction identities commute with evaluation at any $b\in W$,
including $0$. This proves finite freeness and the required base-change
isomorphisms. The finite representative and its zero connection matrix
are already defined over $R_W\otimes_\CC K$ for $K=\CC((u))$;
passing to $R_W[[u]][u^{-1}]$ is extension of this same finite module.
Thus the construction also gives an algebraic rank-five bundle on
$W\times_\CC\Spec K$.

On the product $M\times C$ the ordinary function is zero, so the
Gauss--Manin operator is
\[
 \partial_t-u^{-1}\partial_t0=\partial_t.
\]
It commutes with the Cech, de Rham, and bar differentials. The
comparison maps are independent of $t$, hence the diagram
\begin{equation}\label{eq:v6-horizontal-final}
\begin{tikzcd}[column sep=small]
 \mathcal B^{\mathrm{or}}_{\mathrm{per},W}
   \arrow[r,"\Pi"]\arrow[d,"\partial_t"']&
 \mathcal K^{\mathrm{or}}_W((u))
   \arrow[r,"\mathcal I"]\arrow[d,"\partial_t"]&
 \Omega_{M\times W/W}^\bullet((u))[4]
   \arrow[d,"\partial_t"]\\
 \mathcal B^{\mathrm{or}}_{\mathrm{per},W}
   \arrow[r,"\Pi"']&\mathcal K^{\mathrm{or}}_W((u))
   \arrow[r,"\mathcal I"']&
 \Omega_{M\times W/W}^\bullet((u))[4]
\end{tikzcd}
\end{equation}
commutes as a diagram of complexes. Thus its connection matrix in
$1,\allowbreak H,\allowbreak H^2,\allowbreak H^3,\allowbreak H^4$ is zero. Evaluation of all its arrows at $0$ is
allowed by the base-change calculation, constructing the central
comparison. Parallel transport sends each hyperplane-power basis
class to the corresponding class on the special fiber.

Finally apply the ordinary-function matrix-factorization comparison
\cite{Efimov2018CyclicMF} to the smooth projective pair $(M,0)$.
Its critical locus is $M$ and its critical value is zero, so the
comparison identifies its periodic cyclic homology with the same
de Rham hypercohomology. The normalization by $[4]$ has even parity
and changes none of the periodic parity assignments.
Composing with the inverse of \eqref{eq:v6-global-I-Pi} proves
\eqref{eq:v6-main-HP}. Every space, orientation, differential,
comparison arrow, and basis used in this proof has been constructed
above for the displayed family $E_{\ell,b}=[L_b^{-1}\xrightarrow{\ell}\O_{X_b}]$. No expansion of the
geometric fiber or of the sheaf-moduli problem occurs.
\end{proof}

This theorem concerns the hyperplane-section sheaves on the smooth
quintics and their explicit projective completion by perfect cones.
At the exceptional point the component restrictions are those in
\eqref{eq:v6-exceptional-restrictions}. In this example $f_b=0$
and the local shifted potential of the ambient Darboux model is also
zero; a formula $\nu=f_b g$ is therefore a vacuous factorization and
supplies no additional transport. The regular connection and the
specialization comparison are instead the concrete global constructions
of \eqref{eq:v6-horizontal-final}. For arbitrary polynomial-stable
objects, the finite-model and geometric-comparison hypotheses elsewhere
in the paper have not been verified by this calculation.

\subsection{Other fixed linear systems and line-bundle twists}
\label{ssec:r14-linear-systems}
The preceding calculation concerns the complete hyperplane linear
system. To verify that the finite comparison hypothesis is not
confined to that single class, we carry out the same moduli
construction for degree-two and degree-three divisors and for
line-bundle sheaves on these divisors. The ambient family, the
special fiber, and the restriction functors remain fixed.

\begin{theorem}\label{thm:r14-divisor-examples}
Let $\pi:\cX\to C$ and $L=a^*\O_{\P^4}(1)$ be the resolved
quintic family above. For $m\in\{1,2,3\}$ and $k\in\mathbb Z$,
put $V_m=H^0(\P^4,\O_{\P^4}(m))$, $M_m=\P(V_m)$ and
$d_m=\binom{m+4}{4}-1$. The complexes
\[
 E^{m,k}_{s,b}=[L_b^{k-m}\xrightarrow{s}L_b^k],
      \qquad [s]\in M_m,\quad b\in C,
\]
form a smooth proper locus $M_m\times C$ in the classical
scalar-rigidified moduli space of simple perfect complexes on
the fibers. Its central points are their derived restriction
triples in $\mathfrak X_1\times^h_{\mathfrak S}\mathfrak X_2$.
With the orientation $K_{M_m}$ its DT sheaf is
$\CC_{M_m}[d_m]$. There are horizontal quasi-isomorphisms
\[
 \CyclicChains_{\mathrm{per}}(\MF(M_m,0))
 \simeq R\Gamma(M_m,\Omega_{M_m}^\bullet((u))[d_m])
 \simeq R\Gamma(M_m,\mathcal B_{m,\mathrm{per}}^{\mathrm{or}})
\]
in the derived category of two-periodic complexes, using the
specified zero critical presentation. The right-hand complex
is the oriented bar/BV complex of this central moduli locus.
It is identified with the central fiber by the maps
\eqref{eq:r14-iota0-basis}.

These data satisfy Hypothesis~\ref{hyp:finite-cyclic}. The finite
model has even rank $r_m=d_m+1$, zero odd part and connection
$\partial_t$. The three ranks are $5,15,35$, respectively, and
\begin{equation}\label{eq:r14-divisor-DT-groups}
 \mathrm{DT}_{\mathrm{coh}}^j(M_{m,b})=
 \begin{cases}
  \CC,&j=2a-d_m,\quad 0\le a\le d_m,\\
  0,&\text{otherwise}.
 \end{cases}
\end{equation}
For $b\ne0$ these objects are pure two-dimensional sheaves on
members of $|L_b^m|$, and on the integral-divisor locus they are
Gieseker stable. The stated proper completion parametrizes the
perfect complexes displayed above, rather than imposing a
polynomial stability condition on every central point.
\end{theorem}
\begin{proof}
\emph{Cohomology and the universal object.}
First take $k=0$. For $1\le m\le3$ the hypersurface sequence on
$\P^4$ gives, on a smooth fiber,
\[
 H^0(X_b,L_b^m)=V_m,\qquad H^{>0}(X_b,L_b^m)=0.
\]
The upper bound $m\le3$ ensures in particular that
$H^0(\P^4,\O(m-5))=H^4(\P^4,\O(m-5))=0$.
On the central fiber tensor the normalization sequence by $L_0^m$:
\[
 0\longrightarrow L_0^m\longrightarrow
 j_{1*}L_1^m\oplus j_{2*}L_2^m
 \longrightarrow j_{S*}L_S^m\longrightarrow0.
\]
For the quartic component and the boundary quartic surface the
hypersurface sequences give
\[
 \begin{aligned}
 H^0(X_1,L_1^m)&=V_m,&H^{>0}(X_1,L_1^m)&=0,\\
 H^0(S,L_S^m)&=H^0(\P^3,\O(m)),&H^{>0}(S,L_S^m)&=0.
 \end{aligned}
\]
Here $m<4$ removes the degree-four defining equation from the
spaces of degree-$m$ sections. On
$X_2=\operatorname{Bl}_\Gamma\P^3$, projection formula and
$R\beta_*\O_{X_2}=\O_{\P^3}$ give the same cohomology as on
$\P^3$. The map on sections of the normalization sequence is
\[
 V_m\oplus H^0(\P^3,\O(m))\longrightarrow
 H^0(\P^3,\O(m)),\qquad (s,t)\longmapsto s|_{h=0}-t.
\]
It is surjective with kernel $V_m$. Thus the calculation on the
central fiber agrees with the one on each smooth fiber.
Serre duality using the chosen relative volume form gives,
for every $b$,
\begin{equation}\label{eq:r14-power-cohomology}
 \begin{array}{c|cccc}
 &H^0&H^1&H^2&H^3\\ \hline
 \O_{X_b}&\CC&0&0&\CC\\
 L_b^m&V_m&0&0&0\\
 L_b^{-m}&0&0&0&V_m^\vee.
 \end{array}
\end{equation}
The restriction maps and the trace identify the relative
cohomology sheaves with these constant vector bundles, by
cohomology and base change. In particular
$\pi_*L^m=V_m\otimes\O_C$.

On $\cX\times M_m$ the universal complex is
\[
 \mathcal E_m=
 [L^{-m}\boxtimes\O_{M_m}(-1)\longrightarrow\O_{\cX\times M_m}].
\]
Both terms are locally free and flat over $M_m\times C$.
Consequently this is a relative perfect family at all parameters,
including those at which the differential has a kernel.

\emph{The actual deformation space.}
The sheaf endomorphism complex is
\[
 L_b^{-m}\xrightarrow{w\mapsto(ws,sw)}\O_{X_b}^{\oplus2}
       \xrightarrow{(a,c)\mapsto(a-c)s}L_b^m.
\]
Its hypercohomology has the two rows of
\eqref{eq:v6-two-rows}, with $V$ replaced by $V_m$.
Both displayed nonzero maps have rank one. The columns run
from $-1$ to $1$, so a higher differential joining rows $3$ and
$0$ would require a nonexistent displacement of four columns.
It follows that
\[
 (\dim\ext^0,\dim\ext^1,\dim\ext^2,\dim\ext^3)
                   =(1,d_m,d_m,1),\qquad\ext^{<0}=0.
\]
The middle groups are intrinsically $T_{[s]}M_m$ and
$T_{[s]}^*M_m$, by the Euler sequence and its dual.
Applying derived global sections to
$L_b^{-m}\to\O_{X_b}\to E^{m,0}_{s,b}$ also gives
\[
 H^0(R\Gamma E^{m,0}_{s,b})=\CC,\quad
 H^1(R\Gamma E^{m,0}_{s,b})=0,\quad
 H^2(R\Gamma E^{m,0}_{s,b})=\ker(V_m^\vee\xrightarrow{s^\vee}\CC),
\]
and no other hypercohomology.

Let $A$ be a local Artinian base deformation and $E_A$ a perfect
deformation of this object. A minimal finite free representative
of $R\Gamma E_A$ has terms $A$ and $A^{d_m}$ only in degrees
$0$ and $2$, so its differential is zero. A generator in degree
zero defines $q_A:\O_{X_A}\to E_A$. The homotopy fiber $K_A$
is perfect and its derived residue fiber is $L_b^{-m}$.
At each local ring of $X_A$, a minimal finite free representative
therefore has one term of rank one in degree zero. Hence $K_A$
is a line bundle. The square-zero exact sequence for units and
$H^1(X_b,\O_{X_b})=0$ identify it with the existing lift
$L_A^{-m}$. Thus
\[
 E_A\simeq[L_A^{-m}\xrightarrow{s_A}\O_{X_A}],\qquad
 s_A\in V_m\otimes A,\quad s_A\bmod\mathfrak m_A\ne0.
\]
Changing the generator of $H^0(E_A)$ or the identification of
$K_A$ multiplies $s_A$ by a unit. Conversely, an isomorphism of
the complexes induces a diagram of these triangles and hence
such a rescaling, exactly as in
\eqref{eq:v6-cone-isom-diagram}. The infinitesimal functor after
scalar rigidification is therefore the functor of $M_m\times C$.
Lieblich's simple-complex algebraic space is locally of finite
presentation \cite{Lieblich2006}. The family morphism is
formally etale, and is injective on geometric points by the
same recovery of the section. It is an open immersion. Its
source is projective and smooth over $C$, proving the asserted
proper moduli locus.

\emph{Restriction to the special fiber and the orientation.}
Restricting the displayed universal two-term complex gives
\[
 E_{i,s}=[L_i^{-m}\xrightarrow{s_i}\O_{X_i}],\qquad
 E_{S,s}=[L_S^{-m}\xrightarrow{s_S}\O_S].
\]
Their iterated restrictions are identified by the identity of
this latter complex. All its degrees are retained. Tensoring
the normalization sequence by $E_s^\vee\otimes E_s$ and taking
sections gives
\[
 \RH_{X_0}(E_s,E_s)\xrightarrow{\sim}
 \fib\left(\bigoplus_i\RH_{X_i}(E_{i,s},E_{i,s})
                \longrightarrow\RH_S(E_{S,s},E_{S,s})\right).
\]
The last map is the difference of restrictions. It is the
endomorphism deformation complex of the matching triple.
The restriction functor gives the corresponding morphism of
deformation Lie algebras, and the displayed equivalence proves
that it is a quasi-isomorphism. Thus the central triples have
the same formal moduli locus $M_m$.

The relative orientation and the residue-compatible Lagrangian
construction supply the d-critical structures as in
Theorem~\ref{thm:v6-gluing-locus}. Since $M_m$ is smooth,
$\mathcal S^0_{M_m}=0$ \cite[Example 2.15]{Joyce}; these
structures are zero. Equivalently a minimal critical
presentation has $d_m$ classical variables and a smooth
classical quotient of dimension $d_m$, so its gradient ideal
is zero. Its normalized potential is $0$. The zero critical
charts identify
\[
 K_{M_m,0}=K_{M_m}^{\otimes2},\qquad
 K_{M_m}=\O_{M_m}(-r_m).
\]
Choose $K_{M_m}$ as the orientation. The coordinate Jacobians
and their determinants identify its local trivializations and
satisfy the cocycle equations. The oriented vanishing-cycle
gluing theorem \cite{BBDJS} then gives
$\mathcal P_{m,b}=\CC_{M_m}[d_m]$ with the same identifications
for all $b$.

\emph{The maps and the finite cohomology.}
Use $\mathcal A=\Sym T_{M_m}$ and its two zero augmentations.
The bar complex and the maps are
\[
 \begin{aligned}
 \mathcal B^{-q}&=(\ker(\mathcal A\to\O_{M_m}))^{\otimes q},\\
 I(v_1\wedge\cdots\wedge v_q)
   &=\sum_{\sigma\in\mathfrak S_q}\mathrm{sgn}(\sigma)
                      [v_{\sigma(1)}|\cdots|v_{\sigma(q)}],\\
 \Pi[a_1|\cdots|a_q]&=\frac1{q!}\pi_1(a_1)\wedge\cdots\wedge\pi_1(a_q).
 \end{aligned}
\]
Tensor with $K_{M_m}$ and let $\mathcal I$ be signed contraction
with top forms as in \eqref{eq:v6-contraction}.
The verifications in Proposition~\ref{prop:v6-global-bar} use
only the rank of the vector bundle, not the number four:
$\Pi I=1$, $bI=\Pi b=0$, and the normalized bar--Koszul
homotopy obeys $bh+hb=1-I\Pi$, $hI=\Pi h=0$.
With $\Delta=\mathcal I^{-1}d\mathcal I$ the same equations give
\[
 D_B=-b+uI\Delta\Pi,\qquad D_B^2=0,\quad
 D_BI=Iu\Delta,\quad\Pi D_B=u\Delta\Pi,
\]
\[
 D_B(-h)+(-h)D_B=1-I\Pi.
\]
On each coordinate open, coefficientwise extension of the chosen
homotopy $h$ preserves the lower Laurent bound and proves that the
local comparison is a quasi-isomorphism. The maps $I$, $\Pi$ and
$\mathcal I$, and therefore $D_B$, are globally defined by
symmetric-degree projection, exterior products and contraction of
tensor fields. These maps commute with coordinate changes and their
triple compositions. Thus they give the global chain comparison
\eqref{eq:r14-actual-central-diagram} and its chain section.
Its being a quasi-isomorphism is local on $M_m$ and follows from
the chartwise homotopies just constructed. No gluing of those
chosen local homotopies is needed. Derived global sections of this
sheaf quasi-isomorphism give the asserted hypercohomology comparison.

The exterior Euler sequences on $M_m=\P^{d_m}$ give
\[
 0\longrightarrow\Omega^p_{M_m}\longrightarrow
 \Lambda^pV_m^\vee\otimes\O_{M_m}(-p)
                 \longrightarrow\Omega^{p-1}_{M_m}\longrightarrow0,
 \quad1\le p\le d_m.
\]
Each middle term is acyclic. Induction gives
$H^q(M_m,\Omega^p)=\CC$ for $q=p$ and $0$ otherwise.
The de Rham spectral sequence therefore has no nonzero
higher differentials. The Cech cocycles
\[
 \alpha_{ij}=d\log(s_j/s_i),\qquad
 z_p=\alpha_{i_0i_1}\wedge\cdots\wedge\alpha_{i_{p-1}i_p}
\]
represent its $d_m+1$ nonzero classes: restriction to a
coordinate $\P^p$ has iterated residue one on
$dx_1/x_1\wedge\cdots\wedge dx_p/x_p$.
This also proves \eqref{eq:r14-divisor-DT-groups} with its
integer grading.

On the finite Cech--de Rham total complex choose the $z_p$ as
cohomology representatives and split each term as boundaries,
these representatives, and a complement to cycles. The
differential maps the latter complement isomorphically to the
boundaries in the next term. Its inverse gives a homotopy
$dh+hd=1-ip$. Extend this fixed contraction coefficientwise
over $\Gamma(V,\O_C)$ and in $u$, and conjugate by
multiplication by $u^p$ on $p$-forms. Since $0\le p\le d_m$,
this conjugation changes the Laurent bound by a uniformly
bounded amount. It produces the finite complex
\eqref{eq:r14-explicit-finite-P}, with inclusion $e_p\mapsto
Z_p=u^pz_p$. Evaluation at any $b$, including $0$, commutes
with every map. All $d_m$ here are even, so this finite complex
has only even terms. The maps $I\mathcal I^{-1}$ give exactly
\eqref{eq:r14-iota0-basis}. They are quasi-isomorphisms by the
contraction and the residue basis calculation.

The ambient critical family is $M_m\times C$ with $F=0$.
Consequently $\partial_t$ commutes with the Cech, de Rham and
bar differentials and all the displayed maps. The finite basis
is horizontal, its matrix is zero, and its coefficients are
regular at $0$. Finally the comparison for $\MF(M_m,0)$ in
\cite{Efimov2018CyclicMF} gives the smooth-fiber cyclic
identification. The shift $[d_m]$ has even parity in these
three examples. Over a field the comparison zigzag can be
represented by a quasi-isomorphic chain model; it supplies the
maps $\iota_b$ of Hypothesis~\ref{hyp:finite-cyclic}. Its
central map is already the constructed bar map above.

\emph{The twists and the sheaf locus.}
Tensoring by the globally defined line bundle $L^k$ sends
$\mathcal E_m$ to
\[
 [L^{k-m}\boxtimes\O_{M_m}(-1)\longrightarrow
                                  L^k\boxtimes\O_{M_m}].
\]
It commutes with restriction and base change, and induces
isomorphisms of dg endomorphism algebras preserving trace.
It therefore transports all the deformation, d-critical and
orientation identifications used above. This proves the
assertions for every $k$.
For $b\ne0$, $s$ is a non-zero-divisor on the integral smooth
threefold, so this complex is the indicated line-bundle sheaf
on a Cartier divisor. A Cartier divisor in a smooth threefold
is Cohen--Macaulay of pure dimension two. Tensoring its
structure sheaf by an invertible sheaf preserves purity.
When the divisor is integral, every nonzero subsheaf has
generic rank one; the reduced-Hilbert-polynomial argument
following Proposition~\ref{prop:v6-projective-locus} proves
stability. At the special fiber the definition retains the
two-term complex, including its possible degree-$-1$
cohomology. This completes the stated comparison for these
specific divisor-supported objects and their specified limits.
\end{proof}

\subsection{Complete linear systems of arbitrary positive degree}\label{ssec:all-positive-degrees}
The bound $m\leq3$ in the preceding construction is unnecessary when one uses the relative bundle of sections rather than the full space of homogeneous polynomials of a fixed degree. Let $S_j=H^0(\P^4,\O_{\P^4}(j))$, with $S_j=0$ for $j<0$, and retain the resolved quintic family $\pi:\cX\to C$, the morphism $a:\cX\to\P^4$, and $L=a^*\O_{\P^4}(1)$.

\begin{lemma}\label{lem:all-degree-sections}
For every integer $m>0$, the sheaf $V_m=\pi_*L^m$ is a vector bundle on $C$, commutes with base change, and has rank
\begin{equation}\label{eq:all-degree-rank}
 r_m=\binom{m+4}{4}-\binom{m-1}{4},\qquad d_m=r_m-1,
\end{equation}
where the second binomial coefficient equals zero for $1\leq m\leq4$. It admits the presentation
\begin{equation}\label{eq:all-degree-presentation}
 0\longrightarrow S_{m-5}\otimes\O_C
 \xrightarrow{\,hq_4-tq_5\,}S_m\otimes\O_C
 \longrightarrow V_m\longrightarrow0.
\end{equation}
For every fiber, including the central fiber,
\[
 R\Gamma(X_b,L_b^m)=V_{m,b},\qquad
 R\Gamma(X_b,L_b^{-m})=V_{m,b}^\vee[-3],\qquad
 R\Gamma(X_b,\O_{X_b})\simeq\CC\oplus\CC[-3].
\]
\end{lemma}
\begin{proof}
The hypersurface sequence on a smooth fiber is
\[
0\to\O_{\P^4}(m-5)\xrightarrow{hq_4-tq_5}\O_{\P^4}(m)\to L_b^m\to0.
\]
For $m>0$, the middle cohomology of both line bundles vanishes; moreover $m-5\geq-4$, so $H^4(\P^4,\O(m-5))=0$. Thus the cokernel of multiplication by the degree-five equation is $H^0(X_b,L_b^m)$ and the higher cohomology vanishes. Multiplication is injective because the polynomial ring is a domain.

For the special fiber, write $T_j=H^0(\P^3,\O_{\P^3}(j))$ and take $T_j=0$ for $j<0$. Cohomology of the quartic hypersurface, the boundary quartic surface and the blowup $X_2=\operatorname{Bl}_\Gamma\P^3$ gives
\[
\begin{aligned}
H^0(X_1,L_1^m)&=S_m/q_4S_{m-4},\\
H^0(X_2,L_2^m)&=T_m,\\
H^0(S,L_S^m)&=T_m/(q_4|_{h=0})T_{m-4},
\end{aligned}
\]
and all three objects have zero higher cohomology. For the blowup this follows from $R\beta_*\O_{X_2}=\O_{\P^3}$ and the projection formula. The normalization sequence identifies $H^0(X_0,L_0^m)$ with the kernel of the difference of restrictions from the first two displayed spaces to the third. The map from $T_m$ is surjective. A compatible pair $(a\bmod q_4,b)$ satisfies $a|_{h=0}-b=(q_4|_{h=0})c$; lifting $c$ to $\widetilde c\in S_{m-4}$ and replacing $a$ by $a-q_4\widetilde c$ gives compatible polynomials. The ambiguity is $(h)\cap(q_4)=(hq_4)$. Hence
\[
H^0(X_0,L_0^m)\simeq S_m/(hq_4)S_{m-5},\quad
h^0(X_0,L_0^m)=\dim S_m-\dim S_{m-5}=r_m.
\]
The fiber dimension is constant, giving base change. The cokernel in \eqref{eq:all-degree-presentation} is locally free since multiplication by $hq_4-tq_5$ has constant rank; the fiberwise identifications give the asserted presentation. Finally relative Gorenstein--Serre duality, with the chosen volume form, identifies the negative tensor powers and their base change as stated.
\end{proof}

\begin{proposition}\label{prop:divisor-stability-all-degrees}
Let $X$ be a smooth projective Calabi--Yau threefold with $\operatorname{Pic}(X)=\Z[H]$, $H$ ample. For every effective Cartier divisor $D\in|mH|$, $m>0$, the sheaf $\O_D(kH)$ on $X$ is pure of dimension two and $H$-Gieseker stable, for every integer $k$.
\end{proposition}
\begin{proof}
The divisor $D$ is Cohen--Macaulay and pure of dimension two. Fix a nonzero proper subsheaf $A\subset\O_D(kH)$ and write $D=\sum_i m_iD_i$ as a cycle of prime divisors. At the generic point of $D_i$, the module $A$ has some length $0\leq r_i\leq m_i$. Put $D'=\sum_i r_iD_i$ and $D''=D-D'$. When $D''\ne0$, purity gives an inclusion $A\subset\O_{D'}(kH-D'')$ with cokernel supported in dimension at most one. When $D''=0$ retain the original inclusion. Write $D'\sim aH$ with $1\leq a\leq m$ and $c=H^3>0$. The leading terms of the Hilbert polynomials are
\[
\begin{aligned}
P_{\O_D(kH)}(n)&=\frac{mc}{2}n^2+c(mk-m^2/2)n+O(1),\\
P_{\O_{D'}(kH-D'')}(n)&=\frac{ac}{2}n^2+c(a(k-m+a)-a^2/2)n+O(1).
\end{aligned}
\]
The linear coefficients of their reduced Hilbert polynomials are respectively $2k-m$ and $2k-2m+a$. For $a<m$ the latter is smaller, and passing to $A$ only decreases the reduced Hilbert polynomial. If $a=m$, then the quotient $\O_D(kH)/A$ is nonzero and supported in dimension at most one, and its positive Hilbert polynomial gives the strict inequality. This proves stability even for reducible or nonreduced $D$.
\end{proof}

\begin{theorem}\label{thm:all-positive-degrees}
For $m>0$ and $k\in\Z$, let $\rho_m:M_m=\P_C(V_m)\to C$ parametrize one-dimensional subspaces and define the relative perfect complexes
\begin{equation}\label{eq:all-degree-universal}
\mathcal E^{m,k}=\big[L^{k-m}\boxtimes\O_{M_m}(-1)\longrightarrow L^k\big]
\quad\text{on }\cX\times_C M_m
\end{equation}
in cohomological degrees $-1,0$. This defines an open immersion into the scalar-rigidified classical moduli of simple perfect complexes. The locus is smooth and projective over $C$, and for $b\neq0$ its objects are precisely the $L_b$-Gieseker-stable sheaves $\O_D(kL_b)$ for $D\in|L_b^m|$. The special-fiber objects are their effective full derived matching triples on $X_1$, $X_2$, and $S$.

With the orientation $K_{M_m/C}|_{M_{m,b}}$ the d-critical structure is zero and the oriented DT sheaf is $\CC_{M_{m,b}}[d_m]$. The oriented periodic two-sided bar model has a horizontal quasi-isomorphism
\begin{equation}\label{eq:all-degree-bar}
\mathcal B^{\mathrm{or}}_{m,\mathrm{per}}\xrightarrow{\ I\Pi\ }
\big(\Omega^\bullet_{M_m/C}((u)),u\,d\big)[d_m].
\end{equation}
The corresponding finite complex has rank $r_m$ in parity $d_m\pmod2$, zero in the opposite parity, and regular connection with horizontal basis $1,H,\ldots,H^{d_m}$. In particular
\begin{equation}\label{eq:all-degree-groups}
\mathrm{DT}_{\mathrm{coh}}^j(M_{m,b})=
\begin{cases}\CC,&j=2a-d_m,\ 0\leq a\leq d_m,\\0,&\text{otherwise}.\end{cases}
\end{equation}
The central comparison is the specialization of \eqref{eq:all-degree-bar}; the corresponding ordinary critical realization satisfies
\[
R\Gamma(M_{m,b},\mathcal B^{\mathrm{or}}_{m,\mathrm{per}})
\simeq\CyclicChains_{\mathrm{per}}(\MF(M_{m,b},0))[d_m].
\]
\end{theorem}
\begin{proof}
Tensoring by $L^k$ commutes with derived restriction, base change, endomorphisms and trace; hence first take $k=0$. The endomorphism complex of $E=[L_b^{-m}\xrightarrow{s}\O_{X_b}]$ is
\[
L_b^{-m}\xrightarrow{w\mapsto(ws,sw)}\O_{X_b}^{\oplus2}
\xrightarrow{(a,c)\mapsto(a-c)s} L_b^m
\]
in degrees $-1,0,1$. Lemma~\ref{lem:all-degree-sections} leaves only cohomology rows $0$ and $3$, and each nonzero map in those rows has rank one. Their separation by three degrees and the three-column range exclude higher differentials. Consequently
\[
\dim\operatorname{Ext}^{0,1,2,3}(E,E)=(1,d_m,d_m,1),
\qquad\operatorname{Ext}^{<0}(E,E)=0.
\]
The projective-bundle Euler sequence identifies the middle Ext groups with tangent and cotangent spaces. Taking derived global sections of $L^{-m}\to\O\to E$ gives
\[
H^0(R\Gamma E)=\CC,\quad H^1=0,\quad
H^2=\ker(V_{m,b}^{\vee}\xrightarrow{s^{\vee}}\CC),
\]
with all other cohomology zero. For a perfect deformation $E_A$ over a local Artinian $C$-algebra $A$, cancel invertible differentials in a finite-free model of $R\Gamma E_A$. Its remaining degree-zero and degree-two terms cannot support a differential, so $H^0(R\Gamma E_A)\simeq A$. A generator yields the evaluation triangle
\[
K_A\longrightarrow\O_{X_A}\longrightarrow E_A\longrightarrow K_A[1].
\]
Local cancellation in a perfect representative shows that $K_A$ is a line bundle. The units sequence and $H^1(X_b,\O_{X_b})=0$ show that its lift is, up to isomorphism, the existing $L_A^{-m}$. Hence $E_A$ comes from a section $s_A\in V_m\otimes A$ nonzero on the residue field, modulo a unit. Conversely an isomorphism of evaluation triangles identifies precisely those sections up to rescaling. These are the infinitesimal points of $\P_C(V_m)$ after scalar rigidification. The representability criterion of \cite{Lieblich2006} gives an etale universally injective map, hence an open immersion. The source is a projective bundle, so smooth and proper. A primitive representative over any discrete valuation ring supplies the explicit valuative extension. Stability on $X_b$ follows from Proposition~\ref{prop:divisor-stability-all-degrees} and $\operatorname{Pic}(X_b)=\Z[L_b]$.

The central complex is termwise a pair of vector bundles, so its two restrictions agree on $S$ and furnish an effective matching triple. Applying the normalization triangle to $E^\vee\otimes E$ identifies its deformation algebra with the homotopy fiber of the difference of the two component restriction maps. The d-critical structure on its smooth classical locus is zero, and the orientation $K_{M_m/C}$ gives the shifted constant DT sheaf. The normalized bar--Koszul antisymmetrization $I$ and projection $\Pi$ of Section~\ref{sec:quantization} are functorial under changes of frames of $T_{M_m/C}$. Their local contraction homotopies show that \eqref{eq:all-degree-bar} is a quasi-isomorphism; global gluing requires the maps themselves, not a global choice of the local homotopies. The shift by $[d_m]$ and its differential sign are the normalized orientation conventions.

Over a base open trivializing $V_m$, the projective bundle is $\P^{d_m}\times C$. Its de Rham cohomology has the basis $1,H,\ldots,H^{d_m}$. On the standard cover, representatives can be chosen as Cech products of $d\log(s_j/s_i)$; multiplication by $u^a$ on degree-$a$ forms produces cocycles for $u\,d$. Splitting the finite Cech--de Rham complex into cohomology and contractible summands yields a finite cohomology complex whose basis specializes at $0$. The Chern class $H$ is defined by a relative line bundle and is horizontal. The relative connection is thus $\partial_t$ on this basis, regular at $0$. The comparison maps preserve these classes; the ordinary cyclic identification is given by the matrix-factorization comparison of \cite{Efimov2018CyclicMF}. This proves all the asserted identifications for the specified zero critical model.
\end{proof}

\subsection{From section complexes to point and curve ideals}
\label{ssec:transition-incidence}
The preceding theorem computes the entire generic-to-central comparison
for the section complexes: the generic potential, actual central matching
objects, oriented bar differential, cohomology, and regular connection.
Section~\ref{sec:divisor-sheaves} next treats the different objects
$j_{D*}I_{Z/D}$ in the order one point, two points, a line, a line
with one point, and a line with two point defects. In those examples,
``incidence'' means the specified scheme-theoretic containment $Z\subset D$.
It does not mean that a point of the sheaf moduli lies in $X_i$ rather
than $X_0$. Before passing to them we retain the further perfect-complex
examples and the fixed-fourfold deformation calculations.

\subsection{Several sections and Grassmannian moduli}\label{sec:grassmannian}
Fix $m>0$. Lemma~\ref{lem:all-degree-sections} constructs $V_m=\pi_*L^m$ and proves
\[
 r=\operatorname{rk}V_m=\binom{m+4}{4}-\binom{m-1}{4}.
\]
For $1\le q\le r$, let $G=\operatorname{Gr}_C(q,V_m)$ parametrize rank-$q$ subbundles, with tautological $\mathcal U\subset V_m\otimes\mathcal O_G$. Set $d=q(r-q)$. Dualizing the evaluation $\mathcal U\otimes\mathcal O\to L^m$ gives
\begin{equation}\label{eq:grass-17}
 \mathcal E^{m,k,q}=[L^{k-m}\longrightarrow L^k\otimes\mathcal U^\vee]
       \quad\text{on }\mathcal X\times_C G,
\end{equation}
in degrees $-1,0$. For $q=1$ it differs from the initial universal complex only by a line bundle from the parameter space.

\begin{theorem}\label{thm:grassmannian-moduli}
The family~\eqref{eq:grass-17} gives an open immersion of $G$ into the classical scalar-rigidified moduli space of simple perfect complexes on the fibers. This locus is smooth and projective over $C$. Its fiber objects satisfy
\begin{equation}\label{eq:grass-18}
 \dim\operatorname{Ext}^{0,1,2,3}(E,E)=(1,d,d,1),\qquad
                         \operatorname{Ext}^{<0}(E,E)=0,
\end{equation}
and, with $h_b=c_1(L_b)$,
\begin{equation}\label{eq:grass-19}
 \operatorname{ch}(E_{U,b}^{m,k,q})=q e^{kh_b}-e^{(k-m)h_b}.
\end{equation}
The central objects are their full derived restriction triples.

More generally, the same assertion about the classical moduli holds for an oriented projective flat Gorenstein threefold family with $H^{0,3}(\mathcal O)=\mathbb C$, $H^{1,2}(\mathcal O)=0$, $R\pi_*L^m=V_m$ locally free and $R\pi_*L^{-m}=V_m^\vee[-3]$, compatibly with base change. The resolved quintic satisfies these hypotheses.
\end{theorem}
\begin{proof}
Tensor by $L^{-k}$ and put $W=U^\vee$. The sheaf endomorphism complex has the three terms
\[
 L_b^{-m}\otimes W^\vee\ \longrightarrow\
 \mathcal O\oplus(\mathcal O\otimes\operatorname{End}W)
 \ \longrightarrow\ L_b^m\otimes W.
\]
Their differentials are $v\mapsto(vs,sv)$ and $(a,B)\mapsto sa-Bs$. Their composite is zero. We compute their derived global sections and then reconstruct the full deformation functor.

Here is the full hypercohomology calculation, including its matrices.
Put $V=V_{m,b}$ and choose a basis $v_1,\ldots,v_r$ such that
$U=\langle v_1,\ldots,v_q\rangle$. In the dual basis $e_1,\ldots,e_q$
of $W$, the section column is $s=\sum_{a=1}^q v_a\otimes e_a$.
Represent an element of $V\otimes W$ by a $q\times r$ matrix.
The matrix of $s$ is
\[
 S=\begin{pmatrix}I_q&0_{q\times(r-q)}\end{pmatrix},
 \qquad D_0(a,B)=aS-BS=\begin{pmatrix}aI_q-B&0\end{pmatrix}.
\]
The word ``row'' below refers to the sheaf-cohomology degree $j$ in
$E_1^{p,j}=H^j(X_b,\mathcal Hom^p(E,E))$, not to a row of $S$.
Its nonzero first-page terms are as follows.
\[
\begin{array}{c|ccc}
 j\backslash p&-1&0&1\\\hline
 3&V^\vee\otimes W^\vee\xrightarrow{D_3}
       &\mathbb C\oplus\operatorname{End}W&0\\
 2&0&0&0\\
 1&0&0&0\\
 0&0&\mathbb C\oplus\operatorname{End}W\xrightarrow{D_0}
       &V\otimes W.
\end{array}
\]
For $D_3$, write the dual input as $(C,D)$, where $C$ is the
$q\times q$ block, with the transpose convention from the trace pairing,
and $D$ is the remaining $q\times(r-q)$ block. The differential
$v\mapsto(vs,sv)$ gives
\[
 D_3(C,D)=(\operatorname{tr}C,C).
\]
Thus, listing the $q^2$ entries of the first block before the remaining
$q(r-q)$ entries, the two numerical matrices are
\begin{equation}\label{eq:grass-explicit-matrices}
 D_0=\begin{pmatrix}\operatorname{vec}(I_q)&-I_{q^2}\\
                        0&0\end{pmatrix},\qquad
 D_3=\begin{pmatrix}\operatorname{tr}&0\\ I_{q^2}&0\end{pmatrix}.
\end{equation}
Their sizes are $qr\times(1+q^2)$ and $(1+q^2)\times qr$.
For example, for $q=2,r=5$, their nonzero blocks are
\[
 D_0=\begin{pmatrix}
 1&-1&0&0&0\\0&0&-1&0&0\\0&0&0&-1&0\\1&0&0&0&-1\\
 \multicolumn{5}{c}{0_{6\times5}}
 \end{pmatrix},\qquad
 D_3=\left(\begin{array}{rrrr|c}
 1&0&0&1&\\1&0&0&0&\\0&1&0&0&0_{5\times6}\\
 0&0&1&0&\\0&0&0&1&
 \end{array}\right).
\]
In general $\ker D_0=\mathbb C(1,I_q)$, $\operatorname{rk}D_0=q^2$,
$\ker D_3=\{(0,D)\}$ has dimension $q(r-q)$, and
$\operatorname{coker}D_3$ is one-dimensional, detected by
$(a,B)\mapsto\operatorname{tr}B-a$. The total degree is $p+j$.
Consequently the four surviving total degrees $0,1,2,3$ have dimensions
$1,q(r-q),q(r-q),1$. A differential joining $j=3$ to $j=0$ would
increase $p$ by four; the three columns $-1,0,1$ exclude it.
This proves~\eqref{eq:grass-18}, including on the entire singular fiber
$X_0$. No smoothness of $X_0$ or of the component moduli stacks was
used. The defining triangle also gives
\begin{equation}\label{eq:grass-20}
 H^0(R\Gamma E)=U^\vee,\quad H^2(R\Gamma E)=(V_{m,b}/U)^\vee,
                       \quad H^j=0\ (j\ne0,2).
\end{equation}
For an Artinian base deformation $A$, a minimal finite-free representative of $R\Gamma E_A$ has ranks $q,r-q$ only in degrees~0,2, hence zero differential. Evaluation from a basis of $H^0$ gives
\[
 K_A\longrightarrow\mathcal O_{X_A}^{\oplus q}\longrightarrow E_A
                     \longrightarrow K_A[1].
\]
The residue complex of $K_A$ is $L_b^{-m}$. Local cancellation in a perfect representative shows that $K_A$ is a line bundle. The square-zero units sequence and $H^1(X_b,\mathcal O)=0$ make its lift unique up to isomorphism; the existing lift is $L_A^{-m}$. Thus $E_A$ is recovered from $q$ sections. Their dual map $A^q\to V_m\otimes A$ has a unit $q$-minor, so its image is a rank-$q$ subbundle.

Changing the evaluation basis acts by $\mathrm{GL}_q(A)$; changing the generator of $K_A$ is absorbed by this action. Conversely an isomorphism of the complexes recovers these evaluation triangles and the same subbundle. The residual automorphisms are scalar. Hence the classical infinitesimal moduli after scalar rigidification is precisely $G$.

Lieblich's representability~\cite{Lieblich2006}, the infinitesimal criterion for etaleness and the same reconstruction over algebraically closed fields give an etale universally injective morphism, hence an open immersion. The source is smooth projective. Over a discrete valuation ring its limit is the saturated submodule $U_\eta\cap V_{m,A}$; the quotient is torsion-free and therefore free. The universal perfect complex extends there.

The Kodaira--Spencer identification is
\[
 \operatorname{Ext}^1(E,E)=\operatorname{Hom}(U,V_{m,b})/\operatorname{End}U
                         =\operatorname{Hom}(U,V_{m,b}/U).
\]
Restriction of~\eqref{eq:grass-17} is termwise, and its two restrictions to $S$ agree. Additivity gives~\eqref{eq:grass-19}.
\end{proof}

\begin{theorem}\label{thm:grassmannian-DT}
Equip the classical Grassmannian loci $G_b$ of Theorem~\ref{thm:grassmannian-moduli} with the d-critical realization of Section~\ref{ssec:orientation-input} and the orientation $K_{G/C}|_{G_b}$. Their d-critical structures are zero and their DT sheaves are $\mathbb C_{G_b}[d]$. For partitions $\lambda\subset q\times(r-q)$,
\begin{equation}\label{eq:grass-21}
 \mathrm{DT}_{\rm coh}^j(G_b)
     =\bigoplus_{2|\lambda|-d=j}\mathbb C\,\sigma_\lambda.
\end{equation}
There is a global oriented periodic bar model with a comparison
\begin{equation}\label{eq:grass-22}
 \mathcal B^{\rm or}_{G,\mathrm{per}}\xrightarrow{\mathcal I\Pi}
                     (\Omega^\bullet_{G/C}((u)),u\,d)[d].
\end{equation}
It is a quasi-isomorphism on every fiber and commutes with specialization. Its finite cohomology has rank $\binom rq$ in parity $d$, zero in the other parity, and a regular connection with horizontal Schubert basis. The central comparison is the specialization of this map. For the specified ordinary critical realization,
\[
 R\Gamma(G_b,\mathcal B^{\rm or}_{G,\mathrm{per}})
                  \simeq CC_{\mathrm{per}}(\operatorname{MF}(G_b,0))[d].
\]
\end{theorem}
\begin{proof}
A smooth scheme has $\mathcal S^0=0$~\cite{Joyce}, so the d-critical structures in the statement are zero. The tangent bundle is $\mathcal U^\vee\otimes(V_m/\mathcal U)$; hence
\[
 K_{G/C}=(\det\mathcal U)^r\otimes(\det V_m)^{-q}.
\]
Its square is the virtual canonical bundle of the zero critical structure. This chosen orientation gives the stated DT sheaf~\cite{BBDJS}. No splitting of the entire Artin enhancement is asserted.

Set $\mathcal A=\operatorname{Sym}T_{G/C}$ with the two zero augmentations. Form the two-sided bar resolution and the tensorial maps $I,\Pi$ of Section~\ref{sec:v5-bar-bridge}. Tensor by $K_{G/C}$ and use signed contraction with top forms. With $\Delta=(-1)^d\mathcal I^{-1}d_{G/C}\mathcal I$ and the suspension signs of Theorem~\ref{thm:v5-periodic-bar}, its differential is $-b+uI\Delta\Pi$. The local identities prove~\eqref{eq:grass-22}. The maps and differential glue tensorially. Local homotopies prove that the global sheaf map is a quasi-isomorphism; a global choice of these homotopies is unnecessary.

The Schubert cells of $\operatorname{Gr}(q,r)$ are indexed by partitions in the rectangle $q\times(r-q)$. Their closures give a cohomology basis in even degrees~\cite{Bri}. The classes $\sigma_\lambda=s_\lambda(\mathcal U^\vee)$ have degree $2|\lambda|$, proving~\eqref{eq:grass-21}. Their number is $\binom rq$. They are globally defined on the relative Grassmannian.

Over a base open trivializing $V_m$, choose closed Cech--de Rham representatives $z_\lambda$ in bidegree $(|\lambda|,|\lambda|)$. For differential $u\,d$, put $Z_\lambda=u^{|\lambda|}z_\lambda$. The finite model is
\begin{equation}\label{eq:grass-23}
 P^{\bar d}=\bigoplus_\lambda R e_\lambda,\qquad
 P^{\overline{d+1}}=0,\quad d_P=0,\qquad R=\mathcal O_C\otimes\mathbb C((u)),
\end{equation}
with chain map $e_\lambda\mapsto I\mathcal I^{-1}(Z_\lambda)$ to the bar complex. Split the bounded Cech--de Rham complex over $\mathbb C$ into cohomology and contractible summands and extend coefficientwise. Conjugation by $u^p$ on $p$-forms changes the Laurent bound by at most $d$. Thus the map is a quasi-isomorphism in the stated completion and commutes with evaluation at $b=0$.

Chern classes of $\mathcal U^\vee$ have closed absolute de Rham representatives; their relative classes and Schur polynomials are horizontal. The bases agree on base overlaps and
\begin{equation}\label{eq:grass-24}
 \nabla_{\partial_t}e_\lambda=0.
\end{equation}
In this affine quintic family $V_m$ is free, so a global trivialization supplies the finite chain model as well. All routes in
\[
\begin{tikzcd}[column sep=large,row sep=large]
 H(P_b)\ar[r,"\sim"]\ar[d,"\mathrm{Hol}"']&\mathbb H(G_b,\mathcal B_b^{\rm or})\ar[d,"\mathrm{Hol}"]\\
 H(P_0)\ar[r,"\sim"]&\mathbb H(G_0,\mathcal B_0^{\rm or})
\end{tikzcd}
\]
send $e_\lambda$ to the same Schubert class. Efimov's theorem~\cite{Efimov2018CyclicMF} gives the cyclic realization with suspension $[d]$.
\end{proof}

\textbf{Which derived intersection has been calculated.}
The component morphisms are the derived Lagrangian morphisms
$r_i:\mathfrak{Perf}(X_i)\to\mathfrak{Perf}(S)$ of
Theorem~\ref{thm:sc-lagrangian-restriction}. Their sources are not asserted to be
smooth schemes, nor are they ``smooth derived enhancements.'' The
central universal complex in~\eqref{eq:grass-17} defines the actual map
\[
 G_0\longrightarrow t_0\bigl(\mathfrak{Perf}(X_1)
                  \times^h_{\mathfrak{Perf}(S)}\mathfrak{Perf}(X_2)\bigr),
 \quad U\longmapsto(E_{1,U},E_{2,U},\alpha_U).
\]
The normalization triangle~\eqref{eq:sc-mapping-triangle} identifies
its endomorphism deformation algebra with that of $E_{U,0}$.
The reconstruction in Theorem~\ref{thm:grassmannian-moduli} identifies
the classical scalar-rigidified formal neighborhoods of these matching
objects with those of $G_0$. Thus $G_0$ is the smooth \emph{classical}
moduli locus used in this example. Its reduced critical enhancement still
has tangent and obstruction spaces of dimension $d$, and is locally
$\operatorname{dCrit}_{G_0}(0)$, not the ordinary smooth scheme $G_0$.

There are two bar calculations, with different purposes. For derived
component charts $\operatorname{Spec}A_i\to\operatorname{Spec}A_S$,
$\operatorname{Bar}(A_1,A_S,A_2)$ computes the structure algebra of
that component-moduli intersection. To compute its \emph{oriented DT
realization} on the just identified locus, use its induced reduced
critical chart $(U,0)$ and its canonical orientation. In that chart
$U\times^h_{T^*U}U$ is a presentation of the reduced critical enhancement.
The bar of these two zero sections retracts to
\[
 \mathcal O_U[\xi_1,\ldots,\xi_d],\quad |\xi_i|=-1,\quad
 \delta=0,
 \qquad D_{\mathrm{DT}}=u\sum_i\partial_{x_i}\partial_{\xi_i}.
\]
Contraction with the orientation coefficient gives
$(\Omega_U^\bullet((u)),u d)[d]$. These are the maps in
\eqref{eq:grass-22}. The two zero sections are a local presentation of
the reduced matching problem; they are not replacements of
$\mathfrak{Perf}(X_1)$ and $\mathfrak{Perf}(X_2)$ by smooth submanifolds.
The identification is at the level of the specified oriented critical
realization, and does not assert an isomorphism between the unreduced
Artin stacks and a shifted cotangent scheme.

\noindent\textbf{Examples and the Chern character.}
For $m=1$ one has $r=5$, and the families give
\[
\begin{array}{c|ccccc}
q&1&2&3&4&5\\\hline
\dim G_b&4&6&6&4&0\\
\operatorname{rk}E&0&1&2&3&4\\
\operatorname{rk}\mathrm{DT}^{\rm per}&5&10&10&5&1.
\end{array}
\]
For $q=2$, independent linear forms form a regular sequence on a smooth quintic. The Koszul resolution gives $E^{1,0,2}_{U,b}\simeq I_{Z/X_b}(1)$, where $Z=X_b\cap H_1\cap H_2$. These rank-one torsion-free sheaves are Gieseker stable. For $q=5$, the Euler sequence gives the rigid vector bundle $T_{\mathbb P^4}(-1)|_{X_b}$. For arbitrary $m,q$, the theorem concerns simple perfect complexes; it does not impose sheaf stability on every cokernel.

With $h=c_1(L_b)$ the components of~\eqref{eq:grass-19} are
\[
\begin{split}
\operatorname{ch}_0(E)&=q-1,\\
\operatorname{ch}_1(E)&=((q-1)k+m)h,\\
\operatorname{ch}_2(E)&=(\tfrac12(q-1)k^2+km-\tfrac12m^2)h^2,\\
\operatorname{ch}_3(E)&=(\tfrac16(q-1)k^3+\tfrac12k^2m
                  -\tfrac12km^2+\tfrac16m^3)h^3.
\end{split}
\]
These characters specify the chosen loci but do not classify all objects with that class. A precise rank-one criterion is as follows.

\begin{proposition}\label{prop:rank-one-character}
Let $D\subset X_b$ be a smooth integral divisor in $|L_b^m|$. Let $F$ be torsion-free of rank one on $D$, with a specified $F^{\vee\vee}\simeq L_b^k|_D$. If $\operatorname{ch}(j_*F)=e^{kh}(1-e^{-mh})$, then $F\simeq L_b^k|_D$.
\end{proposition}
\begin{proof}
On a smooth surface $F^{\vee\vee}$ is a line bundle and its quotient by $F$ is zero-dimensional. Thus
\begin{equation}\label{eq:grass-25}
\operatorname{ch}(j_*F)=e^{kh}(1-e^{-mh})-                       \operatorname{length}(F^{\vee\vee}/F)[\mathrm{pt}].
\end{equation}
Equality in degree six makes the nonnegative length zero.
\end{proof}

For sheaves without this double-dual condition, or higher-rank torsion-free sheaves, character equality alone does not give the reconstruction used in Theorem~\ref{thm:grassmannian-moduli}.

\textbf{The same five steps in degree two.}
For $m=2$ the vector space of sections has dimension $15$, and the
classical parameter scheme is $G_b=\operatorname{Gr}(q,15)$, including
$b=0$. Its reduced potential is $0$ by the matrix and reconstruction
proof above. Its specified categorical realization is
$\operatorname{MF}(G_b,0)$ with the DT suspension $[q(15-q)]$.
The actual central objects are
\[
 [L_1^{k-2}\to L_1^k\otimes U^\vee],\quad
 [L_2^{k-2}\to L_2^k\otimes U^\vee],\quad
 [L_S^{k-2}\to L_S^k\otimes U^\vee],
\]
with their tautological matching equivalence; the maps are the same
$q$ quadratic sections restricted termwise. Their oriented bar and
horizontal basis are exactly those of~\eqref{eq:grass-22}--\eqref{eq:grass-24}.
In particular
\[
 \mathrm{DT}_{\mathrm{coh}}^j(G_b)=
 \bigoplus_{\substack{\lambda\subset q\times(15-q)\\
                         2|\lambda|-q(15-q)=j}}\mathbb C\sigma_\lambda,
 \qquad \nabla\sigma_\lambda=0.
\]
The periodic ranks for $q=1,2,3,4$ are respectively $15,105,455,1365$.
All are in even parity since $q(15-q)$ is even. The ranks for $15-q$
are the same. For $q=1$ this is the quadric-section example
$M_b=\mathbb P^{14}$, with one DT class in each even degree from
$-14$ to $14$. These calculations are on the moduli of objects on
$X_b$, and at $0$ on the identified matching locus, not separately on
either component threefold.
\textbf{All points of the Grassmannian, not only general subspaces.}
The proof of Theorem~\ref{thm:grassmannian-moduli} uses only linear
independence in $H^0(X_b,L_b^m)$ and the whole-fiber cohomology
vanishings.  It does not require a base-point-free system on $X_b$,
a regular sequence of section equations, or smooth support of a
cokernel.  Thus the rank, Schubert, potential, and connection
calculations apply on the entire Grassmannian on the generic and
special fibers.  The reduced derived obstruction bundle remains
$\Omega_G$, even though the classical parameter is smooth.
The rank-deficient alternative is analyzed at the end of this
section; it is a different, nonsimple moduli problem.

\subsection{Complexes with several cohomology sheaves}\label{sec:complexes-several-sheaves}
We now construct new perfect objects by applying an actual relative
autoequivalence $\Phi$, for example the spherical twist below, to the
universal Grassmannian family. The reason their calculations agree is the
quasi-isomorphism of enhanced endomorphism algebras
\[
 R\operatorname{Hom}_{X_b}(E,E)\xrightarrow{\ \Phi\ }
 R\operatorname{Hom}_{X_b}(\Phi E,\Phi E),
\]
with its composition operations. Thus their formal moduli functors are
isomorphic on every test algebra. We transport the cyclic pairing,
critical charts, and orientation through this isomorphism. This identifies
the potentials, the oriented DT complexes, and the Gauss--Manin operators;
it is not an assertion based on equality of Chern characters. At the
central fiber we first apply $\Phi$ to the actual perfect complex on
$X_0$, and then restrict that transformed complex to $X_1,X_2,S$.
The following proposition spells out this construction.

\begin{proposition}\label{prop:transport-equivalences}
Suppose enhanced equivalences $\Phi_T:\operatorname{Perf}(\mathcal X_T)\to\operatorname{Perf}(\mathcal X_T)$ have inverses and commute coherently with base change $T\to C$. The image of the universal family of Theorem~\ref{thm:grassmannian-moduli} is an isomorphic smooth proper classical simple-perfect-complex locus $G^\Phi\simeq G$. With the transported zero d-critical structure and orientation, the cohomological, oriented bar, finite-connection, and normalized ordinary cyclic statements of Theorem~\ref{thm:grassmannian-DT} hold for this locus.
\end{proposition}
\begin{proof}
The functor and its inverse act on the relative object-moduli groupoids on every test scheme. In particular,
\[
 R\operatorname{Hom}(E,E)\xrightarrow{\sim}
 R\operatorname{Hom}(\Phi E,\Phi E)
\]
preserves the full deformation algebra, simplicity and negative self-Ext vanishing. The open immersion of Theorem~\ref{thm:grassmannian-moduli} is transported to the image locus. The inverse identifies it with $G$, proving smoothness and properness.

The transformed universal complex is relative perfect. Its restriction to $X_0$ and then $X_i,S$ gives effective matching triples. No equivalence on the individual component categories is needed for this step. It is the restriction of the actual transformed complex.

If $a:G\to G^\Phi$ is the resulting isomorphism, then
\[
 a^*K_{G^\Phi/C}=K_{G/C},\qquad
 a^*\mathcal P_{G_b^\Phi}=\mathbb C_{G_b}[d].
\]
Transporting the tensorial bar maps in~\eqref{eq:grass-22} gives the same finite model~\eqref{eq:grass-23} and horizontal Schubert classes. Pullback by $a_b$ also identifies the specified ordinary categories $\operatorname{MF}(G_b^\Phi,0)$ and $\operatorname{MF}(G_b,0)$.
In particular this applies without change to the $m=2$ Grassmannian
families, with their displayed Schubert groups and connection; no
new critical equation is created by applying the equivalence.
\end{proof}

This is not a statement about every independently constructed intrinsic DT category. If a stability condition $(Z,\mathcal P)$ is supplied, its transport is $(Z\circ\Phi_*^{-1},\Phi\mathcal P)$. No stability condition on the whole derived category is being assumed without input. Tensoring by an extending line bundle $N$ and shifting by $a$ change the charge to $(-1)^a e^{c_1(N_b)}\gamma$.

\begin{proposition}\label{prop:relative-spherical-twist}
The relative spherical twist
\begin{equation}\label{eq:grass-26}
 \mathsf T(E)=\operatorname{Cone}\bigl(L\pi^*R\pi_*E
                                \xrightarrow{\mathrm{ev}}E\bigr)
\end{equation}
is an equivalence of the relative perfect-complex category of the resolved quintic family, compatible with derived base change. Proposition~\ref{prop:transport-equivalences} therefore applies to its iterates, inverse, shifts and extending line-bundle twists.
\end{proposition}
\begin{proof}
Proper perfect pushforward and pullback preserve perfect complexes. The unit and relative trace give $R\pi_*\mathcal O\simeq\mathcal O_C\oplus\mathcal O_C[-3]$. Over the affine smooth curve the extension class lies in $\operatorname{Ext}^4_C(\mathcal O_C,\mathcal O_C)=0$; choose the splitting respecting the unit and trace.

Put $\mathsf S=L\pi^*$ and $\mathsf R=R\pi_*$, and denote the left adjoint of $\mathsf S$ by $\mathsf L$. Gorenstein duality with the chosen relative volume gives $\mathsf L=\mathsf R[3]$. For $\mathsf C=\operatorname{Cone}(\mathrm{id}\to\mathsf R\mathsf S)$,
\[
 \mathsf C\simeq[-3],\qquad
 \mathsf R\longrightarrow\mathsf R\mathsf S\mathsf L\longrightarrow\mathsf C\mathsf L
\]
is the canonical adjoint comparison. On a fiber its value on $E$ is the adjoint of the perfect pairing
\[
 R\Gamma(E)\otimes R\operatorname{Hom}(E,\mathcal O)
                         \longrightarrow\mathbb C[-3].
\]
Under duality it is the identity on $R\Gamma(E)$. The relative comparison is an equivalence because its cone on every perfect $E$ is perfect on the base and has zero fibers. In the convention of~\cite[Theorem~1.1]{AL} the cotwist is $\mathsf C[-1]$; the two verified conditions are its invertibility and $\mathsf R\simeq\mathsf C\mathsf L$. The theorem proves invertibility of~\eqref{eq:grass-26}. Its inverse is
\[
 \mathsf T^{-1}(E)=\operatorname{Cone}\bigl(E\longrightarrow
                          L\pi^*(R\pi_*E[3])\bigr)[-1].
\]
All functors and natural maps are taken in dg enhancements. Their unit, trace, evaluation and inverse identities commute with derived base change. This argument uses Gorenstein duality, so includes $X_0$.
\end{proof}

Take $k=0$ and write $\mathcal Q_{U,b}=\mathsf T(E^{m,0,q}_{U,b})[1]$. We calculate the cohomology and extension class of this object next.

Put $W_U=(V_{m,b}/U)^\vee$. The truncation of~\eqref{eq:grass-20} and the evaluation triangle give, by cone composition,
\begin{equation}\label{eq:grass-27}
 W_U\otimes\mathcal O[-2]\longrightarrow L_b^{-m}[1]
 \longrightarrow\mathsf T(E_{U,b})\longrightarrow W_U\otimes\mathcal O[-1].
\end{equation}
The degree-zero evaluation is $U^\vee\otimes\mathcal O\to E_{U,b}$, whose cone is $L_b^{-m}[1]$. The remaining degree-two evaluation gives the first map of~\eqref{eq:grass-27}. It follows that
\begin{equation}\label{eq:grass-28}
 \mathcal H^{-2}(\mathcal Q_{U,b})=L_b^{-m},\qquad
 \mathcal H^0(\mathcal Q_{U,b})=W_U\otimes\mathcal O,
\end{equation}
with all other cohomology zero. For $q<r$ both terms are nonzero even on smooth fibers. Under
\[
 \operatorname{Ext}^3(W_U\otimes\mathcal O,L_b^{-m})
                   =\operatorname{Hom}(W_U,V_{m,b}^\vee),
\]
the extension class is, up to the cone sign, the inclusion $(V_{m,b}/U)^\vee\hookrightarrow V_{m,b}^\vee$. It is nonzero, so this is not the direct sum of the cohomology sheaves. Its character is
\begin{equation}\label{eq:grass-29}
 \operatorname{ch}(\mathcal Q_{U,b})=r-q+e^{-mh_b}.
\end{equation}
The positive-rank complexes therefore have the same proper Grassmannian moduli and the same oriented cohomological comparison as Theorem~\ref{thm:grassmannian-DT}. Tensoring by $L^k$ multiplies this character by $e^{kh_b}$. Their stability is asserted only when a supplied stability condition is transported; the constructed loci are simple-perfect-complex loci.

\subsection{Self-extensions in the fixed ambient projective space}
Assume that $Y=\P(V)$ is the projective space of one-dimensional subspaces of a five-dimensional vector space $V$. Consider sheaves $\I=\O_Z$, where $Z$ is a complete intersection of a quintic and a hyperplane. Let $q,h$ be their homogeneous equations. The hyperplane equation $h\in V^\vee$ is defined up to a scalar, while $q$ is a nonzero element of $\mathrm{Sym}^5(V^\vee/\CC h)$ up to a scalar. If $Q$ denotes the universal quotient bundle on $\P(V^\vee)$, this family of complete intersections is parametrized by
\[
 \P(\mathrm{Sym}^5Q)\longrightarrow\P(V^\vee).
\]
For a fixed smooth quintic with equation $s$, restriction of $s$ to each hyperplane gives a section of this projective bundle. Its image parametrizes the hyperplane-section structure sheaves on that quintic. This construction describes the indicated family; it does not classify all stable sheaves of the same numerical type, and it does not identify the full derived moduli stack with that projective bundle.

For these objects the determinant on $Y$ is trivial and the determinant on $X$ is $\O_X(1)$. We next compute the full self-extension groups in both ambient varieties. The threefold calculation above is compatible with this fixed-ambient computation.

\subsubsection{Self-extensions of a complete-intersection structure sheaf}
Let $V$ have dimension five, let $H\simeq\P^3$ be a hyperplane in $Y=\P(V)$, and let $X\subset Y$ be a smooth quintic meeting $H$ transversely. Put $Z=X\cap H$, let $Q=H^0(H,\O_H(1))$, and let $q\in\Sym^5Q$ define $Z$ in $H$. Choose a trace identification $H^3(H,\O_H(-4))\simeq\CC$ for the vector-space dualities below. We compute the self-extensions over both $X$ and $Y$.

\begin{proposition}\label{prop:Ext-example}
For $F=\O_Z$, viewed as a sheaf on $Y$, the nonzero self-Ext groups are
\begin{align*}
 \ext_Y^0(F,F)&=\CC,\\
 \ext_Y^1(F,F)&=Q\oplus\bigl(\Sym^5Q/\CC q\bigr),\\
 \ext_Y^2(F,F)&=Q^\vee\oplus\bigl(\Sym^6Q/qQ\bigr),\\
 \ext_Y^3(F,F)&=\CC.
\end{align*}
Their dimensions are respectively $1,59,84,1$. In particular $\ext_Y^4(F,F)=0$. On the quintic threefold $X$ the dimensions are instead $1,4,4,1$ in degrees zero through three.
\end{proposition}
\begin{proof}
To compute $Ext^\bullet_Y(\mathcal{O}_Z, \mathcal{O}_Z)$ we
replace the first copy of $\mathcal{O}_Z$ with its Koszul resolution:
$$
 0 \to \mathcal{O}_Y(-6) \to \mathcal{O}_Y(-5) \oplus \mathcal{O}_Y(-1)  \to
 \mathcal{O}_Y \to 0.
$$
The two defining equations annihilate $\O_Z$, so applying sheaf Hom to the displayed Koszul resolution produces a complex with zero differential. Taking derived global sections gives a split total complex, whose spectral sequence has $E_2$ terms containing 
$H^q(Y, \mathcal{O}_Z(p))$ for $p = 0, 1, 5, 6$ and $q = 0, 1, 2$ (since $Z$ has dimension 2). 
To compute the cohomology we can view $\mathcal{O}_Z$ as a sheaf supported on the 
hyperplane $H \simeq \P^3$. Then we can use the long exact sequence induced by the 
resolution 
$$
0 \to \mathcal{O}_H(-5) \to \mathcal{O}_H \to \mathcal{O}_Z \to 0.
$$
As the intermediate cohomology groups $H^q (H, \mathcal{O}_H(p))$ vanish for 
$q = 1, 2$ and all $p$, the long exact sequence reduces to short exact sequences 
$$
0 \to H^0(H, \mathcal{O}_H(p-5) )\to H^0(H, \mathcal{O}_H(p))\to H^0(H, \mathcal{O}_Z(p) )
\to 0,
$$
and 
$$
0 \to H^2(H, \mathcal{O}_Z(p))  \to H^3(H, \mathcal{O}_H(p-5)) \to H^3(H, \mathcal{O}_H(p))
\to 0,
$$
with $p = 0, 1, 5, 6$.  In particular, $H^1(H, \mathcal{O}_Z(p))$ is zero for all $p$. 
The results for the other cases are given as follows:
\begin{itemize}

\item For $p = 0$, we have  $H^0(H, \mathcal{O}_Z) = \mathbb{C},$ and $
H^2(H, \mathcal{O}_Z) = Q^\vee$;

\item For $p = 1$, we have $H^0(H, \mathcal{O}_Z(1) ) = Q, 
H^2(H, \mathcal{O}_Z(1) ) = \mathbb{C}$;

\item For $p = 5$, we have $H^0(H, \mathcal{O}_Z(5) ) = \mathrm{Sym}^5(Q)/\mathbb{C}q,$ and $ 
H^2(H, \mathcal{O}_Z(5) ) = 0$;

\item For $p = 6$, we have $H^0(H, \mathcal{O}_Z(6) ) = \mathrm{Sym}^6(Q)/Q q,$ and $
H^2(H, \mathcal{O}_Z(6) ) = 0$.
\end{itemize}

Substituting these values in the above spectral sequence with only nonzero terms
$$
\begin{cases}

E_2^{0, 0} = \mathbb{C},
\\
E_2^{1, 0} = Q \oplus \mathrm{Sym}^5(Q)/\mathbb{C}q, 
\\
E_2^{2, 0}
= \mathrm{Sym}^6(Q)/Q q,
\\

E_2^{0, 2}= Q^\vee,
\\
E_2^{1, 2} =\mathbb{C},
\end{cases}
$$
all further differentials vanish because the sheaf-Hom complex has zero differential and is already a direct sum of its terms. Thus the asserted Ext groups follow; degree considerations alone would not exclude every possible higher differential.

Since $\dim Q=4$, $\dim\mathrm{Sym}^5Q=56$ and $\dim\mathrm{Sym}^6Q=84$, the ambient dimensions are $1,59,84,1$. On $X$ the one-step resolution $[\O_X(-1)\to\O_X]$ gives the terms $\O_Z$ and $\O_Z(1)$, so the same computation gives $1,4,4,1$. 
\end{proof}

The nonzero $\ext_Y^3(F,F)$ contributes degree two to the full ambient
tangent complex. This term is retained in the normal deformation
triangle of \cite[Proposition~4.19]{KS}. For the fiber
objects, Section~\ref{sec:v6-hyperplanes} extends the calculation to all
$E_{\ell,b}$, including the central complexes.

\textbf{The fixed-fourfold calculation for a quadric section.}
For comparison take $E=\mathcal O_D$, $D\in|2H|$ on a smooth
quintic. The deformation reconstruction gives
$\dim\operatorname{Ext}^1_X(E,E)=14$ and zero reduced potential.
For the fiber embedding in the total family the normal spaces are
$W_1=\mathbb C$ and $W_2=\mathbb C^{14}$. For the different
embedding $X\hookrightarrow\mathbb P^4$ the normal bundle is
$\mathcal O_X(5)$. Applying sheaf Hom to
$[\mathcal O_X(-2)\to\mathcal O_X]$ with target $\mathcal O_D(5)$
gives $[\mathcal O_D(5)\xrightarrow{0}\mathcal O_D(7)]$.
The intermediate cohomology vanishes by the complete-intersection
resolution, and consequently
\[
 \begin{split}
 \dim W_1&=h^0(\mathcal O_D(5))=r_5-r_3=125-35=90,\\
 \dim W_2&=h^0(\mathcal O_D(7))=r_7-r_5=315-125=190.
 \end{split}
\]
The section complex extends to $\mathbb P^4$, so its Atiyah--conormal
obstruction vanishes. The normal triangle therefore gives
$\dim\operatorname{Ext}^1_{\mathbb P^4}(i_*E,i_*E)=104$ and
$\dim\operatorname{Ext}^2_{\mathbb P^4}(i_*E,i_*E)=204$.
These are fourfold self-extension calculations, not substitutes for
the fiberwise DT moduli or for the two-component central bar.

\subsection{The formal neighborhood of the quintic in the fixed ambient fourfold}
The map $i_*$ in \eqref{eq:normal-deformation} uses the actual
embedding of the threefold. Replacing a neighborhood of $X$ by
$\Tot(N_{X/Y})$ would instead replace its deformation algebra by the
one for the zero section of a line bundle. The latter has a
retraction to $X$, whereas the actual embedding need not have one.
The following calculation tests this replacement in the embedding
$X\hookrightarrow Y=\P^4$ used in
Proposition~\ref{prop:Ext-example}. Its role is to determine whether
one may compute the normal higher operations using that replacement.
The fiber embedding $X_b\hookrightarrow\cX$ has the different normal
line bundle $\O_{X_b}\otimes T_{C,b}$ and is treated by the
Atiyah--Kodaira--Spencer map \cite[Eq.~(4.32)]{KS}.
\begin{proposition}\label{prop:quintic-formal}
For a smooth quintic threefold $X\subset\P^4$, the formal completion of $\P^4$ along $X$ is not isomorphic over $X$ to the formal completion of $\Tot(N_{X/\P^4})$ along its zero section.
\end{proposition}
\begin{proof}
An isomorphism over $X$ would compose with the bundle projection to give a formal retraction to $X$. Its first infinitesimal restriction would split the conormal sequence, equivalently the tangent sequence
\[
 0\longrightarrow T_X\longrightarrow T_{\P^4}|_X
 \longrightarrow\O_X(5)\longrightarrow0.
\]
The Euler sequence restricted to $X$, together with $H^1(X,\O_X)=0$, gives $h^0(T_{\P^4}|_X)=5h^0(\O_X(1))-1=24$. The hypersurface sequence gives $h^0(\O_X(5))=\binom94-1=125$. A splitting would make the map on global sections surjective, a contradiction.
\end{proof}
The divisor ideal $\mathcal I_X$ has associated graded algebra
\[
 \bigoplus_{m\geq0}\mathcal I_X^m/\mathcal I_X^{m+1}
       \simeq\Sym_{\O_X}(N_{X/Y}^\vee).
\]
The extension from this graded algebra to the completed neighborhood
is part of the embedding. For the moduli problem it acts on $E$
through the class
\[
 E\xrightarrow{\At(E)}E\otimes\mathbb L_X[1]
   \xrightarrow{\id_E\otimes\kappa_i[1]}E\otimes N^\vee[2].
\]
This is $\alpha_E$ in \cite[equation~(4.26)]{KS} and produces the
connecting maps in \eqref{eq:normal-deformation}. Proposition~
\ref{prop:quintic-formal} therefore excludes a uniform replacement of
all these objectwise calculations by the zero-section model. For a
particular $E$ the displayed composite may nevertheless vanish.
For example, the hyperplane complex
$[\O_X(-1)\xrightarrow{\ell}\O_X]$ is the restriction of the
same two-term vector-bundle complex on $\P^4$; it extends to the
first infinitesimal neighborhood of $X$. Its class $\alpha_E$ is
zero by the obstruction calculation in \cite{HT}. This is the
objectwise vanishing used in the hyperplane application.

\subsection{Exhaustion of the section and Grassmannian parameter spaces}
\label{ssec:section-locus-exhaustion}
The words ``locus'' and ``smooth'' in this section refer to the
specified moduli functors, not to discarding singular points of those
functors.  For clarity, the parameter spaces are
\[
 \mathbb P_C(V_m),\qquad \operatorname{Gr}_C(q,V_m),
 \qquad V_m=\pi_*L^m.
\]
Every point of each space was included in the deformation
reconstruction, and both are smooth over the base curve.  In
particular their singular classical loci are empty.  On a reduced
critical chart of relative dimension $d$ their algebra still is
\[
 \O_U[\eta_1,\ldots,\eta_d],\quad |\eta_i|=-1,\quad
 \delta\eta_i=0.
\]
The degree-one tangent cohomology is $\Omega_U$ and is not discarded.
The zero normalized d-critical structure is forced by the full
classical deformation calculation.  Replacing it by a genuinely
nonzero d-critical class would change this moduli problem.

One can enlarge the presentation by allowing dependent sections,
but this does not add a singular stratum to the same simple-object
Grassmannian.  Suppose a map
$L_b^{k-m}\to L_b^k\otimes W$, with $\dim W=q$, has section rank
$s<q$.  A constant change of basis gives the direct sum
\[
 E\simeq [L_b^{k-m}\to L_b^k\otimes W']
             \oplus (L_b^k)^{\oplus(q-s)},\qquad\dim W'=s.
\]
When $s>0$ both summands are nonzero, so their projections give
nontrivial idempotents in $\operatorname{End}(E)$; the object is not
simple.  When $s=0$ the entire differential is zero and the source
shift is another nonzero direct summand.  Such objects require a
stack-valued invariant with their larger automorphism groups, not
the scalar-rigidified invariant calculated above.  Rank here is
the rank of the vector-space map to $H^0(X_b,L_b^m)$, not the
rank of evaluation at an individual point of $X_b$.

The full projective and Grassmannian moduli problems are also
unchanged by a special coefficient choice for the quintic, as long
as the whole-fiber cohomology and semistable-family hypotheses
continue to hold.  Special coefficients can instead make a
\emph{different} object-moduli problem singular, as for the curve
ideals in Section~\ref{sec:cubic-geometry}.  A fixed Chern character
alone does not identify every perfect object with a section complex;
the present argument does not make that assertion.

\section{Point and curve ideals on moving surfaces, and general surface sheaves}
\label{sec:divisor-sheaves}
The examples in this section concern sheaves on a smooth quintic
$X_b$ and their specified perfect-complex limits when that quintic
degenerates to $X_1\cup_S X_2$. For each example we specify the actual objects
and classical moduli locus before computing a potential or a
cohomology group. 
We start by considering pure 2 dimensional sheaves with support on divisors moving in a fixed linear system in the threefold, such that the sheaf is defined as pushforward of ideal sheaf of one point lying scheme theoretically on the divisor, then two points, then a line, a line with one point,
and a line with two point.
We then treat the general incidence problem $Z\subset D$ and larger
linear systems. For us the term ``\textit{incidence condition}" always means containment as
closed subschemes. For each identified smooth moduli locus $M$,
the zero critical function and canonical orientation give
$\mathscr P_M=\CC_M[\dim M]$; the two-sided bar below includes its
oriented BV differential. Obviously when the moduli is singular, this
replacement by a constant sheaf is not made and the situation becomes more involved.

\subsection{Primitive support and the full rank-one description}
Let $X\subset\P^4$ be a smooth quintic, let $H=c_1(\O_X(1))$, and
normalize $[\mathrm{pt}]$ to have integral one. Write
\[
 \mathrm{ch}(F)=(0,H,\beta,n[\mathrm{pt}]),
 \qquad \beta\in H^4(X,\mathbb Q),\quad n\in\mathbb Q.
\]
Thus $n$ denotes the degree of $\mathrm{ch}_3(F)$ and not $\chi(F)$.
The Lefschetz theorem gives $\operatorname{Pic}(X)=\mathbb Z[\O_X(1)]$.

The word ``smooth'' in the hypothesis applies to the threefold $X$;
the support surface $D$ is allowed to be singular. This proposition
concerns one fixed smooth quintic, or any smooth fiber $X_b$ of the
degeneration. It is not a stability assertion on $X_0=X_1\cup_S X_2$.
In fact the class $H$ cannot split into two nonzero effective divisor
classes on this fixed $X$, since each would be a positive integral
multiple of the generator $H$. This is precisely the special case
noted in \cite[footnote~4]{GST}. Splitting of the surface when the
threefold itself degenerates, and splitting inside $|mH|$ for $m>1$,
are different situations. We discuss both explicitly below.

\begin{proposition}\label{prop:primitive-support}
Every pure two-dimensional coherent sheaf $F$ on $X$ with
$\mathrm{ch}_1(F)=H$ is scheme-theoretically supported on a unique
integral divisor $D\in|H|$. It has the form $j_*G$, where
$j:D\hookrightarrow X$ and $G$ is torsion-free of rank one on $D$.
It is Gieseker stable. If $D$ is smooth, there are a unique line bundle
$M$ and a unique zero-dimensional subscheme $Z\subset D$ such that
$G=M\otimes I_{Z/D}$, with $M=G^{\vee\vee}$.
\end{proposition}
\begin{proof}
The codimension-one cycle of $F$ is
$\sum_a r_a[D_a]$, where the $D_a$ are prime divisors and $r_a$ are
positive generic lengths. Every effective prime divisor has class
$m_aH$ for a positive integer $m_a$: positivity follows by intersecting
with $H^2$. Therefore $\sum_a r_am_a=1$. There is a single component
$D$, of class $H$, with generic length one. The ambient variety is
smooth, so $D$ is Cartier. A local equation of $D$ kills $F$ at its
generic point. Its image in $F$ is consequently supported in dimension
at most one; purity forces that image to vanish. Thus $I_DF=0$ and
$F=j_*G$. Purity over the integral surface says precisely that $G$ is
torsion-free; its generic rank is one.

A nonzero subsheaf of $G$ has rank one. If it is proper, its quotient
is nonzero and has dimension at most one. The leading quadratic
coefficients of the two Hilbert polynomials agree, whereas the
quotient polynomial is positive for sufficiently large arguments.
Dividing by the common leading coefficient proves the strict
Gieseker inequality. A subsheaf of $j_*G$ on $X$ is annihilated by the ideal of $D$,
so it is the pushforward of a unique subsheaf of $G$ on $D$.\\
Moreover $\chi(X,j_*G(aH))=\chi(D,G(aH|_D))$.
The same reduced-Hilbert-polynomial inequality therefore holds on $X$:
closed-immersion pushforward with the restricted polarization does not
introduce new subsheaves or change this inequality.
On a smooth surface the reflexive rank-one hull $M=G^{\vee\vee}$
is invertible and $M/G$ is zero-dimensional. Tensoring the inclusion
$G\subset M$ with $M^{-1}$ gives the ideal of a unique finite
subscheme $Z$. This proves the asserted uniqueness as well.
\end{proof}

When $X_b$ degenerates, the central support may be $D_1\cup D_2$,
with $D_i\subset X_i$. Neither purity nor the primitive class on the
smooth fiber proves stability of the whole central sheaf, or of its
component restrictions, for the chosen extending polarization.
Our central moduli problem uses the specified perfect matching triples;
component sheaves, when present, and their derived boundary restrictions
are kept separately. In higher-degree systems on a smooth fiber,
\cite[Introduction and Section~2]{GST} instead treats the Joyce--Song
pair $(I_{Z/X},s)$. A nonzero section is sufficient for that pair
because $I_{Z/X}$ is rank-one torsion-free and stable. It is not a
proof that its surface-sheaf cokernel is Gieseker stable on a reducible
support. A concrete destabilizing example is given in
Section~\ref{ssec:general-subschemes}.

The last description is used only on smooth supports. On a singular
integral divisor, the rank-one reflexive hull need not be a line
bundle. Also, even on a smooth $D$, its line bundle $M$ need not be isomorphic to $j^*\O_X(kH)$ for any integer $k$. Such extra line bundles occur on Noether--Lefschetz
loci in $|H|$ and have to be retained when the charge is fixed.

\begin{proposition}\label{prop:general-divisor-ch}
Let $D\in|H|$ be smooth, set $h=H|_D$, $\alpha=c_1(M)$ and
$\ell=\operatorname{length}(Z)$. For $F=j_*(M\otimes I_{Z/D})$ one has
\begin{align}
 \beta&=j_*\alpha-\frac12H^2,\label{eq:general-beta}\\
 n&=\frac12\alpha^2-\frac12\alpha h+\frac56-\ell,\label{eq:general-n}\\
 \chi(F(m))&=\frac52m^2+(H\beta)m+n+\frac{25}{6}\notag\\
 &=\frac52m^2+\left(\alpha h-\frac52\right)m
       +\frac12(\alpha^2-\alpha h)+5-\ell.\label{eq:general-Hilbert}
\end{align}
Products in \eqref{eq:general-n} are intersection numbers on $D$.
In particular, $n+25/6$ is an integer. If $d=\alpha h$, then
\begin{equation}\label{eq:defect-bound}
 0\leq\ell\leq \frac{d^2}{10}-\frac d2+\frac56-n.
\end{equation}
\end{proposition}
\begin{proof}
The normal line bundle is $\O_D(1)$. Grothendieck--Riemann--Roch
for the Cartier inclusion reads
\[
 \mathrm{ch}(j_*G)=j_*\bigl(\mathrm{ch}(G)\operatorname{td}
       (\O_D(1))^{-1}\bigr).
\]
In dimensions at most two,
\[
 \mathrm{ch}(G)=1+\alpha+\tfrac12\alpha^2-\ell[\mathrm{pt}]_D,
 \qquad \operatorname{td}(\O_D(1))^{-1}
       =\frac{1-e^{-h}}h=1-\tfrac12h+\tfrac16h^2.
\]
The degree-zero, degree-two, and degree-four parts give the three
Chern-character components, using $j_*h=H^2$ and $\int_Dh^2=5$.
Adjunction gives $K_D=h$. Hence Riemann--Roch on $D$ shows
$\alpha(\alpha-h)/2\in\mathbb Z$. Alternatively this integrality
follows directly from the resulting Euler characteristic.
The tangent sequence of the quintic gives $c_2(X)=10H^2$ and $H^3=5$.
Multiplying $\mathrm{ch}(F)e^{mH}\operatorname{td}(X)$ in degree six
therefore gives
$\frac52m^2+(H\beta)m+n+\frac1{12}\int_XHc_2(X)$,
which is \eqref{eq:general-Hilbert}.
Finally the Hodge index theorem gives $\alpha^2\leq d^2/5$.
Solving \eqref{eq:general-n} for $\ell$ proves \eqref{eq:defect-bound}.
\end{proof}

For the restricted line bundle $M=j^*\O_X(kH)=\O_D(k)$ these formulas become
\begin{align}
 \mathrm{ch}(F)&=e^{kH}(1-e^{-H})-\ell[\mathrm{pt}],\label{eq:ambient-defect-charge}\\
 \beta&=(k-\tfrac12)H^2,\qquad
 n=\tfrac52k^2-\tfrac52k+\tfrac56-\ell,\notag\\
 \chi(F(m))&=\tfrac52(m+k)^2-\tfrac52(m+k)+5-\ell.\notag
\end{align}
For $\ell=0$ this recovers the hyperplane-section examples in
Section~\ref{sec:v6-hyperplanes}. Positive $\ell$ supplies additional
charges rather than additional presentations of those same objects.

\begin{example}\label{ex:NL-line}
Suppose $X$ contains a line $C_0$, and choose a smooth hyperplane
section $D$ containing it. Such a choice is possible whenever the
restricted linear system has a smooth member along $C_0$; this is the
hypothesis on $D$ here. Adjunction on $D$ gives
$C_0^2+C_0h=-2$, so $C_0h=1$ and $C_0^2=-3$.
For $M=\O_D(C_0)$ and a length-$\ell$ subscheme $Z$, the resulting
stable sheaf satisfies
\[
 \beta=[C_0]-\tfrac12H^2,\qquad n=-\tfrac76-\ell,
 \qquad \chi(F(m))=\tfrac52m^2-\tfrac32m+3-\ell.
\]
Since $C_0h=1$ is not divisible by $5$, this determinant is not
$\O_D(k)$ for any integer $k$. It is a concrete charge outside the
subfamilies obtained by restricting $\O_X(kH)$. The support and determinant deformations
must remain on the locus where this class is algebraic.
This positive line bundle is distinct from the ideal of the line:
$\O_D(-C_0)\otimes I_{Z/D}$ instead has
$\beta=-[C_0]-\tfrac12H^2$ and $n=-\tfrac16-\ell$.
The complete line-ideal calculation is
Theorem~\ref{thm:line-full-comparison}.
\end{example}

\subsection{Relative surface Hilbert schemes and derived point modifications}
Let $U\subset|H|$ parametrize smooth divisors and let
$\mathscr D\to U$ be their universal surface. The relative Hilbert
scheme $\operatorname{Hilb}^{\ell}(\mathscr D/U)$ is smooth over $U$
of relative dimension $2\ell$, by the smooth-surface Hilbert-scheme
theorem \cite{Fogarty}. Its tautological ideal, tensored with
$\O_{\mathscr D}(k)$, gives a flat family of the sheaves in
\eqref{eq:ambient-defect-charge}. Because $U\subset|H|=\P^4$ has dimension four and its surface
Hilbert-scheme fibers have dimension $2\ell$, the scheme
$\operatorname{Hilb}^{\ell}(\mathscr D/U)$ has dimension $4+2\ell$. This family is proper over $U$, but $U$ itself omits
singular hyperplanes. \\ A proper completion over the degeneration curve
is a separate requirement of the specialization theorem.
For a line bundle on $D$ not isomorphic to $j^*\O_X(kH)$, the analogous
family is taken over a specified relative Picard or
Noether--Lefschetz locus carrying that line bundle (locally on its
Picard gerbe when necessary). Its dimension and smoothness depend on
that locus.

The dimension statement just made does not classify the smooth and
singular loci of sheaf moduli over all hyperplanes: it excludes singular
surface supports. For length two the complete generic-fiber sheaf
moduli is instead the smooth eight-dimensional projective bundle in
Theorem~\ref{thm:two-point-incidence}, even though some supporting
hyperplanes are singular. The next proposition is an object-level
construction, related to Proposition~\ref{prop:general-divisor-ch}.

The following construction gives actual perfect-complex limits on a
specified part of the family without assuming that every limit is a
sheaf. It extends the section-complex construction, discussed in previous sections, by point defects.

\begin{proposition}\label{prop:derived-point-modification}
Let $\pi:\mathcal X\to(C,0)$ be the same resolved quintic
degeneration as in Section~\ref{ssec:quintic-construction}; $C$ is its
smooth base curve, not a curve subscheme of a surface. For a finite-type
morphism $\tau:T\to C$, put
$\mathcal X_T=\mathcal X\times_C T$ and let $p_T:\mathcal X_T\to T$
be projection. Pull back $L=a^*\mathcal O_{\mathbb P^4}(1)$ to this
space. Fix $m\geq1$ and $k\in\mathbb Z$, and let $A_T$ be the perfect
complex of an actual section
\[
 A_T=[L^{k-m}\otimes p_T^*\mathcal A\xrightarrow{s}L^k],
 \qquad |L^{k-m}\otimes p_T^*\mathcal A|=-1,
\]
where $\mathcal A$ is a line bundle on $T$. Suppose
$Z\subset\cX_T$ is finite flat of degree $\ell$ over $T$, lies in
the relative smooth locus of $\cX_T/T$, and $s|_Z=0$.
The evaluation chain map $A_T\to Q_T:=L^k|_Z$ defines
\begin{equation}\label{eq:point-modification-triangle}
 E_T\longrightarrow A_T\longrightarrow Q_T\longrightarrow E_T[1].
\end{equation}
Then $E_T$ is perfect and this triangle commutes with every derived
base change in $T$ and every derived restriction to $X_i$ and $S$.
At a point $v\in T$ mapping to $b\in C$, where $s_v$ cuts out a
Cartier divisor $D_v\in|mH_b|$, the fiber is
$j_{D_v*}I_{Z_v/D_v}\otimes L_b^k$ on $X_b$. In particular the assertion
covers both $|H_b|$ and $|2H_b|$. Over $\tau^{-1}(0)$, restricting the
same complex to $X_1,X_2,S$ supplies an effective matching triple,
retaining all cohomology of the three terms. The relative-smooth-support
condition is a sufficient condition for this construction, not a condition
on the full moduli problem: explicit perfect replacements at $S$ are
constructed in Section~\ref{sec:crossing-perfect}.
\end{proposition}
\begin{proof}
A finite flat family $Z/T$ has a perfect structure sheaf relative to
$T$. The local smooth-coordinate resolution makes the relevant
perfectness criterion concrete. Over $T=\Spec R$, let
$A=R[x_1,\ldots,x_d]$ and let $M$ be the finite projective
$R$-module of a point quotient. The actions of $x_j$ are commuting
endomorphisms $X_j$ of $M$. An $A$-projective resolution is
\[
 \left(A\otimes_R M\otimes_R\Lambda^\bullet R^d,\quad
       \sum_{j=1}^d(x_j\otimes1-1\otimes X_j)\iota_{e_j^*}\right)
       \longrightarrow M.
\]
The exterior degree $r$ is in cohomological degree $-r$.
Successive division by the monic operators $x_j-X_j$ identifies
its quotient with $M$ and proves exactness in negative degrees.
The same finite-resolution criterion is local in the etale topology
for a smooth morphism. Thus a coherent $T$-perfect quotient
supported in the relative smooth locus is perfect on the ambient
scheme. It applies to $Q_T=L^k|_Z$. Extending by zero preserves
perfectness here: the support is closed and contained in the
smooth open, and outside that support the complex vanishes.
The source $A_T$ is already a bounded complex of line bundles.

The degree-zero evaluation $L^k\to L^k|_Z$ kills the image of $s$
by the stated incidence condition. It is therefore a chain map.
The category of perfect complexes is closed under fibers, proving
perfectness of $E_T$. Derived pullback is exact and preserves both
this chain map and its fiber, which proves every compatibility in
the statement. If $s$ is a non-zero-divisor, $A_T$ on that fiber
resolves $j_*\O_D(k)$; evaluation onto $Q_T$ is surjective. Its
kernel is exactly $j_*I_{Z/D}(k)$. Finally $E_0$ is a perfect complex
on $X_0$, and associativity of derived restriction identifies its two
restrictions to $S$. This is the matching equivalence of its effective
triple.
\end{proof}

The same fiber construction works with any relative perfect complex
$A_T$ and any relative perfect quotient complex $Q_T$, once an actual
morphism is specified. It imposes no bound of two on the number of
cohomology sheaves of $E_T$. The point construction above establishes
perfectness on its stated open locus; collision with the common surface is not covered by that sufficient
hypothesis. Section~\ref{sec:crossing-perfect} provides explicit
replacement complexes and length-two perfect quotients there.

\textbf{The relative Hilbert scheme over singular supporting surfaces.}
The open $U$ above excludes singular surfaces, but it is not the
parameter space used for the complete one-point and two-point
calculations below.  Those calculations use
$\{(Z,D):Z\subset D\in|mH|,\ \operatorname{length}Z=r\}$ over
\emph{all} divisors, with $r=1,2$ and $m=1,2$.  The evaluation-kernel
proof identifies their whole generic spaces with projective bundles
over $X$ and $X^{[2]}$.  Thus a singular support or a tangent double
point is already present in them.  For curve defects, by contrast,
$\operatorname{Quot}_X(I_\Lambda,2)$ is genuinely singular, and
Section~\ref{sec:global-two-defect} computes its full intrinsic DT
complex rather than using the smooth-surface Hilbert-scheme theorem.
Higher lengths and arbitrary curve classes are governed by their
actual cyclic deformation potentials, not by the dimension formula
$4+2r$ for the smooth-support open.

\subsection{Reconstructing the sheaf or perfect complex from its point data}
\label{ssec:point-reconstruction}
The charge formula of Proposition~\ref{prop:general-divisor-ch} and
Proposition~\ref{prop:derived-point-modification} have different roles.
The former computes the Chern character on a smooth support. The
latter constructs a perfect object from an actual evaluation map,
including on the singular threefold, provided the point quotient is
perfect there. Neither proposition classifies the singular locus
of the full sheaf moduli space. 
For any fiber $Y=X_b$, write $L_Y=L|_Y$, and let $Z\subset Y$ be a
nonempty finite subscheme. Its ideal \emph{in the threefold} is
$I_{Z/Y}$. If a section $s\in H^0(Y,I_{Z/Y}L_Y)$ cuts out a Cartier
divisor $D$, there is an exact sequence
\begin{equation}\label{eq:incidence-ideal-sequence}
 0\longrightarrow L_Y^{-1}\xrightarrow{s}I_{Z/Y}
 \longrightarrow j_{D*}I_{Z/D}\longrightarrow0.
\end{equation}
The last ideal is taken \emph{on the surface}. We write
\begin{equation}\label{eq:incidence-perfect-object}
 E_{Z,s}:=\operatorname{Cone}(L_Y^{-1}\xrightarrow{s}I_{Z/Y}).
\end{equation}
Here $L_Y^{-1}$ is in degree $-1$ and $I_{Z/Y}$ in degree zero.
When $s$ is not a non-zero-divisor, this notation means the whole
complex, not its zeroth cohomology. By the octahedral triangle it
also equals the homotopy fiber of
$[L_Y^{-1}\xrightarrow{s}\O_Y]\to\O_Z$.

\begin{proposition}\label{prop:point-moduli-reconstruction}
Let $Y$ be a connected projective Gorenstein threefold with a chosen
isomorphism $\omega_Y\simeq\O_Y$, with
$H^0(\O_Y)=H^3(\O_Y)=\CC$ and $H^1(\O_Y)=H^2(\O_Y)=0$.
Let $Z\subset Y$ be finite of length $r>0$, with $\O_Z$ perfect
on $Y$. Assume
\[
 H^{>0}(Y,I_{Z/Y}L_Y)=0,
 \qquad W_Z=H^0(Y,I_{Z/Y}L_Y),\qquad 0\ne s\in W_Z.
\]
Then $E_{Z,s}$ is simple and has no negative self-extensions.
On ordinary local Artinian bases, its deformations, modulo scalar
automorphisms, are exactly the deformations of the incidence pair
$(Z,[s])$. The assertion also holds for deformations of $Y$ induced
from the given oriented family, provided $L_Y$ is the restriction
of $L$ and the stated fiberwise vanishing holds. Consequently a
finite-presentation incidence family satisfying these conditions
maps by an open immersion to the classical scalar-rigidified
perfect-complex moduli space.
\end{proposition}
\begin{proof}
Put $I=I_{Z/Y}$ and $E=E_{Z,s}$. Both are perfect, by the triangle
$I\to\O_Y\to\O_Z$ and \eqref{eq:incidence-perfect-object}.
Serre duality and the vanishing hypothesis give
\[
 R\operatorname{Hom}_Y(I,L_Y^{-1})\simeq W_Z^\vee[-3].
\]
Apply $R\operatorname{Hom}_Y(-,L_Y^{-1})$ to
$L_Y^{-1}\to I\to E$. Its long exact sequence gives
\begin{equation}\label{eq:reconstruct-Ext-gap}
 \operatorname{Ext}^1_Y(E,L_Y^{-1})=\CC,
 \qquad \operatorname{Ext}^3_Y(E,L_Y^{-1})=\ker(s^\vee:W_Z^\vee\to\CC),
\end{equation}
and all other groups vanish. The arrow $s^\vee$ is surjective because
$s\ne0$. The degree-one generator is the connecting morphism of
the triangle defining $E$.

Let $A$ be an ordinary local Artinian algebra, let $Y_A$ be the
induced family and $E_A$ a perfect deformation. Proper perfect
pushforward makes $R\operatorname{Hom}_{Y_A}(E_A,L_A^{-1})$ a
perfect $A$-complex. A minimal finite-free representative has terms
only in degrees one and three, by \eqref{eq:reconstruct-Ext-gap}.
There is no differential between them. In particular its degree-one
cohomology is free of rank one and commutes with residue base
change. A generator reconstructs a triangle
\begin{equation}\label{eq:reconstruct-ideal-triangle}
 L_A^{-1}\longrightarrow I_A\longrightarrow E_A
 \longrightarrow L_A^{-1}[1]
\end{equation}
lifting the given triangle.

The sequence $I\to\O_Y\to\O_Z$ gives
$H^0(I)=0$, $H^1(I)=\CC^{r-1}$, $H^2(I)=0$, $H^3(I)=\CC$.
Duality therefore gives
\[
 \operatorname{Hom}_Y(I,\O_Y)=\CC,
 \quad\operatorname{Ext}^1_Y(I,\O_Y)=0,
 \quad\operatorname{Ext}^2_Y(I,\O_Y)=\CC^{r-1},
\]
with no other groups. The same minimal-complex argument, now in
degrees zero and two, lifts the inclusion $I\to\O_Y$ uniquely up
to a unit to $I_A\to\O_{Y_A}$. Its cone $Q_A$ has derived residue
fiber $\O_Z$ in degree zero. The derived flatness criterion over a
nilpotent Artinian base shows that $Q_A$ is an $A$-flat coherent
sheaf in degree zero. Explicitly, the local criterion for flatness,
applied successively to small extensions of $A$, says that a
bounded $A$-perfect complex whose derived residue is in degree
zero is represented by an $A$-flat sheaf. Derived Nakayama makes
its support a nilpotent thickening of the finite support of $Z$;
properness then makes that support finite over $A$.

The map $\O_{Y_A}\to Q_A$ is surjective, since its reduction is
surjective and its cokernel is zero by Nakayama. Thus
$Q_A=\O_{Z_A}$ for a finite flat subscheme of length $r$, and
$I_A=I_{Z_A/Y_A}$. The first arrow of
\eqref{eq:reconstruct-ideal-triangle} is a section of
$I_{Z_A/Y_A}L_A$, reducing to $s$. The two choices of generators
multiply it by units, giving precisely its projective class.
An isomorphism of the complexes respects these one-dimensional
spaces, so recovers the same subscheme and projective section.
Conversely an isomorphism of the incidence data gives the
isomorphism of the displayed cone. This proves the equivalence
of the classical infinitesimal deformation problems, including
base deformations.

For completeness, $\operatorname{End}(I)=\CC$: its inclusion into
$\operatorname{Hom}(I,\O_Y)=\CC$ identifies it with the scalars.
Applying $R\operatorname{Hom}(E,-)$ and
$R\operatorname{Hom}(-,I)$ to the defining triangle, and using
$R\Gamma(I L_Y)=W_Z$ in degree zero, shows
$\operatorname{Ext}^{<0}(E,E)=0$. Moreover
$\operatorname{Hom}(E,I)$ is the kernel of
$\CC\to W_Z$, $1\mapsto s$, and hence is zero.
The next part of the same exact sequence embeds
$\operatorname{End}(E)$ into the copy of $\CC$ in
\eqref{eq:reconstruct-Ext-gap}; the identity of $E$ is nonzero,
so this is an isomorphism. Simplicity and universal gluability
therefore hold. The family morphism to the algebraic space of
\cite{Lieblich2006} is formally etale by the reconstruction, and
universally injective by the same argument over algebraically
closed fields. On finite-presentation loci it is an open immersion.
\end{proof}

On $X_0$, the normalization triangle identifies the endomorphism
algebra of each such $E_{Z,s}$ with the deformation algebra of its
matching triple. Thus the same reconstruction identifies its
classical locus inside
$\RPerf(X_1)\times^h_{\RPerf(S)}\RPerf(X_2)$.
If that locus is smooth of dimension $d$, its d-critical structure
is zero, and the chosen canonical orientation gives
$\CC[d]$. Its oriented two-sided bar realization is
\begin{equation}\label{eq:point-smooth-bar}
 \left(\operatorname{Bar}_{\operatorname{Sym}T_M}
       (\O_M,\O_M)\otimes K_M\right)((u)),
 \qquad D_B=-b+uI\Delta\Pi,
\end{equation}
with the comparison
$\mathcal B^{\rm or}_{M,\mathrm{per}}\simeq
 (\Omega_M^\bullet((u)),u d)[d]$.
This is a critical-chart realization of the actual matching locus;
its underlying structure-algebra bar alone does not have the DT
differential. Compact supports are used when $M$ is not proper.

\textbf{A point on the common divisor is a different case.}
The perfectness hypothesis above cannot be removed. At a point of
$S$ use the completed local ring
\[
 R=\CC[[x,y,z,w]]/(xy),\qquad D=(z=0),\qquad p=(x,y,z,w).
\]
The module $R/(z)$ is perfect, since $z$ is a non-zero-divisor.
The residue field $k(p)$ is not perfect: a residue field has finite
projective dimension only over a regular local ring, whereas
$R$ has dimension three and embedding dimension four. From
\begin{equation}\label{eq:node-ideal-obstruction}
 0\longrightarrow I_{p/D}\longrightarrow R/(z)
 \longrightarrow k(p)\longrightarrow0
\end{equation}
it follows that $I_{p/D}$ is not perfect over $R$. Otherwise the
third term would be perfect as the cone of two perfect objects.
This is the local form of an actual collision: in $xy=t$, after
$t=s^2$, take the point $x=y=s$, $z=w=0$ on the divisor $z=0$.
Its ordinary sheaf limit is \eqref{eq:node-ideal-obstruction}.
Thus extending the projective incidence parameter space does not
by itself extend its universal sheaf into $\RPerf(X_0)$.

There are also perfect length-two subschemes at the crossing.
For $a\ne0$,
\begin{equation}\label{eq:node-two-lci}
 R/(x-a y,z-b y,w-c y)\simeq\CC[[y]]/(y^2)
\end{equation}
is cut out by a regular sequence of length three in $R$ and is
perfect. This example is important for the scope of the calculations:
restricting all point supports to $X_0\setminus S$ gives a definite
open moduli problem, but it does not exhaust every possible
perfect length-two limit. Propositions~\ref{prop:crossing-two-perfect} and
\ref{prop:crossing-perfect-exhaustion} calculate the complete
length-two chart, including its singular parameter stratum.
They distinguish perfect quotients at the crossing from nonperfect
ones.  A global identification of the image of all replacement
complexes is not implied by adjoining these local parameters.

\subsection{One point on a moving hyperplane section}
\label{ssec:one-point-stepwise}
\subsubsection{The moduli problem and its generic potential}\label{ssec:incidence-monodromy}
An incidence pair in this example is a point $x\in X_b$ together
with a hyperplane divisor $D\in|H|$ containing that point. In families,
$Z\subset X_b\times T$ is finite flat of degree one over $T$,
$D\subset X_b\times T$ is a relative hyperplane divisor, and the
containment $Z\subset D$ holds as a closed subscheme. The sheaf is
$j_{D*}I_{Z/D}$, not the ideal $I_{Z/X_b}$ on the threefold.
The parameter space contains singular hyperplane surfaces; its own
smoothness is proved below. Its dimension is six, from three dimensions
for $x$ and three for a hyperplane containing $x$.

The preceding projective loci have globally defined horizontal hyperplane classes. We now examine a different moduli locus, where the cohomology remembers the monodromy of the Calabi--Yau threefold itself.
Let $V_1=H^0(\P^4,\O_{\P^4}(1))$. The evaluation sequence on $\cX$ is
\begin{equation}\label{eq:incidence-bundle}
0\longrightarrow\mathcal E\longrightarrow V_1\otimes\O_{\cX}
\longrightarrow L\longrightarrow0,
\qquad Y=\P_{\cX}(\mathcal E),\ \rho=\pi p.
\end{equation}
On a smooth fiber, $Y_b=\{(x,[\ell]):x\in X_b,\ \ell(x)=0\}$. Associate to this pair the sheaf
\begin{equation}\label{eq:incidence-sheaf}
F_{x,\ell,b}=j_{\ell,b*}I_{x/D_{\ell,b}},\qquad D_{\ell,b}=\{\ell=0\}\cap X_b,
\end{equation}
where the ideal is taken on the divisor $D_{\ell,b}$.

\begin{proposition}\label{prop:incidence-locus}
For $b\ne0$, $Y_b$ is a smooth projective six-dimensional locally complete moduli locus of stable pure two-dimensional sheaves \eqref{eq:incidence-sheaf}. With orientation $K_{Y_b}$ its DT sheaf is $\CC_{Y_b}[6]$ and $\mathrm{DT}^j_{\mathrm{coh}}(Y_b)=H^{j+6}(Y_b,\CC)$. Let $\xi=c_1(\O_Y(1))$. There are isomorphisms of local systems, also horizontal for the Gauss--Manin connection,
\[
 \bigoplus_{a=0}^{3} R^{k-2a}\pi_*\CC\big|_{C^*}
 \xrightarrow{\;\sum\xi^ap^*\;}R^k\rho_*\CC\big|_{C^*},
\qquad C^*=C\setminus\{0\}.
\]
\end{proposition}
\begin{proof}
A hyperplane divisor of a smooth quintic is integral by the Lefschetz Picard theorem: every nonzero effective prime divisor has positive class and $\operatorname{Pic}(X_b)=\Z[L_b]$. The ideal of a point on this Cartier divisor is torsion-free of rank one and Gieseker stable; any proper subsheaf has smaller reduced Hilbert polynomial. The vanishing
\[
 H^i(X_b,I_{x/X_b}\otimes L_b)=0\ (i>0),\qquad
 h^0(X_b,I_{x/X_b}\otimes L_b)=4
\]
follows from evaluation of sections and $H^{>0}(X_b,L_b)=0$. It verifies the positivity assumption in the relative Hilbert-scheme identification of the local sheaf deformation problem \cite{GST}. The projective bundle $Y_b\to X_b$ has fiber $\P^3$, hence is smooth projective of dimension six. Its smooth d-critical structure is zero and the chosen orientation gives the stated DT sheaf. Leray--Hirsch with the globally defined Chern class $\xi$ proves the decomposition; its powers are horizontal since they are Chern classes of a line bundle on the total family.
\end{proof}

\subsubsection{The generic DT groups and the connection}
The normalized potential on the smooth moduli $Y_b$ is $f_b=0$.
With orientation $K_{Y_b}$ the DT sheaf is $\CC_{Y_b}[6]$, and its
specified categorical realization is $\operatorname{MF}(Y_b,0)$,
with the shift $[6]$ in the normalized periodic comparison.
The full generic cohomology is given by
\begin{equation}\label{eq:one-point-generic-polynomial}
 P_{Y_b}(z)=(1+z^2+z^4+z^6)
                (1+z^2+204z^3+z^4+z^6),
 \qquad \mathrm{DT}^j_{\mathrm{coh}}(Y_b)=H^{j+6}(Y_b,\CC).
\end{equation}
Its even and odd periodic dimensions are $16$ and $816$,
respectively. The connection is the direct sum of the quintic
Gauss--Manin systems under the displayed projective-bundle maps;
$1,\xi,\xi^2,\xi^3$ are horizontal. Although the ordinary potential
is zero, $Y_b$ is not a constant family of varieties, so the
connection need not have zero matrix in an algebraic frame.
The horizontal-lift construction of Theorem~\ref{thm:v5-Cech-GM}
with $F=0$ is the relative de Rham Gauss--Manin connection.
It is flat on $C^*$, and its monodromy is computed next.

\begin{theorem}\label{thm:incidence-monodromy}
For the resolved quintic degeneration, let $N_X$ denote the logarithm of unipotent monodromy on $H^3(X_b,\Q)$. Then
\[
 N_X^2=0,\qquad\operatorname{rank}(N_X)=21.
\]
The induced monodromy logarithm on the odd periodic DT cohomology of $Y_b$ has rank $84$. Consequently its natural Gauss--Manin connection cannot extend across $0$ as a finite locally free connection without a pole, analytically or over $\CC((u))[[t]]$.
\end{theorem}
\begin{proof}
The weight spectral sequence of the proper semistable two-component degeneration gives
\[
\begin{aligned}
\operatorname{Gr}^W_2H^3_{\lim}&\simeq
\operatorname{coker}\big(H^2(X_1)\oplus H^2(X_2)\xrightarrow{r_1-r_2}H^2(S)\big),\\
\operatorname{Gr}^W_4H^3_{\lim}&\simeq
\ker\big(H^2(S)(-1)\xrightarrow{(\iota_{1*},-\iota_{2*})}H^4(X_1)\oplus H^4(X_2)\big).
\end{aligned}
\]
Here the first map is restriction and the second is its Poincare-dual Gysin map. By the monodromy--weight theorem the nilpotent operator induces an isomorphism from the weight-four graded piece to the weight-two graded piece with a Tate twist, and $N_X^2=0$. The quartic threefold has $H^2(X_1,\Q)=\Q H_1$, whereas the blowup $X_2$ has $H^2(X_2,\Q)=\Q H_2\oplus\Q E$. Both $H_1$ and $H_2$ restrict to the hyperplane class $H_S$ on the K3 surface $S$, and $E|_S=[\Gamma]=5H_S$. Consequently the image of restriction has rank one. The Gysin map has the same rank. Since $\dim H^2(S,\Q)=22$, both displayed graded pieces have dimension $21$. This gives $\operatorname{rank}(N_X)=21$.

The odd cohomology of a smooth quintic threefold is concentrated in degree three. The projective-bundle decomposition of Proposition~\ref{prop:incidence-locus} contains four copies of this local system, shifted to degrees $3,5,7,9$. The DT shift $[6]$ does not change parity. Hence the logarithm on odd periodic cohomology has rank $4\cdot21=84$.

If a holomorphic connection extended without a pole on a disk, its matrix differential equation would admit a single-valued horizontal frame, forcing trivial local monodromy, a contradiction. Formally, a unipotent regular-singular connection with nonzero nilpotent residue $R$ has meromorphic normal form $\nabla=\partial_t-t^{-1}R$. A Laurent-series horizontal section $v=\sum_jv_jt^j$ must satisfy $(j\,1-R)v_j=0$. For $j\neq0$ the matrix is invertible, and for $j=0$ the solutions lie in $\ker R$, giving fewer horizontal sections than the rank. A pole-free free connection over $\CC((u))[[t]]$ has a full horizontal frame by coefficient recursion, again a contradiction.
\end{proof}

\begin{proposition}\label{prop:incidence-log-extension}
After restricting the base analytically to a disk $\Delta$ about $0$, the projective bundle $\rho:Y\to\Delta$ is semistable and has central fiber
\[
Y_0=\P_{X_1}(\mathcal E|_{X_1})\cup_{\P_S(\mathcal E|_S)}
\P_{X_2}(\mathcal E|_{X_2}).
\]
The logarithmic de Rham complex $\Omega^\bullet_{Y/\Delta}(\log Y_0)$ computes the locally free logarithmic extension of the Gauss--Manin system. Its fiber at $0$ is canonically the limiting nearby-cycle cohomology $H^k_{\lim}(Y_b,\CC)$, and the connection has nilpotent residue $-R_k$ with $R_k=(2\pi i)^{-1}\log T_k$. After periodization and the DT shift $[6]$, its odd residue has rank $84$.
\end{proposition}
\begin{proof}
The projective bundle is locally trivial on the smooth total space $\cX$. Near a crossing its morphism to the base is locally $xy=t$ times smooth projective-chart directions, hence semistable. Relative logarithmic differentials have relation $d\log x+d\log y=d\log t$; modulo base logarithmic forms, $d\log y=-d\log x$, giving a local frame. The logarithmic de Rham theorem for proper semistable maps \cite{Ste} gives the asserted locally free extension and identifies the central fiber with nearby-cycle cohomology. On a multivalued flat frame $v$ over $\Delta^*$, with continuation $v\mapsto vT_k$, the frame $e=v\exp(-R_k\log t)$ is single-valued because $\exp(-2\pi iR_k)=T_k^{-1}$. The connection matrix in this frame is $-R_k\,dt/t$. Its residue is therefore $-R_k$. The projective-bundle classes extend and remain horizontal, so the odd residue is four copies of the quintic residue, as in Theorem~\ref{thm:incidence-monodromy}.
\end{proof}

The logarithmic central fiber in the last proposition is nearby-cycle cohomology, not an identification with an arbitrary unconstructed ordinary DT bar complex. Thus the all-degree section-complex example and the incidence-sheaf monodromy calculation have genuinely different specialization behavior.

\subsubsection{The actual central matching locus with the point away from $S$}
\label{ssec:one-point-central}
We now compute a specified central moduli space, rather than
identifying the whole auxiliary fiber $Y_0$ with perfect objects.
Set
\[
 U_i=X_i\setminus S,\qquad U=X_0\setminus S=U_1\amalg U_2,
 \qquad M_{1,0}^{\circ}=\P_U(\mathcal E|_U).
\]
The superscript $\circ$ is part of the notation: only the point
support is required to avoid $S$. Every nonzero hyperplane section
vanishing at that point is retained, including sections whose
central cone has negative cohomology. The variety
$M_{1,0}^{\circ}$ is a smooth six-dimensional, generally nonproper,
classical moduli locus. Proposition~\ref{prop:point-moduli-reconstruction}
applies because $p\in U$ is a perfect point of $X_0$ and
$V\to L_p$ is surjective.

Here are its component objects. If $p\in U_i$, put
$Q_i=\O_p$, $Q_{3-i}=0$, and $Q_S=0$. Restricting the evaluation
triangle gives
\begin{equation}\label{eq:one-point-central-triples}
 \begin{aligned}
 A_i&=[L_i^{-1}\xrightarrow{s_i}\O_{X_i}],
 &E_i&=\operatorname{fib}(A_i\to Q_i),\\
 A_S&=[L_S^{-1}\xrightarrow{s_S}\O_S],
 &E_S&=A_S.
 \end{aligned}
\end{equation}
The two equivalences $L\iota_i^*E_i\simeq E_S$ are the restrictions
of the same triangle; they specify the matching map, not just its
class in a Grothendieck group. If $s_i$ cuts out a Cartier surface
$D_i$ on both components, these are
$j_{D_i*}I_{p/D_i}$ on the component containing $p$ and
$j_{D_{3-i}*}\O_{D_{3-i}}$ on the other. If in addition $s_S$ is
nonzero, the common object is the sheaf of the curve
$C_D=D_1\cap S=D_2\cap S$. When these non-zero-divisor conditions
fail, \eqref{eq:one-point-central-triples} is still the definition
and the whole boundary complex is used.

The deformation triangle of Section~2 identifies the actual
matching deformation algebra with
$R\operatorname{Hom}_{X_0}(E_{p,s},E_{p,s})$.
The reconstruction proposition therefore makes
$M_{1,0}^{\circ}$ a smooth locus in the scalar-rigidified derived
Lagrangian intersection. Its normalized potential is zero. With
orientation $K_{M_{1,0}^{\circ}}$ its bar/BV model is
\eqref{eq:point-smooth-bar} with $d=6$.

\begin{proposition}\label{prop:central-open-cohomology}
Let $P_c(T;z)=\sum_q\dim H_c^q(T,\CC)z^q$ denote the
compactly supported Poincare polynomial. For the two components
of the resolved quintic,
\begin{equation}\label{eq:component-open-cohomology}
 P_c(U_1;z)=81z^3+z^6,
 \qquad P_c(U_2;z)=z^2+123z^3+z^4+z^6.
\end{equation}
In particular write
\begin{equation}\label{eq:central-U-polynomial}
 P_U(z):=P_c(U;z)=z^2+204z^3+z^4+2z^6.
\end{equation}
The actual central open moduli locus above has
\begin{equation}\label{eq:one-central-DT}
 \mathrm{DT}_{\mathrm{coh}}^j(M_{1,0}^{\circ})
 =H_c^{j+6}(M_{1,0}^{\circ},\CC),\qquad
 P_c(M_{1,0}^{\circ};z)=(1+z^2+z^4+z^6)P_U(z).
\end{equation}
Its normalized periodic even and odd dimensions are $16$ and $816$.
\end{proposition}
\begin{proof}
For $X_1$, the hypersurface Chern-class calculation gives
$\chi(X_1)=-56$ and hence $b_3(X_1)=60$; its other nonzero
Betti numbers are $b_0=b_2=b_4=b_6=1$. The curve $\Gamma$
has degree $20$ and genus $51$, since
$2g(\Gamma)-2=(4+5-4)\cdot20=100$.
The blowup formula gives
$b_2(X_2)=b_4(X_2)=2$, $b_3(X_2)=102$, and
$b_0(X_2)=b_6(X_2)=1$.

In the localization sequence for $U_i\subset X_i$ the maps
$H^0(X_i)\to H^0(S)$ are isomorphisms. The image of
$H^2(X_i)\to H^2(S)$ is the one-dimensional span of $h_S$;
for $X_2$ use $E_\Gamma|_S=5h_S$. The maps
$H^4(X_i)\to H^4(S)$ have rank one. Since $S$ is a K3 surface,
$H^1(S)=H^3(S)=0$ and $b_2(S)=22$.
The long exact sequence therefore gives
\[
\begin{array}{c|rrrrrrr}
 &H_c^0&H_c^1&H_c^2&H_c^3&H_c^4&H_c^5&H_c^6\\\hline
 U_1&0&0&0&\CC^{81}&0&0&\CC\\
 U_2&0&0&\CC&\CC^{123}&\CC&0&\CC.
\end{array}
\]
For example $H_c^3(U_i)$ is an extension of $H^3(X_i)$ by
the cokernel of the restriction on $H^2$, of dimension $21$.
This proves \eqref{eq:component-open-cohomology}.

The moduli identification and zero critical structure were proved
above, so the oriented DT sheaf is the shifted constant sheaf.
The projective bundle $M_{1,0}^{\circ}\to U$ is proper of relative
dimension three. The projective-bundle formula applies with
compact supports and the global class $\xi=c_1(\O(1))$,
giving \eqref{eq:one-central-DT}. Periodization and the even
shift $[6]$ give $4\cdot4=16$ even and $4\cdot204=816$ odd
classes. These are the cohomology groups of the oriented bar/BV
complex of the stated matching locus, not cohomology of one
component stack separately.
\end{proof}

\subsubsection{The connection}
The generic invariant and the central open invariant have the
same parity dimensions, but their integral gradings differ:
\begin{equation}\label{eq:one-central-difference}
 P_c(M_{1,0}^{\circ};z)-P_{Y_b}(z)
  =(1+z^2+z^4+z^6)(z^6-1).
\end{equation}
Dimension agreement after periodization does not define a
horizontal specialization. On $C^*$ the connection is the four
shifted copies of the quintic Gauss--Manin connection displayed
in Proposition~\ref{prop:incidence-locus}. It is flat there.
Its logarithmic extension has the residue of
Proposition~\ref{prop:incidence-log-extension}, which is nonzero.
It cannot be replaced by a pole-free connection with the same
restriction to $C^*$.

There is a natural comparison with the central open locus that
can be examined directly. The proper parameter family
$\rho:Y=\P_{\cX}(\mathcal E)\to\Delta$ contains
$M_{1,0}^{\circ}$ as the complement of
$Z_S=\P_S(\mathcal E|_S)$ in its central fiber. On this complement
$R\psi_\rho\CC$ is the constant complex. Extension by zero gives
\begin{equation}\label{eq:one-open-nearby-map}
 R\Gamma_c(M_{1,0}^{\circ},\CC)[6]
 \longrightarrow R\Gamma(Y_0,R\psi_\rho\CC)[6]
 \longrightarrow R\Gamma(Z_S,(R\psi_\rho\CC)|_{Z_S})[6]
 \xrightarrow{+1}.
\end{equation}
The first term is the central DT realization just constructed;
the second is the limiting generic DT realization. The last term
is not zero. To see this without assuming a DT compactification,
the semistable plumbing model $xy=t$ retracts, over $S$, onto
the circle bundle of $N_{S/X_1}=\O_S(h_S)$.
The Gysin sequence for that circle bundle, whose Euler class is
$h_S$, gives cohomology dimensions
\[
 1,21,21,1\quad\hbox{in degrees }0,2,3,5.
\]
Indeed multiplication by $h_S$ has ranks one from $H^0(S)$ to
$H^2(S)$ and from $H^2(S)$ to $H^4(S)$.
Pullback to $Z_S$ adds the projective-bundle factor
$1+z^2+z^4+z^6$. Thus the third term of
\eqref{eq:one-open-nearby-map} has nonzero cohomology in both
parities, and the first arrow is not a quasi-isomorphism, even
after periodization. This calculation concerns the topological
nearby-cycle complex of the displayed proper parameter family;
it does not supply perfect objects at the excluded points.

Consequently the one-point example proves the generic cohomology,
the bar cohomology of this actual central open locus, and the
flat logarithmic connection, but does not satisfy the regular
comparison theorem for these two endpoints. A different, larger
central perfect-complex moduli problem and its oriented
specialization morphism would have to be constructed and tested
by Proposition~\ref{prop:specialization-cone-test}.

\subsubsection{One point on a member of $|2H|$.}
Let $V_2=H^0(X_b,\mathcal O_{X_b}(2))$, of dimension $15$.
The exact evaluation sequence
\[
 0\to K_{1,2}\to V_2\otimes\mathcal O_{X_b}
            \to\mathcal O_{X_b}(2)\to0
\]
gives a rank-$14$ bundle. The generic sheaf moduli in question is
$Y_{1,2}=\mathbb P_{X_b}(K_{1,2})$, with projective fiber
$\mathbb P^{13}$ and dimension $16$. Its objects are
$\operatorname{Cone}(\mathcal O_{X_b}(-2)\to I_{p/X_b})
=j_{D*}I_{p/D}$, for $D\in|2H_b|$ containing $p$ scheme-theoretically.
The reconstruction of Proposition~\ref{prop:point-moduli-reconstruction}
uses $L^2$ instead of $L$; the evaluation surjection and Kodaira
vanishing verify its hypothesis. These particular sheaves are stable
also for reducible or nonreduced $D$, by the slope calculation in
Lemma~\ref{lem:quadric-small-curve-stability}.

The classical space is smooth, its minimal potential is zero, and
with orientation $K_{Y_{1,2}}$ the generic DT groups are
\[
 \mathrm{DT}_{\mathrm{coh}}^j(Y_{1,2})
       =\bigoplus_{a=0}^{13}H^{j+16-2a}(X_b,\mathbb C).
\]
The normalized periodic dimensions are $(56,2856)$. On $C^*$ its
connection is the direct sum of the fourteen shifted copies of
$\nabla_X$; the powers of $\xi=c_1(\mathcal O_{Y_{1,2}}(1))$ are
horizontal. Its odd monodromy logarithm has rank $14\cdot21=294$.

For the central objects with $p\in U=(X_1\setminus S)\sqcup
(X_2\setminus S)$, the evaluation kernel still has rank fourteen,
since $L^2$ is generated by the pulled-back quadratic sections.
Thus this central open is a $\mathbb P^{13}$-bundle over $U$.
Its derived restrictions are the cone formulas of
\eqref{eq:incidence-perfect-object} with $L^{-2}$; its actual classical
matching locus is identified by the same reconstruction. Its reduced
potential is zero, and its oriented bar computes
$\bigoplus_{a=0}^{13}H_c^{j+16-2a}(U,\mathbb C)$.
The compact-support polynomial of $U$ calculated above gives periodic
dimensions $(56,2856)$. Equality of these two ranks is not a
specialization isomorphism: the natural generic monodromy is nontrivial.
Section~\ref{sec:crossing-perfect} constructs an additional perfect
point limit at $S$, rather than claiming that this open is the full
central moduli space.
\subsection{Two points on a moving hyperplane section}
\label{ssec:two-point-stepwise}
\subsubsection{The moduli problem, generic potential, and cohomology}
Here a ``two-point subscheme'' means a finite subscheme of length two.
It may be two distinct points or a nonreduced double point with a tangent
direction. The incidence pair is $Z\subset D\subset X$, with
$D\in|H|$, and the containment is scheme-theoretic. In particular,
a hyperplane containing a double point must vanish on its infinitesimal
tangent direction, not merely at its underlying closed point.
The sheaf $F_{Z,D}=j_{D*}I_{Z/D}$ is an object on the whole smooth
quintic. The theorem below identifies its actual moduli locus, which
is smooth even over singular hyperplane surfaces. It does not say
that a central matching completion remains smooth.

We next allow the support hyperplane to be singular, while retaining
a smooth proper moduli space of the sheaves themselves. For a
smooth quintic $X$, let $X^{[2]}$ be its Hilbert scheme of length-two
subschemes and put $V=H^0(X,\O_X(1))$. The tautological rank-two
bundle is denoted $\O_X(1)^{[2]}$.

\begin{theorem}\label{thm:two-point-incidence}
The evaluation kernel
\[
 0\longrightarrow\mathcal K_2\longrightarrow V\otimes\O_{X^{[2]}}
     \longrightarrow\O_X(1)^{[2]}\longrightarrow0
\]
is a rank-three vector bundle. Its projective bundle
$Y_2=\P_{X^{[2]}}(\mathcal K_2)$ is a smooth projective
locally complete locus in the scalar-rigidified moduli of stable
sheaves
\[
 F_{Z,D}=j_{D*}I_{Z/D},\qquad Z\subset D\in|H|,
 \qquad\operatorname{length}(Z)=2.
\]
It has dimension eight and charge
$(0,H,-H^2/2,-7[\mathrm{pt}]/6)$.
With orientation $K_{Y_2}$, its DT sheaf is $\CC_{Y_2}[8]$.
Writing $P_X(z)=1+z^2+204z^3+z^4+z^6$, its DT cohomology is
$H^{j+8}(Y_2,\CC)$, where
\begin{equation}\label{eq:two-point-Poincare}
 P_{Y_2}(z)=(1+z^2+z^4)
 \left\{\frac{P_X(z)^2+P_X(-z^2)}2+(z^2+z^4)P_X(z)\right\}.
\end{equation}
In the normalized periodic realization the even and odd dimensions
are $62172$ and $3672$; its Euler characteristic is $58500$.
In the resolved quintic degeneration, the logarithm of monodromy
on its odd periodic DT cohomology has square zero and rank $378$.
\end{theorem}
\begin{proof}
The very ample line bundle $\O_X(1)$ separates every length-two
subscheme, including a tangent vector at a double point. The
restriction $V\to H^0(\O_Z(1))$ is consequently surjective for
every $Z$. The universal finite flat scheme has degree two, so
its tautological pushforward is locally free of rank two; the
fiberwise surjectivity proves the vector-bundle assertion.
The projective bundle parametrizes the nonzero sections vanishing
on $Z$, modulo scalar, which are exactly the hyperplanes $D$
containing $Z$.

For completeness the Hilbert scheme here is smooth even though
Hilbert schemes of larger numbers of points on a threefold may be
singular. In smooth local coordinates $x_1,x_2,x_3$, choose an
open chart where the quotient algebra of a length-two scheme has
basis $1,x_1$. Its equations are
\[
 x_1^2=a x_1+b,\qquad x_2=c_2x_1+d_2,
 \qquad x_3=c_3x_1+d_3.
\]
These six coefficients are free parameters. Such charts cover the
nonreduced locus after changing the distinguished coordinate;
the distinct-points locus is also smooth of dimension six. The
same argument in etale coordinates proves smoothness on $X$.
Thus $Y_2$ is smooth projective of dimension $6+2=8$.

The universal incidence sequence is
\[
 0\longrightarrow I_{Z/D}\longrightarrow\O_D
     \longrightarrow\O_Z\longrightarrow0.
\]
Its last two terms are flat over $Y_2$, hence so is its kernel.
Each hyperplane $D$ on $X$ is integral by the primitive-class
argument, so this is a rank-one torsion-free sheaf on $D$ and its
pushforward is stable by Proposition~\ref{prop:primitive-support}.
Its character is \eqref{eq:ambient-defect-charge} with $k=0$, $\ell=2$.

We verify that this parameter space has the sheaf deformation
theory, not just a map to the sheaf moduli space. Evaluation and
$H^{>0}(X,\O_X(1))=0$ give
\[
 H^i(X,I_{Z/X}(1))=0\quad(i>0),\qquad
 h^0(X,I_{Z/X}(1))=3
\]
for every length-two $Z$. This is exactly the positivity condition
in \cite[Definition~2.2]{GST}. The relative Hilbert scheme of pairs
$Z\subset D$, identified with $Y_2$ by the evaluation kernel, is
therefore locally complete for the sheaves $j_{D*}I_{Z/D}$ by
\cite[Theorem~1, Proposition~2.1 and Corollary~2.3]{GST}.
These hypotheses hold uniformly, including singular hyperplanes.
At the level of points, the sheaf recovers $D$ from its annihilator;
the inclusion $I_{Z/D}\subset\O_D$ is recovered up to scalar from
its double dual. Indeed $D$ is an integral Cohen--Macaulay surface
and $Z$ is zero-dimensional, so
$\mathcal Hom_D(I_Z,\O_D)=\O_D$ and
$I_Z^{\vee\vee}=\O_D$. Thus it recovers $Z$ as well.
This also identifies the projective parameter space with the
specified classical moduli locus. The d-critical structure on its
smooth truncation is zero. The chosen orientation gives
$\CC_{Y_2}[8]$, rather than an unspecified orientation local system.

We compute its cohomology rather than only its Euler number.
The ordered length-two incidence space is
$\widetilde{X^2}=\Bl_\Delta(X\times X)$. The involution exchanging
the two factors lifts, and the quotient is $X^{[2]}$. The
exceptional divisor parametrizes tangent directions and is fixed
by the involution; this also follows from the preceding local
coordinates, where the quotient of the difference coordinate by
its sign records the doubled point. With rational coefficients,
cohomology of a finite group quotient is the invariant cohomology:
apply proper base change to the finite quotient map and take
invariants, an exact functor in characteristic zero. The blowup
formula, with center of codimension three, consequently gives
\begin{equation}\label{eq:Hilb-two-cohomology}
 H^*(X^{[2]},\mathbb Q)\simeq
    \bigl(H^*(X,\mathbb Q)^{\otimes2}\bigr)^{\mathfrak S_2}
       \oplus H^*(X,\mathbb Q)[-2]
       \oplus H^*(X,\mathbb Q)[-4].
\end{equation}
The action on the tensor square includes the graded sign.
The two exceptional summands are invariant: the involution is
$-1$ on the normal vector to the diagonal but is the identity on
its projectivization and on its hyperplane class. The trace on the
tensor square in each total degree is encoded by $P_X(-z^2)$;
averaging with the identity gives
$(P_X(z)^2+P_X(-z^2))/2$. The projective-bundle theorem then gives
\eqref{eq:two-point-Poincare}. The Chern class computation
$c(T_X)=(1+H)^5/(1+5H)$ gives $\chi(X)=-200$, and the Lefschetz
cohomology outside the middle degree gives $b_3(X)=204$.
Evaluating \eqref{eq:two-point-Poincare} at $1$ and $-1$ yields
$65844$ and $58500$. Half their sum and difference are $62172$
and $3672$. The DT shift is eight and preserves parity. The
zero critical presentation has the oriented bar realization
$(\Omega^\bullet_{Y_2}((u)),u\,d)[8]$ and normalized cyclic
realization $\CyclicChains_{\mathrm{per}}(\MF(Y_2,0))[8]$.

Finally, the blowup and projective-bundle decompositions are
relative constructions over the smooth part of the degeneration
curve; their maps are induced by pullback, exceptional Gysin maps,
and global Chern classes. They therefore commute with monodromy.
The even cohomology of a quintic is the constant four-dimensional
span of $1,H,H^2,H^3$. Formula~\eqref{eq:Hilb-two-cohomology}
shows that the odd cohomology of $X^{[2]}$ contains four copies of
$H^3(X)$ from the tensor square and two from the exceptional
summands. The $\P^2$-bundle triples this to eighteen copies.
Theorem~\ref{thm:incidence-monodromy} computes $N_X^2=0$ and
$\operatorname{rank}N_X=21$ for the resolved quintic degeneration.
Thus the odd logarithm on $Y_2$ has square zero and rank
$18\cdot21=378$.
\end{proof}

\subsubsection{The Gauss--Manin connection and its logarithmic limit}
Over the punctured degeneration curve set
$\mathbb V=\bigoplus_j R^j\pi_*\CC$, graded by cohomological degree.
The relative blowup and projective-bundle maps used in the proof
identify the cohomology local system of $Y_2$ with
\[
 \bigoplus_{a=0}^2
 \left(\operatorname{Sym}_{\mathrm{gr}}^2\mathbb V
       \oplus\mathbb V[-2]\oplus\mathbb V[-4]\right)[-2a].
\]
Here $\operatorname{Sym}_{\mathrm{gr}}^2$ denotes invariants for
the graded interchange, and $[-2a]$ raises the cohomological degree
by $2a$. Over $\CC$ the Tate twists are suppressed.
The Gauss--Manin connection on its first summand is the restriction
of $\nabla_X\otimes1+1\otimes\nabla_X$ to the graded invariants;
on each exceptional summand it is $\nabla_X$.
The projective-bundle factor has zero connection matrix in its
basis $1,\xi,\xi^2$, since $\xi$ is a global first Chern class.
Applying the DT suspension $[8]$ and periodization gives the
connection on the specified DT realization.
Thus the odd part is eighteen shifted copies of the degree-three
quintic system, not merely a vector space of the same dimension.
After the semistable logarithmic extension its nilpotent
monodromy on this odd part is eighteen copies of $N_X$.
An identification with an ordinary central matching DT complex
would still require the specialization morphism of
Problem~\ref{prob:logarithmic-categorical-comparison}.

Tensoring these sheaves by an extending $L^k$ transports the
moduli locus and its chosen orientation, and changes the charge
according to \eqref{eq:ambient-defect-charge}. The nonzero
monodromy is retained. Consequently this example requires a
logarithmic limit of its natural cohomological variation; it
cannot satisfy the pole-free extension hypothesis at $0$.
The existence of its ordinary central DT bar comparison is an
additional geometric question, not a consequence of the numerical
Euler characteristic $58500$.

\subsubsection{Central matching objects and an explicit open cohomology calculation}
\label{ssec:two-point-central}
For the central calculation put $U=X_0\setminus S$ as in
\eqref{eq:central-U-polynomial}. Its Hilbert scheme $U^{[2]}$
is smooth of dimension six. A pair of distinct points may lie on
the same component or on different components. A length-two
subscheme supported at one point includes its tangent direction;
this scheme-theoretic datum is essential in the containment
$Z\subset D$.

The map $a:X_0\to\P^4$ is not an embedding on $X_2$. Thus the
argument that $L_b$ separates all length-two schemes for $b\ne0$
does not apply unchanged at $0$. Write
\[
 B_0^{\rm sep}=\{Z\in U^{[2]}:
 V\longrightarrow H^0(Z,L_0|_Z)\ \hbox{is surjective}\}.
\]
On this open set the evaluation kernel $\mathcal K_{2,0}$ has
rank three, and we define
\begin{equation}\label{eq:two-central-moduli}
 M_{2,0}^{\circ}=\P_{B_0^{\rm sep}}(\mathcal K_{2,0}).
\end{equation}
It is a smooth eight-dimensional open locus in the classical
scalar-rigidified matching moduli space, by
Proposition~\ref{prop:point-moduli-reconstruction}.
Its universal object is \eqref{eq:incidence-perfect-object} on
$X_0$. If $Z=Z_1\amalg Z_2$, with $Z_i\subset U_i$, its two
component objects are
\begin{equation}\label{eq:two-central-triples}
 E_i=\operatorname{fib}\left(
 [L_i^{-1}\xrightarrow{s_i}\O_{X_i}]\longrightarrow\O_{Z_i}\right),
 \qquad
 E_S=[L_S^{-1}\xrightarrow{s_S}\O_S].
\end{equation}
An empty $Z_i$ means the zero quotient. Formula
\eqref{eq:two-central-triples} treats both distributions $2+0$
and $1+1$ of the length. On the component-sheaf locus the first
objects are $j_{D_i*}I_{Z_i/D_i}$, and the common derived boundary
is unchanged by these point modifications because $Z\cap S$ is
empty. It is nevertheless the full displayed boundary complex.

\begin{proposition}\label{prop:two-central-open-DT}
Let $\Gamma\subset S$ be the smooth complete-intersection curve
of degrees $(4,5)$, of genus $51$. With the canonical orientation,
the normalized potential of \eqref{eq:two-central-moduli} is zero
and its oriented bar cohomology is
\[
 \mathrm{DT}_{\mathrm{coh}}^j(M_{2,0}^{\circ})
   =H_c^{j+8}(M_{2,0}^{\circ},\CC).
\]
Its compactly supported Poincare polynomial is
\begin{equation}\label{eq:two-central-Poincare}
\begin{split}
 P_c(M_{2,0}^{\circ};z)=(1+z^2+z^4)\biggl\{&
 \frac{P_U(z)^2+P_U(-z^2)}2+(z^2+z^4)P_U(z)\\
 &{}-z^4(1+102z+z^2)\biggr\}.
\end{split}
\end{equation}
In particular its periodic even and odd dimensions are $62166$
and $3366$. This formula is for the open locus in
\eqref{eq:two-central-moduli}, not for a proper central completion.
\end{proposition}
\begin{proof}
Let $E^{\circ}=E_\Gamma\setminus S$ inside $U_2$. Its map to
$\Gamma$ is an affine-line bundle: on every exceptional $\P^1$
the strict transform of $S$ meets one point. The morphism $a$ is
an immersion away from this exceptional divisor and its only
contracted tangent directions there are the affine-line directions.
Consequently the non-surjective evaluation locus is
\[
 B_0^{\rm bad}=\operatorname{Hilb}^2(E^{\circ}/\Gamma),
\]
an affine-plane bundle over $\Gamma$.
This also holds scheme-theoretically. In local blowup coordinates
$(u,v,w)\mapsto(u,uv,w)$, a length-two quotient with basis $1,v$
has equations
\[
 v^2=\alpha v+\beta,\qquad u=cv+d,\qquad w=ev+f.
\]
The evaluation by $1,u,uv,w$ has rank one exactly when
$c=e=c\alpha+d=0$, or $c=e=d=0$. The remaining parameters are
$\alpha,\beta,f$, giving the asserted relative affine plane.
The same description of distinct points says that both lie on
one contracted fiber. No point in $U_1$ has the same image as a
point in $U_2$, since the former image is outside $(h=0)$.
These observations cover the bad locus.

The blowup of $U\times U$ along its diagonal, followed by the
involution exchanging its factors, is $U^{[2]}$. With rational
coefficients, taking invariants under this finite group is exact.
The compact-support blowup formula therefore gives
\begin{equation}\label{eq:Hilb2-compact-formula}
 P_c(U^{[2]};z)
 =\frac{P_U(z)^2+P_U(-z^2)}2+(z^2+z^4)P_U(z).
\end{equation}
The sign in $P_U(-z^2)$ is the graded interchange sign. This
argument works for the disconnected smooth threefold $U$ as well:
the cross-component part is the product $U_1\times U_2$.

To subtract the bad locus at the level of Betti numbers, rather
than merely Euler characteristics, we verify that the restriction
\begin{equation}\label{eq:bad-restriction-surjective}
 H_c^q(U^{[2]},\CC)\longrightarrow H_c^q(B_0^{\rm bad},\CC)
\end{equation}
is surjective in every degree. First the maps
$H_c^q(U_2)\to H_c^q(E^{\circ})$ are surjective.
Indeed $U_2\setminus E^{\circ}=\P^3\setminus S$ has compact
cohomology only in degrees three and six, of dimensions $21$ and
one. The localization sequence, together with
\eqref{eq:component-open-cohomology} and
$H_c^q(E^{\circ})=H^{q-2}(\Gamma)$, proves the assertion; in
degrees two and four the restriction is an isomorphism, and in
degree three its kernel has dimension $21$.

Kunneth and exactness of finite-group invariants now make
$H_c^q(\operatorname{Sym}^2 U_2)\to
H_c^q(\operatorname{Sym}^2 E^{\circ})$ surjective.
Restriction to the fiberwise symmetric square over $\Gamma$ is
surjective as well. On the ordered cover this is restriction
along the diagonal $\Gamma\to\Gamma\times\Gamma$, with the
same complex rank-two affine fibers on both sides. Under the
Thom identifications it is cup product. Every class $a$ of
$H^*(\Gamma)$ is obtained from the invariant class
$(1\otimes a+a\otimes1)/2$, so surjectivity persists after the
finite quotient. The proper Hilbert--Chow pullback to $U^{[2]}$
has the same restriction on $B_0^{\rm bad}$, since the Hilbert
scheme of two points on an affine line equals its symmetric
square. This proves \eqref{eq:bad-restriction-surjective}.

The compact-support localization sequence consequently identifies
$H_c^q(B_0^{\rm sep})$ with the kernel in
\eqref{eq:bad-restriction-surjective}. The bad locus has polynomial
$z^4 P_\Gamma(z)=z^4(1+102z+z^2)$. Subtract it from
\eqref{eq:Hilb2-compact-formula}, and apply the proper
$\P^2$-bundle formula to \eqref{eq:two-central-moduli}.
This proves \eqref{eq:two-central-Poincare}.
The moduli reconstruction and zero critical structure give the
oriented bar model \eqref{eq:point-smooth-bar} with $d=8$.
Evaluating the polynomial at $1$ and $-1$ gives $65532$ and
$58800$; their half-sum and half-difference are $62166$ and
$3366$, as claimed.
\end{proof}

\subsubsection{Connection and the remaining central strata}
On $C^*$ the generic two-point DT connection is the one described
after Theorem~\ref{thm:two-point-incidence}: it is induced from
$\nabla_X\otimes1+1\otimes\nabla_X$ on graded symmetric tensors,
from $\nabla_X$ on the two exceptional summands, and from the
horizontal projective-bundle classes. It is flat. Its odd
monodromy logarithm has rank $378$, so the natural logarithmic
connection does not admit a pole-free extension of the same system.

The central open realization just computed has different parity
ranks from the generic realization, namely $(62166,3366)$ instead
of $(62172,3672)$. It therefore cannot be the central endpoint of
the finite regular comparison for the full generic two-point
moduli space. The omitted objects include the evaluation-rank-drop
locus $B_0^{\rm bad}$, where the space of sections has dimension
four rather than three, and the perfect length-two schemes at
$S$ exemplified by \eqref{eq:node-two-lci}. Other limiting
perfect complexes may be needed for a proper completion.
Their critical structures, orientation identifications, and
specialization maps have not been computed by the open-locus
calculation. Thus this example establishes a generic invariant,
an actual central open bar invariant, and a flat logarithmic
variation; it does not establish the unrestricted ordinary
central-fiber degeneration identity.

\subsubsection{Two points on a member of $|2H|$.}
The universal length-two scheme defines the rank-two tautological
bundle $\mathcal O_{X_b}(2)^{[2]}$ on $X_b^{[2]}$. Separation of every
length-two scheme gives
\[
 0\to K_{2,2}\to H^0(X_b,\mathcal O_{X_b}(2))\otimes\mathcal O_{X_b^{[2]}}
       \to\mathcal O_{X_b}(2)^{[2]}\to0.
\]
Consequently $Y_{2,2}=\mathbb P_{X_b^{[2]}}(K_{2,2})$ has rank-$13$
evaluation kernel, fiber $\mathbb P^{12}$, and dimension $18$.
The same proof as Theorem~\ref{thm:two-point-incidence} identifies its
actual sheaf deformation functor. Its minimal potential is zero, and
\[
 \mathrm{DT}_{\mathrm{coh}}^j(Y_{2,2})
       =\bigoplus_{a=0}^{12}H^{j+18-2a}(X_b^{[2]},\mathbb C).
\]
Equivalently its ordinary Poincare polynomial is
\[
 \left(\sum_{a=0}^{12}z^{2a}\right)
 \left\{\frac{P_X(z)^2+P_X(-z^2)}2+(z^2+z^4)P_X(z)\right\}.
\]
Its normalized periodic dimensions are $(269412,15912)$.
The tensor, exceptional-Gysin, and projective-bundle maps of the
length-two calculation give its flat connection on $C^*$. The odd part
contains $13\cdot6=78$ shifted copies of $H^3(X_b)$, so its monodromy
logarithm has rank $78\cdot21=1638$ and square zero.

On the central separation open of
Proposition~\ref{prop:two-central-open-DT}, the quadratic evaluation
kernel has rank thirteen. The bad pairs on a single exceptional fiber
remain bad: every section of $L^2=a^*\mathcal O(2)$ is constant on that
fiber. Conversely quadrics separate every length-two image in
$\mathbb P^4$, so there are no additional failures on this open.
Thus the base of the central projective bundle is the same separation
open, and its fiber is now $\mathbb P^{12}$. Its oriented bar cohomology
has polynomial $\sum_{a=0}^{12}z^{2a}$ times the base polynomial in
Proposition~\ref{prop:two-central-open-DT}, with DT shift $[18]$.
Its periodic dimensions are $(269386,14586)$, not the generic ones.
Section~\ref{sec:crossing-perfect} computes a separate chart for perfect
double points at $S$ and the matrices of their restrictions. Neither
this extra chart nor the separation open by itself constitutes a proper
central completion.
\subsection{The ideal of a line: a complete generic-to-central comparison}
\label{ssec:line-example}
We next replace the point subscheme by a line $\Lambda\subset X_b$.
This is not the positive line bundle $\O_D(\Lambda)$ in
Example~\ref{ex:NL-line}: on a smooth surface containing the line,
$I_{\Lambda/D}=\O_D(-\Lambda)$. The objects considered here are
\[
 F_{\Lambda,D}=j_{D*}I_{\Lambda/D},\qquad
 \Lambda\subset D\in|H|.
\]
For this example choose the forms in the resolved degeneration
generally enough to satisfy the transversality conditions specified
below. The central limits of lines in the resolved threefold are
not all lines contained in its two components. Identifying those
limits is the main geometric step.

\subsubsection{Generic moduli, potential, and cohomology}
For a general smooth quintic the scheme of lines is reduced of
length $2875$. At every such line the normal bundle is
$\O_{\P^1}(-1)^{\oplus2}$. In particular the line is infinitesimally
rigid. The restriction of linear forms to $\Lambda$ is surjective,
and
\[
 h^0(I_{\Lambda/X_b}(1))=3,
 \qquad H^{>0}(I_{\Lambda/X_b}(1))=0.
\]
The locally complete sheaf-moduli statement of
\cite[Theorem~1 and Definition~2.2]{GST} applies; alternatively
use the reconstruction below. The locus for a fixed line is
$\P^2$, parametrizing all hyperplanes containing it, including
those whose surface section is singular. The rank-one ideal on the
integral surface is stable by Proposition~\ref{prop:primitive-support}.
With its canonical orientation the potential is zero and
\begin{equation}\label{eq:line-generic-DT}
 \mathrm{DT}_{\mathrm{coh}}^j(M_{\Lambda,b})
 =H^{j+2}(\P^2,\CC)
 =\begin{cases}\CC,&j=-2,0,2,\\0,&\text{otherwise}.\end{cases}
\end{equation}
The full line locus therefore has $\CC^{2875}$ in these three
degrees and even periodic dimension $8625$.
We now construct the central objects and verify the comparison;
the number $2875$ alone would not do so.

\subsubsection{Which curves occur on the resolved central fiber}
Let $G=\operatorname{Gr}(2,5)$, and denote its tautological
rank-two bundle by $\mathcal U$. Choose $q_4,h$ so that the
quartic line scheme
\[
 T_4=Z(q_4)\subset G
\]
is a smooth curve and the quartic surface $S=(h=q_4=0)$ contains
no lines. On $T_4$, multiplication by $h$ is an injective bundle
map, giving the line bundle
\[
 Q_4=\operatorname{Sym}^5\mathcal U^\vee/
                   h\operatorname{Sym}^4\mathcal U^\vee.
\]
On $G_h=\operatorname{Gr}(2,4)$, the lines in $(h=0)=\P^3$,
multiplication by $q_4$ is injective, giving the rank-four bundle
\[
 Q_h=\operatorname{Sym}^5\mathcal U^\vee/
                         q_4\mathcal U^\vee.
\]
Choose $q_5$ so that its two quotient sections are transverse to
zero, and so that the resulting lines meet the indicated divisors
and curve $\Gamma$ transversely at distinct points. These are the
open transversality conditions for this example. The quotient
bundles are generated by the restrictions of quintic forms, so
a general choice gives transverse zero schemes. The resulting
limiting-line conditions are also the ones in
\cite[Theorems~2.3 and~4.1]{WuLines}.

A zero of $q_5$ in $Q_4$ is a line $\Lambda_1\subset X_1$ whose
unique intersection $p$ with $S$ lies on $\Gamma$.
Its curve in the resolved union is
\begin{equation}\label{eq:line-type-A}
 C_A=\Lambda_1\cup_p e_p,
\end{equation}
where $e_p\simeq\P^1$ is the exceptional fiber of
$X_2\to\P^3$ over $p$. The second component is necessary:
it is a curve already present in the fixed component $X_2$, not
an extra expanded target. Its $L$-degree is zero, whereas
$L\cdot\Lambda_1=1$.

A zero of $q_5$ in $Q_h$ is a line $\Lambda_2\subset\P^3$ for
which $q_5|_{\Lambda_2}$ is divisible by $q_4|_{\Lambda_2}$.
Its four intersection points with $S$ all lie on $\Gamma$.
Its strict transform
\begin{equation}\label{eq:line-type-B}
 C_B=\widetilde\Lambda_2\subset X_2
\end{equation}
is disjoint from the strict transform of $S$. It meets the four
exceptional fibers at points different from their intersections
with $S$.

\begin{proposition}\label{prop:resolved-line-family}
Under these transversality conditions, after shrinking the base
about $0$, the $2875$ lines on the smooth fibers extend to a finite
etale family $T\to C$ of embedded curves in $\cX/C$.
Its central fiber consists of $1600$ curves of type
\eqref{eq:line-type-A} and $1275$ curves of type
\eqref{eq:line-type-B}. Their universal structure sheaf is
perfect on every fiber, including $X_0$. Each central curve is a
local complete intersection and is infinitesimally rigid and
unobstructed as an embedded curve in $X_0$.
\end{proposition}
\begin{proof}
First we identify the limiting curves and the etale parameters.
The equation for a line in $X'_t=(hq_4=tq_5)$ is the vanishing of
a section of $\operatorname{Sym}^5\mathcal U^\vee$ on $G$.
Near a line in $T_4$, the five normal coordinates to $T_4$ are
eliminated by the five equations in
$h\operatorname{Sym}^4\mathcal U^\vee$. On the closure of the
nonzero-fiber line locus, divide the remaining equation by $t$.
At $t=0$ it is precisely the section of $Q_4$ induced by $q_5$.
Its transverse zero gives a unique etale branch over $t$.
Near $G_h$, the same elimination uses the two normal coordinates
specified by $h$, followed by the four equations of $Q_h$ after
dividing by $t$. Transversality again gives one etale branch per
zero. These neighborhoods are disjoint because $S$ has no lines.
Properness of $G$ ensures that every limiting line is in one of
them. Thus the closure of the generic line scheme, with these
saturated equations, is finite etale after shrinking $C$.
This argument uses the actual equations, not only their zero-cycle
classes.

Pull back the resolution ideal $(h,t)$ to the universal line over
such a branch. For a type-$B$ branch, the line equation satisfies
$h=t a$ on that line, so the ideal is principal and its lift to the
resolution is the strict transform. The four points where
$q_4=q_5=0$ give the exceptional intersections described above;
their finite ratio $q_5/q_4$ is different from the boundary
section of the exceptional bundle.
For a type-$A$ branch, $h$ has one simple zero on the central
line. Blowing up $(h,t)$ on the smooth surface formed by the
family of lines is a blowup at that point. Its central fiber is
the original line together with its exceptional $\P^1$, giving
\eqref{eq:line-type-A}. The maps are embeddings in the resolved
fourfold. In the chart $t=a\tau$, $b=\tau c$ of
\eqref{eq:sc-quintic-charts}, the strict transform has equations
\[
 c=a\,g(a,a\tau),\qquad w=w(a,a\tau)
\]
for regular local functions $g,w$. These are two independent
regular equations, and the base map on the curve family is
$t=a\tau$. The other blowup chart gives the second affine chart
of the exceptional $\P^1$. This proves relative flatness and the
local complete-intersection assertion, including at the node of
$C_A$. Type $B$ is a smooth curve contained in the smooth part
of $X_0$, and is locally a complete intersection there.

We check the embedded deformation theory separately. On $C_A$,
its rank-two normal bundle restricts to
$N_{\Lambda_1/X_1}$ and $N_{e_p/X_2}$ on the two branches.
The first has $h^0=1,h^1=0$ because $T_4$ is smooth of the
expected dimension one. The second is
$\O_{\P^1}\oplus\O_{\P^1}(-1)$: motion of $p$ along $\Gamma$
gives the trivial summand and the exceptional normal direction
gives $\O(-1)$. The normal matching sequence gives
\[
 0\to H^0(N_{C_A/X_0})\to\CC\oplus\CC
 \longrightarrow T_p S\to H^1(N_{C_A/X_0})\to0.
\]
The two maps to $T_pS$ are respectively the derivative of
$\Lambda_1\mapsto\Lambda_1\cap S$ and the tangent of $\Gamma$.
Their difference is an isomorphism by the transverse zero
condition on $Q_4$. Thus both outer groups vanish.

For $C_B$, the normal bundle is the elementary transform of
$N_{\Lambda_2/\P^3}=\O(1)^{\oplus2}$ at its four transverse
intersections with $\Gamma$:
\[
 0\to N_{C_B/X_2}\to\O(1)^{\oplus2}
 \to\bigoplus_{p\in\Lambda_2\cap\Gamma}\CC_p\to0.
\]
The map on global sections is the derivative of the four
conditions imposed by $Q_h$ and is an isomorphism at a transverse
zero. Hence $H^0(N_{C_B/X_0})=H^1(N_{C_B/X_0})=0$ as well.
The ideal and structure sheaf of either curve are perfect, since
the immersion into $X_0$ is regular of codimension two.
The normal-complex criterion for deformations of a regular
immersion now verifies that the displayed relative curve family
is etale in its Hilbert moduli problem; there are no unaccounted
embedded first-order directions.

Finally the counts can be checked by Chern roots, independently
of the deformation argument. Let $x,y$ be the roots of
$\mathcal U^\vee$ and put $\sigma_1=x+y$. On $T_4$,
$c_1(Q_4)=5\sigma_1$ and
\[
 \int_G\sigma_1 c_5(\operatorname{Sym}^4\mathcal U^\vee)=320.
\]
Thus the first zero scheme has length $1600$. The second has
length
\[
 \int_{G_h}c_4\bigl(\operatorname{Sym}^5\mathcal U^\vee-
                         \mathcal U^\vee\bigr)=1275.
\]
For clarity these integrations mean the coefficient of
$x^{n-1}y^{n-2}$ in $(x-y)$ times the degree-$2(n-2)$
symmetric polynomial for $\operatorname{Gr}(2,n)$. The first
uses $c_5=\prod_{a=0}^4((4-a)x+ay)$ and $n=5$; the second
uses the degree-four term of
$\prod_{a=0}^5(1+((5-a)x+ay))/(1+x)(1+y)$ and $n=4$.
The sum is $2875$, also equal to
$\int_G c_6(\operatorname{Sym}^5\mathcal U^\vee)$.
The sections were chosen transverse, so the counts are reduced
points. They agree with the classical limiting-line calculation
\cite{WuLines}; the perfect resolved curves and their normal
calculations above supply the additional data used here.
\end{proof}

\subsubsection{The central Lagrangian intersection and its bar cohomology}
Let $\mathcal C_T\subset\cX\times_C T$ be the curve family just
constructed. The bundle
\[
 \mathcal W=\pi_{T*}(I_{\mathcal C_T}L)
\]
has rank three and commutes with base change. Indeed every curve
has $H^0(\O_{C_t})=\CC$, $H^1(\O_{C_t})=0$,
$h^0(L|_{C_t})=2$, and $H^1(L|_{C_t})=0$.
For $C_A$, these statements follow from normalization, using
degrees $(1,0)$ on its two branches; for $C_B$ they are the
usual statements on $\P^1$. Linear forms restrict surjectively
onto the two sections. The ideal sequence and the whole-fiber
cohomology table prove the asserted vanishing and base change.
Set $M_{\rm line}=\P_T(\mathcal W)$ and form its universal complex
\[
 E_{C,s}=\operatorname{Cone}(L^{-1}\xrightarrow{s}I_C).
\]
It is relatively perfect. The reconstruction proof of
Proposition~\ref{prop:point-moduli-reconstruction} applies with
$Z=C$: all it uses after the first duality calculation is
$H^1(\O_C)=0$, $H^0(\O_C)=\CC$, and perfectness. Its reconstructed
quotient is then a flat curve subscheme rather than a finite
subscheme. The embedded rigidity in
Proposition~\ref{prop:resolved-line-family} fixes that curve
infinitesimally. Thus the remaining object deformations are
exactly the projective sections in $\mathcal W$.
This proves that $M_{\rm line}$ is a smooth projective relative
classical moduli locus, of relative dimension two, including at $0$.

For type $A$, write $C_1=\Lambda_1$, $C_2=e_p$.
The relative regular equations above restrict to regular
equations on $X_i$ and on $S$, so the central matching objects are
\begin{equation}\label{eq:line-central-matching}
 E_i=\operatorname{Cone}(L_i^{-1}\xrightarrow{s_i}I_{C_i/X_i}),
 \qquad
 E_S=\operatorname{Cone}(L_S^{-1}\xrightarrow{s_S}I_{p/S}).
\end{equation}
For type $B$, replace the three ideals on the right by
$\O_{X_1},I_{C_B/X_2},\O_S$, respectively.
The matching equivalence is induced by restriction of the
universal cone. On the component-sheaf locus, these cones are
the pushforwards of the appropriate curve ideals on their
surface supports; the formulas also retain all exceptional
perfect-complex limits. The normalization triangle identifies
these matching deformation algebras with the actual global
endomorphism algebras. Consequently their classical moduli locus
is precisely the central fiber of $M_{\rm line}$.

\begin{theorem}\label{thm:line-full-comparison}
For the resolved family and the transversality conditions of
Proposition~\ref{prop:resolved-line-family}, the moduli locus
$M_{\rm line}=\P_T(\mathcal W)$ has zero normalized critical
potential on every fiber and canonical relative orientation.
Its generic and central cohomological DT groups are both
\begin{equation}\label{eq:line-full-DT}
 \mathrm{DT}_{\mathrm{coh}}^j(M_{{\rm line},b})
 =\begin{cases}\CC^{2875},&j=-2,0,2,\\0,&\text{otherwise}.\end{cases}
\end{equation}
At $0$ these are the cohomology groups of the oriented two-sided
bar model on the actual matching locus
\eqref{eq:line-central-matching}. The periodic connection is
regular and flat. On a disk on which the finite etale curve
family is labeled, its matrix is zero in the basis
$1,\xi,\xi^2$ on each of the $2875$ projective planes.
This locus satisfies the hypotheses and conclusion of
Theorem~\ref{thm:smooth-proper-comparison}.
\end{theorem}
\begin{proof}
The preceding constructions verify smoothness, projectivity,
perfectness, and the central moduli identification, rather than
assuming them. Smoothness forces the d-critical structure to be
zero. The relative canonical bundle gives the compatible
orientation and hence the shifted constant DT sheaves $\CC[2]$.
On each projective plane use the zero-section bar algebra,
orientation tensor, and mixed differential
\eqref{eq:point-smooth-bar} with $d=2$.
The intrinsic antisymmetrization and projection give the
horizontal comparison to $(\Omega^\bullet((u)),u d)[2]$.
It is the oriented DT realization of the matching locus just
identified. The three projective-plane cohomology classes give
\eqref{eq:line-full-DT}. Their Chern-class representatives are
horizontal, and the finite etale base has a locally constant
basis of its points. Thus the connection is the direct sum of
these three constant classes with the finite etale local system;
it extends regularly at $0$. On a disk that local system is
constant, giving the displayed zero matrix and the specified
horizontal generic-to-central identification.
\end{proof}

\subsubsection{The same complete line comparison for $|2H|$.}
Keep exactly the relative perfect curve family
$\mathcal C_{\mathcal T}\subset\mathcal X\times_C\mathcal T$ of
Proposition~\ref{prop:resolved-line-family}. Put
\[
 V_{\mathrm{line},2}=\pi_{\mathcal T*}
                  (I_{\mathcal C_{\mathcal T}}L^2),\qquad
 M_{\mathrm{line},2}=\mathbb P_{\mathcal T}(V_{\mathrm{line},2}).
\]
On a smooth line $L^2$ has three sections and no $H^1$. On a type-A
central curve its degrees are $(2,0)$; the normalization sequence gives
$3+1-1=3$ sections and $H^1=0$. Type B has degree two on
$\mathbb P^1$ and the same answer. Quadratic forms restrict surjectively,
so the ideal sequence and whole-fiber vanishings give
\[
 \operatorname{rk}V_{\mathrm{line},2}=15-3=12,\qquad
 R^{>0}\pi_{\mathcal T*}(I_{\mathcal C_{\mathcal T}}L^2)=0.
\]
The reconstruction proof used in Theorem~\ref{thm:line-full-comparison}
therefore applies with $L^{-2}$, giving the actual smooth classical
relative moduli locus with fiber a disjoint union of $2875$ copies of
$\mathbb P^{11}$. Its reduced critical potential is zero. Its generic
and central DT groups, with the relative canonical orientation, are
\[
 \mathrm{DT}_{\mathrm{coh}}^j(M_{\mathrm{line},2,b})
 =\begin{cases}\mathbb C^{2875},&j=-11,-9,\ldots,9,11,\\
                 0,&\text{otherwise}.
   \end{cases}
\]
In particular its normalized periodic realization has dimension $34500$
in odd parity and zero in even parity. The shift by eleven is essential.
The actual central objects in~\eqref{eq:line-central-matching} are now
$\operatorname{Cone}(L_i^{-2}\to I_{C_i/X_i})$ and their common derived
restriction. Their reduced zero-critical oriented bar is the same
bar--Koszul construction in dimension eleven. Its hyperplane-power
basis is horizontal; its connection extends across $0$ and has zero
matrix after labeling the finite etale curve family. This proves all
five comparison steps for this degree-two line locus, not just its
generic numerical Euler characteristic.
\textbf{Non-general quintics and obstructed lines.}
The word ``general'' in the reduced-line calculation is used to
ensure a reduced finite line scheme, not to discard singular
surface supports.  When $H^0(N_{\Lambda/X})$ is nonzero, the
curve parameter has additional deformation and obstruction
terms and the zero-potential proof does not apply.  The full
nonreduced curve component in Section~\ref{sec:cubic-geometry}
is an explicit such case, with a cubic potential and a complete
central comparison.  Point defects on an obstructed curve
also require the mixed cyclic operations; a product of that
cubic with an independently calculated defect invariant is
not asserted without a product identification of the
oriented deformation problem.

\subsection{A line and one point, including an embedded point}
\label{ssec:line-one-point}
The phrase ``a line and a point'' denotes a subscheme $Z$ with
\begin{equation}\label{eq:line-point-quotient}
 0\longrightarrow I_{Z/X_b}\longrightarrow I_{\Lambda/X_b}
 \longrightarrow Q\longrightarrow0,
 \qquad\operatorname{length}Q=1.
\end{equation}
If the quotient is supported off $\Lambda$, this is the disjoint
union of the line and a point. If it is supported on $\Lambda$,
it specifies an embedded point and its quotient direction.
Retaining these quotients is necessary for a proper generic
moduli space; just choosing a point of $D\setminus\Lambda$
would describe only an open subset.

\subsubsection{Generic moduli and its zero potential}
For a fixed infinitesimally rigid line in a smooth quintic, put
$B_{\Lambda,1}=\operatorname{Quot}_{X_b}(I_\Lambda,1)$.
There is an isomorphism
\[
 B_{\Lambda,1}\simeq\operatorname{Bl}_\Lambda X_b.
\]
Indeed a length-one quotient is a one-dimensional quotient of
the ideal at a point, so the Quot scheme is its projectivization.
The ideal of a smooth codimension-two curve is locally a regular
sequence, hence of linear type; its projectivization is the
blowup. Off the line the quotient is unique; on the line the
fiber is a projective line of conormal quotient directions.
The Hilbert--Chow identification and its local triviality for
a smooth curve identify this with the corresponding Hilbert
moduli locus when the curve is infinitesimally rigid
\cite[Theorems~2 and~3]{RicolfiSmooth}.

The ideal $I_\Lambda(1)$ is generated by the three linear forms
vanishing on the line. Thus the map from its three global
sections to every length-one quotient $Q(1)$ is surjective.
The kernel bundle $\mathcal V_{\Lambda,1}$ on $B_{\Lambda,1}$
has rank two, and
\[
 M_{\Lambda,1}=\P_{B_{\Lambda,1}}(\mathcal V_{\Lambda,1})
\]
parametrizes exactly $Z\subset D$, with $Z$ as in
\eqref{eq:line-point-quotient}. In particular
$H^{>0}(I_Z(1))=0$ everywhere on this locus.
Theorem~1 and Definition~2.2 of \cite{GST} make this a locally
complete moduli locus of the sheaves $j_{D*}I_{Z/D}$, not merely
a space of chosen presentations. The surface $D$ is integral,
so these sheaves are Gieseker stable. The space is smooth and
projective of dimension $3+1=4$. Its normalized potential is
zero and its canonical orientation gives the DT sheaf $\CC[4]$.

\begin{proposition}\label{prop:line-one-generic}
For a fixed rigid line on a smooth quintic,
\begin{equation}\label{eq:line-one-Poincare}
 \mathrm{DT}_{\mathrm{coh}}^j(M_{\Lambda,1})
 =H^{j+4}(M_{\Lambda,1},\CC),\qquad
 P_{M_{\Lambda,1}}(z)
 =(1+z^2)\{P_X(z)+z^2+z^4\}.
\end{equation}
Its periodic even and odd dimensions are $12$ and $408$.
In the resolved quintic family its generic Gauss--Manin system,
on each labeled branch of the line family, has two copies of
$H^3(X_b)$ in odd parity, and its odd monodromy logarithm has
rank $42$ and square zero. Its connection is flat and admits the
corresponding logarithmic extension.
\end{proposition}
\begin{proof}
The blowup formula adds $H^{*-2}(\Lambda)$ to $H^*(X_b)$.
The projective-bundle formula adds its copy shifted by two.
Both maps are given by the pullback and Gysin maps of the
relative blowup and powers of the projective hyperplane class,
so they are horizontal in families. This proves the polynomial
and the connection description, not just the dimensions.
The only odd cohomology of the quintic is its $204$-dimensional
$H^3$, whereas the line and the added hyperplane classes are
even. Thus the parity dimensions are $2(4+2)=12$ and
$2\cdot204=408$. The residue calculation is
$2\operatorname{rank}N_X=42$, by
Theorem~\ref{thm:incidence-monodromy}. Flatness is the ordinary
Gauss--Manin flatness of these smooth proper moduli over $C^*$;
the nonzero logarithm prevents a pole-free extension of that
same variation.
\end{proof}

\subsubsection{Central perfect objects and their bar cohomology}
Fix one of the resolved central curves $C_A$ or $C_B$ in
Proposition~\ref{prop:resolved-line-family}, and put
$C^{\circ}=C\cap U$. Only point quotients supported in
$U=X_0\setminus S$ are used in this calculation. Their quotient
space is
\[
 B_{C,1}^{\circ}=\operatorname{Bl}_{C^{\circ}} U.
\]
The universal length-one quotient is perfect on $X_0$, so its
kernel $I_Z$ and the cone
$E_{Z,s}=\operatorname{Cone}(L_0^{-1}\to I_Z)$ are perfect.
Restrict to the open subset where
\begin{equation}\label{eq:line-one-evaluation}
 H^0(I_C L_0)\longrightarrow H^0(Q L_0)
\end{equation}
is surjective, and form the projective bundle of its rank-two
kernel. Denote this four-dimensional locus by
$M_{C,1}^{\circ}$.

These are actual smooth moduli loci. To verify this at a central
curve with two branches, cover it by a neighborhood of its node
and neighborhoods of the finite quotient support. Near the node
the quotient vanishes, so the ideal remains a local
complete-intersection curve ideal in an Artinian deformation.
Near the quotient support the threefold and the curve branch
are smooth; the local Hilbert--Chow construction for a smooth
curve recovers the Cohen--Macaulay curve and the finite quotient
\cite{RicolfiSmooth}. These recovered curves agree away from
the finite support and hence glue. Their global embedded
normal cohomology vanishes by
Proposition~\ref{prop:resolved-line-family}, so the curve itself
has no additional infinitesimal deformation on $X_0$.
The remaining deformations are exactly the Quot deformations.
Finally the proof of
Proposition~\ref{prop:point-moduli-reconstruction} reconstructs
$I_Z$ and its section from $E_{Z,s}$: here $H^1(\O_Z)=0$
and $h^0(\O_Z)=2$, so the same degree gaps used there apply.
This verifies the moduli assertion, including embedded points
on the smooth part of $C$.

The component restrictions are explicitly
\begin{equation}\label{eq:line-one-central-triples}
 I_{Z_i}=\operatorname{fib}(I_{C_i/X_i}\to Q_i),\qquad
 E_i=\operatorname{Cone}(L_i^{-1}\to I_{Z_i}),
\end{equation}
where the point quotient is assigned to its supporting component.
An empty $C_i$ has ideal $\O_{X_i}$ and an absent quotient is zero.
The common boundary object is unchanged from
\eqref{eq:line-central-matching}, since $Q_S=0$.
The normalization triangle proves that these are the actual
matching deformation problems. Thus the central oriented DT
complex is $\CC_{M_{C,1}^{\circ}}[4]$ and its two-sided bar is
\eqref{eq:point-smooth-bar} with $d=4$.

\begin{proposition}\label{prop:line-one-central}
The compactly supported polynomials of these specified central
open matching loci are
\begin{align}
 P_c(M_{C_A,1}^{\circ};z)
  &=(1+z^2)(z^2+204z^3+3z^4+2z^6),
     \label{eq:line-one-central-A}\\
 P_c(M_{C_B,1}^{\circ};z)
  &=(1+z^2)(z^2+207z^3+2z^4+2z^6).
     \label{eq:line-one-central-B}
\end{align}
Their DT cohomology is $H_c^{j+4}$ of the respective variety.
The periodic parity dimensions are $(12,408)$ for type $A$
and $(10,414)$ for type $B$.
\end{proposition}
\begin{proof}
In type $A$, $C_A^{\circ}$ is two disjoint affine lines.
The ideal $I_{C_A}L_0$ is generated by the three linear equations
of its image line: their pullback cuts out the line together
with the exceptional fiber, with the regular equations used
in Proposition~\ref{prop:resolved-line-family}.
Thus \eqref{eq:line-one-evaluation} is surjective for every
point quotient on $U$. The compact-support blowup formula
adds $z^2P_c(C_A^{\circ};z)=2z^4$ to $P_U(z)$, and the
$\P^1$-bundle gives \eqref{eq:line-one-central-A}.

In type $B$, the inverse image of the line in $\P^3$ contains
$C_B$ and the four exceptional fibers over its intersections
with $\Gamma$. The three sections of $I_{C_B}L_0$ fail to
generate exactly in those additional directions. On
$B_{C_B,1}^{\circ}=\operatorname{Bl}_{C_B}U$, the locus where
\eqref{eq:line-one-evaluation} is zero is the union of the
strict transforms of these four affine exceptional fibers.
They are disjoint copies of $\mathbb A^1$. This also includes
the missing quotient direction above each intersection with
$C_B$. In local coordinates for the blowup of $\Gamma$, the
line equations have the form $u(v-a),w-bu$ and its strict
transform has ideal $(v-a,w-bu)$. After blowing up this ideal,
the common zero of the divided generating sections is exactly
the strict transform of $u=w=0$. This proves the claimed
scheme-theoretic description of the bad quotient locus.

The unremoved blowup has polynomial
$P_U(z)+z^2+z^4$, since $C_B=\P^1$.
Restriction from its $H_c^2$ to the four copies of
$H_c^2(\mathbb A^1)$ has rank one. One way to check the rank
is to use its Poincare-dual Gysin map: the four curves have
identical nonzero classes. Before the blowup they are fibers
over four points of the connected curve $\Gamma$, and after
the blowup each class is reduced by one exceptional fiber
class over $C_B$. The latter classes are equal because
$C_B\simeq\P^1$. The original fiber class is nonzero, as is
also seen by restricting the Thom class of
$E^{\circ}\to\Gamma$ to a fiber. Hence their common strict
transform class is nonzero. All other compact cohomology
of the four bad affine lines vanishes.

The localization exact sequence therefore changes
$h_c^2=2$ to $1$ and $h_c^3=204$ to $207$, leaving
$h_c^4=2$ and $h_c^6=2$ unchanged. Applying the projective
bundle formula on the complement gives
\eqref{eq:line-one-central-B}. The zero potential and
canonical orientation already identified above give the
DT grading and the stated parity dimensions.
\end{proof}

\subsubsection{What the comparison proves in this example}
For all $2875$ generic lines the periodic dimensions are
$(34500,1173000)$. The union of the central open loci just
computed has dimensions $(31950,1180650)$, using $1600$
and $1275$ of the two types. Thus it is not the central
endpoint of the regular comparison for the full generic
line-and-point space. It omits point quotients on $S$ and,
in type $B$, the rank-drop locus in
\eqref{eq:line-one-evaluation}. At a point of $S$, a residue
field quotient is not perfect on $X_0$, so the ordinary
kernel construction cannot fill this boundary within
$\Perf(X_0)$. On $C^*$ the connection is flat as in
Proposition~\ref{prop:line-one-generic}, with a nonzero
logarithmic residue. Additional central perfect objects and
their oriented specialization morphism are needed for a
full logarithmic comparison. In contrast, the example with
the line and \emph{no} point quotient satisfies the regular
comparison completely, by Theorem~\ref{thm:line-full-comparison}.

\subsubsection{A line and one point on a member of $|2H|$.}
For a labeled rigid line, the quotient scheme
$\operatorname{Quot}_X(I_\Lambda,1)=\operatorname{Bl}_\Lambda X$
is unchanged. The three generators of $I_\Lambda(1)$ generate its
length-one quotient. Multiplying by a linear form nonzero at the
quotient point proves surjectivity of
$H^0(I_\Lambda(2))\to H^0(Q(2))$. Since
$h^0(I_\Lambda(2))=12$ and $H^{>0}(I_\Lambda(2))=0$, the new section
kernel has rank eleven and has no higher cohomology. Thus the generic
moduli is a $\mathbb P^{10}$-bundle over $\operatorname{Bl}_\Lambda X$,
including the embedded-point quotient directions. Its dimension is
$13$, its minimal potential is zero, and
\[
 \mathrm{DT}_{\mathrm{coh}}^j(M_{\Lambda,1,2})
    =\bigoplus_{a=0}^{10}
             H^{j+13-2a}(\operatorname{Bl}_\Lambda X,\mathbb C).
\]
The blowup has even and odd Betti totals $6$ and $204$. Because the
DT shift is odd, the normalized periodic dimensions are $(2244,66)$
per labeled line. The nonconstant connection is eleven shifted copies
of the quintic connection together with the constant exceptional
classes. Its nontrivial monodromy now lies in \emph{even} DT parity,
with logarithm of rank $11\cdot21=231$ per line.

For the two central open quotient loci of
Proposition~\ref{prop:line-one-central}, the same degree-two
surjectivity holds precisely away from the indicated bad curves.
The bad type-B curves remain bad because the pullback quadrics vanish
at the point of the contracted line. Therefore the projective fibers
are $\mathbb P^{10}$ over the very same two bases. The oriented bar
cohomology is their compact cohomology, multiplied by
$1+z^2+\cdots+z^{20}$ and shifted by thirteen. The periodic dimensions
per curve are $(2244,66)$ for type A and $(2277,55)$ for type B.
These are statements about the two actual open matching loci, not
about the missing full central boundary. The perfect replacement at
$S$ in Section~\ref{sec:crossing-perfect} illustrates why extra
cohomology and quotient directions must be retained there.
\subsection{A line and two points: the first singular curve-defect example}
\label{ssec:line-two-points}
Fix a rigid line $\Lambda$ on a smooth quintic $X$. The compact
parameter for two point defects is
\[
 B_{\Lambda,2}=\operatorname{Quot}_X(I_\Lambda,2).
\]
A point is a quotient $I_\Lambda\twoheadrightarrow Q$ of length
two, with kernel $I_Z$. It includes two isolated points, double
points off the line, and embedded point structures on the line.
As in \cite{RicolfiSmooth}, rigidity of the smooth line identifies
this with the corresponding Hilbert--Chow component of
one-dimensional subschemes. The actual incidence parameter is
\begin{equation}\label{eq:line-two-incidence}
 M_{\Lambda,2}=
 \{(I_\Lambda\twoheadrightarrow Q,[s]):
 [s]\in\P H^0(I_\Lambda(1)),\ s\mapsto0\text{ in }H^0(Q(1))\}.
\end{equation}
Scheme-theoretically it is the zero locus of the tautological
map from $\O_{\P^2}(-1)$ to the rank-two tautological quotient
bundle on $B_{\Lambda,2}\times\P^2$.
The associated sheaf is $j_{D*}I_{Z/D}$, with $D=(s=0)$.
It is stable because $D$ is an integral hyperplane section.
Stability alone does not make \eqref{eq:line-two-incidence}
a smooth space or a locally complete sheaf-moduli locus everywhere.

\subsubsection{The section-space rank really jumps}
From $0\to I_Z(1)\to I_\Lambda(1)\to Q(1)\to0$ one obtains
\begin{equation}\label{eq:line-two-rank-sequence}
 0\to H^0(I_Z(1))\to\CC^3\xrightarrow{\mathrm{ev}_Q}\CC^2
 \to H^1(I_Z(1))\to0.
\end{equation}
On the rank-two open set the fiber of
\eqref{eq:line-two-incidence} is $\P^0$, and the vanishing
hypothesis of \cite[Definition~2.2]{GST} holds. On the rank-one
locus the fiber is $\P^1$ and $H^1(I_Z(1))=\CC$.
This latter locus is nonempty even for two distinct isolated
points: choose $p,q\in X\setminus\Lambda$ with the same image
under the linear projection from $\Lambda$ to $\P^2$.
The evaluations of the three linear forms defining $\Lambda$
are then proportional. Such points can be taken on the residual
plane quartic in a plane containing $\Lambda$.
Thus positivity at degree one has not been established for the
whole line-and-two-point component; in fact it fails on the
specified locus. The locus where it holds is open and its
locally complete sheaf-moduli interpretation follows from the
local deformation argument of \cite{GST}.

\subsubsection{A concrete nonzero critical potential}
There is a singular point even on that positivity open set.
In smooth local coordinates write $I_\Lambda=(x,y)$ in
$R=\CC[x,y,z]$ and let
\begin{equation}\label{eq:line-two-singular-ideal}
 J=(x^2,xy,y^2,xz,yz)=\mathfrak m I_\Lambda,
 \qquad Q=I_\Lambda/J\simeq\CC^2.
\end{equation}
The tangent space of the Quot scheme is
$\operatorname{Hom}_R(J,Q)$, of dimension ten, because $J$
has five minimal generators and $Q$ is killed by $\mathfrak m$.
Its local dimension is six. The three global generators of
$I_\Lambda(1)$ surject onto $Q$, so the incidence space agrees
locally with this singular Quot scheme. In particular the actual
sheaf moduli here cannot have a zero potential on a smooth
six-dimensional chart.

For explicit local calculations one has the following critical
presentation of the affine Quot model
\cite[Theorem~2.2 and Proposition~2.10]{DavisonRicolfi}.
For $n$ point defects take three endomorphisms $A,B,C$ of
$\CC^n$, two vectors $v_1,v_2$, and a covector $w$, impose the
open condition that $v_1,v_2$ generate under the endomorphisms,
and divide by $\operatorname{GL}_n$. This smooth space has
dimension $2n^2+3n$. Its ordinary potential, with one consistent
choice of arrow signs, is
\begin{equation}\label{eq:line-Quot-potential}
 W=\operatorname{Tr}(C[A,B])+w(Av_2-Bv_1).
\end{equation}
Its critical locus is the affine Quot scheme of the line ideal.
This is a specified local critical realization, not an inference
from an Euler characteristic.

\begin{proposition}\label{prop:line-two-determinant-potential}
For $n=2$, on the chart where $v_1,v_2$ form a basis, the
potential \eqref{eq:line-Quot-potential} is, after removal of two
nondegenerate quadratic pairs, the determinant of a $3\times3$
matrix, independent of one additional smooth coordinate.
Its critical scheme is
$\mathbb A^1\times\{M\in\operatorname{Mat}_3:\operatorname{rank}M\le1\}$.
At the point \eqref{eq:line-two-singular-ideal}, the perverse
vanishing-cycle complex of this determinant realization has
one-dimensional stalk cohomology in degrees $-6,-4,-1$ and
zero stalk cohomology in other degrees.
\end{proposition}
\begin{proof}
Use the frame to set $v_1=e_1$, $v_2=e_2$. The equation paired
with $w$ is $Ae_2-Be_1=0$. Its two entries and the two entries
of $w$ form nondegenerate quadratic pairs. Eliminating them
sets
\[
 A=\begin{pmatrix}a&b\\c&d\end{pmatrix},\quad
 B=\begin{pmatrix}b&e\\d&f\end{pmatrix},\quad
 C=\begin{pmatrix}g&h\\i&j\end{pmatrix}.
\]
A direct multiplication gives
\[
 \operatorname{Tr}(C[A,B])
 =(g-j)(bd-ce)+h\{c(b-f)+d(d-a)\}
                 +i\{e(a-d)+bf-b^2\}.
\]
With $\gamma=g-j$, $u=a-d$, $v=f-b$, this is
\begin{equation}\label{eq:determinant-local-potential}
 f=\det\begin{pmatrix}
 \gamma&h&i\\ u&b&c\\-v&e&d
 \end{pmatrix}.
\end{equation}
The remaining coordinate $g+j$ is free. The gradient consists
of the $2\times2$ minors, proving the critical-scheme assertion
and its dimension six. At the zero matrix all gradient equations
have order two, so the tangent dimension is ten, in agreement
with the direct Quot calculation.

The determinant is homogeneous. Its Milnor fiber at the vertex
has the homotopy type of the affine fiber $\det M=1$, namely
$\operatorname{SL}_3(\CC)$, which retracts by polar decomposition
to $\operatorname{SU}(3)$. The latter has cohomology the exterior
algebra on generators in degrees three and five, as follows
also from its fibration over $S^5$ with fiber $S^3$.
The reduced cohomology is therefore in degrees $3,5,8$.
The free coordinate does not change the Milnor fiber, and the
smooth chart has dimension ten after eliminating the quadratic
pairs. Perverse normalization shifts these degrees by $-9$,
giving $-6,-4,-1$. Nondegenerate quadratic stabilization preserves
this normalized statement with the chosen orientation.
\end{proof}

The corresponding explicit local bar is the two-sided bar of
the zero section and $df$ in $T^*\mathbb A^{10}$. Its Koszul
retract and periodic differential are
\[
 K_f=\CC[x_1,\ldots,x_{10}][\eta_1,\ldots,\eta_{10}],
 \qquad D=-\sum_a f_a\partial_{\eta_a}
                      +u\sum_a\partial_{x_a}\partial_{\eta_a}.
\]
The variable corresponding to the free coordinate has $f_a=0$.
The maps and the bar connection are those of Section~5.
For a specified family $f(x,t)$ the operator is
$\partial_t-u^{-1}\partial_t f$, with its overlap homotopies;
it commutes with $D$ and is flat on cohomology.

Theorem~\ref{thm:determinant-intrinsic} below identifies this
particular determinant d-critical class with the intrinsic class of the
sheaf-moduli germ. It does so by calculating the kernel which defines
that d-critical class, not by appealing to a numerical DT/PT identity.
Theorem~\ref{thm:determinant-relative} gives the corresponding actual
central local comparison. A global invariant still requires the other
Quot charts and their oriented descent. The rank jump in
\eqref{eq:line-two-rank-sequence} also persists for $|H|$; in $|2H|$
it is removed by Lemma~\ref{lem:quadric-quot-evaluation}, without
removing the determinant singularity.

\subsubsection{Central restrictions and the comparison to be constructed}
For a relative curve ideal $I_C$, an actual perfect quotient
$I_C\to Q$ of length two gives $I_Z=\operatorname{fib}(I_C\to Q)$
and $E=\operatorname{Cone}(L^{-1}\to I_Z)$.
On each central component and on the common surface one must
use the two triangles
\begin{equation}\label{eq:line-two-central-triangles}
 I_{Z,i}\to I_{C,i}\to Q_i\xrightarrow{+1},
 \qquad L_i^{-1}\to I_{Z,i}\to E_i\xrightarrow{+1},
\end{equation}
and the same formulas with subscript $S$, all by derived
restriction. They define the actual central matching triple
when the indicated relative perfect objects exist.
On its derived charts the structure algebra is
$A_1\otimes_A^L A_2$; the oriented critical model and BV
operator must then be computed from that deformation problem.
The determinant calculation above is not automatically its
potential at every central point.

The global generic calculation for the full singular Quot
component and its degree-two surface incidence family is given in
Section~\ref{sec:global-two-defect}.  It is not a smooth
parameter-space calculation: the extensions in the determinant DT
sheaf are computed there.  The corresponding full ordinary
central perfect-complex completion is not identified by the local
triangles alone.
Theorem~\ref{thm:determinant-intrinsic} proves the intrinsic comparison
at the determinant germ, and Theorem~\ref{thm:determinant-relative}
identifies its actual local central matching problem. The additional
global calculations are the oriented descent across the other Quot
and evaluation-rank strata, and a central perfect-complex completion
with its specialization map. The flat local chain operator does not
replace these global data.

\subsubsection{A line and two points on a member of $|2H|$.}
In degree two the evaluation behavior improves, but the critical
structure does not become zero. Lemma~\ref{lem:quadric-quot-evaluation}
proves that for every length-two quotient $I_\Lambda\twoheadrightarrow Q$
(including nonreduced and embedded quotients)
\[
 H^0(I_\Lambda(2))\longrightarrow H^0(Q(2))
 \text{ is onto},\qquad
 H^{>0}(I_Z(2))=0,\qquad h^0(I_Z(2))=10.
\]
Consequently the complete degree-two incidence scheme is
$\mathbb P_{B_{\Lambda,2}}(V_{10})$, a projective bundle with fiber
$\mathbb P^9$ over the \emph{singular} scheme
$B_{\Lambda,2}=\operatorname{Quot}_X(I_\Lambda,2)$.
The positivity hypothesis of~\cite{GST} holds on this whole locus,
so it has the sheaf deformation theory. A projective bundle over
$B_{\Lambda,2}$ is not a smooth moduli space.

At the embedded quotient~\eqref{eq:line-two-singular-ideal}, its
classical analytic germ is
\[
 \mathbb A^{10}\times
     \{M\in\operatorname{Mat}_3:\operatorname{rank}M\leq1\},
\]
where nine of the ten smooth coordinates are projective-bundle
coordinates. Theorem~\ref{thm:determinant-intrinsic} identifies the
intrinsic local d-critical potential with $\det M$, not $0$.
Its normalized local DT stalk has dimension one in each of degrees
$-15,-13,-10$. The corresponding degree-one stalk has degrees
$-6,-4,-1$. These are local vanishing-cycle groups, not global Betti
totals of the proper incidence scheme.

Theorem~\ref{thm:determinant-relative} identifies this nonzero chart
with the actual central matching deformation problem for defects
supported at a relative smooth point of a limiting perfect curve.
It gives the two-sided bar, the twisted de Rham comparison, and the
flat connection explicitly for both degrees. Section~\ref{sec:cubic-geometry}
gives a separate proper, genuinely singular comparison whose limiting
curve itself meets $S$ nontransversely in the moduli problem. A global
DT group for all of $B_{\Lambda,2}$ would additionally require the
oriented descent and the other Quot charts; the local stalk calculation
is not used to replace those data.
\subsection{Arbitrary curve and point subschemes and higher-degree divisors}
\label{ssec:general-subschemes}
Let $X$ now be a fixed smooth quintic, put $H=c_1(\O_X(1))$,
and let $Z\subset X$ have dimension at most one. We allow isolated
and embedded points. Write
\[
 \gamma=[Z]_{\mathrm{one\text{-}dimensional}},\qquad
 e=\chi(\O_Z).
\]
The letter $e$ is not the degree $n$ of $\ch_3$ of the surface
sheaf. For $D\in|mH|$ containing $Z$ scheme-theoretically, the
sequence defining the sheaf is
\begin{equation}\label{eq:general-Z-sequence}
 0\to\O_X(-mH)\xrightarrow{s}I_{Z/X}
       \to F_{Z,D}:=j_{D*}I_{Z/D}\to0.
\end{equation}
The two ideals in this expression belong to different schemes.
The condition $s\in H^0(I_{Z/X}(mH))$ is precisely the
scheme-theoretic containment $Z\subset D$. In particular it
includes the equations of every embedded point; containment of
the sets of closed points would not suffice.

\textbf{The charge and the generic deformation problem.}
On a Calabi--Yau threefold, Riemann--Roch for $\O_Z$ gives
$\ch(\O_Z)=(0,0,\gamma,e[\mathrm{pt}])$. Additivity in
\eqref{eq:general-Z-sequence} therefore gives
\begin{equation}\label{eq:general-Z-character}
 \ch(F_{Z,D})=
 \left(0,mH,-\frac{m^2}{2}H^2-\gamma,
             \left(\frac{5m^3}{6}-e\right)[\mathrm{pt}]\right).
\end{equation}
In degree one, a line with $r$ point defects has $e=1+r$ and
hence
\[
 \beta=-\tfrac12H^2-[\Lambda],\quad
 n=-\tfrac16-r,\quad
 \chi(F_{Z,D}(a))=\tfrac52a^2-\tfrac72a+4-r.
\]
This identifies the charges of the three line examples, including
embedded defects. It does not assert that every sheaf of that
charge is an ideal sheaf on its support.

Fix a bounded Hilbert moduli problem $B=I_e(X,\gamma)$ with its
universal ideal $\mathscr I$ and projection $q:X\times B\to B$.
The positivity condition relevant to the result of Gholampour,
Sheshmani and Thomas is exactly
\begin{equation}\label{eq:GST-positivity-exact}
 H^i(X,I_{Z/X}(mH))=0\quad(i>0)\quad\hbox{for every }[Z]\in B.
\end{equation}
This is \cite[Definition~2.2]{GST}, not the assertion that $m$
is always large, and not a stability condition on $F_{Z,D}$.
For a fixed bounded Hilbert problem, uniform Serre vanishing
ensures it for sufficiently large $m$. When it holds, the
incidence moduli scheme is the projective bundle
\begin{equation}\label{eq:general-incidence-projective}
 M_{m,\gamma,e}=\P_B(\mathscr V),\qquad
 \mathscr V=q_*\mathscr I(mH),\qquad
 c_m:=\operatorname{rk}\mathscr V
       =\frac{5m^3+25m}{6}-mH\gamma-e,
\end{equation}
provided $c_m>0$. The quotient convention is not used here:
$\P_B(\mathscr V)$ parametrizes one-dimensional subspaces.
The universal inclusion gives \eqref{eq:general-Z-sequence} in
families. Conversely a section defining $Z\subset D$ gives that
subspace, proving the projective-bundle identification, as in
\cite[Proposition~2.1 and Lemma~3.2]{GST}. Its locally complete
sheaf deformation theory is \cite[Theorem~1 and Corollary~2.3]{GST}.
These assertions include singular Hilbert moduli $B$; a projective
bundle over a singular scheme is not thereby smooth.

\textbf{Stability is a separate question in higher degree.}
For $m=1$, Proposition~\ref{prop:primitive-support} proves
Gieseker stability for every ideal $I_{Z/D}$ of rank one on its
integral support, whether or not $Z$ is pure. For $m>1$, the
support can split inside the fixed smooth quintic and the same
argument no longer applies. There are actual destabilizing
subsheaves. For example, take distinct smooth transverse
$D_1,D_2\in|H|$ and a smooth curve
$C=D_2\cap Q$ cut by a quadric, with $C$ not a component of
$D_1\cap D_2$. Put $D=D_1+D_2$ and $F=I_{C/D}$ on $X$.
There is a subsheaf
\[
 A=\O_{D_1}(-H)\hookrightarrow F.
\]
It is the kernel of $\O_D\to\O_{D_2}$ and vanishes on $C$.
Using $[C]=2H^2$ and \eqref{eq:general-Z-character}, the leading
terms of the two Hilbert polynomials with respect to $H$ are
\[
 P_F(a)=5a^2-20a+O(1),\qquad
 P_A(a)=\tfrac52a^2-\tfrac{15}{2}a+O(1).
\]
Their reduced linear coefficients are $-4$ and $-3$, respectively.
Thus $A$ destabilizes $F$. This example occurs on a smooth
quintic; it is independent of degeneration of the threefold.
The Joyce--Song construction in \cite{GST} instead uses the
stable rank-one ideal $I_{Z/X}$ and its nonzero section. Its
complex is quasi-isomorphic to $F$, but the stability condition
on the pair is not Gieseker stability of $F$.

The stability assertion for the particular sheaves
$\O_D(kH)$ in Proposition~\ref{prop:divisor-stability-all-degrees}
has a different proof, using the explicit sub-divisor slopes.
It must not be extended to $I_{Z/D}$ with an arbitrary curve
part by analogy. For a finite point defect, those slopes have
the same linear coefficients as for $\O_D(kH)$, so that argument
still proves stability in Picard rank one: sub-divisor supports
have strictly smaller reduced linear coefficient, and a proper
subsheaf with full generic support has a quotient of dimension
at most one and hence strictly smaller reduced polynomial.
The preceding curve example changes that linear coefficient.

\textbf{A cohomological projective-bundle statement and its hypothesis.}
The numerical formula in \cite[Corollary~3.3]{GST} follows from
a relation of Behrend functions. It does not specify the
vanishing-cycle sheaf or its monodromy. The following is the
precise additional hypothesis under which its cohomological
analogue is available.

\begin{proposition}\label{prop:cohomological-projective-bundle}
Let $p:M=\P_B(\mathscr V)\to B$, with $\operatorname{rk}\mathscr V=c>0$.
Suppose the two spaces have specified oriented d-critical
realizations and their DT sheaves satisfy an isomorphism
\begin{equation}\label{eq:oriented-smooth-pullback}
 \mathscr P_M\simeq p^*\mathscr P_B[c-1].
\end{equation}
Then
\begin{equation}\label{eq:cohomological-projective-formula}
 \mathrm{DT}_{\mathrm{coh}}^j(M)
 \simeq\bigoplus_{a=0}^{c-1}
       \mathrm{DT}_{\mathrm{coh}}^{j+c-1-2a}(B).
\end{equation}
For a family with the same compatible oriented realizations,
these maps are horizontal: its Gauss--Manin system is the direct
sum of these shifted copies of the system on $B$.
The hypothesis holds for smooth $B,M$ with their zero critical
structures and canonical orientations. More generally it holds
when compatible local potentials for $M$ are the smooth
pullbacks of those for $B$, with the orientation comparisons
in \eqref{eq:oriented-smooth-pullback}.
\end{proposition}
\begin{proof}
Let $\xi=c_1(\O_M(1))$. The projective-bundle maps defined by
$1,\xi,\ldots,\xi^{c-1}$ give
$Rp_*\CC_M\simeq\bigoplus_{a=0}^{c-1}\CC_B[-2a]$.
Use the projection formula with \eqref{eq:oriented-smooth-pullback}
and then compactly supported pushforward to a point; $p$ is
proper. The shifts give exactly
\eqref{eq:cohomological-projective-formula}.
In a relative family the same powers are Chern classes of a
line bundle on the total space and are horizontal. Hence the
comparison commutes with the Gauss--Manin operators.
For smooth zero-critical spaces, both sides of
\eqref{eq:oriented-smooth-pullback} are the constant sheaf shifted
by $\dim M=\dim B+c-1$. For compatible smooth pullback charts,
the smooth base-change property of vanishing cycles, with the
perverse shift and the specified orientations, gives that
isomorphism locally and its prescribed overlap maps glue it.
\end{proof}

For singular $B$, neither smoothness of the map $p$ as an
ordinary morphism nor equality of symmetric obstruction-theory
ranks proves \eqref{eq:oriented-smooth-pullback}. That oriented
comparison is an additional theorem to establish for the
particular sheaf deformation problem. In particular, replacing
its DT sheaf by a constant sheaf would erase the nonzero
potential in Proposition~\ref{prop:line-two-determinant-potential}.
The construction of an arbitrary local potential uses the
cyclic brackets of Proposition~\ref{prop:general-cyclic-reduction}, not only the Hilbert
polynomial of $Z$.

\textbf{Explicit higher-degree point examples.}
Put $r_m=h^0(X,\O_X(mH))=(5m^3+25m)/6$ for $m>0$.
For length one and length two, $\O_X(mH)$ separates the relevant
subscheme, and \eqref{eq:GST-positivity-exact} holds for every
$m\ge1$. Thus the generic moduli are projective bundles of
ranks $r_m-1$ and $r_m-2$ over $X$ and $X^{[2]}$, respectively.
They are smooth of dimensions $r_m+1$ and $r_m+3$, and their
potentials are zero with canonical orientations. In the notation
$P_X(z)=1+z^2+204z^3+z^4+z^6$, their Poincare polynomials are
\begin{align}
 P_{1,m}(z)&=(1+z^2+\cdots+z^{2(r_m-2)})P_X(z),\label{eq:higher-point-one}\\
 P_{2,m}(z)&=(1+z^2+\cdots+z^{2(r_m-3)})
 \left\{\frac{P_X(z)^2+P_X(-z^2)}2+(z^2+z^4)P_X(z)\right\}.
 \label{eq:higher-point-two}
\end{align}
The coefficient of $z^{j+r_m+1}$ in the first, or
$z^{j+r_m+3}$ in the second, is the corresponding generic DT
Betti number. These follow from the same reconstruction and
projective-bundle proofs as the degree-one cases. The connection
is the quintic Gauss--Manin connection on each $X$-summand,
and its tensor-induced connection on the Hilbert-square
summands, with horizontal projective hyperplane powers.
It is flat on the punctured degeneration curve. The nonzero
quintic monodromy is retained; increasing the divisor degree
does not remove the need for a logarithmic limit or resolve
perfectness of point collisions at $S$.

\textbf{Central objects for a general subscheme.}
Suppose an actual family of perfect ideal objects $I_{Z_T}$ and
an actual morphism $L_T^{-m}\to I_{Z_T}$ have been constructed
across $0$. Its cone is perfect and its central restrictions
are, without suppressing any degrees,
\begin{equation}\label{eq:general-central-Z-cones}
 E_i=\operatorname{Cone}(L_i^{-m}\to Lj_i^*I_{Z_0}),\qquad
 E_S=\operatorname{Cone}(L_S^{-m}\to Lj_S^*I_{Z_0}).
\end{equation}
They define a point of the Lagrangian intersection by the
canonical equality of the iterated restrictions. Its deformation
algebra is the homotopy fiber of the two restriction maps on
endomorphism algebras. Its oriented local critical model gives
the DT bar differential, and Section~5 gives the horizontal
operator once a compatible family of such ordinary potentials
has been identified. These assertions are functorial for perfect
objects of arbitrary amplitude.

To obtain the full central DT group for a specified sheaf
problem, one must also choose its bounded central moduli locus
and construct its orientation and specialization map. Equations
\eqref{eq:general-Z-character} and
\eqref{eq:general-central-Z-cones} alone do not choose such a
completion. The nonperfect point limit
\eqref{eq:node-ideal-obstruction} gives a concrete obstruction
to taking the ordinary relative Hilbert scheme as that
completion without modification. Theorem~\ref{thm:line-full-comparison}
overcomes this obstruction for the line locus by constructing
its actual perfect central curves. For more general $Z$, the
criterion in Proposition~\ref{prop:specialization-cone-test}
states exactly what additional vanishing is required for a
full specialization isomorphism.

\subsection{General pure surface sheaves and the scope of deformation invariance}
\label{ssec:general-pure-sheaves-scope}
On a smooth quintic, a pure sheaf of primitive first Chern class
$H$ is of the form $j_*G$ on an integral hyperplane section.
On a smooth such section its description is
$G=M\otimes I_{W/D}$, with $M$ an arbitrary line bundle and
$W$ finite. In particular, $M$ need not be a restriction of
$\O_X(kH)$ and need not be the inverse of an effective curve.
The ideal-sheaf examples above therefore do not classify all
pure two-dimensional sheaves. On a singular integral surface
its rank-one reflexive hull can also fail to be a line bundle.
These are genuine additional moduli problems, not changes of
notation in a section complex.

For any specified oriented finite-type locus of such objects,
the formal potential is computed from its actual cyclic
endomorphism algebra as in Proposition~\ref{prop:general-cyclic-reduction}. A family of
compatible ordinary charts has the chain connection of
Theorem~\ref{thm:v5-Cech-GM}, and the local central bar comparison of
Theorem~\ref{thm:v5-periodic-bar}. These constructions are general in the perfect
object. Finiteness and the identification of the same central
geometric complex are separate hypotheses of
Theorem~\ref{thm:finite-cyclic}. On a smooth proper relative
moduli locus with its constructed central matching
identification, Theorem~\ref{thm:smooth-proper-comparison}
verifies them and proves the corresponding cohomological
specialization. With logarithmic variation one instead uses the
nearby-cycle triangle of
Proposition~\ref{prop:specialization-cone-test}; its cone must be
calculated rather than set to zero by flatness.

The examples now have the following precise division. For
section complexes in every degree and for the rigid-line
families of Theorem~\ref{thm:line-full-comparison}, the generic
and actual full central groups, the bar comparisons, and the
regular flat connections are calculated. For one point, two
points, and a line with one point, the generic DT groups and
the specified nonproper central matching opens are calculated;
the logarithmic systems retain nontrivial monodromy. These
opens are not declared to be a proper completion. For a line
with two point defects the intrinsic nonzero determinant structure
and its local central comparison are proved in
Section~\ref{sec:singular-comparison}. Section~\ref{sec:crossing-perfect}
constructs additional perfect point limits at $S$, and
Section~\ref{sec:cubic-geometry} gives a proper nonreduced example with
its complete nonzero-potential comparison. These results do not yet
give a global central comparison for the entire line-and-two-point
moduli space. The
higher-degree projective-bundle reduction is a cohomological
formula only with the explicit orientation and critical-chart
hypothesis of Proposition~\ref{prop:cohomological-projective-bundle}.

Thus flatness of the constructed connection is proved in its
stated setting, and the full deformation-invariance theorem is
verified on the specified complete examples. An unrestricted
pole-free identification for all surface sheaves would
contradict the nontrivial monodromy already computed in the
one-point example. The extension to the other examples must
retain that monodromy and prove the actual central
specialization morphism; it cannot be obtained by treating a
numerical degeneration identity as an identity of DT complexes.

\textbf{The invariant on all strata of an arbitrary surface-sheaf problem.}
For a specified global oriented moduli space, the local potentials of
Proposition~\ref{prop:general-cyclic-reduction} and the critical
bar maps apply at its singular as well as smooth points.  Global
cohomology must use the restricted DT complexes and connecting maps
of Proposition~\ref{prop:boundary-localization}.  On the singular
component computed in Section~\ref{sec:global-two-defect}, replacing
the perverse complex by a constant sheaf or by its intersection
complex gives the wrong answer; the proof calculates exactly the
missing extensions.  On an as yet unidentified central image,
there is no canonical global cohomology group to which a formal
local calculation can be assigned.  This is why the global ordinary
central comparison for general charges is stated as
Problem~\ref{prob:logarithmic-categorical-comparison}, rather than
as a consequence of a perfect obstruction theory alone.

\subsection{The component surface and the boundary restriction}
The resolved degeneration of Section~\ref{ssec:quintic-construction} has component
$X_2=\Bl_\Gamma\P^3$, where $\Gamma$ is a complete intersection of
degrees $4$ and $5$. A generic plane meets $\Gamma$ transversely in
$20$ distinct points. Its strict transform is a smooth surface
$B=\Bl_{20}\P^2$. Write $h$ for the plane class and $e_a$ for the
exceptional curve classes. The calculations needed for the categorical
boundary problem are as follows.

\begin{proposition}\label{prop:component-surface-charge}
For the inclusion $j_B:B\hookrightarrow X_2$, one has
\[
 N_{B/X_2}=\O_B(h),\qquad K_B=-3h+\sum_{a=1}^{20}e_a,
 \qquad C_B=B\cap S\sim4h-\sum_{a=1}^{20}e_a.
\]
For a line bundle $M$ on $B$ with $\alpha=c_1(M)$ and a finite
subscheme $Z\subset B$ of length $\ell$, the sheaf
$F_B=j_{B*}(M\otimes I_Z)$ has
\begin{align}
 \mathrm{ch}_1(F_B)&=H_2,\notag\\
 \mathrm{ch}_2(F_B)&=j_{B*}(\alpha-\tfrac12h),\label{eq:component-charge}\\
 \int_{X_2}\mathrm{ch}_3(F_B)
  &=\tfrac12\alpha^2-\tfrac12\alpha h+\tfrac16-\ell.\notag
\end{align}
If $Z\cap C_B=\varnothing$, its derived restriction to $S$ is the
pushforward of $M|_{C_B}$. In particular a factorization through the
boundary family $\O_{C_B}(k)$ requires the isomorphism
$M|_{C_B}\simeq\O_{C_B}(k)$, not just a numerical positivity condition.
\end{proposition}
\begin{proof}
The divisor $B$ is a member of $|H_2|$, so its normal line bundle is
$\O_B(h)$. The usual blowup formula gives $K_B$, and the restriction
of $S\sim4H_2-E$ gives $C_B$. For its character use GRR as in
Proposition~\ref{prop:general-divisor-ch}, now with $h^2=1$.
The inverse Todd factor is again
$1-h/2+h^2/6$, giving \eqref{eq:component-charge}.
If the finite defect avoids the boundary, $I_Z$ is the structure
sheaf near $C_B$. The two smooth divisors meet transversely there,
so derived restriction reduces to the ordinary restriction of $M$.
The last assertion is the resulting identification of the boundary
object, not an inference from its Chern character.
\end{proof}

The boundary curve is a plane quartic before blowing up, which
accounts for the coefficient $4$ in its strict-transform class.
The coefficient $1/6$ is the degree-two coefficient of the inverse
Todd factor $(1-e^{-h})/h$. Both enter the charge and boundary
computations independently of the choice of a moduli compactification.

Even the vanishing of the first Tor sheaf does not force the ordinary
boundary restriction to be pure. This distinction matters when one
tries to replace $\Perf(S)$ by a stable-sheaf moduli space.

\begin{example}\label{ex:boundary-defect}
On a smooth surface with local coordinates $x,y$, let the boundary
curve be $x=0$ and let $I_p=(x,y)$. Multiplication by $x$ on $I_p$
is injective, so its derived restriction has no degree-$-1$
cohomology. Nevertheless
\[
 I_p\otimes\CC[x,y]/(x)=I_p/xI_p
\]
has a nonzero element $[x]$ annihilated by $y$. In a neighborhood
of $p$ on the boundary there is an exact sequence
\[
 0\longrightarrow \CC_p\longrightarrow I_p|_{(x=0)}
     \longrightarrow I_{p/(x=0)}\longrightarrow0.
\]
The class $[y]$ maps to the generator of the last ideal. Thus the
restriction has a zero-dimensional torsion subsheaf. Applying this
local model to a smooth support divisor transverse to $S$ gives a
pure surface sheaf normal to $S$ whose restriction is not a pure
curve sheaf. The full target $\Perf(S)$ retains that restriction;
a compactified-Jacobian target would require a further restriction
of the moduli problem.
\end{example}

\subsection{A categorical realization of the boundary coefficient}
Assume now $\operatorname{Pic}(S)=\mathbb Z h_S$, with $S$ a smooth
quartic K3 surface. Every member of $|h_S|=\P^3$ is integral: a
nontrivial decomposition would express the primitive ample class as
a sum of positive integral multiples of itself. Write
$\mathcal J_k=M_S(0,h_S,4k-2)$ for the moduli space of stable sheaves
with this Mukai vector and polarization $h_S$. We use the
symplectic form defined by the trace pairing \cite{mukai1984symplectic}.

\begin{theorem}\label{thm:boundary-P3}
The assignment $C\mapsto j_{C*}\O_C(k)$ defines a closed Lagrangian
embedding $s_k:\P^3\hookrightarrow\mathcal J_k$ and a section of
the support morphism $\mathcal J_k\to|h_S|$. Its normal bundle is
$T^*\P^3$. Let
\[
 \mathfrak B_k=\P^3\timesh_{\mathcal J_k}\P^3
\]
be its derived self-intersection, oriented by $K_{\P^3}$ on the
classical truncation. Its cohomological DT groups are
\begin{equation}\label{eq:boundary-DT}
 \mathrm{DT}^j_{\mathrm{coh}}(\mathfrak B_k)=
 \begin{cases}\CC,&j=-3,-1,1,3,\\0,&\text{otherwise}.\end{cases}
\end{equation}
The normalized periodic realization has rank four in odd parity and
zero in even parity. Moreover
\[
 \chi\bigl(\mathrm{DT}^*_{\mathrm{coh}}(\mathfrak B_k)\bigr)
     =\int_{\P^3}c_3(T^*\P^3)=-4.
\]
\end{theorem}
\begin{proof}
For an integral genus-three curve $C$, $\O_C(k)$ is torsion-free of
rank one and stable. Adjunction gives genus three and
$\chi(\O_C(k))=4k-2$. The universal Cartier curve is flat over
$\P^3$; its twisted structure sheaf therefore defines $s_k$.
The support of this sheaf is $C$, proving that its composite with
the support morphism is the identity. A section of this separated
morphism is a closed immersion.

The primitive Mukai vector has square $h_S^2=4$. Stability and
Serre duality on the K3 surface identify scalar endomorphisms with
the dual top Ext, and the trace-free obstruction group vanishes.
The K3 sheaf-moduli theorem \cite{mukai1984symplectic} gives a smooth symplectic
moduli space of dimension $v^2+2=6$. Its form pulls back to zero
on $\P^3$, since $H^0(\P^3,\Omega^2)=0$. The embedding has half
the ambient dimension, and hence is Lagrangian. The symplectic
pairing identifies its normal bundle with its cotangent bundle.

The derived intersection is $(-1)$-shifted symplectic and has
smooth classical truncation $\P^3$. The associated d-critical
structure is zero, and its virtual canonical bundle is
$K_{\P^3}^{\otimes2}$. The specified square root $K_{\P^3}$ gives
the vanishing-cycle sheaf $\CC_{\P^3}[3]$ by the oriented critical
construction \cite{Joyce,BBDJS}. Thus its hypercohomology is
$H^{j+3}(\P^3,\CC)$, proving \eqref{eq:boundary-DT}. The zero
critical presentation and the bar comparison of
Section~\ref{sec:v5-bar-bridge} give the normalized periodic model
$(\Omega^\bullet_{\P^3}((u)),u\,d)[3]$ and its stated parity.
Finally the Euler sequence gives $c_3(T\P^3)=4H^3$; dualizing
changes the top Chern class by $(-1)^3$. Both displayed Euler
calculations are consequently $-4$.
\end{proof}

This identifies the full graded cohomology of the specified
boundary self-intersection before taking its Euler characteristic.
It applies to this particular boundary embedding. A product formula
for the full component invariants would additionally require a
comparison for their actual restriction morphisms and their
oriented DT complexes. In particular a value of an intersection
number alone does not specify that comparison.

\textbf{The boundary Gauss--Manin calculation.}
Take a smooth family of polarized quartic K3 surfaces over a smooth
curve $T$, on an open set where the primitive sheaf-moduli spaces
and the embeddings $s_k$ in the theorem form families. The relative
linear system is the projective bundle $\P_T(p_*\O_S(H_S))$ of
relative dimension three. Its central self-intersection is formed
from the two copies of the actual embedding $s_k$, with the same
relative canonical orientation. The proof above applies fiberwise and identifies its DT complex with the relative zero-critical
realization. The classes $1,H,H^2,H^3$ are powers of the first
Chern class of the global tautological line bundle, so are
horizontal for the relative de Rham connection. After the shift
$[3]$ the finite periodic model is
\[
 P^{\overline0}=0,\qquad
 P^{\overline1}=\bigoplus_{a=0}^3\O_T((u))e_a,
 \qquad d_P=0,\qquad\nabla e_a=0.
\]
This follows as a chain comparison from
Theorem~\ref{thm:smooth-proper-comparison}; over a coordinate open
its connection matrix in this basis is zero. In particular it
applies to the constant family of the boundary embedding used
in the resolved-quintic calculation. This computes the boundary
self-intersection realization, not the entire central quintic
moduli space. The latter requires its actual component restriction
maps and the compatibility of their orientations.

\textbf{The boundary self-intersection for $|2h_S|$.}
The same boundary construction can be calculated for the degree-two
linear system, with its stability checked rather than inferred from
primitivity. Assume $\operatorname{Pic}(S)=\mathbb Z h_S$,
$h_S^2=4$. Then $|2h_S|=\mathbb P^9$ and its curves have arithmetic
genus nine. The sheaf $\mathcal O_C(kh_S)$, $C\in|2h_S|$, has Mukai
vector
\[
 v=(0,2h_S,8k-8),\qquad v^2=16.
\]
It is stable even when $C=C_1+C_2$ or $C=2C_1$ with
$C_i\in|h_S|$. Indeed its reduced linear Hilbert polynomial has
constant term $k-1$. A subsheaf supported with generic length one
on a proper component injects into
$\mathcal O_{C_i}((k-1)h_S)$, whose reduced constant term is
$k-3/2$. A proper subsheaf with full generic length has a nonzero
zero-dimensional quotient and hence a smaller reduced polynomial.
An integral $C$ is the usual rank-one torsion-free case.
Thus every member of this projective linear system defines a stable
sheaf; no smoothness of the curve support was imposed.

The \emph{stable locus} $J^{\mathrm{st}}_v$ of K3 sheaves is smooth
symplectic of dimension $v^2+2=18$: the trace-free obstruction space
vanishes by Serre duality for each stable sheaf. Primitivity is not
needed for this assertion about stable points; we do not assert
that the whole semistable compactification is smooth. The section
\[
 |2h_S|=\mathbb P^9\longrightarrow J^{\mathrm{st}}_v,
              \qquad C\longmapsto\mathcal O_C(kh_S)
\]
is a closed embedding, since the support map is its retraction.
Its pulled-back two-form vanishes on projective space and it has
half the dimension, so it is Lagrangian with normal bundle
$T^*\mathbb P^9$. Its derived self-intersection has the zero
reduced critical presentation on $\mathbb P^9$, with canonical
orientation. The two-sided bar/BV calculation therefore gives
\[
 \mathrm{DT}_{\mathrm{coh}}^j
     =\begin{cases}\mathbb C,&j=-9,-7,\ldots,7,9,\\0,&\text{otherwise},\end{cases}
 \qquad\chi(\mathrm{DT}_{\mathrm{coh}})=-10.
\]
Its normalized periodic rank is ten in odd parity, and the ten
hyperplane powers are horizontal in a relative K3 linear-system
family. This is the degree-two boundary self-intersection invariant,
not the entire central quintic sheaf-moduli invariant.

\section{A K3-pencil degeneration with an explicit central comparison}
\label{sec:second-degeneration}
There is another geometric construction related to the quartic component
of the quintic degeneration. It has a nonprojective Calabi--Yau open
threefold with proper surface supports. The extra boundary is essential:
the projective threefold and one of its projective central components
are not themselves Calabi--Yau pairs with only the gluing surface.
We first construct the family, then identify a moduli locus of perfect
objects on every fiber, including the derived matching locus, and
finally calculate both its DT cohomology and its base connection.

\subsection{The pencil, the two blowups, and the removed boundary}
\label{ssec:second-geometry}
Let $Y\subset\P^4$ be a smooth quartic threefold, and let $H$
be its hyperplane class. Choose independent sections
$s,c\in H^0(Y,\O_Y(H))$ such that $S=(s=0)$ is a smooth K3
surface and $C_0=(s=c=0)$ is a smooth plane quartic. Thus
$g(C_0)=3$, and its normal bundle in $Y$ is
$\O_{C_0}(H)\oplus\O_{C_0}(H)$.
The blowup
\[
 Z=\Bl_{C_0}Y
\]
resolves the pencil $[s:c]$, giving $f:Z\to\P^1$;
we use $\lambda=-c/s$ on its finite coordinate chart.
If $E_0$ is its exceptional divisor and $\widetilde H$ is the
hyperplane pullback, the fiber class is
$F=\widetilde H-E_0$. The blowup formula gives
\begin{equation}\label{eq:second-canonical}
 K_Z=-\widetilde H+E_0=-F,
 \qquad K_Z\simeq f^*\O_{\P^1}(-1).
\end{equation}
The fiber $F_\infty$ corresponding to $s=0$ is the strict transform
of $S$, isomorphic to $S$ because $C_0$ is Cartier in it.
It follows that $U=Z\setminus F_\infty$ has trivial canonical
bundle after a choice of frame of $K_Z(F_\infty)$.
Its fibers over $\A^1_\lambda\subset\P^1$ are proper surfaces
isomorphic to $D_\lambda=(c+\lambda s=0)\subset Y$.
We may choose a general pencil so that its singular members have
only ordinary double points. The arguments below use only that
the members are integral quartic surfaces, smooth along $C_0$;
these properties hold for such a choice.

Here is a semistable degeneration of this open threefold. Begin
with the blowup
\[
 \mathscr Y=\Bl_{S\times\{0\}}(Y\times\A^1_t).
\]
Its central fiber is $Y\cup_S P$, where
\[
 P=\P_S(\O_S\oplus\O_S(-H))
\]
uses one-dimensional quotients. Denote by $S_\infty$ the section
with normal bundle $\O_S(-H)$ and by $S_0$ the other section,
with normal bundle $\O_S(H)$. The first is glued to $Y$ along
$S$. The strict transform of $C_0\times\A^1_t$ is smooth and
has central fiber in $S_0$, away from $S_\infty$. Blow up this
strict transform, and write
\[
 \overline{\mathscr U}\longrightarrow\A^1_t,
 \qquad \overline{\mathscr U}_t=Z\quad(t\ne0),
 \qquad \overline{\mathscr U}_0=Y\cup_{S_\infty}V,
 \quad V=\Bl_{C_0}P.
\]
The strict transform $B_\infty\subset V$ of $S_0$ is the central
fiber of the strict transform $\mathscr B$ of
$S\times\A^1_t$. Remove this relative divisor and set
\begin{equation}\label{eq:second-open-family}
 \mathscr U=\overline{\mathscr U}\setminus\mathscr B,
 \qquad U_t=U\ (t\ne0),\qquad
 U_0=Y\cup_S(V\setminus B_\infty).
\end{equation}

\begin{proposition}\label{prop:second-geometry}
The total space $\overline{\mathscr U}$ is smooth and its central
fiber is simple normal crossing. There is a relative logarithmic
volume form
\[
 \omega_{\overline{\mathscr U}/\A^1_t}(\mathscr B)
       \simeq\O_{\overline{\mathscr U}}.
\]
Its restriction gives $\omega_{\mathscr U/\A^1_t}\simeq\O_{\mathscr U}$.
On the two projective central components the identities are
\begin{equation}\label{eq:second-log-pairs}
 K_Y+S=0,\qquad K_V+S_\infty+B_\infty=0,
 \qquad S_\infty\cap B_\infty=\varnothing.
\end{equation}
Let $b:V\to P$, let $E$ be its exceptional divisor, and set
$G=b^*p^*H-E$. Then $G$ is the strict transform of
$p^{-1}(C_0)$, and
\begin{equation}\label{eq:second-ruled-support}
 G\simeq\P_{C_0}(\O\oplus\O(-H)),\qquad
 G\cap B_\infty=\varnothing,\qquad
 G\cap S_\infty=C_0,\qquad N_{G/V}=\O_G(-s_\infty).
\end{equation}
Here $s_\infty$ denotes the section of the ruled surface with
normal bundle $\O_{C_0}(-H)$.
\end{proposition}
\begin{proof}
Near $C_0$, take local equations $u=s$, $v=c$, and a coordinate
along $C_0$. The first blowup has charts
$t=u\tau$ and $u=t\alpha$. These are smooth, with the two
central components crossing normally in the first chart.
The strict transform of $C_0\times\A^1$ lies in the second
chart and is cut out by $\alpha=v=0$. Blowing up this regular
codimension-two center is smooth; since $t$ is an independent
coordinate there, the center meets the central component smoothly
and away from its crossing with $Y$. This proves the assertion
on the total space and the indicated central components.

Before either blowup, $K_{Y\times\A^1/\A^1}+(S\times\A^1)=0$.
For each blowup the center has codimension two and is contained
with multiplicity one in the strict transform of that boundary.
The discrepancy of the canonical divisor is $+E$, whereas the
strict transform of the boundary subtracts $E$. Their sum
therefore remains trivial. This proves the relative logarithmic
formula, including the claim on the complement.
On the projective bundle $K_P=-S_0-S_\infty$.
After the second blowup $B_\infty=b^*S_0-E$ and
$K_V=b^*K_P+E$, giving \eqref{eq:second-log-pairs}.
The two sections were disjoint and the center is contained only
in $S_0$, so their strict transforms remain disjoint.

The surface $p^{-1}(C_0)$ meets $S_0$ transversely along the
center. Its strict transform is disjoint from $B_\infty$.
Blowing up that surface along its Cartier section changes no
surface, giving the identification in
\eqref{eq:second-ruled-support}. Its intersection with
$S_\infty$ is still $C_0$.
The relation $S_0=S_\infty+p^*H$ gives
$B_\infty=S_\infty+G$ after the blowup. Restriction to $G$,
which avoids $B_\infty$, consequently gives
$\O_G(G)=\O_G(-s_\infty)$. This proves the normal-bundle formula.
\end{proof}

In particular $S_\infty$ alone is not anticanonical in the
projective $V$. The source of the second restriction morphism is
the moduli of perfect complexes with proper support avoiding
$B_\infty$. Proposition~\ref{prop:sc-two-boundaries}, applied to
\eqref{eq:second-log-pairs}, gives its Lagrangian structure.
Together with the first component this defines
\begin{equation}\label{eq:second-matching-stack}
 \mathfrak G^{(2)}=
 \mathfrak{Perf}(Y)\times^{\mathrm h}_{\mathfrak{Perf}(S)}
 \mathfrak{Perf}_{\mathrm{supp}\cap B_\infty=\varnothing}(V).
\end{equation}
The subscript specifies the indicated open substack; no
transversality of a sheaf to $S$ is imposed. This is a moduli
problem on the fixed components, not on expanded targets.

\subsection{A family of perfect complexes and its deformation theory}
\label{ssec:second-objects}
Let $\rho:\overline{\mathscr U}\to Y\times\A^1$ denote the
composition of the two blowups, and let $\mathscr E$ be the
exceptional divisor of the second one. The line bundle
\[
 \mathscr L=\rho^*\O_Y(H)\otimes\O(-\mathscr E)
\]
has restrictions $\O_Z(F)$ on a generic fiber,
$\O_Y(H)$ on the first central component, and $\O_V(G)$ on
the second. The pullbacks of $s,c$, divided by the exceptional
factor of the second blowup, define sections $\widetilde s$ and
$\widetilde c$ of $\mathscr L$.
On $Y$ these restrict to $s,c$. On $V$, $\widetilde s$ is zero
and $\widetilde c$ cuts out $G$.
In the local blowup chart with $v$ as exceptional parameter,
$u=t\alpha=tva$ and the divided sections are $ta,1$.
Thus $\widetilde c+\lambda\widetilde s$ is nonvanishing on
$\mathscr B$ for finite $\lambda$.

On $\mathscr U\times\A^1_\lambda$ form the universal complex
\begin{equation}\label{eq:second-universal-complex}
 E_{\lambda,t}=
 [\mathscr L_t^{-1}\xrightarrow{\widetilde c+\lambda\widetilde s}
                          \O_{U_t}],
 \qquad\hbox{in degrees }-1,0.
\end{equation}
Its section is a non-zero-divisor on every fiber: on each central
component its restriction is a nonzero section. It is therefore
the structure sheaf of a relative Cartier divisor, flat over
$(\lambda,t)$. Its support is proper over that parameter space,
being closed in the projective family
$\overline{\mathscr U}\times\A^1_\lambda$ and disjoint from
$\mathscr B\times\A^1_\lambda$.
For $t\ne0$ it is the structure sheaf of the pencil fiber
$D_\lambda\simeq f^{-1}(\lambda)$ in $Z$.
At $t=0$ its support is
\begin{equation}\label{eq:second-central-surface}
 W_\lambda=D_\lambda\cup_{C_0}G,
 \qquad
 (E_{1},E_{2},\alpha)=(\O_{D_\lambda},\O_G,\id_{\O_{C_0}})
       \in\mathfrak G^{(2)}.
\end{equation}
These component restrictions and the identity are derived
restrictions of the actual complex
\eqref{eq:second-universal-complex}. They happen to be sheaves
in this particular example.

\begin{theorem}\label{thm:second-moduli}
The family \eqref{eq:second-universal-complex} defines an open
immersion
$\A^1_\lambda\times\A^1_t\to M^{\mathrm{simp}}$
into the classical scalar-rigidified moduli space of simple
perfect complexes on the fibers, on its specified finite-type
neighborhoods. At every point its self-extension dimensions are
\begin{equation}\label{eq:second-exts}
 \dim\operatorname{Ext}^{0,1,2,3}(E_{\lambda,t},E_{\lambda,t})=(1,1,1,1),
 \qquad \operatorname{Ext}^{<0}=0.
\end{equation}
Its central locus identifies with a locally complete classical
locus of the effective triples \eqref{eq:second-central-surface}
in \eqref{eq:second-matching-stack}.
\end{theorem}
\begin{proof}
On a noncentral fiber, the support is a quartic surface in
$\P^3$, smooth K3 or nodal according to the pencil value. The hypersurface sequence
in $\P^3$ gives $H^{0,1,2}(\O_{D_\lambda})=(\CC,0,\CC)$.
Its normal line bundle as a fiber of $f$ is trivial.
The Cartier resolution of its structure sheaf gives
\[
 \RH_{U_t}(\O_{D_\lambda},\O_{D_\lambda})
  \simeq R\Gamma(\O_{D_\lambda})
           \oplus R\Gamma(\O_{D_\lambda})[-1],
\]
proving \eqref{eq:second-exts} for $t\ne0$.
This splitting is the direct sheaf-Hom calculation for this
structure sheaf with its Cartier resolution, not a splitting
asserted for arbitrary perfect objects.

We compute the central terms explicitly. On the ruled surface
$p_G:G\to C_0$, $Rp_{G*}\O_G=\O_{C_0}$, and restriction to its
section is the identity on this pushforward. The normalization
sequence for $W_\lambda$ is
\[
 0\longrightarrow\O_{W_\lambda}\longrightarrow
       \O_{D_\lambda}\oplus\O_G
       \longrightarrow\O_{C_0}\longrightarrow0.
\]
The map $H^1(G,\O_G)\to H^1(C_0,\O_{C_0})$ is an isomorphism.
It follows that
$H^{0,1,2}(W_\lambda,\O_{W_\lambda})=(\CC,0,\CC)$.

Let $N=\mathscr L_0|_{W_\lambda}$ be the normal line bundle of
this Cartier support in $U_0$. Its restrictions are
$\O_{D_\lambda}(H)$ and $\O_G(-s_\infty)$, both restricting
to $\O_{C_0}(H)$ on the double curve.
All cohomology of $\O_G(-s_\infty)$ vanishes, since its
restriction to each $\P^1$ fiber is $\O_{\P^1}(-1)$.
On $D_\lambda$, $C_0$ is the divisor $H$, so
\[
 0\longrightarrow\O_{D_\lambda}\xrightarrow{s}
 \O_{D_\lambda}(H)\longrightarrow\O_{C_0}(H)\longrightarrow0.
\]
In particular the restriction on global sections is surjective
with one-dimensional kernel. Tensoring the normalization sequence
by $N$ and using these facts gives
$H^{0,1,2}(W_\lambda,N)=(\CC,0,\CC)$; the last copy is
$H^1(C_0,\O_{C_0}(H))$, one-dimensional because
$K_{C_0}=\O_{C_0}(H)$.
The Cartier resolution on $U_0$ now gives
\[
 \RH_{U_0}(\O_{W_\lambda},\O_{W_\lambda})
  \simeq R\Gamma(\O_{W_\lambda})\oplus R\Gamma(N)[-1].
\]
This proves \eqref{eq:second-exts} also at $t=0$, and simplicity
and negative self-Ext vanishing follow.

The relative Kodaira--Spencer map of the parameter $\lambda$
is an isomorphism onto $\operatorname{Ext}^1$: at a central point it is the
nonzero section $s|_{D_\lambda}$ generating the just computed
kernel, and on a noncentral fiber it moves its support in the
pencil. We explain why this tangent calculation gives the full
classical deformation statement here. The simple universally
gluable perfect-object moduli functor is locally of finite
presentation on a neighborhood of each point, by
\cite{Lieblich2006}, applied on the projective compactification
and then restricted to supports avoiding $\mathscr B$.\\
Let $B$ be its complete local ring after scalar rigidification,
and $R=\CC[[t-t_0]]$ the base ring. The relative tangent space
has dimension one, so the complete local algebra has a
presentation $R[[x]]\twoheadrightarrow B$. The given family
induces $B\to R[[\lambda-\lambda_0]]$. Choose $x$ so that its
linear coefficient in $\lambda-\lambda_0$ is nonzero. The
composite is an isomorphism of power-series rings over $R$, by
formal inversion in that variable. It follows that the kernel
of $R[[x]]\to B$ is zero and that the family identifies the two
completed local rings. This proves formal etaleness, and hence
etaleness by local finite presentation.

The supports on $Y$ distinguish the central parameters, and the
pencil fibers distinguish the noncentral parameters, over every
algebraically closed extension. Thus the map is universally
injective, and the etale map is an open immersion.
Finally \eqref{eq:sc-mapping-triangle}, which is also valid for
these proper supports on the open union, identifies its central
formal deformation problem with the matching-triple problem.
This gives the asserted locally complete locus in
\eqref{eq:second-matching-stack}.
\end{proof}

\textbf{The component restriction maps near the central locus.}
The matching in \eqref{eq:second-central-surface} can also be seen
directly on the component sheaves. Let $\mathcal H$ be the open
subset of $\P H^0(Y,\O_Y(H))\setminus\{[s]\}$ where the
intersection with $S$ is a smooth curve. Restriction of sections
gives
\[
 \rho_H:\mathcal H\longrightarrow |H|_S^{\mathrm{sm}}.
\]
Its fiber over $C_0$ consists precisely of the hyperplanes
$D_\lambda=(c+\lambda s=0)$, hence is $\A^1_\lambda$.
Near their structure sheaves, $\mathcal H$ is the classical
component moduli locus: the Cartier resolution gives
\[
 \begin{aligned}
 \operatorname{Ext}^1_Y(\O_{D_\lambda},\O_{D_\lambda})
       &=H^0(D_\lambda,\O_{D_\lambda}(H))\simeq\CC^4,\\
 \operatorname{Ext}^2_Y(\O_{D_\lambda},\O_{D_\lambda})
       &=H^2(D_\lambda,\O_{D_\lambda})\simeq\CC.
 \end{aligned}
\]
The universal hyperplane family induces the isomorphism on the
first group; the same complete-local-ring argument as in the
proof of Theorem~\ref{thm:second-moduli} gives its local
completeness. The displayed obstruction direction remains in
the derived component moduli stack.

On the other component, $N_{G/V}=\O_G(-s_\infty)$ has zero
cohomology. The Cartier calculation at $\O_G$ gives
\[
 \operatorname{Ext}^1_V(\O_G,\O_G)=H^1(G,\O_G)\simeq\CC^3,
 \qquad\operatorname{Ext}^2_V(\O_G,\O_G)=0.
\]
Consequently the smooth family of line bundles
$\Pic^0(G)\simeq\Pic^0(C_0)$ is locally complete at this
component sheaf. The restriction of a line bundle from $G$ to
its section $C_0$ is the inverse of pullback from $C_0$.
Let $J$ be the smooth six-dimensional K3 sheaf-moduli space
near $\O_{C_0}$, as in \cite{mukai1984symplectic}. Over the smooth-curve open it is the relative
Jacobian of degree-zero line bundles, with section $s_0$
assigning $C'\mapsto\O_{C'}$. The classical restriction maps
are therefore
\[
 \begin{tikzcd}[column sep=huge]
 \mathcal H\ar[r,"s_0\circ\rho_H"]&J&
 \Pic^0(C_0)\ar[l,hook',"\text{fiber over }C_0"'] .
 \end{tikzcd}
\]
Matching forces the curve to be $C_0$ and the boundary line
bundle to be $\O_{C_0}$. Since $\Pic^0(G)\to\Pic^0(C_0)$
is an isomorphism, it also forces the component line bundle to
be $\O_G$. After scalar rigidification, the remaining parameter
is exactly $\rho_H^{-1}(C_0)=\A^1_\lambda$.

The support directions and Jacobian directions span the tangent
of $J$ at $\O_{C_0}$. Thus the component tangent map
$\CC^4\oplus\CC^3\to\CC^6$ is surjective with
one-dimensional kernel. The full derived restriction maps,
with the obstruction direction on the first component retained,
are the Lagrangian maps in \eqref{eq:second-matching-stack}.
Their intersection has the deformation complex computed in
\eqref{eq:second-exts}. This identifies the central parameter
and its extra derived direction through the actual component
restriction maps.

\subsection{Generic and special DT complexes and the Gauss--Manin operator}
\label{ssec:second-DT}
Give the smooth moduli family of
Theorem~\ref{thm:second-moduli} its zero d-critical
realization and the relative orientation
$K_{\A^1_\lambda\times\A^1_t/\A^1_t}$, with frame $d\lambda$.
Use the same orientation on the actual central matching locus.
The relative logarithmic volume supplies the shifted forms on
the proper-support object moduli. On a smooth classical scheme
the zero d-critical structure is the unique one; consequently
any compatible reduced critical realization has this structure.
Here the chosen realization is on the actual classical moduli
locus identified in Theorem~\ref{thm:second-moduli}.
The normalization comparison identifies the central locus, and
we transport the stated orientation through that identification.
The full Artin stack is still distinguished from this reduced
critical realization.

\begin{theorem}\label{thm:second-DT-comparison}
For every $t$, including $0$, the DT sheaf on the moduli locus
$M_t=\A^1_\lambda$ of Theorem~\ref{thm:second-moduli} is
$\CC_{\A^1}[1]$. With compact supports,
\begin{equation}\label{eq:second-DT-groups}
 \mathrm{DT}_{\mathrm{coh}}^j(M_t)=
 \begin{cases}\CC,&j=1,\\0,&j\ne1.\end{cases}
\end{equation}
The oriented bar realization on the actual central matching
locus has the same compactly supported groups. Its periodic
compact-support realization and that on the smooth fibers
admit the common finite model
\[
 P^{\overline0}=0,\qquad
 P^{\overline1}=\O_{\A^1_t}\otimes_\CC\CC((u)),\qquad
 d_P=0,\qquad\nabla_{\partial_t}=\partial_t.
\]
The bar comparisons commute with specialization and are
horizontal. The chosen ordinary zero-critical categories on
all $M_t$ are identified with $\MF(\A^1,0)$; their normalized
ordinary periodic cyclic homology is likewise one-dimensional
in odd parity, with trivial variation in $t$.
\end{theorem}
\begin{proof}
A smooth one-dimensional zero-critical locus, oriented by its
canonical line bundle, has vanishing-cycle sheaf $\CC[1]$.
This applies to both endpoints by the moduli identification
and orientation above. Since $H_c^2(\A^1,\CC)=\CC$ and all its
other compactly supported cohomology is zero, the shift gives
\eqref{eq:second-DT-groups}.

For the chain calculation use the global coordinate $\lambda$
and the zero potential. The bar-to-Koszul comparison and
contraction with $d\lambda$ give
\[
 \mathcal B_{M,\mathrm{DT}}
     \longrightarrow (\Omega^\bullet_{\A^1_\lambda}((u)),u\,d)[1].
\]
These formulas are independent of $t$ and are formed on the same
moduli family, so the operator on them is $\partial_t$ and all
maps specialize to the central bar map. A compact-support model
is obtained on $\P^1_\lambda$ using
\[
 \Omega^\bullet_{\P^1}(\log\{\infty\})(-\infty).
\]
Its degree-zero sheaf is $\O(-1)$ and its degree-one sheaf is
$\O(-2)$. Its hypercohomology is a single copy of $\CC$ in degree
two. This can be seen directly from their sheaf cohomology;
there is no possible differential hitting or leaving
$H^1(\O(-2))$. After the critical shift and periodization, this
is the displayed odd rank-one model. Choose its generator once
on $\P^1_\lambda$ and tensor it with $\O_{\A^1_t}$.
The finite comparison and its contracting homotopies can be
chosen over $\CC$, so are horizontal after this tensor product.
No properness of $M_t\to\Spec\CC$ is used: compact support is
implemented by this specified compactification.

For ordinary rather than compactly supported cohomology,
$H^*(\A^1,\CC)$ is a single copy in degree zero. The normalized
zero-critical cyclic comparison of
\cite{Efimov2018CyclicMF} therefore also has odd rank one after
$[1]$. The isomorphism $M_t=\A^1$ identifies the specified
zero-critical categories themselves. This last statement concerns
these chosen ordinary categories. It is not an identification
of ordinary and compactly supported cohomology in their
unperiodized degrees, nor a claim about every construction of
an intrinsic category of the full Artin moduli stack.
\end{proof}

The parameter of this connection is $t$, the degeneration
parameter. The variation of $H^2(D_\lambda)$ as $\lambda$ moves
in the K3 pencil is a different Gauss--Manin system, which may
have nontrivial monodromy. It does not contradict the constant
DT variation just computed for the one-dimensional object-moduli
locus along $t$.

Tensoring \eqref{eq:second-universal-complex} by the pullback of
$\O_Y(kH)$, for any $k\in\mathbb Z$, gives further perfect-object
examples with identical deformation complexes and transported
orientations. On a generic fiber $Z$ their characters are
\[
 \ch(\O_{f^{-1}(\lambda)}(k\widetilde H))
   =\bigl(0,F,k\widetilde HF,2k^2[\mathrm{pt}]\bigr),
 \qquad\chi=2+2k^2.
\]
Here $F^2=0$ and $\widetilde H^2F=4$.
The matching identity on the central fiber is
\[
 [\O_{W_\lambda}(k)]=[j_{1*}\O_{D_\lambda}(k)]
      +[j_{2*}\O_G(k)]-[j_{S*}\O_{C_0}(k)]
\]
in the coherent Grothendieck group of the union.
These examples have proper object supports but a nonproper
parameter space; the compact-support convention is part of
their computed invariant. Modifications by points or by
line bundles on the support not obtained by restriction from the threefold give additional moduli problems whose
central completion and orientation have to be constructed
separately. The general cone construction of
Proposition~\ref{prop:general-perfect-fiber} applies to such
objects once its perfectness hypotheses are met.

\textbf{Nodal members of the pencil.}
The structure-sheaf calculation in this section includes the nodal
quartic surfaces in the pencil.  Their Cartier resolutions, the
vanishings of $H^1(\O)$ and $H^0(N)$ used in the proof, and the
one-parameter deformation reconstruction hold at these members
as well.  Hence the specified classical object-moduli space remains
$\A^1_\lambda$ there; its singular classical locus is empty.
Its reduced derived enhancement still has the dual obstruction
direction.  A Hilbert scheme of points on a fixed nodal surface is
a different moduli problem and cannot be substituted for this
structure-sheaf family.  The defect calculations of
Sections~\ref{sec:singular-comparison} and
\ref{sec:crossing-perfect} show explicitly how such extra data can
alter the deformation problem.

\section{Singular critical charts in degrees one and two}
\label{sec:singular-comparison}
The examples in Section~\ref{sec:v6-hyperplanes} have smooth classical
parameter schemes. Here we keep a singular classical scheme and determine its
nonzero d-critical structure. We first verify the degree-two positivity
and stability assertions used above. We then compare the singular chart
with actual perfect objects on the generic and special threefolds.
The calculation is local in the moduli space; the separate proper
singular comparison in Section~\ref{sec:cubic-geometry} will have a
complete global DT group.

\subsection{Positivity and stability for the quadratic linear system}
\begin{lemma}\label{lem:quadric-small-curve-stability}
Let $X$ be a smooth quintic, $D\in|2H|$, and $Z\subset D$ have dimension
at most one. Suppose the degree of its one-dimensional cycle is $e_1<5$.
Then $j_{D*}I_{Z/D}$ is Gieseker stable for polarization $H$.
In particular this applies to finite subschemes and to a line with any
finite number of isolated or embedded point defects.
\end{lemma}
\begin{proof}
The Cartier divisor $D$ is Cohen--Macaulay of pure dimension two, and
$I_{Z/D}$ is a subsheaf of $\mathcal O_D$ with the same generic lengths.
It is therefore pure. Its Hilbert polynomial has leading terms
\[
 P_F(a)=5a^2-(10+e_1)a+O(1),\qquad
 p_F(a)=a^2-\left(2+\frac{e_1}{5}\right)a+O(1).
\]
If $D$ is integral, any nonzero subsheaf has the same rank-one leading
term. A proper quotient has a positive polynomial of lower degree,
which gives strict stability.

Otherwise Picard rank one implies that $D=D_1+D_2$, with
$D_i\in|H|$ integral; equality $D_1=D_2$ is allowed. A subsheaf with
full generic length has the same leading quadratic coefficient as $F$,
and a nonzero quotient again gives the strict inequality. A subsheaf
with smaller generic length is supported on one reduced $D_i$. Its
inclusion into $\mathcal O_D$ factors through
\[
 \ker(\mathcal O_D\to\mathcal O_{D_{3-i}})
       =\mathcal O_{D_i}(-H)
\]
when the components differ. When $D=2D_1$, the corresponding maximal
submodule annihilated by the ideal of $D_1$ is likewise
$\mathcal O_{D_1}(-H)$. Its reduced linear coefficient is $-3$;
a subsheaf can only decrease this coefficient, or leave it equal and
decrease the constant coefficient. Since
$-3<-2-e_1/5$ for $e_1<5$, such a subsheaf cannot destabilize $F$.
This treats every possible generic length on the divisor.
\end{proof}

\begin{lemma}\label{lem:quadric-quot-evaluation}
Let $\Lambda\subset X$ be a line on a smooth quintic, and let
$I_\Lambda\twoheadrightarrow Q$ be any length-two quotient with kernel
$I_Z$. Then
\[
 H^0(I_\Lambda(2))\twoheadrightarrow H^0(Q(2)),\quad
 H^{>0}(I_Z(2))=0,\quad h^0(I_Z(2))=10.
\]
Consequently the scheme of $(Z,D)$ with $D\in|2H|$ is a
$\mathbb P^9$-bundle over $\operatorname{Quot}_X(I_\Lambda,2)$.
\end{lemma}
\begin{proof}
The ideal $I_\Lambda(1)$ is generated by the three linear equations of
$\Lambda$. Denote their images in $H^0(Q(1))$ by $v_1,v_2,v_3$.
Choose a linear form $\ell_0$ nonvanishing at the finite support of $Q$.
Multiplication by $\ell_0$ identifies $H^0(Q(1))$ with $H^0(Q(2))$.
If the $v_i$ span a two-dimensional vector space, their products with
$\ell_0$ prove surjectivity.

If their span has dimension one, a nonzero vector $v$ in that span
must generate $Q(1)$ as an $\mathcal O_X$-module: the $v_i$ generate
it as a sheaf. Hence $Q(1)$ is a cyclic length-two module, isomorphic
locally to the structure sheaf of a length-two subscheme $W$.
The very ample bundle $\mathcal O_X(1)$ separates $W$, including a
tangent double point. Thus the products $\ell v$, for linear forms
$\ell$, span $H^0(Q(2))$. Each product is the image of a quadratic
section of $I_\Lambda(2)$. A length-two quotient with two generators
at the same reduced point cannot occur in this one-dimensional case;
its sheaf generators already require two independent $v_i$.
These alternatives exhaust the quotient modules.

The sequence of $\Lambda\simeq\mathbb P^1$ and Kodaira vanishing give
$H^{>0}(I_\Lambda(2))=0$ and $h^0(I_\Lambda(2))=15-3=12$.
The quotient sequence therefore gives the asserted vanishings and
$h^0(I_Z(2))=12-2=10$. The universal finite flat quotient has a
rank-two tautological bundle. Surjectivity just proved makes the
kernel a rank-ten vector bundle, commuting with base change.
Its projectivization parametrizes exactly the sections vanishing on
$Z$, modulo scalar. This also verifies the precise positivity
condition of~\cite[Definition~2.2]{GST}; smoothness of the Quot scheme
has not been used.
\end{proof}

The reconstruction in Proposition~\ref{prop:point-moduli-reconstruction}
also works with $L^m$ and a perfect proper subscheme $Z$ of dimension
at most one with $H^1(\mathcal O_Z)=0$. Indeed its first duality
calculation only requires $R\Gamma(I_ZL^m)$ in degree zero.
For the second reconstruction,
\[
 H^0(I_Z)=0,\quad H^2(I_Z)=H^1(\mathcal O_Z)=0,\quad H^3(I_Z)=\mathbb C
\]
again place $R\operatorname{Hom}(I_Z,\mathcal O_X)$ in degrees zero
and two. The lifted quotient is flat and proper, with the same Hilbert
polynomial; its reduction need not be zero-dimensional for that
flatness argument. Thus it recovers a curve subscheme in this case.
For a line with point defects, the exact sequence
$0\to Q\to\mathcal O_Z\to\mathcal O_\Lambda\to0$ proves the required
$H^1$ vanishing. This explains why the singular incidence scheme in the
preceding lemma has the sheaf deformation theory, not just a map to it.

\subsection{Identifying the intrinsic determinant critical structure}
The potential in Proposition~\ref{prop:line-two-determinant-potential}
was obtained from the explicit affine Quot presentation of
\cite[Proposition~2.10]{DavisonRicolfi}. Identifying a classical critical
scheme does not in general identify its d-critical structure. In this
particular determinantal example that extra identification can be proved.
The following calculation supplies it.

\begin{lemma}\label{lem:determinantal-dcritical-class}
Let $R=\mathbb C\{z_1,\ldots,z_s,x_{11},\ldots,x_{33}\}$, or its
formal completion, and let $I=I_2(X)$ be the ideal of the $2\times2$
minors of the matrix $X=(x_{ij})$. The classes $f\bmod I^2$ satisfying
$f\in I$ and $df\in I\Omega_R^1$ are precisely
\[
 a(z)\det X\pmod{I^2},\qquad a(z)\in\mathbb C\{z_1,\ldots,z_s\}.
\]
If $f$ is a minimal critical presentation whose gradient ideal is $I$,
then $a(0)\ne0$.
\end{lemma}
\begin{proof}
In characteristic zero the condition $f\in I$, $df\in I\Omega_R^1$
characterizes the second symbolic power $I^{(2)}$. One can check this
at the smooth rank-one points, where $I$ is generated by independent
normal coordinates: vanishing of the first normal derivatives is
exactly order at least two. Taking the symbolic-power contraction
from that dense smooth locus gives the stated condition. For a generic
matrix the determinantal straightening calculation gives
\[
 I^{(2)}=I^2+I_3(X)=I^2+(\det X)
\]
\cite[Proposition~2.2 and the displayed special case $t=k=2$]{BrunsConca}.
The same equality holds after adjoining the $z$ coordinates and after
analytic or formal completion. Its module generator is killed by every
$x_{ij}$ modulo $I^2$. For example,
\begin{align*}
 x_{11}\det X={}&(x_{11}x_{22}-x_{12}x_{21})
                    (x_{11}x_{33}-x_{13}x_{31})\\
 &-(x_{11}x_{23}-x_{13}x_{21})
                    (x_{11}x_{32}-x_{12}x_{31}),
\end{align*}
and row and column permutations give all nine identities.
The determinant itself is not in $I^2$, since it has degree three
whereas $I^2$ has degree at least four in the $x$ variables. This proves
that the quotient is the free rank-one module generated by its class
over the smooth-factor ring.

If $a(0)=0$, the representative $a(z)\det X+I^2$ has no nonzero
cubic term at the origin. Its gradient ideal is then contained in the
cube of the maximal ideal. It cannot equal $I$, whose nine generating
minors have nonzero quadratic terms. Thus $a(0)$ is nonzero.
\end{proof}

\begin{theorem}\label{thm:determinant-intrinsic}
At the quotient $I_\Lambda/\mathfrak m_p I_\Lambda$ of length two,
let the degree-one incidence locus be restricted to its evaluation
rank-two neighborhood. For degree two use the whole projective-bundle
locus of Lemma~\ref{lem:quadric-quot-evaluation}. In each case the
intrinsic reduced d-critical structure of the sheaf moduli has a local
critical presentation
\[
 f(X,z)=\det X,\qquad
 U=\operatorname{Mat}_3\times\mathbb A^s,
 \qquad s=1\ (|H|),\quad s=10\ (|2H|).
\]
With an orientation identified with the canonical orientation of this
chart, the local perverse DT stalk has one-dimensional cohomology in
\[
 \begin{array}{c|c}
 |H|&-6,-4,-1\\
 |2H|&-15,-13,-10
 \end{array}
\]
and zero stalk cohomology in other degrees. These are stalk groups of
the sheaf DT realization, not the global DT groups of the proper Quot
moduli space.
\end{theorem}
\begin{proof}
The affine Quot calculation identifies its classical germ with
$\mathbb A^1\times\{\operatorname{rank}X\leq1\}$.
At the specified quotient the three linear sections generate its
length-two vector space, so evaluation rank two is an open condition
there. Degree one adds a zero-dimensional projective fiber; degree two
adds a smooth projective fiber of dimension nine. The positivity
vanishing and the reconstruction argument identify these with the
actual classical sheaf-moduli germs, including all Artinian
nilpotents. The ambient minimal embedding dimensions are ten and
nineteen, respectively.

The reduced CY deformation problem supplies a minimal analytic or
formal critical presentation, as in Section~\ref{ssec:sc-oriented-moduli}.
After removal of nondegenerate quadratic pairs, an isomorphism of
these classical germs lifts to an isomorphism of their minimal smooth
ambient germs: lift a basis of the cotangent space and use the inverse
function theorem (or its formal version). Thus its gradient ideal is
$I_2(X)$. Its normalized d-critical class is an element of the kernel
calculated in Lemma~\ref{lem:determinantal-dcritical-class}; it equals
$a(z)\det X\bmod I^2$ with $a$ a unit. Multiplying the first matrix row
by $a(z)$ is an invertible coordinate change preserving the rank-one
ideal and carrying the representative $a(z)\det X$ to $\det X$.
The latter has precisely that critical scheme and d-critical class.
This identifies the intrinsic d-critical chart, rather than assuming
that the critical presentation of the Quot scheme has the right
orientation. The orientation comparison is fixed locally as stated;
a different global orientation may still have transition data.

The fiber $\det X=1$ is $\mathrm{SL}_3(\mathbb C)$ and has the homotopy
type of $\mathrm{SU}(3)$. Its reduced cohomology is one-dimensional in
degrees $3,5,8$. The free $z$ variables do not change this Milnor
fiber. For an ambient dimension $N=9+s$, the normalized perverse
vanishing-cycle stalk has degrees $3-(N-1),5-(N-1),8-(N-1)$.
This gives the two displayed rows. In particular the degree-two
projective directions shift the degree-one answer by $-9$; they do
not remove the singular potential.
\end{proof}

\subsection{The actual central local comparison and its connection}
\begin{theorem}\label{thm:determinant-relative}
Let $\mathcal C\subset\mathcal X$ be a labeled rigid perfect-curve
family of Proposition~\ref{prop:resolved-line-family}. Choose a section
$p(t)$ of its relative smooth locus such that $p(0)$ lies outside
$S$ and outside the exceptional locus of the contraction
$a:\mathcal X\to\mathbb P^4$. Thus this contraction is a local
embedding at $p(0)$. The local sheaf-moduli problems with two defects at $p(t)$, and with
surface class $H$ or $2H$, have the determinant charts of
Theorem~\ref{thm:determinant-intrinsic} uniformly near $t=0$.
Their central objects are the actual derived restrictions of perfect
cones, and their reduced critical realizations in the central
Lagrangian intersection have the same nonzero determinant potential.
On these local charts the periodic cyclic and oriented two-sided bar
models are horizontally quasi-isomorphic, with connection $\partial_t$
in the determinant coordinates. This identifies the displayed local
DT stalks at generic and central fibers. It does not compute global
cohomology by taking a single stalk.
\end{theorem}
\begin{proof}
Relative analytic smooth coordinates near $p(t)$ identify the family
with coordinates $(x,y,z,t)$ and its smooth curve with $(x,y)=0$.
A finite quotient supported in this neighborhood is completely
determined by the actions of these three coordinates and the two
images of the generators of $(x,y)$. These are exactly the matrices
and vectors in~\eqref{eq:line-Quot-potential}, independently of $t$.
Embedded rigidity fixes the one-dimensional support curve; no
unrecorded support-curve deformation is omitted. Relative evaluation
rank two holds in degree one near the chosen quotient and in degree
two by Lemma~\ref{lem:quadric-quot-evaluation}.
The reconstruction therefore gives the same determinantal classical
moduli germ over the base, with the indicated projective coordinates.
Here evaluation rank two at $0$ follows because the contraction is a
local embedding at the chosen point and the linear equations of the
image line generate its two-dimensional normal ideal there. Their
images span $I_{C,p}/\mathfrak m_p I_{C,p}$. Multiplication by a linear
form nonzero at $p$ gives the same assertion for quadratic sections.

All point quotients here have support in the relative smooth locus,
so they are perfect, as are $I_{\mathcal C}$ and their kernels $I_Z$.
The actual sheaf complexes are
\[
 E=\operatorname{Cone}(L^{-m}\to I_Z),\quad m=1,2.
\]
Restriction to $X_i$ and $S$ gives the two triangles
\eqref{eq:line-two-central-triangles}, with $L^{-m}$, and the
canonical matching isomorphism. Since $p(0)$ lies on exactly one
component, the quotient on the other component and on $S$ is zero;
the curve ideal itself retains its actual boundary restriction.
The normalization endomorphism triangle identifies the full formal
deformation algebra of this perfect central cone with that of its
matching triple. Thus the determinantal germ is the actual central
classical germ of that matching problem.

The proof of Lemma~\ref{lem:determinantal-dcritical-class} works with
coefficients holomorphic in $t$. It identifies the relative induced
d-critical class as $a(z,t)\det X$, with $a$ invertible. Rescale one
row by $a$ and transport the orientation. This produces the common
representative $f=\det X$ independent of $t$. Its reduced model is
\[
 K_f=\mathcal O_U[\eta_1,\ldots,\eta_N],\quad
 \delta\eta_j=f_j,\quad
 D_K=-\sum_j f_j\partial_{\eta_j}
                  +u\sum_j\partial_{x_j}\partial_{\eta_j}.
\]
The coordinates $x_j$ include all the matrix and free coordinates;
the derivatives in the free coordinates are zero. The Hessian entries
in the tangent differential are the derivatives of the $2\times2$
minors, and are not discarded.

Put $A=\mathcal O_U[p_1,\ldots,p_N]$, with augmentations
$p_j\mapsto f_j$ and $p_j\mapsto0$. Its actual two-sided bar
$\operatorname{Bar}(\mathcal O_{U,df},A,\mathcal O_{U,0})$ retracts
onto $K_f$ by the maps $I,\Pi$ of Section~\ref{sec:v5-bar-bridge}.
Equipping this bar with the transferred mixed operator gives
\[
 \mathcal B_{f,\mathrm{per}}^{\mathrm{or}}
       \simeq(\Omega_U^\bullet((u)),u d-df\wedge)[N]
       \simeq CC_{\mathrm{per}}(\operatorname{MF}(U,f))[N].
\]
This is a critical presentation of the reduced actual matching locus
identified above, not an assertion that the unreduced component
stacks are $U$. In these coordinates $\partial_t f=0$, all bar and
Koszul maps are independent of $t$, and
$[\partial_t,D_K]=0$. The connection and the chart comparison are
therefore horizontal and flat. On a common small analytic Milnor
neighborhood the vanishing-cycle stalk calculation is constant in
$t$ and gives the claimed local specialization. These local
comparisons require their orientation and overlap maps before they
can be used to compute a global invariant.
\end{proof}

The theorem treats a genuinely singular moduli germ on both fibers.
Its curve may be of type A and meet $S$, but the two quotient defects
in this theorem were placed at a relative smooth point. The next
section treats defects at $S$ explicitly. Neither locality condition
is silently replaced by a properness assertion.

\section{Perfect complexes at the common surface and logarithmic transport}
\label{sec:crossing-perfect}
A nontransverse derived intersection retains derived matching data.
It does not turn a nonperfect coherent module on $X_0$ into a perfect
one. The distinction is concrete in~\eqref{eq:node-ideal-obstruction}.
This section constructs replacement perfect complexes at the crossing,
computes their component and boundary restrictions, and determines
what kind of connection can be compatible with the generic point
invariants. None of these calculations requires transversality of the
component-moduli intersection.

\subsection{A perfect replacement for a point reaching the common surface}
Work in the completed local model of a semistable family, and make the
base change $t=s^2$. Put
\[
 \begin{gathered}
 \mathscr R=\mathbb C[[s,x,y,z,w]]/(xy-s^2),\qquad R_0=\mathscr R/(s),\\
 D_s=(z=0),\qquad p_s=(x=y=s,\ z=w=0).
 \end{gathered}
\]
This is a local calculation on the given degeneration, not an added
component of the target. The ramified base change permits a section
$p_s$ whose limit lies at $S$.

\begin{proposition}\label{prop:crossing-one-perfect}
The following complex of finite free $\mathscr R$-modules, in degrees
$-2,-1,0$, is perfect and has generic fiber $I_{p_s/D_s}$ pushed
forward into the threefold:
\begin{equation}\label{eq:crossing-one-matrices}
 \mathscr E=\left[
 \mathscr R\xrightarrow{\ d^{-2}\ }\mathscr R^3
       \xrightarrow{\ d^{-1}\ }\mathscr R^2\right],\qquad
 d^{-2}=\begin{pmatrix}w\\-(x-s)\\z\end{pmatrix},\quad
 d^{-1}=\begin{pmatrix}z&0&-w\\0&z&x-s\end{pmatrix}.
\end{equation}
It is a flat coherent sheaf over $\mathbb C[[s]]$, but its special
fiber is not required to be pure. Its full derived restrictions to
$X_1=(x=0)$ and $X_2=(y=0)$ are, locally,
\[
 E_1\simeq\mathcal O_{D_1}\oplus\mathcal O_{C_1},\qquad
 E_2\simeq I_{p/D_2},
\]
where $D_i=(z=0)\subset X_i$ and $C_1=(z=w=0)\subset X_1$.
Both have the same derived restriction to $S$:
\[
 E_S\simeq\mathcal O_{(z=0)\subset S}\oplus\mathbb C_p.
\]
The matching equivalence is induced by the identical matrices on $S$.
The splitting displayed here is local; the matrix complex, not a choice
of splitting, is the intrinsic datum.
\end{proposition}
\begin{proof}
Direct multiplication gives
\[
 d^{-1}d^{-2}
   =\begin{pmatrix}zw-wz\\-z(x-s)+(x-s)z\end{pmatrix}=0.
\]
Write $K(z)=[\mathscr R\xrightarrow{z}\mathscr R]$ and
$P=[\mathscr R\xrightarrow{(-w,x-s)^T}\mathscr R^2]$.
The displayed complex is their total tensor product, with its Koszul
sign. On $\mathscr R/(z)$ the sequence $(x-s,w)$ is regular: $x-s$
is a non-zero-divisor, and after its quotient $w$ remains a free
power-series variable. The single syzygy between its two generators
is $(-w,x-s)$. Hence $\mathscr E$ resolves the ideal
$(x-s,w)$ in $\mathscr R/(z)$. This ring is a domain and multiplication
by $s$ is injective on the ideal. A torsion-free module over the discrete
valuation ring $\mathbb C[[s]]$ is flat. This proves the flat-sheaf
assertion as well as perfectness of the displayed resolution.
For $s\ne0$, the relation $xy=s^2$ implies $y=s$ whenever $x=s$;
thus the ideal cuts out precisely $p_s$ in $D_s$.

For the central restrictions one substitutes $s=0$ in the matrices
and then $x=0$ or $y=0$. On $X_1$ the factor $P$ becomes
$[\mathcal O\xrightarrow{(-w,0)^T}\mathcal O^2]$, whose cokernel
is $\mathcal O/(w)\oplus\mathcal O$, with no negative cohomology.
Tensor with $K(z)$ to obtain the asserted $E_1$.
On $X_2$, the pair $(x,w)$ is regular on the surface $z=0$, so the
same two-generator presentation resolves $I_{p/D_2}$.
Finally on $S$ the factor is again $(-w,0)^T$, followed by $K(z)$.
This gives $\mathbb C_p\oplus\mathcal O_{(z=0)}$. These are actual
derived restrictions because the original terms are free. In particular
all coherence homotopies of the matching are inherited from one
complex on $R_0$.
\end{proof}

This replacement is different from the ordinary flat family of
\emph{pure} point ideals whose central object is
$I_{p/D_0}$. The latter is not perfect. The new family permits a
one-dimensional torsion contribution on one component in its special
fiber. Equivalently, let
\[
 A=K(z),\qquad Q=K(x-s,z,w).
\]
There is the natural augmentation-compatible map $A\to Q$, and
$\mathscr E\simeq\operatorname{fib}(A\to Q)$.
Although $Q$ is perfect, it is not a flat family of length-one
quotient sheaves at $s=0$: $Q_0$ has two cohomology sheaves supported
on $(x=z=w=0)$. Thus allowing a \emph{perfect quotient complex}
rather than insisting on a finite flat point quotient supplies this
boundary point. This is an explicit instance of the arrow construction
of Proposition~\ref{prop:general-perfect-fiber}.

For the original global problem, a precise larger parameter functor
is the functor of triples $(A_T,Q_T,\epsilon_T)$, where $A_T$ is a
section complex for $|mH|$, $Q_T\in\operatorname{Perf}(\mathcal X_T)$,
and $\epsilon_T:A_T\to Q_T$ is an actual derived morphism, whose
restriction over $C^*$ is the usual point evaluation. One imposes
specified numerical classes and Tor-amplitude bounds on finite-type
subfunctors. Taking the fiber gives a morphism to
$\mathfrak{Perf}(\mathcal X/C)$. Proposition~\ref{prop:crossing-one-perfect}
constructs an additional point of this functor in the local crossing
model and gives all its matching matrices. \\ The functor includes
framing data; forgetting these data need not be an equivalence of
moduli problems. Properness of its image, global simplicity, and an
orientation on a selected image are not consequences of the local
matrix calculation. In particular it has not been identified with a
proper global compactification of the entire one-point moduli space.

For either $m=1$ or $m=2$, choose a divisor equation with nonzero
differential transverse to the displayed point path. It is a local
coordinate $z$ after an analytic coordinate change and trivialization
of $L^m$. Thus the same complex treats point collisions on both linear
systems. Global degrees change the parameter spaces and charges, not
this local resolution.

\subsection{A length-two perfect quotient supported at the crossing}
The situation is different for certain double points. Use
$\mathscr R=\mathbb C[[t,x,y,z,w]]/(xy-t)$ and take a Hilbert chart
with quotient basis $(1,x)$. Write its multiplication rules as
\begin{equation}\label{eq:doublepoint-rules}
 x^2=ax+b,\quad y=cx+d,\quad z=ex+f,\quad w=gx+h.
\end{equation}
The equation $xy=t$ in that two-dimensional algebra is equivalent to
\[
 d=-ca,\qquad cb=t.
\]
These equations retain the base parameter and the scheme structure,
not just the two support points.

\begin{proposition}\label{prop:crossing-two-perfect}
On the open set $c\ne0$ the Hilbert chart~\eqref{eq:doublepoint-rules}
is smooth of relative dimension six over $t$. Its universal quotient
is a perfect relative length-two sheaf, resolved by the Koszul complex
of the regular sequence
\begin{equation}\label{eq:doublepoint-koszul}
 y-cx+ca,\qquad z-ex-f,\qquad w-gx-h.
\end{equation}
At $a=b=e=f=g=h=0$, $c=1$, $t=0$, its restrictions to $X_1$ and
$X_2$ are the residue-field quotients at $p$, whereas its derived
restriction to $S$ is $\mathbb C_p\oplus\mathbb C_p[1]$.
It therefore gives actual perfect matching objects at $S$, with a
nontrivial negative-degree boundary term.
\end{proposition}
\begin{proof}
Eliminate $d$ and $b$ by $d=-ca$ and $b=t/c$. The parameters
$(a,c,e,f,g,h)$, with $c$ invertible, are free. Quotienting the
threefold ring by~\eqref{eq:doublepoint-koszul} gives
\[
 \mathcal O_T[x]/(x^2-ax-t/c),
\]
free over the parameter ring with basis $(1,x)$. The three displayed
equations are a regular sequence: first eliminate $y$, then $z,w$;
the remaining monic equation in $x$ gives the relative finite flat
quotient. Alternatively their independent eliminated variables and
the Cohen--Macaulay threefold fibers verify the regular-sequence
criterion. The Koszul complex is therefore a finite free resolution.

At the stated central point the sequence is $(y-x,z,w)$.
On $X_1$ it becomes $(y,z,w)$ and on $X_2$ it becomes $(-x,z,w)$.
These are the regular point ideals of the two smooth components.
On $S$ it is $(0,z,w)$, whose Koszul complex is
$K(z,w)\otimes[\mathcal O_S\xrightarrow{0}\mathcal O_S]$.
Its cohomology is exactly $\mathbb C_p$ in degrees $0$ and $-1$.
Termwise restriction of one finite free complex supplies the matching
isomorphism and proves the asserted derived, not merely numerical,
identifications.
\end{proof}

The ideals are recovered as the fibers of the unit maps into these
quotient complexes. For a surface section containing the length-two
scheme, the actual surface objects are
\[
 E_i=\operatorname{Cone}(L_i^{-m}\to I_{p/X_i}),\qquad
 E_S=\operatorname{Cone}\!\bigl(L_S^{-m}\to
               \operatorname{fib}(\mathcal O_S\to
                         \mathbb C_p\oplus\mathbb C_p[1])\bigr).
\]
The arrows are the restrictions of the chosen section and quotient
map; replacing $E_S$ by a pure curve sheaf would lose part of this
object. Near a point of $S$ away from the blowup center $\Gamma$,
the map to $\mathbb P^4$ is an embedding and both $L$ and $L^2$
separate the displayed length-two quotient. The incidence moduli is
then locally a projective bundle over the smooth six-dimensional
relative Hilbert chart. Its classical dimensions are $8$ and $18$,
respectively. Its reduced potential is zero by the minimal-coordinate
argument in Section~\ref{sec:v6-hyperplanes}; its derived tangent
complex still has the dual obstruction space. The local oriented bar
and connection are the zero-critical calculation in these dimensions.
Thus being supported at $S$ does not by itself force a singular
classical moduli germ. At $b=c=0$, in contrast, the same Hilbert
parameter equation is $bc=t$ and the quotient is not covered by the
regular-sequence argument. The calculation does not erase that stratum.

The normalization identity also checks the charge of this
crossing quotient:
\[
 [Q_0]=[j_{1*}k(p)]+[j_{2*}k(p)]
       -[j_{S*}(k(p)\oplus k(p)[1])]=2[k(p)]
       \quad\text{in }G_0(X_0).
\]
Thus the two point objects on the components, matched through
$k(p)\oplus k(p)[1]$, represent a length-two limit, not a
length-one limit.  Retaining the Tor term is essential, but
does not allow the length-one charge to be assigned to this
particular triple.

\subsection{What enlarging the central fiber can and cannot change}
\begin{proposition}\label{prop:no-polefree-point-completion}
For the natural generic one-point and two-point DT systems in either
$|H|$ or $|2H|$, no finite-rank connection defined without a pole on a
disk about $0$ can restrict to the computed Gauss--Manin system on
the punctured disk. The same statement holds for the line-with-one-point
system. It remains true after a finite ramified base change.
\end{proposition}
\begin{proof}
If $\nabla=d+A(t)dt$ is holomorphic on the disk, the fundamental matrix
solution of $G'=-AG$, $G(0)=I$, is holomorphic and single-valued.
The restriction of this connection consequently has identity monodromy.
The natural systems calculated above have nonzero nilpotent logarithms:
for one point their ranks are $84$ and $294$; for two points their
odd ranks are $378$ and $1638$; for a labeled line and one point the
nonzero ranks are $42$ and $231$ (the latter in even DT parity).
Thus their monodromy is not the identity. Under $t=s^e$ its logarithm
is multiplied by $e$, and remains nonzero in characteristic zero.
This proves the assertion. Changing only the central objects leaves
this generic local system unchanged and cannot remove the obstruction.
\end{proof}

The appropriate finite extension for these systems is logarithmic.
Here it can be written explicitly, rather than inferred from a bare
flatness assertion. Let $N_X=\log T_X$ on $H^3(X_b)$; its square is
zero and its rank is $21$. In a Jordan basis,
\begin{equation}\label{eq:quintic-log-matrix}
 N_X=J_2^{\oplus21}\oplus0_{162},\qquad
 J_2=\begin{pmatrix}0&1\\0&0\end{pmatrix}.
\end{equation}
For a multivalued horizontal basis $v$ of any of the computed systems,
put $e=\exp(-\log(t)N/(2\pi i))v$. This basis is single-valued and
its extension defines
\begin{equation}\label{eq:computed-log-connection}
 \nabla e=-\frac{N}{2\pi i}\frac{dt}{t}\,e.
\end{equation}
This sign convention gives monodromy $\exp N$ for horizontal sections.
The connection is flat, and its residue and all tensor operations are
specified. Its fiber at $0$ is a limiting realization with monodromy,
not automatically the DT group of an ordinary central moduli space.

For one point, $N$ is $r_m-1$ copies of $N_X$ on the odd part and
zero on the even part. For two points, let
$V=H^*(X_b,\mathbb C)$ with its cohomological grading and take the
operator $N_V$ equal to $N_X$ on $H^3$ and zero otherwise. On the
Hilbert-square cohomology the operator is
\[
 (N_V\otimes1+1\otimes N_V)|_{\operatorname{Sym}_{\mathrm{gr}}^2V}
       \ \oplus\ N_V[-2]\ \oplus\ N_V[-4].
\]
Take $r_m-2$ copies for the projective-bundle classes. This supplies
the whole logarithmic matrix, not just its odd rank. For example,
on $\bigwedge^2H^3(X_b)$ it has Jordan blocks
\[
 J_3^{\oplus210}\oplus J_2^{\oplus3402}\oplus0_{13272}.
\]
Indeed write $H^3=\mathbb C^{21}\otimes\mathbb C^2\oplus\mathbb C^{162}$,
with the usual rank-one nilpotent on $\mathbb C^2$.
The decomposition of its exterior square has
$\binom{21}{2}$ copies of $\operatorname{Sym}^2\mathbb C^2$,
$21\cdot162$ copies of $\mathbb C^2$, and
$\binom{22}{2}+\binom{162}{2}$ trivial factors.
The remaining eighteen even Hilbert-square classes are constant.
Thus the even ranks of $N$ for two points are $11466$ and $49686$,
and the ranks of $N^2$ are $630$ and $2730$, in degrees one and two,
respectively. These statements refer to the two choices of the
supporting linear system, not to cohomological degree.

The derived matching construction and the point replacements above
supply actual central objects. To identify the logarithmic limiting
realization~\eqref{eq:computed-log-connection} with an intrinsic central
DT complex on a global chosen moduli space, one still needs the
oriented specialization morphism and its vanishing-cycle cone in
Proposition~\ref{prop:specialization-cone-test}. That cone is the mathematical
quantity to calculate. Flatness of~\eqref{eq:computed-log-connection}
and a local perfect resolution cannot set it to zero. This limitation
follows from the displayed monodromy matrices, not from excluding
nontransverse intersections from the theory.

\subsection{The singular stratum of the length-two chart}
\label{ssec:crossing-exhaustion}
The equation $bc=t$ in \eqref{eq:doublepoint-rules} is a
parameter-space equation.  Its singular central locus cannot be
replaced by a smooth open without examining its universal object.
The next calculation examines the whole chart.

\begin{proposition}\label{prop:crossing-perfect-exhaustion}
At a geometric point of the central length-two chart
\eqref{eq:doublepoint-rules}, the quotient is perfect on
$X_0$ if and only if $(b,c)\ne(0,0)$.  Thus the two perfect opens
are $c\ne0$ and $b\ne0$.  The central parameter stratum $b=c=0$
consists of nonperfect quotients; it is not a deleted singular
stratum of the moduli of perfect finite quotient sheaves.
A surface-section cone formed from one of those nonperfect ideals
is likewise nonperfect.  This holds in both divisor degrees.
\end{proposition}
\begin{proof}
For $c\ne0$, Proposition~\ref{prop:crossing-two-perfect} supplies the
finite Koszul resolution.  If $b\ne0$, the equation
$x^2-ax-b=0$ makes $x$ invertible on the quotient.  Its support is
therefore contained in the smooth branch $x\ne0$ of $X_0$.
A coherent finite-length module on this smooth threefold has a
finite locally free resolution, so is perfect.

Now let $b=c=0$, so $d=0$.  Eliminate the two independent equations
$z-ex-f$, $w-gx-h$.  After translating the remaining coordinates the threefold ring
is $R'[[z',w']]$, with $R'=\CC[[x,y]]/(xy)$, and the module is
annihilated by $z',w'$.  Restricting a hypothetical finite free
resolution to $R'$ would give a bounded flat resolution, since
$R'[[z',w']]$ is flat over $R'$.  The module is finite over the
noetherian local ring $R'$, so finite flat dimension would imply
finite projective dimension there.  At a crossing support point the module is
\[
 M=R'/(y,x(x-a)).
\]
If $a\ne0$, localizing at $x=0$ gives the residue field of the
nodal curve.  If $a=0$, it gives $R'/(y,x^2)$.
Resolve $R'/(y)$ by the infinite two-periodic sequence with
alternating differentials $y,x$.  After tensoring with either
module, the $y$ maps are zero and the $x$ maps are nilpotent.
For $a=0$ their kernel and cokernel each have dimension one;
for the residue field both differentials are zero.  Thus
$\operatorname{Tor}^{R'}_j(R'/(y),M)$ is nonzero for arbitrarily
large $j$.  A module of finite projective dimension would make
these Tor groups eventually zero.  This contradiction proves
nonperfectness at the crossing; the other support point, when
$a\ne0$, does not change it.

The triangle $I_Z\to\O_{X_0}\to\O_Z$ shows that $I_Z$ is
nonperfect as well.  If the cone of $L^{-m}\to I_Z$ were perfect,
then $I_Z$ would be perfect as an extension of two perfect
objects.  This is impossible.  Changing $m$ changes the space of
sections but not this local obstruction.
\end{proof}

This classifies the ordinary quotient sheaves on this entire chart.
It does not exclude perfect \emph{complex} limits with additional
cohomology: Proposition~\ref{prop:crossing-one-perfect} constructs
one explicitly.  Such an object has a different local deformation
problem from the nonperfect quotient above.  Replacing its complex
by the nonperfect degree-zero quotient would not compute a point of
$\RPerf(X_0)$.

\begin{lemma}\label{lem:node-not-critical-truncation}
The reduced germ $\{bc=0\}\times\A^r$ at a point of $b=c=0$
is not the classical critical locus of a normalized function on a
minimal smooth ambient germ.  In particular, the equation of the
ordinary length-two parameter space does not by itself provide the
missing $(-1)$-shifted critical enhancement.
\end{lemma}
\begin{proof}
Its minimal embedding ideal is the principal ideal $(bc)$ in the ring
\[
 A=\CC\{b,c,z_1,\ldots,z_r\},
\]
or in its formal completion.
An isomorphism with another minimal ambient presentation lifts by
lifting a cotangent-space basis and applying the inverse function
theorem.  A normalized critical function must vanish on the
reduced germ, so write $f=bc\,h$.  If $df\in(bc)\Omega_A^1$,
then reduction of $\partial_bf$ modulo $b$ gives $ch|_{b=0}=0$;
hence $b\mid h$.  Similarly $c\mid h$.  Thus $f\in(bc)^2$.
At the crossing every partial derivative belongs to
$\mathfrak m(bc)$, where $\mathfrak m$ is the maximal ideal.
By Nakayama these derivatives cannot
generate $(bc)$.  This contradicts the defining gradient ideal
of a critical chart.  Removing nondegenerate quadratic pairs
reduces any nonminimal critical chart to such a minimal chart,
so adding those pairs cannot repair the contradiction.
\end{proof}

For comparison, the nonzero function $f=(bc)^2/2$ has
\[
 \operatorname{Crit}(f)=\operatorname{Spec}
       A/(bc^2,b^2c),
\]
whose reduced support is the node, but whose classical scheme has
additional nilpotents at the crossing.  It is not the reduced
Hilbert parameter $\operatorname{Spec}A/(bc)$.  Its derived critical
enhancement cannot be assigned to that different parameter space
without changing the moduli problem.

The distinction is substantial.  A symmetric perfect obstruction
theory must be induced from the specified object deformation
problem; it cannot be supplied merely by naming a singular
parameter equation.  The genuine singular models in
Sections~\ref{sec:singular-comparison}, \ref{sec:cubic-geometry}
and~\ref{sec:global-two-defect} have nonzero critical functions and
are not subject to this obstruction.

\subsection{Why a bounded completion must specify its objects}
\begin{proposition}\label{prop:all-perfect-limits-unbounded}
Even a fixed generic perfect object, a fixed fiberwise numerical
class, and a fixed Tor-amplitude interval do not define a bounded
space of all perfect limits.  More precisely, over a coordinate disk
there are relative perfect complexes with the same generic fiber
and the same fiberwise $K$-class whose central endomorphism dimensions
are arbitrarily large.
\end{proposition}
\begin{proof}
Let $E$ be a given relative perfect complex, and put
\[
 T=[\O_{\mathcal X}\xrightarrow{t}\O_{\mathcal X}]
       \quad\text{in degrees }-1,0,\qquad
 E^{(n)}=E\oplus T^{\oplus n}.
\]
The complex $T$ is acyclic on every nonzero fiber.  Its derived
central restriction is $\O_{X_0}[1]\oplus\O_{X_0}$, so
\[
 E^{(n)}_\eta\simeq E_\eta,\quad
 [E^{(n)}_b]=[E_b]\text{ in }K_0(X_b),\quad
 E^{(n)}_0\simeq E_0\oplus
       (\O_{X_0}[1]\oplus\O_{X_0})^{\oplus n}.
\]
The amplitude lies in the fixed union of the amplitude of $E$ and
$[-1,0]$, independently of $n$.  The last display gives at least
$n^2$ independent central endomorphisms.  In a finite-type family
of perfect objects on the proper $X_0$, the derived pushforward of
the universal endomorphism complex is perfect locally on the
parameter.  On a finite affine cover its finite resolutions give
a uniform bound on the dimensions of its degree-zero fibers.
The displayed unbounded dimensions contradict such a bound.
\end{proof}

The proposition does not rule out a proper completion with a
specified simplicity or stability condition.  It shows why the
phrase ``all perfect limits'' without that condition cannot replace
a moduli definition.  The extra summands in the proof are not
simple and hence do not belong to the simple-object components
calculated in Sections~\ref{sec:v6-hyperplanes} and
\ref{sec:cubic-geometry}.  A proposed enlargement of the point
problem must specify which replacement complexes it includes and
prove boundedness and orientation compatibility for that selection.
This is a mathematical restriction on a finite DT invariant, not a
transversality assumption on the derived Lagrangian intersection.

\section{A proper singular comparison from a tangential central curve}
\label{sec:cubic-geometry}
We construct an example in the same quintic degeneration for which the
\emph{generic and central classical moduli are both nonreduced}, the
intrinsic potential is cubic rather than zero, and the entire DT group
of a specified proper moduli locus is calculated. The degeneration is
chosen along a discriminant in the coefficient space. Its generic
quintics are smooth, but their line scheme has a chosen double point.
Thus this example is not subject to the transverse, reduced-line
hypothesis of Proposition~\ref{prop:resolved-line-family}. Only the
specified local component of the line scheme is used, not all lines
on every smooth quintic.

\subsection{An explicit tangency in the quartic and blowup geometry}
Use homogeneous coordinates $[u:v:x_2:x_3:x_4]$ on $\mathbb P^4$,
and put
\[
 h=u,\qquad \Lambda=(x_2=x_3=x_4=0),\qquad
 p=[0:1:0:0:0].
\]
Fix $c\in\mathbb C^*$ and choose forms
\begin{align}
 q_4&=x_2u^3+x_3u^2v+x_4v^3+x_2^2uv+R_4,
                                      \label{eq:cubic-quartic}\\
 q_5&=x_2u^4+x_3v^4-x_4v^4+c x_2^2v^3+R_5,
                                      \label{eq:cubic-quintic}
\end{align}
where $R_i\in I_\Lambda^3(i)$. Take the remaining coefficients in
these two linear systems sufficiently general that $X_1=(q_4=0)$,
$S=(u=q_4=0)$, and $\Gamma=(u=q_4=q_5=0)$ are smooth. Such choices
exist: outside $\Lambda$ the systems with these fixed jets are generated
without a base point, so Bertini applies. Along $\Lambda$ the three
normal derivatives of $q_4$ are $u^3,u^2v,v^3$, with no common zero.
The only point of $\Lambda$ on $u=0$ is $p$, and there $q_4$ and
$q_5$ have independent linear normal terms $x_4$ and $x_3-x_4$.
These checks prove smoothness at the constrained points as well.

The differential of the quartic line equation at $\Lambda$ is
\[
 H^0(\mathcal O_{\mathbb P^1}(1))^{\oplus3}
       \longrightarrow H^0(\mathcal O_{\mathbb P^1}(4)),
 \quad(\ell_2,\ell_3,\ell_4)
       \longmapsto u^3\ell_2+u^2v\ell_3+v^3\ell_4.
\]
Its image contains all five monomials of degree four. Its kernel
is generated by $(v,-u,0)$. Thus the local quartic line scheme
$T_4$ is a smooth curve. A parameter $\lambda$ on it can be chosen
so that its universal line has equations, to second order,
\begin{equation}\label{eq:cubic-line-jet}
 x_2=\lambda v+O(\lambda^3),\qquad
 x_3=-\lambda u+O(\lambda^3),\qquad
 x_4=-\lambda^2u+O(\lambda^3).
\end{equation}
Indeed substitution into $q_4$ cancels the first-order terms and
cancels $\lambda^2uv^3$ against the displayed $x_4v^3$ term;
$R_4$ has order at least three in $\lambda$.

On $S$ near $p$, use $\xi=x_2/v$, $\eta=x_3/v$.
The equation $q_4=0$ eliminates $x_4/v$ to order at least three.
The evaluation map $T_4\to S$, $\Lambda_\lambda\mapsto
\Lambda_\lambda\cap S$, therefore has expansion
\[
 \xi=\lambda+O(\lambda^3),\qquad \eta=O(\lambda^3).
\]
The curve $\Gamma\subset S$ has equation
\[
 \eta=-c\xi^2+O(\xi^3).
\]
The two smooth curve germs have intersection multiplicity exactly
two at $p$. Both are Lagrangian in the K3 surface $S$ with its
holomorphic symplectic form. They are the actual restrictions of
ideal sheaves of a quartic line and of an exceptional fiber in
$X_2=\operatorname{Bl}_\Gamma\mathbb P^3$, respectively.

\begin{proposition}\label{prop:cubic-relative-curve}
After choosing an analytic disk $\Delta$ about $0$, there is a
holomorphic function $\varepsilon(t)$, with $\varepsilon(0)=0$, such
that the resolved family
\begin{equation}\label{eq:cubic-resolved-family}
 \mathcal X=\operatorname{Bl}_{(u,t)}
       \{u q_4=t(q_5+\varepsilon(t)v^5)\}\longrightarrow\Delta
\end{equation}
is semistable with smooth quintic fibers for $t\ne0$, and has an
isolated proper family of perfect curves with classical parameter
scheme
\[
 \mathcal T\simeq\operatorname{Spec}_\Delta
                  \bigl(\mathcal O_\Delta[\lambda]/(\lambda^2)\bigr).
\]
The reduced central curve is $C_0=\Lambda\cup_p e_p$, where $e_p$ is
the exceptional fiber in $X_2$. It is a regular codimension-two
immersion into $X_0$. Its actual local Hilbert deformation problem
has one tangent and one obstruction direction, and completed local
ring $\mathbb C[[\lambda]]/(\lambda^2)$.
\end{proposition}
\begin{proof}
We construct the discriminant path using the equations of lines, so
that its double scheme structure is controlled. On the Grassmannian
chart near $\Lambda$ write a line as
\[
 x_2=A_2u+B_2v,\qquad x_3=A_3u+B_3v,\qquad
 x_4=A_4u+B_4v.
\]
Its six coordinates are the six displayed coefficients. The line
condition for the unresolved quintic is the vanishing of the six
coefficients of $u q_4-t(q_5+\varepsilon v^5)$ on this graph.
At $t=0$, the five coefficients involving a positive power of $u$
have rank five, by the preceding quartic calculation. Solve them
for five of the Grassmannian coordinates, keeping $B_2=\lambda$
as the sixth. The analytic implicit function theorem applies.
These five equations do not involve $\varepsilon$, since
$\varepsilon v^5$ changes only the last coefficient.

After this elimination the remaining equation is
$-t(g(\lambda,t)+\varepsilon)=0$. On the closure of the nonzero-fiber
line locus it is $g(\lambda,t)+\varepsilon=0$. By
\eqref{eq:cubic-line-jet},
\[
 g(\lambda,0)=c\lambda^2+O(\lambda^3),\quad
 g(0,0)=g_\lambda(0,0)=0,\quad g_{\lambda\lambda}(0,0)=2c\ne0.
\]
Solve $g_\lambda(\lambda_c(t),t)=0$ and set
$\varepsilon(t)=-g(\lambda_c(t),t)$.
The parameter-dependent one-variable Morse lemma, or division of
$g-g(\lambda_c(t),t)$ by $(\lambda-\lambda_c(t))^2$, gives a unit
$u_0(\lambda,t)$ times that square. After taking its analytic square
root and recentering, the saturated line equation is exactly
$\lambda^2=0$ on every fiber. This proves the asserted finite flat
parameter scheme and its double, not reduced, structure.

Pull back the ideal $(u,t)$ to its universal line. Blowing it up
lifts that line family to~\eqref{eq:cubic-resolved-family}. In the
chart $t=u\tau$ its surface of curves has base map $t=u\tau$;
the other two transverse coordinates of the ambient resolved
threefold are functions of $u,\tau,\lambda$. Their graph equations
are two independent regular equations. The other blowup chart has
the same description and supplies the second chart of $e_p$.
This proves relative flatness and the codimension-two regular
immersion, including after base change to the nonreduced
$\mathcal T$. The central curve has two branches $\Lambda$ and
$e_p$, meeting at $p$; their $L$-degrees are $1$ and $0$.
The resolution has the same local charts and canonical-bundle
identity as Section~\ref{ssec:quintic-construction}, because
$q_5+\varepsilon(0)v^5=q_5$. Smoothness of $X_1,S,\Gamma$ makes
the total resolved space smooth along its central fiber. The local
base equation is $xy=t$. Properness and shrinking $\Delta$ then
ensure that all its nonzero fibers are smooth.

It remains to check that the curve parameter is its Hilbert
moduli, not just a double covering of an unexplained image. At $C_0$
the normal matching sequence is
\[
 0\to H^0(N_{C_0/X_0})\to
 H^0(N_{\Lambda/X_1})\oplus H^0(N_{e_p/X_2})
 \longrightarrow T_pS\to H^1(N_{C_0/X_0})\to0.
\]
The two middle spaces are one-dimensional. Their maps are the
derivatives of the evaluation curve $T_4\to S$ and the inclusion
$\Gamma\hookrightarrow S$. In bases $(\partial_\xi,\partial_\eta)$
their difference has matrix
\[
 \begin{pmatrix}1&-1\\0&0\end{pmatrix}.
\]
It has one-dimensional kernel and cokernel. More strongly, the
embedded deformation functor of this regular curve is the matching
of its two component curve deformations: locally the two graph
equations above patch if and only if their values on $S$ agree.
Conversely restriction of a regular-immersion deformation recovers
these graphs and this equality, on every Artinian base. \\
The component curve schemes are the smooth germs $T_4$ and $\Gamma$;
their actual equality equations on $S$ have order-two contact as
computed above. Their fiber product has completed ring
$\mathbb C[[\lambda]]/(\lambda^2)$. The lifted universal line
induces this same equation, not merely the same tangent dimension.
Thus it is the complete local embedded deformation family.
For $t\ne0$, the identical statement follows directly from the six
Grassmannian equations and their five eliminated coordinates.
The regular curve families and these local identifications identify
the selected relative Hilbert germ over $\Delta$.

An isolated finite part of the proper line scheme is open and closed
in a sufficiently small analytic neighborhood of this family. Its
support has no additional branch by the eliminated equation.
Taking this part and then its lifted curves makes the specified
parameter family proper over $\Delta$.
\end{proof}

The proposition can also be read in the henselian or formal
neighborhood of the coefficient point. The analytic disk formulation
is used here so that its ordinary critical charts and parallel
transport have their usual analytic meaning. It changes coefficients
of the quintic family, not the central quartic and blowup geometry.

\subsection{Surface sheaves and their nonzero potential}
For $m=1,2$, put $r_1=5$, $r_2=15$ and
\[
 c_m=r_m-(m+1),\qquad d_m=c_m-1;
 \quad(c_1,d_1)=(3,2),\quad(c_2,d_2)=(12,11).
\]
Let $\mathcal C\subset\mathcal X\times_\Delta\mathcal T$ be the
universal curve of Proposition~\ref{prop:cubic-relative-curve}.

\begin{theorem}\label{thm:cubic-actual-moduli}
The bundle $W_m=\pi_{\mathcal T*}(I_{\mathcal C}L^m)$ is locally
free of rank $c_m$, commutes with base change, and has no higher
direct images. The universal cones
\[
 E_{C,s}=\operatorname{Cone}(L^{-m}\xrightarrow{s}I_C)
\]
define an open proper classical scalar-rigidified perfect-complex
moduli locus $M_m=\mathbb P_{\mathcal T}(W_m)$ near these objects.
It is locally over $\Delta$ the scheme
$\mathbb P^{d_m}\times\operatorname{Spec}\mathcal O_\Delta[\lambda]/(\lambda^2)$,
and is singular and nonreduced on every fiber.
Its induced reduced d-critical structure has a global presentation
\begin{equation}\label{eq:cubic-global-potential}
 U_m=\mathbb P^{d_m}\times\mathbb A^1_\lambda\times\Delta,
 \qquad f_m=\frac{a(t)}3\lambda^3,\qquad a(0)\ne0.
\end{equation}
After shrinking the disk, $a(t)$ is a unit. The canonical bundle of
$U_m/\Delta$ gives a compatible orientation. In the central fiber the
same presentation is the reduced critical realization of the actual
component-moduli Lagrangian intersection at the restrictions of
$E_{C,s}$. The potential is nonzero modulo the square of the critical
ideal.
\end{theorem}
\begin{proof}
On a smooth line, $L^m$ has $m+1$ sections and no $H^1$.
On $C_0=\Lambda\cup e_p$, normalization gives
\[
 0\to H^0(L^m|_{C_0})\to
 H^0(\mathcal O_\Lambda(m))\oplus H^0(\mathcal O_{e_p})
       \to\mathbb C_p\to H^1(L^m|_{C_0})\to0.
\]
The last evaluation is onto. Thus $h^0=m+1$ and $H^1=0$.
The restrictions of degree-$m$ polynomials supply every section on
$\Lambda$, with the forced constant value on $e_p$. The ideal
sequence and whole-fiber cohomology give $H^{>0}(I_CL^m)=0$ and
$h^0(I_CL^m)=r_m-(m+1)$ on every fiber. Cohomology and base change
apply also to the infinitesimal curve parameters, giving $W_m$.
Moreover $H^0(\mathcal O_C)=\mathbb C$ and
$H^1(\mathcal O_C)=0$, again by normalization at $0$.

Apply the reconstruction following
Lemma~\ref{lem:quadric-quot-evaluation} to a deformation of $E_{C,s}$.
The one-dimensional space $\operatorname{Ext}^1(E,L^{-m})$ recovers
$L^{-m}\to I_C\to E$; the one-dimensional
$\operatorname{Hom}(I_C,\mathcal O)$, with zero
$\operatorname{Ext}^1$, recovers $I_C\to\mathcal O\to\mathcal O_C$.
The quotient is a flat curve deformation. Proposition~\ref{prop:cubic-relative-curve}
identifies its full formal parameter as $\mathcal T$, including its
nilpotents. The remaining choice is the projective section in $W_m$.
This proves the asserted moduli identification and simplicity, not
just the existence of a parameter map. Properness follows from the
projective bundle over the finite flat proper $\mathcal T$.
Over an analytic disk, a vector bundle on this trivial first-order
thickening is trivial after shrinking: its restriction to the disk
is trivial and the lifting obstruction and ambiguity are coherent
$H^1$ on the Stein disk. Thus the displayed product description of
$M_m$ is available.

We determine its intrinsic potential. At a point, minimal ambient
coordinates are the $d_m$ projective chart coordinates and
$\lambda$. Its critical ideal is $(\lambda^2)$, not $(0)$.
A normalized local d-critical representative satisfies
$f|_{\lambda=0}=0$ and $df\in(\lambda^2)\Omega^1$.
Expanding in $\lambda$ shows that its class modulo
$(\lambda^2)^2=(\lambda^4)$ is
\[
 a(z,t)\lambda^3/3\pmod{\lambda^4}.
\]
The coefficient is a unit: otherwise a minimal critical potential
would have no cubic term at that point, and its derivatives could
not generate the quadratic ideal $(\lambda^2)$.
These classes are induced by the CY cyclic deformation pairing,
not a prescribed zero structure. Their cubic coefficient is the
pairing of the quadratic obstruction with the remaining tangent
coordinate. Composition, homotopy transfer to this finite tangent
space, and the relative trace make this coefficient holomorphic in
$t$ on each chart. Only this cubic coefficient is needed modulo
$(\lambda^4)$; no convergence assertion for an arbitrary infinite
formal potential is used here. On the global product
embedding of $M_m$ in $U_m$, their coefficients glue to a holomorphic
function on $\mathbb P^{d_m}\times\Delta$. Such a function is
constant along the projective factor, hence equals $a(t)$.
The representative~\eqref{eq:cubic-global-potential} has exactly the
critical scheme $M_m$ and this intrinsic d-critical class, proving
the global critical presentation. Its class is nonzero because
$\lambda^3\notin(\lambda^4)$.

Finally the central ideals have restrictions
\[
 I_{C_0}|_{X_1}^L=I_{\Lambda/X_1},\qquad
 I_{C_0}|_{X_2}^L=I_{e_p/X_2},\qquad
 I_{C_0}|_S^L=I_{p/S}.
\]
The equalities follow from the two relative regular graph equations;
no Tor terms have been deleted. The cones restrict to
\begin{equation}\label{eq:cubic-actual-matching}
 E_i=\operatorname{Cone}(L_i^{-m}\to I_{C_i/X_i}),\quad
 E_S=\operatorname{Cone}(L_S^{-m}\to I_{p/S}),\quad
 (C_1,C_2)=(\Lambda,e_p).
\end{equation}
Their homotopy matching is induced by one universal complex.
Normalization identifies its endomorphism deformation algebra with
that of~\eqref{eq:cubic-actual-matching}. In the curve directions the
actual component ideal-moduli are $T_4$ and $\Gamma$, restricting
to $S$ as above. In Darboux coordinates on $S$, straighten the first
curve to the zero section. The second is the graph of a one-form
whose first nonzero term is a nonzero multiple of $\lambda^2d\lambda$.
Its primitive has a cubic leading term. The actual component-chart bar can be written before passing to
any critical normal form. In Darboux coordinates let the two
restriction maps be
\[
 \lambda\longmapsto(\xi,\eta)=(\lambda,\alpha(\lambda)),\qquad
 \mu\longmapsto(\xi,\eta)=(\mu,\beta(\mu)).
\]
Their difference has $\alpha(\lambda)-\beta(\lambda)
=c_0\lambda^2+O(\lambda^3)$, $c_0\ne0$. The structural two-sided
bar over $\mathbb C[[\xi,\eta]]$ has the Koszul representative
\begin{equation}\label{eq:cubic-component-Koszul}
 \begin{gathered}
 \mathbb C[[\lambda,\mu]][\epsilon_\xi,\epsilon_\eta],
 \qquad|\epsilon_\xi|=|\epsilon_\eta|=-1,\\
 d\epsilon_\xi=\lambda-\mu,\qquad
 d\epsilon_\eta=\alpha(\lambda)-\beta(\mu).
 \end{gathered}
\end{equation}
For example this follows by resolving the diagonal in
$S\times S$ by its two coordinate differences and pulling back
along the two displayed restriction maps. Eliminate the regular
pair $(\lambda-\mu,\epsilon_\xi)$; the remaining algebra is
\[
 \mathbb C[[\lambda]][\epsilon],\qquad
 d\epsilon=\alpha(\lambda)-\beta(\lambda).
\]
Its cyclic primitive is the actual Lagrangian-intersection potential
$\int_0^\lambda(\alpha(v)-\beta(v))\,dv$, with leading term
$c_0\lambda^3/3$. This computes a nontransverse component-moduli
intersection, not an unrelated zero-section example. The
surface-section reconstruction adds the projective coordinates;
the preceding intrinsic-class calculation identifies its reduced
critical structure with~\eqref{eq:cubic-global-potential}.
This also fixes the central orientation through the common canonical
orientation of $U_m$. The smooth-Lagrangian graph calculation is
compatible with the quantization/perverse-sheaf comparison of
\cite{GunninghamSafronov}; it is being applied to these actual local
curve restriction maps, not to an arbitrary smoothness assertion
about the full component stacks.
\end{proof}

The classical dimension of a fiber is $d_m$, whereas its tangent
space has dimension $d_m+1$. The reduced derived tangent complex
has this tangent space in degree zero and its dual in degree one,
with the differential given by the Hessian of $f_m$.
At $\lambda=0$ that Hessian vanishes, but the quadratic obstruction
map $a(t)\lambda^2$ is nonzero. This is a genuinely singular
quasi-smooth critical enhancement, not a smooth scheme with an
artificial quadratic stabilization.

\subsection{Cyclic cohomology, central bar, and the finite connection}
\begin{theorem}\label{thm:cubic-comparison}
For the proper loci $M_m$ of Theorem~\ref{thm:cubic-actual-moduli},
$m=1,2$, the generic and central DT groups are
\begin{equation}\label{eq:cubic-DT-groups}
 \mathrm{DT}_{\mathrm{coh}}^j(M_{m,b})
   =\begin{cases}
       \mathbb C^2,&j=2a-d_m,\quad 0\leq a\leq d_m,\\
       0,&\text{otherwise}.
     \end{cases}
\end{equation}
The periodic realization has rank six in even parity for $|H|$, and
rank twenty-four in odd parity for $|2H|$. These are the groups of
the chosen double-line component, not of the entire line scheme.
The actual central matching locus has an oriented two-sided bar/BV
realization horizontally quasi-isomorphic to the normalized periodic
cyclic complex of $\operatorname{MF}(U_{m,b},f_{m,b})$.
Its finite Gauss--Manin connection is regular and flat across $0$.
In the basis of the two one-variable classes and projective
hyperplane powers, its nonconstant block is
\begin{equation}\label{eq:cubic-connection-matrix}
 \nabla_{\partial_t}
    =\partial_t-\frac{a'(t)}{3a(t)}
                    \begin{pmatrix}1&0\\0&2\end{pmatrix},
\end{equation}
repeated $d_m+1$ times. This proves the finite regular comparison
with both geometric endpoints identified, in a nonzero critical
example.
\end{theorem}
\begin{proof}
The vanishing-cycle sheaf of $\lambda^3/3$ at its unique critical
point is the reduced cohomology of three points: it is
$\mathbb C^2$ in perverse degree zero. A choice of cube root of the
unit $a(t)$ identifies this family with the same function for every
$t$. Smooth pullback to the projective factor and the canonical
orientation give
\[
 \mathcal P_{M_{m,b}}=\mathbb C_{\mathbb P^{d_m}}[d_m]
                           \otimes\mathbb C^2.
\]
The underlying topological critical support is proper. Taking its
hypercohomology proves~\eqref{eq:cubic-DT-groups} and the asserted
parities. In particular the cubic critical structure contributes
two classes at every projective cohomology degree; replacing it by
zero would give only one and would be incorrect.

The cyclic calculation is equally explicit. In the $\lambda$ direction
use
\[
 \mathbb C[\lambda]((u))\xrightarrow{\ u\partial_\lambda-a(t)\lambda^2\ }
                        \mathbb C[\lambda]((u))d\lambda.
\]
Its cohomology has basis $e_0=[d\lambda]$, $e_1=[\lambda d\lambda]$
in degree one. One proof works first $u$-adically: modulo $u$ the
map is multiplication by $-a\lambda^2$, injective with cokernel
basis $1,\lambda$. Successive division by $\lambda^2$ lifts this
basis and gives no kernel at any order. Completeness and inversion
of $u$ give the displayed answer. The explicit reduction identity is
\begin{equation}\label{eq:cubic-reduction}
 [\lambda^{k+3}d\lambda]
       =\frac{(k+1)u}{a(t)}[\lambda^k d\lambda],\qquad k\geq0,
 \qquad[\lambda^2d\lambda]=0.
\end{equation}
Tensor with the projective de Rham complex. Its hyperplane powers
supply $d_m+1$ independent even classes. The DT normalization is
the ambient shift $[d_m+1]$, giving parity $d_m$ and the computed
ranks. Efimov's comparison~\cite{Efimov2018CyclicMF} identifies
this with the specified normalized periodic cyclic realization.

For the central bar, on a projective coordinate chart with variables
$z_1,\ldots,z_{d_m}$ set
\[
 R=\mathcal O_\Delta[z_1,\ldots,z_{d_m},\lambda],\qquad
 A=R[p_1,\ldots,p_{d_m},p_\lambda].
\]
The two augmentations are $p_j,p_\lambda\mapsto0$, and
$p_j\mapsto0$, $p_\lambda\mapsto a(t)\lambda^2$.
Their two-sided bar retracts onto
\[
 K=R[\eta_1,\ldots,\eta_{d_m},\eta_\lambda],\quad
 \delta\eta_j=0,\quad\delta\eta_\lambda=a(t)\lambda^2.
\]
Thus the central structure differential already records the nonzero
quadratic matching equation. Its oriented periodic differential is
\begin{equation}\label{eq:cubic-BV}
 D_K=-a(t)\lambda^2\partial_{\eta_\lambda}
       +u\left(\sum_j\partial_{z_j}\partial_{\eta_j}
                        +\partial_\lambda\partial_{\eta_\lambda}\right).
\end{equation}
The bar projection, antisymmetrization, and contraction by top forms
of Section~\ref{sec:v5-bar-bridge} compare this with
$(\Omega_{U_m/\Delta}^\bullet((u)),u d-df_m\wedge)[d_m+1]$.
The maps respect the orientation and the projective coordinate
changes. They therefore glue and specialize. Theorem~\ref{thm:cubic-actual-moduli}
identified this critical chart with the induced realization on the
actual matching locus~\eqref{eq:cubic-actual-matching}; this is the
geometric input which would be absent in an unrelated polynomial
singularity example.

The chain connection before taking cohomology is
\[
 \nabla_{\partial_t}=\partial_t-\frac{a'(t)}{3u}\lambda^3.
\]
Its commutator with $u d-a\lambda^2d\lambda\wedge$ is zero:
the two terms are $-a'\lambda^2d\lambda\wedge$ and
$+a'\lambda^2d\lambda\wedge$. The transferred bar connection of
Proposition~\ref{prop:v5-horizontal-bar} satisfies the same chain
identity and intertwines the comparison maps up to the specified
homotopy. Applying~\eqref{eq:cubic-reduction} for $k=0,1$ gives
\[
 \nabla e_k=-\frac{k+1}{3}\frac{a'}a e_k,
\]
which is~\eqref{eq:cubic-connection-matrix}. Projective hyperplane
classes are horizontal, so no further connection block appears.
The entries are holomorphic because $a$ is a unit. On the simply
connected disk the classes $a(t)^{(k+1)/3}e_k$ are horizontal.
Their ratios to the values at $0$ give the actual generic-to-central
transport; the same representatives lift through the central bar
maps. This constructs a finite free complex with regular chain
comparison, rather than deducing invariance merely from a zero
curvature statement. It verifies the hypotheses of
Theorem~\ref{thm:finite-cyclic} for this singular locus.
\end{proof}

\subsection{An identified higher-contact component}
\label{ssec:higher-contact}
The distinction between a general and a special coefficient choice
can also be recorded for contacts of order greater than two.
The following statement computes the invariant when the entire
classical component and its intrinsic critical structure have been
identified.  Proposition~\ref{prop:cubic-relative-curve} verifies that
geometric input for contact order two; an arbitrary higher contact
is not inferred just from a first-order tangent calculation.

\begin{proposition}\label{prop:higher-contact-critical}
Let $\mu\geq2$, let $\Delta$ be a sufficiently small analytic disk,
and suppose an actual relative perfect-complex moduli component,
including its central matching identification, has classical space
\[
 M=\P^d\times\operatorname{Spec}_\Delta
          \bigl(\O_\Delta[\lambda]/(\lambda^\mu)\bigr).
\]
Assume it has the compatible relative reduced d-critical structure
induced by the cyclic deformation problem.  Then its normalized
critical structure has a global presentation, near $M$,
\[
 U=\P^d\times\A^1\times\Delta,\qquad
 f=\frac{a(t)}{\mu+1}\lambda^{\mu+1},\qquad a(t)\ne0,
\]
after a change of the transverse coordinate preserving $M$.
For the canonical orientation of $U/\Delta$, both endpoint DT
groups are $\CC^\mu$ in degrees $2j-d$, $0\leq j\leq d$, and
zero in the other degrees.  The periodic connection has blocks
\[
 \nabla_{\partial_t}=\partial_t-
     \frac{a'(t)}{(\mu+1)a(t)}
        \operatorname{diag}(1,2,\ldots,\mu).
\]
Its bar comparison includes the actual equation $\lambda^\mu=0$
and is regular and horizontal on the whole component.
\end{proposition}
\begin{proof}
In the fixed product embedding the critical ideal is $(\lambda^\mu)$.
The condition $df\in(\lambda^\mu)\Omega^1_{U/\Delta}$ and the
normalization at $\lambda=0$ show that its d-critical class is
represented modulo $\lambda^{2\mu}$ by
\[
 f=\sum_{k=\mu+1}^{2\mu-1}a_k(z,t)\lambda^k.
\]
Its first coefficient is a unit: otherwise at a zero of that
coefficient the gradients could not generate the degree-$\mu$
initial term of the ideal.  The product embedding makes all the
coefficients global on $\P^d\times\Delta$.  They are independent
of $z$, because a holomorphic function on $\P^d$ is constant.
The displayed representative itself has gradient ideal
$(\lambda^\mu)$ near $M$.  Write it as
$a(t)\lambda^{\mu+1}v(\lambda,t)/(\mu+1)$, with $v(0,t)=1$.
The change $\lambda'=\lambda v(\lambda,t)^{1/(\mu+1)}$ is
analytic and invertible after shrinking.  It preserves the
fat-point ideal and gives the claimed normal form.  This proves
more than an equality of tangent dimensions.

The Milnor fiber in the transverse direction consists of
$\mu+1$ points; its reduced cohomology has dimension $\mu$.
The perverse DT sheaf is therefore
$\CC_{\P^d}[d]\otimes\CC^\mu$, with the chosen orientation.
For the periodic calculation the one-variable differential is
$u\partial_\lambda-a\lambda^\mu$.  Its cokernel is free on
$e_k=[\lambda^k d\lambda]$, $0\leq k<\mu$, and its kernel is zero.
This follows by division modulo $u$, where the cokernel has precisely
that basis, followed by complete successive lifting.  The reduction
formula is
\[
 [\lambda^{k+\mu+1}d\lambda]
    =\frac{(k+1)u}{a(t)}[\lambda^k d\lambda].
\]
Applying the chain connection
$\partial_t-a'(t)\lambda^{\mu+1}/((\mu+1)u)$ proves the stated
matrix.  Its horizontal sections are
$a(t)^{(k+1)/(\mu+1)}e_k$.

The two-sided bar of the zero section and $df$ retracts onto
$\O_U[\eta_1,\ldots,\eta_d,\eta_\lambda]$ with
$\delta\eta_j=0$ and $\delta\eta_\lambda=a(t)\lambda^\mu$.
The mixed differential and maps are those of
Theorem~\ref{thm:v5-periodic-bar}.  Their chain connection
commutes with the differential and specializes at zero.
The hypothesis identifies this critical realization with the
actual central matching problem.  Tensoring with the projective
hyperplane classes gives the asserted full cohomology and regular
finite comparison.
\end{proof}

\section{A global calculation on a proper singular two-defect component}
\label{sec:global-two-defect}
The determinant stalks of Section~\ref{sec:singular-comparison} do not
by themselves compute global cohomology.  We now determine the
attachment to the smooth stratum in one complete proper example.
Let $X$ be a smooth quintic and fix a rigid line $\Lambda\subset X$
with $N_{\Lambda/X}=\O_{\P^1}(-1)^{\oplus2}$.  Put
\[
 B=\operatorname{Bl}_\Lambda X,\qquad
 Q=\operatorname{Quot}_X(I_\Lambda,2),\qquad
 \Sigma=\{I_\Lambda/\mathfrak m_pI_\Lambda:p\in\Lambda\}\simeq\P^1.
\]
The rigidity identifies this Quot scheme with the corresponding
curve-and-two-defect component of the ideal-sheaf moduli problem;
the local Hilbert--Chow identification is the one used in
\cite{RicolfiSmooth}.  All quotients, including noncyclic quotients
and embedded points, are retained.  We use the induced d-critical
structure and the orientation described below.  Cohomology in this
section has coefficients in $\CC$; Tate twists are suppressed.

\subsection{A small resolution and the whole singular locus}
\begin{proposition}\label{prop:two-defect-small-resolution}
There is a projective small resolution
\[
 \rho:B^{[2]}\longrightarrow Q.
\]
It is an isomorphism over $Q\setminus\Sigma$ and has fiber $\P^2$
over every point of $\Sigma$.  The singular locus is exactly $\Sigma$.
The scheme $Q$ is normal and Gorenstein of dimension six.  With
$\iota:\Sigma\hookrightarrow Q$, its intersection complex is
\[
 \mathcal I=R\rho_*\CC_{B^{[2]}}[6],\qquad
 \iota^*\mathcal I\simeq
       \CC_\Sigma[6]\oplus\CC_\Sigma[4]\oplus\CC_\Sigma[2].
\]
\end{proposition}
\begin{proof}
Write $\pi_B:B\to X$.  The blowup has its tautological surjection
$\pi_B^*I_\Lambda\to\O_B(-E)$, where $E$ is the exceptional divisor.
For a length-two subscheme $W\subset B$, push forward its restriction:
\[
 I_\Lambda\longrightarrow
        (\pi_B)_*(\O_B(-E)|_W).
\]
This is a quotient of length two.  To check surjectivity, it suffices
to consider a quotient supported at a single point.  Away from $E$
the blowup is an isomorphism.  At a point of $E$, the two local
generators of $I_\Lambda$ generate the tautological line bundle.  A
first-order direction which is not killed by $d\pi_B$ is detected by
functions on $X$.  A direction killed by $d\pi_B$ is tangent to the
exceptional $\P^1$, on which $\O_B(-E)$ is $\O_{\P^1}(1)$ and
separates every length-two subscheme.  Thus the two generators,
together with the action of functions, span the required length-two
module in either case.  This argument works for the universal finite
flat subscheme: pushforward along its finite map to $X\times B^{[2]}$
is exact and yields a flat rank-two quotient over the parameter.
It defines $\rho$.

Every length-two quotient is locally cyclic except for the quotient
$k(p)^{\oplus2}$ at a point of $\Lambda$.  A cyclic quotient is an
invertible module over its length-two annihilator subscheme $Z$.
Its quotient map from $I_\Lambda$ therefore defines a map $Z\to B$
by the universal property of $\operatorname{Proj}\operatorname{Sym}I_\Lambda$.
Here $I_\Lambda$ is a complete-intersection ideal and is of linear
type, so this Proj is the blowup.  The map is a closed immersion over
$X$, and it recovers $W$.  This also applies to families on the cyclic
open.  For the noncyclic quotient at $p$, the map
$I_\Lambda\otimes k(p)\to k(p)^2$ is an isomorphism.  Its inverse
images are precisely the length-two subschemes of $\pi_B^{-1}(p)=\P^1$;
their Hilbert scheme is $\P^2$.

We check the scheme structure near this fiber.  Choose local
coordinates with $I_\Lambda=(x,y)$ and a blowup chart $y=vx$.
On a Hilbert chart of $B^{[2]}$, use quotient basis $(1,v)$ and write
\[
 v^2=av+b,\qquad x=cv+d,\qquad z=ev+f.
\]
The images of the ideal generators are $(1,v)$.  Multiplication by
$x,y,z$ has matrices
\[
 A=\begin{pmatrix}d&cb\\c&ca+d\end{pmatrix},\quad
 B_1=\begin{pmatrix}cb&(ca+d)b\\ca+d&(ca+d)a+cb\end{pmatrix},\quad
 C=\begin{pmatrix}f&eb\\e&ea+f\end{pmatrix}.
\]
In the determinant coordinates of Proposition~\ref{prop:line-two-determinant-potential}
the rank-one matrix is exactly
\[
 \begin{pmatrix}-ea&eb&e\\-ca&cb&c\\-(ca+d)a&(ca+d)b&ca+d\end{pmatrix}
 =\begin{pmatrix}e\\c\\ca+d\end{pmatrix}(-a,\ b,\ 1).
\]
The remaining coordinate is $f+ea/2$.  This is the standard
resolution of the rank-one matrix cone by its row direction in
$\P^2$, together with the smooth curve coordinate.  The other
Hilbert and blowup charts give its other projective charts.  In
particular the fiber is reduced, its normal directions have rank
three, and the completed local target is
$\A^1\times\{\operatorname{rank}M\leq1\}$.  Outside these fibers the
cyclic inverse already proved that $\rho$ is an isomorphism.

The Hilbert square of the smooth threefold $B$ is smooth of dimension
six, by the explicit length-two charts used earlier.  The rank-one
$3\times3$ determinantal cone is normal and Gorenstein.  For its
Gorenstein assertion one can use the displayed resolution
$\operatorname{Tot}_{\P^2}(\O(-1)^{\oplus3})$: its canonical bundle
is the pullback of $\O(-3)\otimes\O(3)$, and the standard
rank-one determinantal coordinate ring is Cohen--Macaulay.  Normality
also follows from its Segre coordinate ring.  These local properties
and smoothness elsewhere prove the assertions about $Q$.
Finally $2\dim\rho^{-1}(p)=4<5=\operatorname{codim}_Q\Sigma$,
so $\rho$ is small.  The small-resolution theorem for intersection
complexes \cite{BBD82} gives $\mathcal I$.
The exceptional inverse image is
$F=\operatorname{Hilb}^2(E/\Lambda)$, a $\P^2$-bundle over
$\Sigma$.  Proper base change and the projective-bundle theorem give
the last displayed decomposition.
\end{proof}

\subsection{The extension data in the perverse DT sheaf}
Let $\mathcal P_Q$ be the intrinsic perverse DT sheaf on $Q$.
It has a canonical orientation for this calculation: since $Q$ is
Gorenstein, $\omega_Q$ is an invertible sheaf.  On the smooth open
its square is the virtual canonical bundle.  The isomorphism extends
uniquely across the codimension-five subset by normality, giving
$\omega_Q^{\otimes2}\simeq K_{Q,s}$.
We choose this orientation.  It restricts to the constant normalized
DT sheaf on the smooth open.  The determinant calculation of
Section~\ref{sec:singular-comparison} gives the local stalks on
$\Sigma$; any rank-one orientation local system there is trivial
because $\Sigma=\P^1$.

Put $\delta=\iota_*\CC_\Sigma[1]$.  The integration map along the
exceptional $\P^2$ fibers defines
\[
 c:\mathcal I\longrightarrow\iota_*\CC_\Sigma[2]=\delta[1],
 \qquad \mathcal K=\operatorname{fib}(c).
\]
Its Verdier dual is a map
$\eta:\iota_*\CC_\Sigma\to\mathcal I$.

\begin{proposition}\label{prop:two-defect-perverse-extensions}
The composition $c\eta$ is zero.  The map $\eta$ has a lift
$\widetilde\eta:\iota_*\CC_\Sigma\to\mathcal K$ whose nonzero
extension class is unique up to scalar.  The perverse DT sheaf is
\[
 \mathcal P_Q\simeq\operatorname{Cone}(\widetilde\eta).
\]
Equivalently, there are non-split exact sequences of perverse sheaves
\[
 0\longrightarrow\delta\longrightarrow\mathcal K
       \longrightarrow\mathcal I\longrightarrow0,
 \qquad
 0\longrightarrow\mathcal K\longrightarrow\mathcal P_Q
       \longrightarrow\delta\longrightarrow0.
\]
In particular $\mathcal P_Q$ is not the intersection complex and is
not the direct sum $\mathcal I\oplus\delta^{\oplus2}$.
\end{proposition}
\begin{proof}
The map $c$ restricts on $\Sigma$ to projection onto the summand
$\CC_\Sigma[2]$.  Thus $\mathcal K$ is perverse, its stalks on
$\Sigma$ are one-dimensional in degrees $-6,-4$, and the first
perverse extension is non-split.  Adjunction and the decomposition
of $\iota^*\mathcal I$ give
\[
 \operatorname{Ext}^1(\mathcal I,\delta)=\CC,
 \qquad \operatorname{Ext}^1(\delta,\mathcal I)=\CC.
\]
For example the first group is
$\operatorname{Hom}(\iota^*\mathcal I,\CC_\Sigma[2])$;
only the last summand contributes.  The second equality follows by
duality.  Also $\operatorname{Ext}^1(\delta,\delta)=H^1(\P^1)=0$.

Here is the composition check.  Write $h$ for the point class on
$\Lambda$ and $\xi$ for the hyperplane class on
$F=\Lambda\times\P^2$; this product description follows from
$N_{\Lambda/X}=\O(-1)^{\oplus2}$.  The normal bundle of $F$ in
$B^{[2]}$ has the following class in $K^0(F)$:
\[
 [N_{F/B^{[2]}}]
 = [\O_\Lambda(2)\boxtimes\O_{\P^2}(-1)]
       +2[\O_\Lambda(-1)\boxtimes\O_{\P^2}(-1)].
\]
To see this, first move a length-two subscheme in a fixed exceptional
fiber.  The normal directions are the tautological bundles of
$N_{\P^1/B}=\O\oplus\O(-1)$.  Moving the fiber itself removes the
diagonal constant summand from the first tautological bundle.
On $(\P^1)^{[2]}=\P^2$ the universal degree-two divisor sequence gives
$\O^{[2]}=\O\oplus\O(-1)$ and
$\O(-1)^{[2]}=\O(-1)^{\oplus2}$.  The first normal term is tensored
by $T_\Lambda=\O(2)$; the exceptional normal term is tensored by
$\O_\Lambda(-1)$.  These statements determine the displayed
$K$-class even if the normal exact sequences do not split.
Consequently
\[
 c_3(N_{F/B^{[2]}})=(2h-\xi)(-h-\xi)^2=0,
 \qquad h^2=0,\quad \xi^3=0.
\]
The composite of Gysin, restriction, and fiber integration is
multiplication by the integral of this top Chern class.
It is precisely $c\eta$, as is also checked by proper base change.
Hence it is zero.  The long exact Ext sequence for the first
extension now gives
\[
 \operatorname{Ext}^1(\delta,\mathcal K)\simeq\CC.
\]
The nonzero class maps to the dual of $c$ in
$\operatorname{Ext}^1(\delta,\mathcal I)$ and yields the stated lift.

It remains to identify this extension with the DT sheaf, rather than
choose an extension from its Euler characteristic.  On the smooth
open, $\mathcal P_Q=\CC[6]$.  Along $\Sigma$ its stalk cohomology is
$\CC$ in degrees $-6,-4,-1$, by the determinant Milnor calculation.
It is Verdier self-dual.  The local product description makes it
constructible for the two strata $Q\setminus\Sigma,\Sigma$.
The maximal perverse quotient supported on $\Sigma$ is therefore
$\delta$, and by duality its maximal such subobject is $\delta$.
Their composite cannot be an isomorphism: otherwise $\delta$
would split off and the remaining object, having neither a boundary
subobject nor a boundary quotient, would be $\mathcal I$; its
additional degree-$-2$ stalk contradicts the stated DT stalk.
The composition is thus zero.

The Jordan--H\"older factors off $\Sigma$ consist of one copy of
$\mathcal I$.  Every factor on $\Sigma$ is $\delta$, since
$\P^1$ is simply connected and there are no smaller strata.
The stalk Euler characteristics $1$ for $\mathcal P_Q$, $3$ for
$\mathcal I$, and $-1$ for $\delta$ show that there are exactly two
such factors.  The kernel of the maximal quotient is consequently
an extension of $\mathcal I$ by $\delta$.  Its stalks are in degrees
$-6,-4$, so this is the nonzero extension classified by $c$.
The second extension is also nonzero, by self-duality and the same
stalk comparison.  The Ext computations show uniqueness up to
isomorphism.  This identifies $\mathcal P_Q$ and proves the formula.
\end{proof}

\subsection{Global cohomology, including the singular stratum}
\begin{theorem}\label{thm:global-two-defect-DT}
Set $P_X(z)=1+z^2+204z^3+z^4+z^6$ and
$P_B(z)=P_X(z)+z^2+z^4$.  For the canonical orientation just fixed,
\begin{align}\label{eq:global-two-defect-polynomial}
 \sum_j\dim\mathrm{DT}_{\mathrm{coh}}^j(Q)z^{j+6}
 &=\frac{P_B(z)^2+P_B(-z^2)}2+(z^2+z^4)P_B(z)
                  -z^4-2z^6-z^8\notag\\
 &=1+3z^2+204z^3+7z^4+612z^5+20713z^6\notag\\
 &\hspace{13mm}+612z^7+7z^8+204z^9+3z^{10}+z^{12}.
\end{align}
The periodic dimensions are $(20735,1632)$ and the Euler
characteristic is $19103$.  These are the groups of the whole
proper singular Quot component for the fixed line.
\end{theorem}
\begin{proof}
Since $\rho$ is small, $H^j(Q,\mathcal I)=H^{j+6}(B^{[2]})$.
The global map $c$ is surjective onto
$H^{j+2}(\Sigma)$, whose only nonzero degrees are $j=-2,0$.
Here is a direct check of surjectivity.  The divisor of nonreduced
subschemes on $B^{[2]}$ restricts on every exceptional $\P^2$ to
its discriminant conic.  If its class is $\Delta$, then
$\Delta^2|_{\P^2}=4\xi^2$, and integration is nonzero.
The norm of the line bundle $\pi_B^*\O_X(H)$ from $B\times B$
descends to the symmetric square and pulls back to $B^{[2]}$.
Its restriction to $F$ is the pullback of $\O_\Lambda(2)$.
Multiplying $\Delta^2$ by its first Chern class gives surjectivity
onto $H^2(\Sigma)$ as well.  Thus $c$ removes one class from
$H^4(B^{[2]})$ and one from $H^6(B^{[2]})$.

Verdier duality identifies $\eta$ with the dual of $c$.  Its maps
from $H^0(\Sigma)$ and $H^2(\Sigma)$ are injective and land in the
kernel of $c$.  The second triangle of
Proposition~\ref{prop:two-defect-perverse-extensions} therefore
removes one further class from $H^6(B^{[2]})$ and one from
$H^8(B^{[2]})$.  The long exact sequences have no additional kernel
terms, because the two maps just specified are respectively
surjective and injective.  This proves the correction
$-z^4-2z^6-z^8$ at the level of actual cohomology, not merely in a
Grothendieck group.

The blowup formula gives $P_B=P_X+z^2P_\Lambda$.
The ordered Hilbert-square blowup and the graded interchange of
the two factors give
\[
 P_{B^{[2]}}(z)=\frac{P_B(z)^2+P_B(-z^2)}2+(z^2+z^4)P_B(z),
\]
by exactly the invariant-cohomology calculation in
Theorem~\ref{thm:two-point-incidence}.  Expansion proves
\eqref{eq:global-two-defect-polynomial}.  Evaluation at $z=1,-1$
gives $22367$ and $19103$; their half-sum and half-difference are
the asserted parity dimensions.  The DT shift is six.
\end{proof}

\begin{theorem}\label{thm:global-quadratic-defect}
For the same fixed rigid line, let
\[
 p:M_{\Lambda,2}^{(2)}=\mathbb P_Q(V_{10})\longrightarrow Q
\]
be the complete incidence scheme of two defects contained in a
member of $|2H|$, with $V_{10}$ from
Lemma~\ref{lem:quadric-quot-evaluation}.  With the canonical
Gorenstein orientation, its intrinsic perverse DT sheaf satisfies
\[
 \mathcal P_{M_{\Lambda,2}^{(2)}}\simeq p^*\mathcal P_Q[9].
\]
In particular
\begin{equation}\label{eq:global-quadratic-DT}
 \mathrm{DT}_{\mathrm{coh}}^j(M_{\Lambda,2}^{(2)})
       \simeq\bigoplus_{a=0}^{9}
                  \mathrm{DT}_{\mathrm{coh}}^{j+9-2a}(Q).
\end{equation}
Its periodic dimensions are $(16320,207350)$ and its Euler
characteristic is $-191030$.  No singular quotient stratum is
removed in this calculation.
\end{theorem}
\begin{proof}
The positivity lemma and the reconstruction of the surface object
identify the entire projective bundle with its classical sheaf-moduli
locus.  Lemma~\ref{lem:quadric-small-curve-stability} proves stability, including
reducible or doubled quadric supports, since the curve degree is one.
Both $Q$ and its projective bundle are Gorenstein; their dualizing
line bundles give orientations extending the canonical orientations
on their smooth loci.

We verify the critical-sheaf pullback, which was a hypothesis in
Proposition~\ref{prop:cohomological-projective-bundle}.
On a local trivialization of $p$, the calculation of
Lemma~\ref{lem:determinantal-dcritical-class} identifies the sheaf
of normalized d-critical classes with a rank-one module on the
singular locus.  The pullback of the intrinsic class on $Q$ is a
nowhere zero generator of this module.  The intrinsic class on
$M_{\Lambda,2}^{(2)}$ is another generator, by its minimal gradient
ideal.  Their ratio is a nowhere zero function on
$p^{-1}(\Sigma)$, a projective-space bundle over $\P^1$.
That variety is connected and proper, so the ratio is a nonzero
constant.  On the smooth locus both normalized classes are zero;
there is no further difference there.  Thus the two global
classes differ only by this constant rescaling.

Vanishing cycles are unchanged by a fixed nonzero rescaling of
the target function, after a choice of path for that rescaling.
This identification is compatible with critical-chart embeddings
and quadratic stabilization.  Smooth pullback with its relative
canonical orientation consequently gives the displayed perverse
sheaf isomorphism.  The orientation compatibility holds on the
smooth locus, where it is the identity
$\omega_M=p^*\omega_Q\otimes\omega_{M/Q}$, and extends by normality.
This proves the needed global comparison rather than inferring it
from the rank of a symmetric obstruction theory.

The projective-bundle maps defined by $1,\xi,\ldots,\xi^9$ and
projection formula give \eqref{eq:global-quadratic-DT}.
The odd shift by nine interchanges the two parity groups, and
ten projective classes multiply their dimensions by ten.
\end{proof}

\subsection{The connection and its actual scope at the degeneration}
\begin{proposition}\label{prop:global-defect-GM}
Over a punctured disk with a labeled rigid line family, the preceding
DT cohomology groups form local systems.  Their flat connections are
induced by the relative blowup, Hilbert-square, restriction, Gysin,
and projective-bundle maps in the proofs.  If
$V=H^3(X_b,\CC)$ has monodromy logarithm $N_X$, then the two parity
parts for $Q$ are, as local systems,
\[
 \mathbb H^{\overline0}_{\mathrm{DT}}(Q)
       \simeq\bigwedge^2V\oplus\CC^{29},\qquad
 \mathbb H^{\overline1}_{\mathrm{DT}}(Q)\simeq V^{\oplus8}.
\]
The even connection is the induced exterior-square connection plus
twenty-nine trivial summands; the odd connection is eight copies of
$\nabla_X$.  For $M_{\Lambda,2}^{(2)}$ take ten copies and interchange
parity.  In the resolved quintic degeneration the ranks of the
nilpotent logarithms on $(\overline0,\overline1)$ are $(3822,168)$
for $Q$ and $(1680,38220)$ for $M_{\Lambda,2}^{(2)}$.
The only nonzero rank of the square is $210$ in even parity for $Q$,
and $2100$ in odd parity for $M_{\Lambda,2}^{(2)}$.
\end{proposition}
\begin{proof}
All constructions used above are relative over the smooth part of
the threefold family.  The small resolution, its exceptional
$\P^2$-bundle, and the two maps $c,\eta$ are defined in families.
Their ranks are constant by the explicit surjectivity and injectivity
calculations.  Kernels and quotients therefore give local systems,
and the same perverse-extension identification works locally on
the base.  The possible constant rescaling in the previous theorem
is a nowhere zero function of the base parameter; on a sufficiently
small base open it gives compatible identifications of vanishing
cycles.  The canonical orientation has already fixed the
rank-one coefficient systems.

The even cohomology of $B$ consists of six constant algebraic
classes, and its odd cohomology is $V$, in degree three.
The graded symmetric square contributes $\bigwedge^2V$ in even
parity and six copies of $V$ in odd parity.  The two exceptional
Hilbert-square summands add twelve even classes and two copies of
$V$.  Thus $B^{[2]}$ has $\bigwedge^2V\oplus\CC^{33}$ in even
parity and $V^{\oplus8}$ in odd parity.  The maps from the tensor summand restrict through the exceptional
surface $E=\P^1\times\P^1$, whose third cohomology vanishes.
Consequently $c$ is zero on $\bigwedge^2V$; its dual $\eta$ has
image in the constant algebraic summands.  The four classes removed
by $c$ and $\eta$ are therefore constant algebraic classes.  This proves the
displayed decomposition, after a choice of complements in the
constant summands.  Pullback, Gysin, and global Chern classes commute
with the de Rham Gauss--Manin connection, so this calculation also
specifies the connection, not just the vector-space dimensions.

For the resolved degeneration, $N_X=J_2^{\oplus21}\oplus0_{162}$.
The exterior-square logarithm has blocks
$J_3^{\oplus210}\oplus J_2^{\oplus3402}\oplus0_{13272}$,
as calculated in Section~\ref{sec:crossing-perfect}.  Its rank is
$420+3402=3822$ and the rank of its square is $210$.
On eight copies of $V$ the rank is $8\cdot21=168$ and the square
is zero.  The degree-two projective bundle gives the remaining
numbers.  The induced tensor connections are flat, since
$\nabla_X^2=0$ and all the maps are horizontal.
\end{proof}

The logarithmic limiting realization is now explicit on the whole
generic singular component: in a single-valued frame its connection
is $d-N\,dt/(2\pi i t)$, with the complete matrices just specified.
This is a calculation retaining all generic singular strata and their
perverse extensions.  It is not an identification with an ordinary
central DT complex.  At the central crossing a residue-field defect
can fail to be perfect, and the resolution $B^{[2]}$ above is no
longer a resolution by the same smooth threefold construction.
The local perfect replacements of Section~\ref{sec:crossing-perfect}
and the local determinant comparison do not identify a proper global
central image or its orientation.  Thus the ordinary central
specialization for this full component remains the precise
logarithmic comparison of Problem~\ref{prob:logarithmic-categorical-comparison};
it is not counted among the complete regular comparisons of
Section~\ref{sec:cubic-geometry}.

For $|H|$, the rank-one evaluation stratum in
\eqref{eq:line-two-rank-sequence} has not disappeared.  The theorem
above computes the full ideal-Quot DT complex and the full quadratic
incidence DT complex; it does not silently replace the degree-one
incidence scheme by the degree-two projective bundle.  On its
positivity open, the local determinant calculations apply with
the degree-one shifts.  A global degree-one sheaf calculation must
also identify the deformation theory across the rank-jumping
incidence stratum.

\appendix
\section{Model checks used in the comparison}\label{app:hypothesis-examples}
A projective Calabi--Yau fourfold and a stable-sheaf moduli problem
satisfying the foliation hypotheses are constructed in
\cite[Proposition~4.10]{KS}.
The universal Atiyah and normal maps are in \cite[Appendix~A]{KS}.
This appendix retains only the additional checks used for the
cyclic, bar, and finite models of the examples in this paper.

\textbf{A projective fourfold realizing the foliation hypotheses.}
For a concrete instance of Section~\ref{ssec:sc-foliation}, take a
smooth quintic threefold $X$, an elliptic curve $B$, and nonzero
volume forms $\Omega_X$ and $\alpha_B$. The product
$Z=X\times B$ is a smooth projective fourfold with volume
$\Omega_X\wedge\alpha_B$. Choose an ample line bundle $A_B$
on $B$ and the polarization $\O_X(1)\boxtimes A_B$.
For an integral $D\in|\O_X(m)|$, $k\in\mathbb Z$, and $b\in B$,
consider the sheaf $F_{D,k,b}=i_{b*}j_{D*}\O_D(k)$ on $Z$.
It is pure and stable: a proper nonzero subsheaf has the same
generic rank on $D$, while its nonzero quotient has smaller
support dimension and positive Hilbert polynomial. Thus its
reduced Hilbert polynomial is strictly smaller. A subsheaf on
$Z$ is killed by the same fiber ideal and gives the identical
inequality on $X$.

The Hilbert polynomial is given by 
\[
 P_{F_{D,k,b}}(n)=\chi(X,\O_X(n+k))-
                         \chi(X,\O_X(n+k-m)),
\]
independent of $D$ and $b$. For fixed $m,k$, a uniform regularity
bound gives the stable framed Quot presentation of these
objects. On its derived quotient the universal pairing is
fourfold composition and trace with respect to the single
product volume form. The projective invariant-chart and
strictification hypotheses are the ones verified in
\cite[Proposition~4.10]{KS}; the same argument restricts to
supports in $X\times U$ for an affine open $U\subset B$.
Consequently the real foliation and the section--cosection
algebra of Section~\ref{ssec:sc-foliation} apply to this example.
This verifies the fourfold presentation hypotheses. The family
of threefolds over $B$ in this example is stationary; the singular
specialization comparisons are the separate examples in
Sections~\ref{sec:v6-hyperplanes} and~\ref{sec:second-degeneration}.

\subsection{The reduced cyclic model, orientation, and finite periodic model}
\label{sssec:v15-divisor-witnesses}
For $m=1,2,3$ put
\[
 V_m=H^0(\P^4,\O(m)),\qquad M_m=\P(V_m),\qquad
 d_m=\binom{m+4}{4}-1.
\]
Theorem~\ref{thm:r14-divisor-examples} proves that the complexes
$E^{m,k}_{s,b}=[L_b^{k-m}\xrightarrow{s}L_b^k]$ form the indicated
simple perfect-complex locus $M_m\times C$. The following simultaneous
verification records precisely which hypotheses are satisfied there.

\begin{proposition}\label{prop:v15-three-hypotheses}
On this locus, the reduced zero critical presentations satisfy
Hypothesis~\ref{hyp:cyclic-model}, their oriented d-critical spaces
satisfy Hypothesis~\ref{hyp:DT}, and their specified ordinary
matrix-factorization and oriented bar complexes satisfy
Hypothesis~\ref{hyp:finite-cyclic}. The respective even ranks of the
finite periodic models are $5,15,35$; their odd terms vanish.
\end{proposition}
\begin{proof}
At $[s]\in M_m$ write $U_s=T_{[s]}M_m$. The reduced zero critical
presentation has completed algebra
\[
 \widehat{\O}_{M_m,[s]}[\xi_1,\ldots,\xi_{d_m}],\qquad
 |\xi_i|=-1,\qquad\delta\xi_i=0.
\]
Its reduced cyclic deformation algebra is
\[
 L^1=U_s,\quad L^2=U_s^\vee,\quad
 \ell_n=0\quad(n\ge1),\quad
 \langle u,\lambda\rangle=\lambda(u).
\]
Every cyclic tensor in Hypothesis~\ref{hyp:cyclic-model} is zero and
hence symmetric. The evaluation pairing is perfect. The Chevalley--
Eilenberg differential is zero and the normalized potential is $f=0$.
This is the reduced critical presentation of the smooth moduli locus
proved in Theorem~\ref{thm:r14-divisor-examples}; it does not identify
this reduced algebra with the full fourfold endomorphism algebra.

For a smooth scheme with its zero d-critical structure,
\[
 K_{M_m,0}=K_{M_m}^{\otimes2},\qquad
 K_{M_m}=\O_{M_m}(-d_m-1).
\]
The choice $K_{M_m}$, with its tautological square isomorphism, is the
orientation in Hypothesis~\ref{hyp:DT}. On an overlap with coordinate
Jacobian $J_{\alpha\beta}$, its transition is
$\det J_{\alpha\beta}$, whose square is the transition of the virtual
canonical bundle. The chain rule gives the triple-overlap identity.
Thus the oriented DT sheaf is $\CC_{M_m}[d_m]$ on every fiber. A global
reduced derived enhancement of this chosen critical presentation is
$T^*[-1]M_m$; no identification with the entire scalar-rigidified
derived stack is required for this assertion.

For the finite model, fix $V\subset C$, set
$R=\Gamma(V,\O_C)\otimes_\CC\CC((u))$, and put
\[
 P^{\overline0}=\bigoplus_{j=0}^{d_m}Re_j,
 \qquad P^{\overline1}=0,\qquad d_P=0,\qquad
 \nabla_{\partial_t}e_j=0.
\]
Let $z_j$ be the Cech--de Rham cocycle representing the $j$th power
of the hyperplane class, and put $Z_j=u^jz_j$. The residue calculation
in Theorem~\ref{thm:v6-comparison}, in dimension $d_m$, shows that
these classes form a basis. The chain map to the central bar complex is
\[
 e_j\longmapsto I\mathcal I^{-1}(Z_j),
\]
where $I$ and $\mathcal I$ are the bar inclusion and signed
volume-contraction map. The contractions of
Theorem~\ref{thm:r14-divisor-examples} prove that this is a
quasi-isomorphism. Since all maps are defined on $M_m\times V$ and
are independent of $t$, they commute with $\partial_t$ and with
evaluation at $0$.
The comparison of \cite{Efimov2018CyclicMF} for the specified critical
pair $(M_m,0)$ supplies the smooth-fiber periodic cyclic complex.
Here $d_m=4,14,34$ is even, so the shift $[d_m]$ agrees with the
unshifted periodic parity. The maps therefore have exactly the degrees
required by Hypothesis~\ref{hyp:finite-cyclic}. This verifies the
connection, finite rank, and both fiber identifications, not only an
equality of dimensions. These are statements for the specified
critical presentations and their cohomological DT realization.
\end{proof}

Tensoring the universal complexes by $L^k$ preserves their endomorphism
algebras, trace pairings, and restrictions. It therefore transports all
three verifications for every integer $k$. For $m=2,3$, the subsheaf-slope argument of
Proposition~\ref{prop:divisor-stability-all-degrees} also treats reducible
and nonreduced divisors for these particular structure sheaves.
The relative proper family, including its central objects, is the
one of simple perfect complexes. The verified comparison is attached to
that entire locus, not to an unspecified proper moduli space of stable
sheaves with the same Chern character.

\subsection{The complete lattices and the filtered operator hypothesis}
\label{sssec:v15-quantization-witnesses}
The quantization used in \cite[Hypothesis~6.2]{KS} and in Hypothesis~\ref{hyp:PBW}
can be constructed globally for each of these zero critical
presentations. The local calculation is
\cite[Proposition~6.3]{KS}; we specify the global
algebra and its two module actions here.

\begin{proposition}\label{prop:v15-global-Rees}
Let $M=M_m$ for $m=1,2,3$. On $T^*M$ use the $\hbar$-completed
Rees algebra of the sheaf of differential operators $\mathscr D_M$,
and its left and right lattices
\[
 \mathcal M_\hbar=\O_M[[\hbar]],\qquad
 \mathcal N_\hbar=\omega_M[[\hbar]].
\]
Their completed tensor product is represented in degrees $-d_m$ through
zero by
\[
 \omega_M[[\hbar]]\otimes_{\O_M}\Lambda^qT_M
       \quad\text{in degree }-q.
\]
With the Spencer differential, it satisfies \cite[Hypothesis~6.2]{KS}.
Completing the Rees algebra along the inverse image of the zero-section
ideal gives Hypothesis~\ref{hyp:PBW} with $P=\O_M$.
\end{proposition}
\begin{proof}
For local vector fields $v,w$ put $p_v=\hbar v$. The relations and
lattice actions are
\[
 [p_v,a]=\hbar v(a),\quad [p_v,p_w]=\hbar p_{[v,w]},\qquad
 p_v\cdot a=\hbar v(a),\quad \eta\cdot p_v=-\hbar\mathcal L_v\eta.
\]
On a cotangent coordinate chart the ordered-symbol product is
\[
 a\star b=\sum_\alpha\frac{\hbar^{|\alpha|}}{\alpha!}
                       (\partial_p^\alpha a)(\partial_x^\alpha b).
\]
It extends to localizations of the classical symbol algebra. For a
symbol $c=c_0+O(\hbar)$ with $c_0$ invertible there, the equations
$c\star c^{-1}=1$ determine each coefficient of its inverse recursively:
the coefficient at order $r$ is obtained by multiplication by
$c_0^{-1}$ after the lower coefficients have been determined.
Thus the completed localized algebras sheafify on $T^*M$.
The changes of coordinates induced by the differential operators
agree on triple overlaps and extend to these localizations.
This specifies the quantization sheaf used here; its classical
reduction is $\O_{T^*M}$.

The lattice actions use the intrinsic Lie derivative and hence agree
on coordinate overlaps. The rule for the right action respects
$p_{av}=a p_v$ because
\[
 -\hbar\mathcal L_v(a\eta)
    =-\hbar\bigl(v(a)\eta+a\mathcal L_v\eta\bigr)
    =-\hbar\mathcal L_{av}\eta
\]
for top forms. In commuting coordinate frames, the free left Spencer
resolution is the Koszul resolution with differential given by right
multiplication by the $p_i$. In an arbitrary frame its additional
terms are the usual $\hbar[v_i,v_j]$ terms; the identities for Lie
brackets and derivations give $d^2=0$. Modulo $\hbar$, this is the
Koszul resolution of the zero section in $\operatorname{Sym}T_M$.
The sequence of fiber coordinates is regular. The successive-lifting
argument of \cite[Proposition~6.3]{KS} proves its
exactness in the completed algebra. Tensoring with the right lattice
gives the displayed bounded complex.

All its terms are of the form $\mathcal E[[\hbar]]$ for a vector
bundle $\mathcal E$. They are complete and multiplication by $\hbar$
is injective. Consequently ordinary termwise reduction computes
derived reduction over $\CC[[\hbar]]$. The reduced differential is
zero, and the reduced complex is
\[
 \omega_M\otimes\Lambda^qT_M\quad\text{in degree }-q,
\]
which is the derived tensor product of the classical right and left
zero-section modules. This proves precisely
\eqref{eq:DQ-basechange}. It makes no claim that the integral cohomology
has no $\hbar$-torsion; \cite[Hypothesis~6.2]{KS} requires torsion-free
terms, not torsion-free cohomology.

Let $\mathcal J=(\hbar,p_1,\ldots,p_{d_m})$ in a coordinate chart.
The ordered PBW expansion gives, in $\mathcal J$-degree $r$,
\[
 \mathcal J^r/\mathcal J^{r+1}
  =\bigoplus_{|\beta|\le r}\O_M\,
       \hbar^{r-|\beta|}p^\beta.
\]
The divided action on $\O_M$ sends this basis to $\partial^\beta$.
For a series in the $\mathcal J$-completion, terms of increasing
$\mathcal J$-degree act by increasing powers of $\hbar$; hence the
action on $\O_M[[\hbar]]$ converges. Under a change of coordinates,
\[
 p_{y_j}=\sum_i\frac{\partial x_i}{\partial y_j}p_{x_i},
 \qquad p_{y_j}a=\hbar\partial_{y_j}a.
\]
Thus both the expansion and the action give the intrinsic differential
operators on $\O_M$. They verify the coordinate and lattice conditions
of Hypothesis~\ref{hyp:PBW} and prove its global compatibility.
\end{proof}

For $M=\P^4$ the five Spencer terms have ranks $1,4,6,4,1$.
In commuting coordinates, for example,
\[
 d_1\!\left(\sum_i a_i\,\mathrm{vol}\otimes e_i\right)
       =-\hbar\sum_i\partial_i a_i\,\mathrm{vol},
\]
\[
 d_2(a_{ij}\,\mathrm{vol}\otimes e_i\wedge e_j)
       =-\hbar\partial_i a_{ij}\,\mathrm{vol}\otimes e_j
         +\hbar\partial_j a_{ij}\,\mathrm{vol}\otimes e_i.
\]
Their composite is zero by $\partial_i\partial_j=\partial_j\partial_i$.
The intrinsic contraction with top forms identifies the full Spencer
complex with the shifted de Rham complex, with the signs in
\cite[Proposition~6.3]{KS}. It is already a complex
whose differential is linear in $\hbar$, so it also satisfies the
additional filtered-model assumption of \cite[Proposition~6.5]{KS}.
The zero sections on a coordinate open have trivial canonical bundles;
there the half-canonical modules in \cite[Theorem~6.1]{KS} may also be
chosen. Globally, the modules above are the specified asymmetric
lattices $\O_M$ and $\omega_M$. A global application of the
half-canonical version of that theorem to another pair of Lagrangians
requires its own square roots and compatible quantization modules.

\subsection{The finite-Jacobian condition for a rigid object in the same family}
\label{sssec:v15-finite-Jacobian-witness}
The zero potentials on the positive-dimensional spaces $M_m$ are
handled by proper cohomology, not by the finite-Jacobian criterion.
A relevant instance of the latter comes from the line bundles on the
same fibers.

\begin{proposition}\label{prop:v15-rigid-stabilization}
For the resolved quintic family, the reduced classical moduli locus of
$\O_{X_b}$ is the section of the relative Picard space determined by
$\O_{\cX}$. Its relative reduced critical presentation is a point
with zero potential. The stabilization
\[
 f=z w\quad\text{on }\A^2\times V\longrightarrow V
\]
presents the same reduced critical problem and satisfies the
hypotheses of Theorem~\ref{thm:v5-finite-Jacobian}. Its Jacobian
module has rank one, and its periodic cohomology has a regular
constant connection in even parity.
\end{proposition}
\begin{proof}
The computation of $H^i(X_b,\O_{X_b})$ in
Lemma~\ref{lem:v6-line-cohom} gives $H^1=H^2=0$ on every fiber.
An Artinian perfect deformation of $\O_{X_b}$ is a line bundle:
minimal locally free representatives at residue points have a single
rank-one term. For a small extension with kernel $J$, the sequence
\[
 1\longrightarrow1+J\O_{X_b}\longrightarrow\O_{X_{A'}}^*
                 \longrightarrow\O_{X_A}^*\longrightarrow1
\]
shows uniqueness of the deformation, since $H^1(X_b,\O_{X_b})=0$.
The bundle $\O_{X_{A'}}$ gives existence. Thus after removing scalar
automorphisms the local classical moduli problem is the base itself.
The reduced tangent and obstruction spaces are both zero.

The stabilized critical algebra is
\[
 R_0[z,w][\xi_z,\xi_w],\quad
 \delta\xi_z=w,\quad\delta\xi_w=z,
 \qquad R_0=\Gamma(V,\O_C).
\]
Its augmentation to $R_0$ is a quasi-isomorphism: the sequence
$(w,z)$ is regular, and the contraction in Corollary~\ref{cor:quadratic}
contracts all positive monomial weights. Its Jacobian module is
$R_0[z,w]/(w,z)=R_0$, so the finite-Jacobian hypotheses hold with
$n=2$. On top forms the relation for a monomial $A=z^a w^b$ is
\[
 [z^{a+1}w^b\,dz\wedge dw]
       =u b[z^a w^{b-1}\,dz\wedge dw],\qquad
 [z^a w^{b+1}\,dz\wedge dw]
       =u a[z^{a-1}w^b\,dz\wedge dw].
\]
They follow by applying $u\,d-d(zw)\wedge$ to $A\,dz$ and $A\,dw$.
They reduce a polynomial top form to a multiple of $dz\wedge dw$;
the contraction of Theorem~\ref{thm:v5-finite-Jacobian} proves
nonvanishing and handles Laurent completion. The rank is one and
the parity is even because $n=2$. Since $\partial_t f=0$, the
connection is $\partial_t$ with zero matrix in this basis.
This is an auxiliary critical presentation of an actual rigid
object on the same threefold family.
\end{proof}

\subsection{The scope of the product and splitting assumptions}
\label{sssec:v15-remaining-inputs}
For the divisor loci, the cyclic potentials and the local reduced
fourfold equations vanish. In a zero-equation shifted chart with an
odd generator $\theta_1$, the identity $\nu=0=0\theta_1$ has primitive
factor $\theta_1$ and yields the constant connection. It is an example
of the extraction formula, but it does not satisfy the extra
non-zero-divisor assumption used for uniqueness in
Proposition~\ref{prop:v5-factor-test}.
For a nonzero fiberwise potential in the complex polarization
\eqref{eq:r13-fiber-compatibility}, Proposition~\ref{prop:v15-fiber-product-obstruction}
proves an obstruction to a regular product. The separate assertions
about arbitrary real BKSY quotients still require a comparison with
holomorphic critical presentations; a choice of degrees does not
provide that comparison.

For the charge decomposition in
Theorem~\ref{thm:categorical-degeneration-conditional}, a verified
instance is the single chosen locus $M_m$ and its single corresponding
central DT complex. The direct-sum decomposition then has one summand
and is the identity. Finite disjoint unions of such selected moduli
problems give finite direct sums. These examples do not supply a
decomposition of the full moduli space by arbitrary component classes,
nor a product of independent invariants on $X_1$ and $X_2$.
The closed-form comparisons, orientation comparisons, and finiteness
conditions on any additional locus must be checked for that locus.

\renewcommand{\thesection}{\arabic{section}}
\renewcommand{\sectionname}{}
\renewcommand{\theHsection}{\arabic{section}}
\setcounter{section}{12}
\bibliographystyle{amsalpha}
\bibliography{Bibliography}

\end{document}